\documentclass[a4paper,12pt]{article}  
\usepackage{amsmath}
\usepackage{amssymb}
\usepackage{graphicx}
\usepackage{amsthm}       \newtheorem{thm}{Theorem}[section]
\newtheorem{prop}[thm]{Proposition}
\newtheorem{lem}[thm]{Lemma}
\newtheorem{cor}[thm]{Corollary}  \theoremstyle{definition}
\newtheorem{df}[thm]{Definition}   \theoremstyle{definition}
\newtheorem{ques}[thm]{Question}
\newtheorem{prob}[thm]{Problem}
\newtheorem{rem}[thm]{Remark}                \theoremstyle{plain}
 \theoremstyle{definition}
\newtheorem{ex}[thm]{Example}   \def\CC{\Bbb{C}}
\def\RR{\Bbb{R}}  
\def\CCI{\hat{\CC}}        \def\NN{\Bbb{N}} 

\def\g{\gamma}
\def\G{\Gamma}
\def\GN{\Gamma ^{\NN }}
\def\B1{{\rm\kern.32em\vrule    width.12em       height1.4ex
depth-.05ex\kern-.28em 1}}
\def\pc{\pi _{\CCI }}

\def\GNCR{\GN \times \CCI \rightarrow \GN \times \CCI }
\begin{document}
\title
{Dynamics of 
postcritically bounded\\  
polynomial 
semigroups}
\author{Hiroki Sumi\\ 
 Department of Mathematics \\ 
Graduate School of Science\\ 
Osaka University \\ 1-1, \ Machikaneyama,\ Toyonaka,\ Osaka,\ 560-0043,\ 
Japan\\ E-mail: sumi@math.sci.osaka-u.ac.jp\\ 
http://www.math.sci.osaka-u.ac.jp/$\sim $sumi/welcomeou-e.html}
\date{November 26, 2007}
\maketitle 

\begin{abstract}
We investigate the dynamics of 
semigroups generated by polynomial maps on the Riemann sphere 
such that the postcritical set in the complex plane is bounded. 
Moreover, we investigate the associated random dynamics of polynomials.
We show that for such a 
polynomial semigroup, 
if $A$ and $B$ are two connected components of the Julia set, 
then one of $A$ and $B$ surrounds the other. 
A criterion for the Julia set 
to be connected is given.  
Moreover, we show that for any $n\in \Bbb{N} \cup \{ \aleph _{0}\} ,$  
there exists a finitely generated polynomial semigroup with bounded 
planar postcritical set 
such that the cardinality of the set of all connected components of the 
Julia set is equal to $n.$ 
%
 Furthermore, we investigate the fiberwise dynamics of skew products 
 related to  
polynomial semigroups with bounded planar postcritical set. 
Using uniform 
fiberwise quasiconformal surgery 
on a fiber bundle, 
we show that 
if the Julia set of such a semigroup 
is disconnected, then 
there exist families of uncountably many mutually disjoint quasicircles 
with uniform dilatation which are parameterized by the Cantor set, densely 
inside the Julia set of the semigroup. 
Moreover, 
we show that 
under a certain condition, 
a random Julia set is almost surely a Jordan curve, but not a quasicircle.
Furthermore, we give a classification of polynomial semigroups 
$G$ such that $G$ is generated by a compact family, 
 the planar postcritical set of $G$ is bounded, and  
$G$ is (semi-) hyperbolic.  
Many new phenomena of polynomial semigroups and random dynamics of 
polynomials that do not occur in the usual dynamics of polynomials  
are found and systematically investigated.
\end{abstract}

\tableofcontents 
\section{Introduction}
 The theory of complex dynamical systems, which has 
 its origin in the important 
 work of Fatou and Julia in the 1910s, 
 has been investigated by many people and discussed in depth.  
In particular, since D. Sullivan showed the famous 
``no wandering domain theorem'' using 
Teichm\"{u}ller theory 
in the 1980s, 
this subject has 
attracted 
many researchers 
from a
wide area. 
For a general 
reference on complex dynamical systems, 
see Milnor's textbook \cite{M}.   
 
 There are several 
areas
 in which we deal with 
  generalized notions of 
classical iteration theory of rational functions.   
One of them is the theory of 
dynamics of rational semigroups 
  (semigroups generated by holomorphic maps on the 
  Riemann sphere $\CCI $), and another one is 
 the theory of   
 random dynamics of holomorphic maps on the Riemann sphere. 

In this paper, we will discuss 
these subjects.
 A {\bf rational semigroup} is a semigroup 
generated by a family of non-constant rational maps on 
$\CCI $, where $\CCI $ denotes the Riemann sphere,
 with the semigroup operation being  
functional composition (\cite{HM1}). A 
{\bf polynomial semigroup} is a 
semigroup generated by a family of non-constant 
polynomial maps.
Research on the dynamics of
rational semigroups was initiated by
A. Hinkkanen and G. J. Martin (\cite{HM1,HM2}),
who were interested in the role of the
dynamics of polynomial semigroups while studying
various one-complex-dimensional
moduli spaces for discrete groups,
and
by F. Ren's group(\cite{ZR, GR}), 
 who studied 
such semigroups from the perspective of random dynamical systems.
Moreover, the research 
on 
rational semigroups is related to 
that 
on 
``iterated function systems" in 
fractal geometry.  
In fact, 
the Julia set of a rational semigroup generated by a 
compact family has 
`` backward self-similarity" 
(cf. Lemma~\ref{hmslem}-\ref{bss}). 
For 
other 
research 
on rational semigroups, see 
\cite{Sta1, Sta2, Sta3, SY, SSS, 
SS, SU1, SU2}, and \cite{S1}--\cite{S10}. 

 The research 
on the  
dynamics of rational semigroups is also directly related to 
that 
on the 
random dynamics of holomorphic maps. 
The first 
study 
in this 
direction was 
by Fornaess and Sibony (\cite{FS}), and 
much research has followed. 
(See \cite{Br, Bu1, Bu2, 
BBR,GQL}.)   

 We remark that the complex dynamical systems 
 can be used to describe some mathematical models. For 
 example, the behavior of the population 
 of a certain species can be described as the 
 dynamical system of a polynomial 
 $f(z)= az(1-z)$ 
 such that $f$ preserves the unit interval and 
 the postcritical set in the plane is bounded 
 (cf. \cite{D}). From this point of view, 
 it is very important to consider the random 
 dynamics of such polynomials (see also Example~\ref{realpcbex}). 
For the random dynamics of polynomials on the unit interval, 
see \cite{Steins}. 
 
 We shall give some definitions 
for the 
dynamics of rational semigroups: 
\begin{df}[\cite{HM1,GR}] 
Let $G$ be a rational semigroup. We set
\[ F(G) = \{ z\in \CCI \mid G \mbox{ is normal in a neighborhood of  $z$} \} ,
\ J(G)  = \CCI \setminus  F(G) .\] \(  F(G)\) is  called the
{\bf Fatou set}  of  $G$ and \( J(G)\)  is  called the {\bf 
Julia set} of $G$. 
We 
let 
$\langle h_{1},h_{2},\ldots \rangle $ 
denote 
the 
rational semigroup generated by the family $\{ h_{i}\} .$
The Julia set of the semigroup generated by 
a single map $g$ is denoted by 
$J(g).$ 
\end{df}

\begin{df}\ 
\begin{enumerate}
\item 
For each rational map $g:\CCI \rightarrow \CCI $, 
we set 
$CV(g):= \{ \mbox{all critical values of }
g: \CCI \rightarrow \CCI \} .$ 
Moreover, for each polynomial map $g:\CCI \rightarrow \CCI $, 
we set $CV^{\ast }(g):= CV(g)\setminus \{ \infty \} .$ 
\item 
Let $G$ be a rational semigroup.
We set 
$$ P(G):=
\overline{\bigcup _{g\in G} CV(g)} \ (\subset \CCI ). 
$$ 
This is called the {\bf postcritical set} of $G.$
Furthermore, for a polynomial semigroup $G$,\ we set 
$P^{\ast }(G):= P(G)\setminus \{ \infty \} .$ This is 
called the {\bf planar postcritical set}
(or {\bf finite postcritical set}) 
 of $G.$
We say that a polynomial semigroup $G$ is 
{\bf postcritically bounded} if 
$P^{\ast }(G)$ is bounded in $\CC .$ 
\end{enumerate}
\end{df}
\begin{rem}
\label{pcbrem}
Let $G$ be a rational semigroup 
generated by a family $\Lambda $ of rational maps. 
Then, we have that 
$P(G)=\overline{\cup _{g\in G\cup \{ Id\} }\ g(\cup _{h\in \Lambda }CV(h))}$, 
where Id denotes the identity map on $\CCI $,  
and that $g(P(G))\subset P(G)$ for each $g\in G.$  
From this formula, one can figure out how the set 
$P(G)$ (resp. $P^{\ast }(G)$) spreads in $\CCI $ (resp. $\CC $). 
In fact, in Section~\ref{Const}, using the above formula, 
we present a way to construct examples of postcritically bounded 
polynomial semigroups (with some additional properties). Moreover, 
from the above formula, one may, in the finitely generated case, 
use a computer to see if a polynomial semigroup $G$ is postcritically bounded much in the same way 
as one verifies the boundedness of the critical orbit for the maps $f_{c}(z)=z^{2}+c.$   
\end{rem}
\begin{ex}
\label{realpcbex}
Let 
$\Lambda := \{ h(z)=cz^{a}(1-z)^{b}\mid 
a,b\in \NN  ,\ c>0,\  
c(\frac{a}{a+b})^{a}(\frac{b}{a+b})^{b}$ $\leq 1\} $ 
and let $G$ be the polynomial semigroup generated by 
$\Lambda .$ 
Since for each $h\in \Lambda $, 
$h([0,1])\subset [0,1]$ and 
$CV^{\ast }(h)\subset [0,1]$, 
it follows that each subsemigroup $H$ of $G$ is postcritically 
bounded. 
\end{ex}
\begin{rem}
\label{pcbound}
It is well-known that for a polynomial $g$ with 
$\deg (g)\geq 2$, 
$P^{\ast }(\langle g\rangle )$ is bounded in $\CC $ if and only if 
$J(g)$ is connected (\cite[Theorem 9.5]{M}).
\end{rem}
As mentioned in Remark~\ref{pcbound}, 
 the planar postcritical set is one 
piece of important information 
regarding the dynamics of polynomials. 
Concerning 
the theory of iteration of quadratic polynomials, 
we have been investigating the famous ``Mandelbrot set''.
    
When investigating the dynamics of polynomial semigroups, 
it is natural for us to 
discuss the relationship between 
  the planar postcritical set and the 
  figure of the Julia set.
The first question in this 
regard 
is: 
\begin{ques}
Let $G$ be a polynomial semigroup such that each 
element $g\in G$ is of degree at least two.
Is $J(G)$ necessarily connected when $P^{\ast }(G)$ is 
bounded in $\CC $?
\end{ques}
The answer is {\bf NO.}
\begin{ex}[\cite{SY}]
Let $G=\langle z^{3}, \frac{z^{2}}{4}\rangle .$ 
Then $P^{\ast }(G) =\{ 0\} $ 
(which is bounded in $\CC $)
and $J(G)$ is disconnected ($J(G)$ is a Cantor set 
of round circles). Furthermore,\ 
according to 
\cite[Theorem 2.4.1]{S5},  
it can be shown that 
a small					 
perturbation $H$ of $G$ 
 still satisfies that 
 $P^{\ast }(H) $ is 
 bounded in $\CC $  and that $J(H)$ is disconnected. 
 ($J(H)$ is a 
Cantor set of quasi-circles with uniform dilatation.)
\end{ex}
\begin{ques}
What happens if $P^{\ast }(G) $ is bounded in $\CC $ 
and $J(G)$ is disconnected? 
\end{ques}
\begin{prob}
Classify postcritically bounded polynomial semigroups.
\end{prob}
In this paper, we show that if $G$ is a postcritically bounded 
polynomial semigroup with disconnected Julia set, then 
$\infty \in F(G)$ (cf. Theorem~\ref{mainth2}-\ref{mainth2-2}), and 
for any two connected components of $J(G)$, one of them surrounds 
the other. This implies that 
there exists an intrinsic total order $``\leq "$ 
(called the 
``surrounding 
order") 
in the space ${\cal J}_{G}$ of 
connected components of $J(G)$, and that 
every connected component of $F(G)$ is either simply 
or doubly connected (cf. Theorem~\ref{mainth1}). 
 Moreover, for such a semigroup $G$, we show  
that the interior of ``the smallest filled-in Julia set'' $\hat{K}(G)$ 
is not empty, and that there exists a maximal element and a 
minimal element 
in the space ${\cal J}_{G}$ endowed with the order $\leq $  
(cf. Theorem~\ref{mainth2}).  
From these results, we obtain the result 
 that for a postcritically bounded polynomial semigroup $G$, 
 the Julia set $J(G)$ is uniformly perfect, 
 even if $G$ is not generated by a compact family of polynomials 
 (cf. Theorem~\ref{mainupthm}). 

 Moreover, 
 we utilize Green's functions with pole at infinity
to show that for a postcritically bounded 
 polynomial semigroup $G$, the cardinality of 
the set of all connected components of $J(G)$ is less than or equal to 
 that of $J(H)$, where $H$ is the ``real affine semigroup'' 
 associated with $G$ (cf. Theorem~\ref{polyandrathm1}). 
From this result, we obtain a sufficient condition for the Julia set 
of a postcritically bounded polynomial semigroup to be connected 
(cf. Theorem~\ref{polyandrathm2}). 
In particular, 
we show that if a postcritically bounded polynomial semigroup $G$ 
is generated by a family of quadratic polynomials, then $J(G)$ is connected 
(cf. Theorem~\ref{polyandrathm3}). 
 The proofs of the results in 
this and the previous paragraphs 
 are not straightforward. 
In fact, we first prove 
(1) that 
for any two connected components of $J(G)$ that 
 are included in $\CC $, one of them surrounds the other; 
 next, using (1) and the theory of Green's functions, we prove 
(2) that 
the cardinality of the set of all connected components of 
 $J(G)$ is less than or equal to that of $J(H)$, where 
 $H$ is the associated real affine semigroup; and finally, 
 using (2) and (1), we prove 
(3) that 
$\infty \in F(G)$, 
 int$(\hat{K}(G))\neq \emptyset $, and other results in the previous 
 paragraph.     

 Moreover, we show that for any $n\in \NN \cup \{ \aleph _{0}\} $, 
 there exists a finitely generated, postcritically bounded,  
 polynomial semigroup $G$ such that the cardinality of the set of 
 all connected components of $J(G)$ is equal to $n$ 
 (cf. Proposition~\ref{fincomp}, Proposition~\ref{countprop} and 
 Proposition~\ref{countcomp}). 
A sufficient condition for the cardinality of the set of all connected components 
of 
a Julia set 
to be equal to $\aleph _{0}$ is also given 
(cf. Theorem~\ref{countthm}). To obtain these results, we 
use 
the fact that 
the map induced by any element of a semigroup on the space of connected components of the Julia set preserves the order $\leq $ (cf. Theorem~\ref{mainth1}).   
  Note that 
this is in contrast to the dynamics of 
a single rational map $h$ or a non-elementary 
 Kleinian group, where it is known that either the Julia set is connected, or 
 the Julia set has uncountably many connected components.  

 Applying the previous results, 
 we investigate the dynamics of 
every sequence, or fiberwise dynamics of the skew product 
associated with the generator system (cf. Section \ref{fibsec}). 
Moreover, we investigate 
the random dynamics of polynomials 
acting on 
the Riemann 
sphere. 
Let us consider a polynomial semigroup $G$ generated by a compact 
family $\G $ of polynomials. For each sequence 
$\g =(\g _{1},\g _{2},\g _{3},\ldots )\in \GN $, 
we examine the dynamics along the sequence $\g $, 
that is, the dynamics of the family of maps 
$\{ \g _{n}\circ \cdots \circ \g _{1}\} _{n=1}^{\infty }$. 
We note that this corresponds to the fiberwise dynamics of the skew 
product (see Section \ref{fibsec}) associated with the generator 
system $\G .$ 
We show that 
if $G$ is postcritically bounded, $J(G)$ is disconnected, 
and $G$ is generated by a compact family $\G $ of 
polynomials; 
then, for almost every sequence $\g \in \GN $, 
there exists exactly one bounded component $U_{\g }$ 
of 
the Fatou set 
of $\g $, the Julia set of $\g $ has Lebesgue measure zero, 
there exists no non-constant limit function in $U_{\g }$ for the 
sequence $\g $, 
and for any point $z\in U_{\g }$, 
the orbit along $\g $ 
tends to the interior of the smallest filled-in Julia set 
$\hat{K}(G)$ of $G$ (cf. Theorem~\ref{mainth3}-\ref{mainth3-1}, 
Corollary~\ref{rancor1}). 
 Moreover, using the uniform fiberwise quasiconformal surgery 
 (cf. Theorem~\ref{hypskewqc}), 
 we find 
sub-skew products $\overline{f}$ 
 such that $\overline{f}$ is hyperbolic and 
 such that every fiberwise Julia 
set 
of $\overline{f}$ is a $K$-quasicircle, 
 where $K$ is a constant not depending 
 on the fibers (cf. Theorem~\ref{mainth3}-\ref{mainth3-3}). 
Reusing 
the uniform fiberwise quasiconformal surgery,
we show that if $G$ is a postcritically bounded polynomial semigroup 
with disconnected Julia set, then for any non-empty 
open subset $V$ 
of $J(G)$, 
there exists a 
$2$-generator subsemigroup $H$ of $G$ such that 
$J(H)$ is the disjoint union of ``Cantor family of quasicircles" 
(a family of quasicircles parameterized by a Cantor set) with uniform 
distortion, and such that $J(H)\cap V\neq \emptyset $ 
(cf. Theorem~\ref{cantorqc}). 
Note that the uniform fiberwise quasiconformal surgery 
is based on solving uncountably many Beltrami  
equations (a kind of partial differential equations). 

 We also investigate (semi-)hyperbolic, postcritically bounded, 
 polynomial semigroups generated by a compact family $\G $ of 
 polynomials.   
We show that if $G$ is such a semigroup with disconnected Julia set, 
and if there exists an element $g\in G$ such that 
$J(g)$ is not a Jordan curve, then, 
for almost 
every sequence $\g \in \GN $, the Julia set of $\g $ is a 
Jordan curve but not a quasicircle, the basin of infinity $A_{\g }$ is a 
John domain, and the bounded component $U_{\g }$ of 
the Fatou set is not a John domain (cf. Theorem~\ref{mainthjbnq}). 
Moreover, we 
classify the 
semi-hyperbolic, postcritically bounded, 
polynomial semigroups generated by a compact family $\G $ of 
polynomials.  
We show that such a semigroup $G$ satisfies either  
(I) every fiberwise Julia 
set is a quasicircle 
with uniform distortion, 
or (II) for almost 
every sequence $\g \in \GN $, the Julia set $J_{\g }$ is a 
Jordan curve but not a quasicircle, the basin of infinity 
$A_{\g }$ is a John domain, and the bounded component $U_{\g }$ of the Fatou 
set is not a John domain, or (III) for every $\alpha ,\beta \in \GN $, 
the intersection of the Julia sets $J_{\alpha }$ and $J_{\beta }$ is not empty, 
and $J(G)$ is arcwise connected (cf. Theorem~\ref{mainthran1}). 
Furthermore, we also 
classify the 
hyperbolic, postcritically bounded, 
polynomial semigroups generated by a compact family $\G $ of 
polynomials.  
We show that such a semigroup $G$ satisfies either 
(I) above, or (II) above, or (III)' for every $\alpha ,\beta \in \GN $, 
the intersection of the Julia sets $J_{\alpha }$ and $J_{\beta }$ is not empty, $J(G)$ is arcwise connected, and for every sequence $\g \in \GN $, 
there exist infinitely many bounded components of $F_{\g }$ 
(cf. Theorem~\ref{mainthran2}). We give some examples 
of 
situation (II) above (cf. Example~\ref{jbnqexfirst}, 
Example~\ref{jbnqex} and 
Section~\ref{Const}). 
  Note that 
situation (II) above is 
a 
special and new phenomenon of 
  random dynamics of polynomials that  
  does not occur in the usual dynamics of polynomials. 

The key to 
investigating the 
dynamics of 
postcritically bounded polynomial semigroups is 
the density of repelling fixed points in the Julia set (cf. 
Theorem~\ref{repdense}), which 
can be shown by 
an application of 
the Ahlfors five island theorem, and the lower semi-continuity of 
$\g \mapsto J_{\g }$ (Lemma~\ref{fibfundlem}-\ref{fibfundlem2}), which is a consequence of potential theory.  
Moreover, one of the keys to 
investigating the 
fiberwise dynamics of 
skew products is, the observation of non-constant limit functions
(cf. Lemma~\ref{nclimlem} and \cite{S1}). 
 The key to 
investigating the 
dynamics of semi-hyperbolic polynomial semigroups is, 
the continuity of the map $\g \mapsto J_{\g }$ 
(this is highly nontrivial; see \cite{S1}) 
and the 
Johnness of the basin $A_{\g }$ of infinity (cf. \cite{S4}). 
Note that the continuity of the map $\g \mapsto J_{\g }$ 
does not hold in general, if we do not assume semi-hyperbolicity. 
 Moreover, one of the 
original aspects 
of this paper 
is the idea of 
 ``combining both 
the theory of rational semigroups and 
 that of random complex dynamics". It is quite natural to 
investigate both fields simultaneously. However, 
no study 
thus far has done so.

 Furthermore, in Section~\ref{Const} and Section~\ref{Poly}, we 
provide 
a way of constructing examples of 
 postcritically bounded polynomial semigroups 
 with 
some additional properties (disconnectedness of Julia set, 
semi-hyperbolicity, hyperbolicity, etc.) 
(cf. Proposition~\ref{Constprop}, Theorem~\ref{shshfinprop}, 
Theorem~\ref{sphypopen}). 
For example, by Proposition~\ref{Constprop}, 
there exists a $2$-generator 
postcritically bounded polynomial semigroup $G=\langle h_{1},h_{2}\rangle $ 
with disconnected 
Julia set such that $h_{1}$ has a Siegel disk.
        
As wee see in Example~\ref{realpcbex} and Section~\ref{Const}, 
it is not difficult to construct many examples,  
it is not difficult to verify the hypothesis ``postcritically 
bounded'', 
and the class of postcritically bounded polynomial semigroups is 
very wide.  

 Throughout the paper, we will see many new phenomena in polynomial 
 semigroups or random dynamics of polynomials that do not occur in 
 the usual dynamics of polynomials. Moreover, these new phenomena are 
 systematically investigated.
 
 In Section~\ref{Main}, we present the main results 
 of this paper. We give some tools in Section~\ref{Tools}. 
 The proofs of the main results are given in Section~\ref{Proofs}. 

\ 

There are many applications of the results of postcritically 
bounded polynomial semigroups in many directions. 
 In the sequel \cite{S8}, we 
 will investigate 
Markov process on $\CCI $ associated with 
the random dynamics of polynomials and 
we will consider the probability $T_{\infty }(z)$ 
of tending to $\infty \in \CCI $ 
starting with the initial value $z\in \CCI .$ 
Applying many results of this paper, 
it will be shown in \cite{S8} that 
if the associated polynomial semigroup $G$ 
is postcritically bounded and the Julia set is 
disconnected, then the function $T_{\infty }$ defined on $\CCI $ 
has many interesting properties which are 
similar to those of the Cantor function. 
Such a kind of ``singular functions in the 
complex plane'' appear very naturally in 
random dynamics of polynomials and the 
results of this paper (for example, 
the results on the space of all connected 
components of a Julia set) are the keys to 
investigating that. 
(The above results have been announced in \cite{S9, S10}.) 

 Moreover, as illustrated before, 
 it is very important for us to recall that 
 the complex dynamics can be applied to describe some mathematical models. 
 For example, the behavior of the population of a 
 certain species can be described as the dynamical systems
  of a polynomial $h$ such that $h$ preserves the unit interval  
 and the postcritical set in the plane is bounded. 
 When one considers such a model, 
 it is very natural to consider the random dynamics of 
 polynomial with bounded postcritical set in the plane 
 (see Example~\ref{realpcbex}).    

 In the sequel \cite{SS}, we will give some 
further results on postcritically 
 bounded polynomial semigroups, based on this paper. 
 Moreover, in the sequel \cite{S7}, 
 we will define a new kind of cohomology theory, in order to 
 investigate the action of finitely generated semigroups, and 
we will apply it to the study of 
the dynamics of postcritically bounded polynomial semigroups.  

\ 

\noindent {\bf Acknowledgement:} 
The author 
thanks 
R. Stankewitz for many 
valuable comments.  
  
\section{Main results}
\label{Main}
In this section we present the statements of the main results.
Throughout this paper, we deal with semigroups $G$ that 
might not 
be generated by a compact family of polynomials.
The proofs are given in Section~\ref{Proofs}.

\subsection{Space of connected components of a Julia set, surrounding order}
\label{concompsec}
We present some results 
concerning the connected components of the 
Julia set of a postcritically bounded polynomial semigroup.
 The proofs are 
given in Section~\ref{pfconcompsec}.
\begin{thm}
\label{mainth0}
Let $G$ be a rational semigroup 
generated by a family $\{ h_{\lambda }\} _{\lambda \in 
\Lambda }.$ 
Suppose 
that  
there exists a connected component 
$A$ of $J(G)$ such that $\sharp A> 1$ and 
$\cup _{\lambda \in \Lambda }
J(h_{\lambda })\subset A.$ Moreover, 
suppose that for any $\lambda \in \Lambda $ such that 
$h_{\lambda }$ is a M\"{o}bius transformation of finite order,\ 
we have $h_{\lambda }^{-1}(A)\subset A.$
Then,\ 
$J(G)$ is connected.
\end{thm} 
\begin{df}
We set 
Rat : $=\{ h:\CCI \rightarrow \CCI \mid 
h \mbox { is a non-constant rational map}\} $
endowed with the topology induced by uniform convergence on $\CCI $ 
with respect to the spherical distance.   
We set 
Poly :$=\{ h:\CCI \rightarrow \CCI 
\mid h \mbox{ is a non-constant polynomial}\} $ endowed with 
the relative topology from Rat.   
Moreover, we set 
Poly$_{\deg \geq 2}
:= \{ g\in \mbox{Poly}\mid \deg (g)\geq 2\} $ 
endowed with the relative topology from 
Rat.  
\end{df}
\begin{rem}
Let $d\geq 1$,  $\{ p_{n}\} _{n\in \NN }$ a 
sequence of polynomials of degree $d$, 
and $p$ a polynomial.  
Then, $p_{n}\rightarrow p$ in Poly if and only if 
the coefficients converge appropriately and $p$ is of degree $d.$ 
\end{rem}
\begin{df} 
Let ${\cal G} $ be the set of all polynomial semigroups 
$G$ with the following 
properties:
\begin{itemize}
\item 
 each element of $G$ is of degree 
at least two, and  
\item  $P^{\ast }(G)$ is 
bounded in $\CC $, i.e., $G$ is postcritically bounded.
\end{itemize}   
Furthermore, we set 
${\cal G}_{con}=
\{ G\in {\cal G}\mid 
J(G)\mbox{ is connected}\} $ and 
${\cal G}_{dis}=
\{ G\in {\cal G}\mid 
J(G)\mbox{ is disconnected}\}.$ 
\end{df}  
\noindent {\bf Notation:}
For a polynomial semigroup $G$,\ 
we denote by 
${\cal J}={\cal J}_{G}$ the set of all 
connected components $J$ 
of $J(G)$ such that $J\subset \CC .$   
Moreover, we denote by 
$\hat{{\cal J}}=\hat{{\cal J}}_{G}$ the set of all connected components 
of $J(G).$ 
\begin{rem}
\label{hatjcptrem}
If a polynomial semigroup $G$ is generated by a compact set 
in Poly$_{\deg \geq 2}$, then 
$\infty \in F(G)$ and thus ${\cal J}=\hat{{\cal J}}.$ 
\end{rem}
\begin{df}
For any connected sets $K_{1}$ and 
$K_{2}$ in $\CC ,\ $  ``$K_{1}\leq K_{2}$'' indicates that 
$K_{1}=K_{2}$, or $K_{1}$ is included in 
a bounded component of $\CC \setminus K_{2}.$ Furthermore, 
``$K_{1}<K_{2}$'' indicates $K_{1}\leq K_{2}$ 
and $K_{1}\neq K_{2}.$ Note that 
``$ \leq $'' is a partial order in 
the space of all non-empty compact connected 
sets in $\CC .$ This ``$\leq $" is called 
the {\bf surrounding order.} 
\end{df}

\begin{thm}
\label{mainth1}
Let $G\in {\cal G}$ (possibly generated by a non-compact 
family). 
Then we have all of the following. 
\begin{enumerate}
\item \label{mainth1-1}
$({\cal J},\ \leq )$ is totally ordered. 
\item \label{mainth1-2}
Each connected component of 
$F(G)$ is either simply or doubly connected. 
\item \label{mainth1-3}
For any $g\in G$ and any connected component 
$J$ of $J(G)$,\ we have that  
$g^{-1}(J)$ is connected. 
Let $g^{\ast }(J)$ be the connected component of 
$J(G)$ containing $g^{-1}(J).$ 
If $J\in {\cal J}$, then 
$g^{\ast }(J)\in {\cal J}.$    
If $J_{1},J_{2}\in {\cal J} $ and $J_{1}\leq J_{2},\ $ then 
$g^{-1}(J_{1})\leq g^{-1}(J_{2})$ 
and $g^{\ast }(J_{1})\leq g^{\ast }(J_{2}).$
\end{enumerate}
\end{thm} 
For the figures of the Julia sets of semigroups $G\in {\cal G}_{dis}$, 
see figure~\ref{fig:3mapcountjulia2} and figure~\ref{fig:dcgraph}. 
\subsection{Upper estimates of $\sharp (\hat{{\cal J}})$}
\label{Upper}
Next, we present some results on the space $\hat{{\cal J}}$
and some results on upper estimates of 
$\sharp (\hat{{\cal J}}).$ The proofs are given in 
Section~\ref{Proof of Upper} and Section~\ref{Proof of Properties}.
\begin{df}
\ 
\begin{enumerate}
\item For a polynomial $g$, we denote by 
$a(g)\in \CC $ the coefficient of 
the highest degree term of $g.$ 
\item 
We set 
RA $:=\{ ax+b\in \RR [x]\mid a,b\in \RR ,a\neq 0\} $ 
endowed with the topology such that, 
$a_{n}x+b_{n}\rightarrow ax+b$ if and only if 
$a_{n}\rightarrow a$ and $b_{n}\rightarrow b.$ 
 The space RA is a semigroup 
with the semigroup operation being functional composition.
Any subsemigroup of RA will be called a {\em real affine semigroup}.
We define a map $\Psi :$ Poly $\rightarrow $ RA as follows: 
For a polynomial $g\in $ Poly,   
we set $\Psi (g)(x):= \deg (g)x+\log | a(g)|.$ 

 Moreover, for a polynomial semigroup $G$, 
 we set $\Psi (G):= \{ \Psi (g)\mid g\in G\}  $ ($\subset ${\em RA}). 

\item 
We set $\hat{\RR }:= \RR \cup \{ \pm \infty \} $ endowed with 
the topology such that 
$\{ (r,+\infty ]\} _{r\in \RR }$ makes a fundamental neighborhood system of 
$+\infty $, and such that $\{ [-\infty ,r)\} _{r\in \RR }$ makes a 
fundamental neighborhood system of $-\infty .$ 
For a real affine semigroup $H$, we set 
$$M(H):= \overline{ \{ x\in \RR \mid \exists h\in H,  h(x)=x, |h'(x)|>1\} } 
\ (\subset \hat{\RR })
,$$
where the closure is taken in the space $\hat{\RR }.$   
Moreover, we denote by ${\cal M}_{H}$ the set of all connected components of 
$M(H).$ 
\item 
We denote by $\eta : $ RA $\rightarrow $ Poly the natural embedding 
defined by $\eta (x\mapsto ax+b)=(z\mapsto az+b)$, where 
$x\in \RR $ and $z\in \CC .$  
\item 
We define a map $\Theta :$ Poly $\rightarrow $ Poly  as follows. 
For a polynomial $g$, we set 
$\Theta (g)(z)=a(g)z^{\deg (g)}.$ 
Moreover, for a polynomial semigroup $G$, 
we set $\Theta (G):= \{ \Theta (g)\mid g\in G\} .$
\end{enumerate} 

\end{df}
\begin{rem}
\ 
\begin{enumerate}
\item 
The map $\Psi :$ Poly $\rightarrow $ RA is a semigroup homomorphism. 
That is, we have $\Psi (g\circ h)=\Psi (g)\circ \Psi (h).$
Hence, for a polynomial semigroup $G$, 
the image $\Psi (G)$ is a real affine semigroup. 
Similarly, the map $\Theta :$ Poly $\rightarrow $ Poly 
is a semigroup homomorphism. Hence, 
for a polynomial semigroup $G$, 
the image $\Theta (G)$ is a polynomial semigroup.  
\item 
The maps $\Psi :$ Poly $\rightarrow $ RA,  
$\eta :$ RA $\rightarrow $ Poly, 
and $\Theta :$ Poly $\rightarrow $ Poly are continuous.
\end{enumerate}
\end{rem}
\begin{df}
For any connected sets $M_{1}$ and $M_{2}$ in $\hat{\RR  }$, 
``$M_{1}\leq _{r} M_{2}$'' indicates that $M_{1}=M_{2}$, or 
each $(x,y)\in M_{1}\times M_{2}$ satisfies 
$x<y.$ Furthermore, 
``$M_{1}< _{r} M_{2}$'' indicates $M_{1}\leq _{r}M_{2}$ and 
$M_{1}\neq M_{2}.$ 
\end{df}
\begin{rem}
The above ``$\leq _{r}$" is a partial order in the space of 
non-empty connected subsets of $\hat{\RR }.$ 
Moreover, for each real affine semigroup $H$, 
$({\cal M}_{H},\leq _{r})$ is totally ordered.
\end{rem}
\begin{thm}
\label{polyandrathm1}
\ 
\begin{enumerate}
\item \label{polyandrathm1-1}
Let $G$ be a polynomial semigroup in ${\cal G}.$ 
Then, we have 
$\sharp (\hat{{\cal J}}_{G})\leq \sharp ({\cal M}_{\Psi (G)}).$ 
More precisely, there exists an injective map 
$\tilde{\Psi }: \hat{{\cal J}}_{G}\rightarrow 
{\cal M}_{\Psi (G)}$ such that if 
$J_{1},J_{2}\in {\cal J}_{G}$ and $J_{1}<J_{2}$, 
then 
$\tilde{\Psi }(J_{1})<_{r}\tilde{\Psi }(J_{2}).$ 
\item \label{polyandrathm1-2}
If $G\in {\cal G}_{dis}$, then we have that 
$M(\Psi (G))\subset \RR $ and $M(\Psi (G))=J(\eta (\Psi (G))).$  
\item \label{polyandrathm1-3}
Let $G$ be a polynomial semigroup in ${\cal G}.$ Then, 
$\sharp (\hat{{\cal J}}_{G})\leq 
\sharp ({\hat{\cal J}}_{\eta (\Psi (G))}).$
\end{enumerate}
\end{thm}
\begin{cor}
\label{polyandracor}
Let $G$ be a polynomial semigroup in ${\cal G}.$ 
Then, we have 
$\sharp ({\hat{\cal J}}_{G})\leq \sharp ({\hat{\cal J}}_{\Theta (G)}).$ 
More precisely, there exists an injective map 
$\tilde{\Theta }: \hat{{\cal J}}_{G}\rightarrow 
\hat{{\cal J}}_{\Theta (G)}$ such that 
if $J_{1},J_{2}\in {\cal J}_{G}$ and $J_{1}<J_{2}$, 
then $\tilde{\Theta }(J_{1})\in {\cal J}_{\Theta (G)}$, 
 $\tilde{\Theta }(J_{2})\in {\cal J}_{\Theta (G)}$, and 
 $\tilde{\Theta }(J_{1})<\tilde{\Theta }(J_{2}).$  
\end{cor}
\begin{thm}
\label{polyandrathm2}
Let $G=\langle h_{1},\ldots ,h_{m}\rangle $ be a finitely generated 
polynomial semigroup in ${\cal G}.$ 
For each $j=1,\ldots ,m$, 
let $a_{j}$ be the coefficient of the highest degree term 
of polynomial $h_{j}.$ Let $\alpha := 
\min  _{j=1,\ldots ,m}\{ \frac{-1}{\deg (h_{j})-1}\log |a_{j}|\} $ 
and $\beta := 
\max _{j=1,\ldots ,m}\{ \frac{-1}{\deg (h_{j})-1}\log |a_{j}|\} .$ 
We set $[\alpha ,\beta ]:= \{ x\in \RR \mid \alpha\leq x\leq \beta \} .$  
If 
$[\alpha ,\beta ]\subset \cup _{j=1}^{m}\Psi (h_{j})^{-1}([\alpha ,\beta ])$,  
then $J(G)$ is connected. 
\end{thm}
\begin{thm}
\label{polyandrathm3}
Let $G$ be a 
polynomial semigroup in ${\cal G}$ generated by a 
(possibly non-compact) 
family of 
polynomials of degree two. Then, $J(G)$ is connected.
\end{thm}
\begin{thm}
\label{polyandrathm4}
Let $G$ be a polynomial semigroup in ${\cal G}$ generated by 
a (possibly non-compact) family $\{ h_{\lambda }\} _{\lambda \in \Lambda }$ of 
polynomials. 
Let $a_{\lambda }$ be the coefficient of the highest degree term of 
the polynomial $h_{\lambda }.$ 
Suppose that for any  $\lambda ,\xi \in \Lambda $, 
we have $(\deg (h_{\xi })-1)\log |a_{\lambda }|=
(\deg (h_{\lambda })-1)\log |a_{\xi }|.$ Then, 
$J(G)$ is connected.
\end{thm}

\subsection{Properties of ${\cal J}$}
\label{Properties}
In this section, we present some results on 
${\cal J}.$ The proofs are given in Section~\ref{Proof of Properties}.
\begin{df}
 For a polynomial semigroup $G$,\ we set 
$$ \hat{K}(G):=\{ z\in \CC 
\mid \bigcup _{g\in G}\{ g(z)\} \mbox{ is bounded in }\CC \} $$ 
and call $\hat{K}(G)$ the {\bf smallest filled-in Julia set} of 
$G.$ 
For a polynomial $g$, we set $K(g):= \hat{K}(\langle g\rangle ).$ 

\end{df}
\noindent {\bf Notation:} 
For a set $A\subset \CCI $, we denote by int$(A)$ the set of 
all  
interior points of $A.$ 
\begin{prop}
\label{fcprop}
Let $G\in {\cal G}.$  If $U$ is a connected component 
of $F(G)$ such that $U\cap \hat{K}(G)\neq \emptyset $, 
then $U\subset $ {\em int}$(\hat{K}(G))$ and $U$ is 
simply connected. Furthermore, 
we have $\hat{K}(G)\cap F(G)=$ {\em int}$(\hat{K}(G)).$

\end{prop}
%
%
%
%
%
%
%
%
%
\noindent {\bf Notation:} 
For a polynomial semigroup $G$ with 
$\infty \in F(G)$, we denote by 
$F_{\infty }(G)$ the connected component of $F(G)$ containing 
$\infty .$ Moreover, for a polynomial $g$ with 
$\deg (g)\geq 2$, we set $F_{\infty }(g):= 
F_{\infty }(\langle g\rangle ).$ 

\ 

The following theorem is the key to obtaining further results of postcritically bounded 
polynomial semigroups and related random dynamics of polynomials. 
\begin{thm}
\label{mainth2}
Let $G\in {\cal G} _{dis}$ (possibly generated by a non-compact family). Then, under the above notation,\ 
we have the following.
\begin{enumerate}
\item \label{mainth2-2}
We have that $\infty \in F(G)$ and the  
 connected component 
$F_{\infty }(G)$ of $F(G)$ containing $\infty $ 
is simply connected. 
Furthermore,\ 
the element $J_{\max }=J_{\max}(G)\in {\cal J}$  
containing $\partial F_{\infty }(G)$ 
is the unique element of ${\cal J}$ satisfying that 
$J\leq J_{\max }$ for each 
$J\in {\cal J}.$  
\item 
\label{mainth2-3}
There exists a unique element 
$J_{\min }=J_{\min }(G)\in {\cal J}$ such that 
$J_{\min }\leq J$ for 
  each element $J\in {\cal J}. $
  Furthermore, let $D$ be the unbounded 
  component of $\CC \setminus J_{\min }. $ 
Then,     
$ P^{\ast }(G) \subset \hat{K}(G)\subset 
  \CC \setminus D $ and 
  $\partial \hat{K}(G)\subset J_{\min }. $
\item \label{mainth2-3b}
If $G$ is generated by a family 
$\{ h_{\lambda }\} _{\lambda \in \Lambda },$ 
then there exist two elements $\lambda _{1}$ 
and $\lambda _{2}$ of $\Lambda $ satisfying: 
  \begin{enumerate}
  \item there exist two elements $J_{1}$ and 
  $J_{2}$ of ${\cal J} $ with  
  $J_{1}\neq J_{2}$ such that $J(h_{\lambda _{i}})
  \subset J_{i}$ for each $i=1,2$;  
  \item $J(h_{\lambda _{1}})\cap J_{\min }=\emptyset $; 
  \item for each $n\in \NN $, 
  we have $h_{\lambda _{1}}^{-n}(J(h_{\lambda _{2}}))
  \cap J(h_{\lambda _{2}})=\emptyset $ and 
  $h_{\lambda _{2}}^{-n}(J(h_{\lambda _{1}}))
  \cap J(h_{\lambda _{1}})=\emptyset $; and 
  \item $h_{\lambda _{1}}$ has an attracting 
  fixed point $z_{1}$ in $\CC $,    
{\em int}$(K(h_{\lambda _{1}}))$ consists of 
  only one immediate attracting basin 
  for $z_{1}$,   
  and 
  $K(h_{\lambda _{2}})\subset $ {\em int}$(K(h_{\lambda _{1}})).$ 
  Furthermore, $z_{1}\in $ {\em int}$(K(h_{\lambda _{2}})).$ 
  \end{enumerate}
\item 
\label{mainth2ast1}
For each $g\in G$ with 
  $J(g)\cap J_{\min }=\emptyset $,\ 
  we have that $g$ has an attracting 
  fixed point $z_{g}$ in $\CC $, {\em int}$(K(g))$ consists of 
  only one immediate attracting basin for $z_{g}$,   
  and  $J_{\min }\subset $ {\em int}$(K(g)).$ 
  Note that 
it is not necessarily true that $z_{g}=z_{f}$ when 
$g,f\in G$ are such that $J(g)\cap J_{\min }=\emptyset $ 
and $J(f)\cap J_{\min }=\emptyset $ (see Proposition~\ref{fincomp}).  
\item 
\label{mainth2-4}
We have that  $\mbox{{\em int}}(\hat{K}(G))\neq 
\emptyset .$ Moreover,  

  \begin{enumerate}
  \item \label{mainth2-4-1}
  $\CC \setminus J_{\min }$ is 
  disconnected, $\sharp J\geq 2$ for each 
  $J\in \hat{{\cal J}}$, and  
  
  \item \label{mainth2-4-2}
  for each $g\in G$ with 
  $J(g)\cap J_{\min }=\emptyset $,\ 
  we have that $J_{\min }<g^{\ast }(J_{\min })$, 
  $g^{-1}(J(G))\cap J_{\min }=\emptyset $,  
  $g(\hat{K}(G)\cup J_{\min })\subset $ {\em int}$(\hat{K}(G))$,    
  and   
  the unique attracting fixed point $z_{g}$  
  of $g$ in $\CC $ belongs to 
  $\mbox{{\em int}}(\hat{K}(G)).$ 
  \end{enumerate} 
  
\item 
\label{mainth2ast2}
Let ${\cal A} $ be the set of all doubly connected components 
of $F(G).$ Then, $\cup _{A\in {\cal A}}A\subset \CC $ and 
$({\cal A},\leq )$ is totally ordered.    
\end{enumerate}
\end{thm}
We present a result on uniform perfectness of the Julia sets 
of semigroups in ${\cal G}.$ 
\begin{df}
A compact set $K$ in $\CCI $ is said to be 
uniformly perfect if 
$\sharp K\geq 2$ and there exists a constant $C>0$ such that 
each annulus $A$ that separates $K$ satisfies that 
mod $A<C$, where mod $A$ denotes the modulus of $A$ 
(See the definition in \cite{LV}). 
\end{df}  
\begin{thm}
\label{mainupthm}
\ 
\begin{enumerate}
  \item \label{mainupthm1}
  Let $G$ be a polynomial semigroup in $ {\cal G}.$ Then,  
  $J(G) $ is uniformly perfect. Moreover, if $z_{0}\in J(G)$ is a 
  superattracting fixed point of an element of $G$, then 
  $z_{0}\in $ {\em int}$(J(G)).$ 
\item 
\label{mainupthm2} 
If $G\in {\cal G}$ and $\infty \in J(G)$, then 
$G\in {\cal G}_{con}$ and $\infty \in $ {\em int}$(J(G)).$  
  \item \label{mainupthm3}
  Suppose that $G\in {\cal G}_{dis}.$   
Let $z_{1}\in J(G)\cap \CC $ be a superattracting 
fixed point of $g\in G.$   
Then 
$z_{1}\in $ {\em int}$(J_{\min })$ and 
  $J(g)\subset J_{\min }.$  
\end{enumerate}

\end{thm}
We remark that in \cite{HM2}, it was shown that there exists a rational semigroup $G$ such that 
$J(G)$ is not uniformly perfect. 

 We now present results on the Julia sets of subsemigroups of an element of ${\cal G}_{dis}.$ 
 \begin{prop}
\label{orderjprop}
Let $G\in {\cal G}_{dis}$ 
and let $J_{1}, J_{2}\in {\cal J}={\cal J}_{G}$  
with $J_{1}\leq J_{2}.$ 
%
%
Let $A_{i}$ be the unbounded component of 
$\CC \setminus J_{i}$ for each $i=1,2.$   
Then, we have the following.
\begin{enumerate}
\item \label{orderjprop1}
Let $Q_{1}=\{ g\in G\mid \exists J \in {\cal J} \mbox{with }
J_{1}\leq J,\ J(g)\subset J\} $ 
and let $H_{1}$ be the subsemigroup of $G$ generated by 
$Q_{1}.$ Then 
$J(H_{1})\subset J_{1}\cup A_{1}.$ 
\item \label{orderjprop2}
Let $Q_{2}=\{ g\in G\mid \exists J \in {\cal J} \mbox{with }
J\leq J_{2},\ J(g)\subset J\} $ 
and let $H_{2}$ be the subsemigroup of $G$ generated by 
$Q_{2}.$ Then 
$J(H_{2})\subset \CC \setminus A_{2}.$ 
\item \label{orderjprop3} 
Let $Q=\{ g\in G\mid \exists J \in {\cal J} \mbox{with }
J_{1}\leq J\leq J_{2},\ J(g)\subset J\} $ 
and let $H$ be the subsemigroup of $G$ generated by 
$Q.$ Then 
$J(H)\subset J_{1}\cup (A_{1}\setminus A_{2}).$
\end{enumerate}
\end{prop}

\begin{prop}
\label{bminprop}
Let $G$ be a 
polynomial 
semigroup generated by a compact subset 
$\G $ of {\em Poly}$_{\deg \geq 2}.$ Suppose that $G\in {\cal G}_{dis}.$ Then,    
there exists an element 
$h_{1}\in \G $ with 
 $J(h_{1})\subset J_{\max } $ and 
there exists an element 
$h_{2}\in \G $ with 
$J(h_{2})\subset J_{\min }.$ 
\end{prop}

\subsection{Finitely generated polynomial semigroups $G\in {\cal G}_{dis}$ such that 
$2\leq \sharp (\hat{{\cal J }}_{G})\leq \aleph _{0}$}
\label{Poly}
In this section, we present some results on various 
finitely generated polynomial semigroups $G\in {\cal G}_{dis}$ such that 
$2\leq \sharp (\hat{{\cal J}}_{G})\leq \aleph _{0}.$ 
The proofs are given in Section~\ref{Proof of Poly}.

It is well-known that for a rational map $g$ 
with $\deg (g)\geq 2$, if $J(g)$ is 
disconnected, then $J(g)$ has uncountably 
many connected components 
(See \cite{M}). 
Moreover, if $G$ is a non-elementary 
Kleinian group with disconnected 
Julia set (limit set), then  
$J(G)$ has uncountably many connected components. 
However,  
for general rational semigroups, 
we have the following examples.
\begin{thm}
\label{fcthm}
Let $G$ be a polynomial semigroup in ${\cal G}$ 
generated by a (possibly non-compact) family $\G $ 
in {\em Poly}$_{\deg \geq 2}.$  
Suppose that there exist mutually distinct 
elements $J_{1},\ldots , J_{n}\in \hat{{\cal J}}_{G}$ 
such that for each $h\in \G $ and each $j\in \{ 1,\ldots ,n\} $,  
there exists an element $k\in \{ 1,\ldots ,n\} $ with 
$h^{-1}(J_{j})\cap J_{k}\neq \emptyset .$ 
Then, we have $\sharp (\hat{{\cal J}}_{G})=n.$   
\end{thm}

\begin{prop}
\label{fincomp}
For any $n\in \NN  $ with $n>1$, there exists a 
finitely generated 
polynomial semigroup $G_{n}=\langle h_{1},\ldots ,h_{2n}\rangle $ 
in ${\cal G}$ satisfying 
$\sharp (\hat{{\cal J}}_{G_{n}})=n.$ 
In fact, let $0<\epsilon <\frac{1}{2}$ and 
we set for each $j=1,\ldots ,n,\ a_{j}(z):=\frac{1}{j} z^{2}$ 
and $\beta _{j}(z):=\frac{1}{j}(z-\epsilon )^{2}+\epsilon .$ 
Then, for any sufficiently large $l\in \NN $, there exists an open neighborhood 
$V$ of $(\alpha _{1}^{l},\ldots ,\alpha _{n}^{l}, 
\beta _{1}^{l},\ldots ,\beta _{n}^{l})$ in  
 {\em (Poly)}$^{2n}$ such that for any 
 $(h_{1},\ldots ,h_{2n})\in V$,  
 the semigroup $G=\langle h_{1},\ldots ,h_{2n}\rangle $ satisfies that 
 $G\in {\cal G}$ and 
 $\sharp (\hat{{\cal J}}_{G})=n.$    

\end{prop}
\begin{thm}
\label{countthm}
Let $G=\langle h_{1},\ldots ,h_{m}\rangle \in {\cal G}_{dis}$ be a 
polynomial semigroup with $m\geq 3.$ Suppose 
that there exists an element $J_{0}\in \hat{{\cal J}}$ 
such that $\cup _{j=1}^{m-1}J(h_{j})\subset J_{0}$, and 
such that for each $j=1,\ldots ,m-1, $ 
we have $h_{j}^{-1}(J(h_{m}))\cap J_{0}\neq \emptyset .$ 
Then, 
we have all of the following.
\begin{enumerate}
\item $\sharp (\hat{{\cal J}})=\aleph _{0}.$ 
\item $J_{0}=J_{\min }$, or $J_{0}=J_{\max }.$ 
\item If $J_{0}=J_{\min }$, then $J_{\max }=J(h_{m})$,  
$J(G)=J_{\max}\cup \bigcup _{n\in \NN \cup \{ 0\} }
(h_{m})^{-n}(J_{\min })$, 
and for any $J\in \hat{{\cal J}}$ with $J\neq J_{\max}$, 
there exists no sequence $\{ C_{j}\} _{j\in \NN } $ 
of mutually distinct elements of $\hat{{\cal J}} $ such that 
$\min _{z\in C_{j}}d(z,J)\rightarrow 0$ as $j\rightarrow \infty .$ 

\item 
 If 
$J_{0}=J_{\max }$, then $J_{\min }=J(h_{m})$,  
$J(G)=J_{\min}\cup \bigcup _{n\in \NN \cup \{ 0\} }
(h_{m})^{-n}(J_{\max })$, and 
for any $J\in \hat{{\cal J}}$ with $J\neq J_{\min }$, 
there exists no sequence $\{ C_{j}\} _{j\in \NN }$ 
of mutually distinct elements of $\hat{{\cal J}} $ such that 
$\min _{z\in C_{j}}d(z,J)\rightarrow 0$ as $j\rightarrow \infty .$  
\end{enumerate} 
\end{thm}
\begin{prop}
\label{countprop}
There exists an open set $V$ in {\em (Poly}$_{\deg \geq 2})^{3}$ such that for any $(h_{1},h_{2},h_{3})\in V$, 
$G=\langle h_{1},h_{2},h_{3}\rangle $ satisfies that 
$G\in {\cal G}_{dis}$,\ 
$\cup _{j=1}^{2}J(h_{j})\subset J_{\min }(G)$,
 $J_{\max }(G)=J(h_{3})$,   
$h_{j}^{-1}(J(h_{3}))\cap J_{\min }(G)\neq \emptyset $ for each 
$j=1,2$, 
and $\sharp (\hat{{\cal J}}_{G})=\aleph _{0}.$  
\end{prop}
\begin{prop}
\label{countcomp}
There exists a $3$-generator polynomial 
semigroup $G=$ \\ $ \langle h_{1},h_{2},h_{3}\rangle $ 
in ${\cal G}_{dis }$ such that   
$\cup _{j=1}^{2}(h_{j})^{-1}(J_{\max }(G))\subset J_{\min }(G)$,
 $J_{\max }(G)=J(h_{3})$,  
$ \sharp (\hat{{\cal J}}_{G})=\aleph _{0} $,\ there exists 
a superattracting fixed point $z_{0}$ of some element 
of $G$ with  $z_{0}\in J(G)$,  and {\em int}$(J_{\min }(G))\neq \emptyset .$
\end{prop}
As mentioned before, these results illustrate new phenomena which can hold in the rational semigroups,  
but cannot hold in the dynamics of a single rational map or Kleinian groups. 

\begin{figure}[htbp]
\caption{The Julia set of a $3$-generator polynomial semigroup $G\in {\cal G}_{dis}$ with 
$\sharp (\hat{{\cal J}}_{G})=\aleph _{0}.$}    
\ \ \ \ \ \ \ \ \ \ \ \ \ \ \ \ \ \ \ \ \ \ \ \ \ \ \ \ \ \ \ \ 
\includegraphics[width=4.9cm,width=4.9cm]{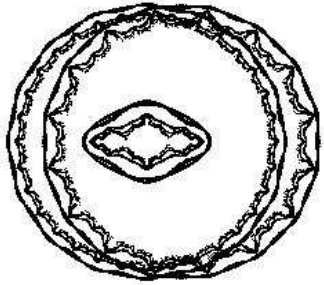}
\label{fig:3mapcountjulia2}
\end{figure}

 For the figure of the Julia set of a $3$-generator polynomial semigroup $G\in {\cal G}_{dis}$ with 
 $\sharp \hat{{\cal J}}_{G}=\aleph _{0}$, see figure~\ref{fig:3mapcountjulia2}.

\subsection{Fiberwise dynamics and Julia sets}
\label{fibsec}
We present some results on the fiberwise 
dynamics of the skew product related to a postcritically bounded 
polynomial semigroup with disconnected Julia set.
 In particular, using the uniform 
 fiberwise quasiconformal surgery on a 
 fiber bundle, 
 we show the existence 
 of a family of quasicircles parameterized 
 by a Cantor set with uniform distortion in 
 the Julia set of such a semigroup. The proofs 
 are given in Section~\ref{pffibsec}. 
 
\begin{df}[\cite{S1,S4}]
\ 
\begin{enumerate}
\item  Let $X$ be a compact metric space, 
$g:X\rightarrow X$ a continuous 
map, and $f:X\times \CCI \rightarrow 
X\times \CCI $ a continuous map. 
We say that $f$ is a rational skew 
product (or fibered rational map on 
trivial bundle $X\times \CCI $) 
over $g:X\rightarrow X$, if 
$\pi \circ f=g\circ \pi $ where 
$\pi :X\times \CCI \rightarrow X$ denotes 
the canonical projection,\ and 
if for each $x\in X$, the 
restriction 
$f_{x}:= f|_{\pi ^{-1}(\{ x\} )}:\pi ^{-1}(\{ x\}) \rightarrow 
\pi ^{-1}(\{ g(x)\} )$ of $f$ is a 
non-constant rational map,\ 
under the canonical identification 
$\pi ^{-1}(\{ x'\} )\cong \CCI $ for each 
$x'\in X.$ Let $d(x)=\deg (f_{x})$, for each 
$x\in X.$
Let $f_{x,n}$ be 
the rational map 
defined by: $f_{x,n}(y)=\pi _{\CCI }(f^{n}(x,y))$, 
for each $n\in \NN ,x\in X$ and $y\in \CCI $,   
where    
$\pi _{\CCI }:X\times \CCI 
\rightarrow  \CCI $ is the projection map.

  Moreover, if $f_{x,1}$ is a polynomial for each $x\in X$, 
 then we say that $f:X\times \CCI \rightarrow X\times \CCI $ is a 
 polynomial skew product over $g:X\rightarrow X.$  
\item  
Let $\G $ be a compact subset of Rat. 
We set  $\GN := \{ \g =(\g _{1}, \g_{2},\ldots )\mid \forall j,\g _{j}\in \G\} $ endowed with the product topology. This is a compact metric space.
Let  
$\sigma :\GN \rightarrow \GN $ be the shift map, which is defined by  
$\sigma (\g _{1},\g _{2},\ldots ):=(\g _{2},\g _{3},\ldots )
.$ Moreover,\ 
we define a map
$f:\GN \times \CCI \rightarrow 
\GN \times \CCI $ by:
$(\g ,y) \mapsto (\sigma (\g ),\g _{1}(y)),\ $
where $\g =(\g_{1},\g_{2},\ldots ).$
This is called 
{\bf the skew product associated with 
the family $\G $ of rational maps.} 
Note that $f_{\g ,n}(y)=\g _{n}\circ \cdots \circ \g _{1}(y).$ 
\end{enumerate}

\end{df}
\begin{rem}
Let $f:X\times \CCI \rightarrow X\times \CCI  $ 
be a rational skew product over 
$g:X\rightarrow X.$ Then, the function 
$x\mapsto d(x)$ is continuous in $X.$ 
\end{rem}
\begin{df}[\cite{S1, S4}]
Let $f:X\times \CCI 
\rightarrow X\times \CCI $ be a 
rational skew product over 
$g:X\rightarrow X.$ Then, we use the following notation. 
\begin{enumerate}
\item For each $x\in X$ and $n\in \NN $, we set 
$f_{x}^{n}:=
f^{n}|_{\pi ^{-1}(\{ x\} )}:\pi ^{-1}(\{ x\} )\rightarrow 
\pi ^{-1}(\{ g ^{n}(x)\} )\subset X\times \CCI .$
\item For each $x\in X$,  
we denote by $F_{x}(f)$ 
the set of points 
$y\in \CCI $ which has a neighborhood $U$ 
in $\CCI $ such that 
$\{ f_{x,n}:U\rightarrow 
\CCI \} _{n\in \NN }$
is normal. Moreover, we set 
$F^{x}(f):= \{ x\} \times F_{x}(f)\ (\subset X\times \CCI ).$  
\item For each $x\in X$, 
we set 
$J_{x}(f):=\CCI \setminus 
F_{x}(f).$ Moreover, we set 
$J^{x}(f):= \{ x\} \times J_{x}(f)$ $ (\subset X\times \CCI ).$ 
These sets $J^{x}(f)$ and $J_{x}(f)$ are called the 
fiberwise Julia sets.
\item We set 
$\tilde{J}(f):=
\overline {\bigcup _{x\in X}J^{x}
(f)}$, where the closure is taken in the product space $X\times \CCI .$
\item For each $x\in X$, we set 
$\hat{J}^{x}(f):=\pi ^{-1}(\{ x\} )\cap \tilde{J}(f).$ 
Moreover, we set $\hat{J}_{x}(f):= 
\pc (\hat{J}^{x}(f)).$ 
\item We set $\tilde{F}(f):=(X\times \CCI)\setminus 
\tilde{J}(f).$
\end{enumerate}

\end{df}
\begin{rem}
We have $\hat{J}^{x}(f)\supset J^{x}(f)$ and 
$\hat{J}_{x}(f)\supset J_{x}(f).$ 
However, 
strict containment can occur. 
For example, let $h_{1}$ be a polynomial having a Siegel disk 
with center $z_{1}\in \CC .$ 
Let $h_{2}$ be a polynomial such that 
$z_{1}$ is a repelling fixed point of $h_{2}.$ 
Let $\G =\{ h_{1},h_{2}\} .$  
Let $f:\G \times \CCI \rightarrow \G \times \CCI $ be 
the skew product associated with the family $\G .$ 
Let $x =(h_{1},h_{1},h_{1},\ldots )\in \GN .$ 
Then, $(x,z_{1})\in \hat{J}^{x}(f)\setminus  J^{x}(f)$ and 
$z_{1}\in \hat{J}_{x}(f)\setminus J_{x}(f).$ 
\end{rem}
\begin{df}
Let $f:X\times \CCI \rightarrow X\times \CCI $ be a 
polynomial skew product over $g:X\rightarrow X.$ 
Then for each $x\in X$, we set
$K_{x}(f):=
\{ y\in \CCI  \mid 
\{ f_{x,n}(y)\} _{n\in \NN }
\mbox{ is bounded } $ in $\CC \}  $, 
 and 
 $A_{x}(f):=\{ y\in \CCI 
\mid f_{x,n}(y)\rightarrow \infty 
,\ n\rightarrow \infty \} .$
Moreover, we set 
$K^{x}(f):= \{ x\} \times K_{x}(f) \ (\subset 
X\times \CCI ) $ and 
$A^{x}(f):= \{ x\} \times A_{x}(f)\ (\subset X\times \CCI ).$ 
\end{df}

\begin{df}
Let $G$ 
be a 
polynomial semigroup generated by a 
subset $\G $ of Poly$_{\deg \geq 2}.$ 
Suppose 
$G\in {\cal G}_{dis}.$ Then 
we set 
$$ \G_{\min }:=\{ h\in \G \mid 
J(h)\subset J_{\min }\} ,$$
where $J_{\min }$ denotes the 
unique minimal element in $({\cal J},\ \leq )$ 
in Theorem~\ref{mainth2}-\ref{mainth2-3}. 
Furthermore, if $\G _{\min }\neq \emptyset $, 
let $G_{\min ,\G }$ be the subsemigroup 
of $G$ that is generated by 
$\G _{\min }.$ 
\end{df}
\begin{rem}
\label{jminrem}
Let $G$ be a polynomial semigroup generated by a compact subset 
$\G $ of Poly$_{\deg \geq 2}.$ Suppose $G\in {\cal G}_{dis}.$ Then, 
 by Proposition~\ref{bminprop},  
 we have $\G _{\min }\neq \emptyset $ and 
 $\G \setminus \G _{\min }\neq \emptyset .$ 
 Moreover, $\G _{\min }$ is a compact subset of $\G .$ For, 
 if $\{ h_{n}\} _{n\in \NN }\subset \G _{\min }$ and 
 $h_{n}\rightarrow h_{\infty }$ in $\G $, 
 then for a repelling periodic point $z_{0}\in 
 J(h_{\infty })$ of $h_{\infty }$, 
we have that $d(z_{0}, J(h_{n}))\rightarrow 0$ as $n\rightarrow \infty $, 
which implies that $z_{0}\in J_{\min }$ and thus $h_{\infty }\in \G _{\min }.$ 
\end{rem}
\noindent {\bf Notation:} 
Let ${\cal F}:= \{ \varphi _{n}\} _{n\in \NN }$  
be a sequence  
of meromorphic functions 
in a domain $V.$ We say that a meromorphic function 
$\psi $ is a limit function of ${\cal F}$ 
if there exists a strictly increasing sequence $\{ n_{j}\} _{j\in \NN }$ of 
positive integers 
such that 
$\varphi _{n_{j}}\rightarrow \psi $ locally uniformly on 
$V$, as $j\rightarrow \infty .$   
\begin{df}
Let $G$ be a rational semigroup. 
\begin{enumerate}
\item 
We say that 
$G$ is hyperbolic if $P(G)\subset F(G).$ 
\item We say that $G$ is semi-hyperbolic if 
there exists a number $\delta >0$ and a 
number $N\in \NN $ such that 
for each $y\in J(G)$ and each $g\in G$, 
we have $\deg (g:V\rightarrow B(y,\delta ))\leq N$ for 
each connected component $V$ of $g^{-1}(B(y,\delta ))$, 
where $B(y,\delta )$ denotes the ball of radius $\delta $ 
with center $y$ with respect to the spherical distance, 
and $\deg (g:\cdot \rightarrow \cdot )$ denotes the 
degree of finite branched covering. 
(For background of semi-hyperbolicity, see \cite{S1} and \cite{S4}.) 
\end{enumerate}
\end{df}
The following Proposition (\ref{nonminnoncpt}-\ref{nonminnoncpt1}
and \ref{nonminnoncpt}-\ref{nonminnoncpt2}) 
means that for a polynomial semigroup 
$G\in {\cal G}_{dis}$ generated by a 
compact subset $\G $ of Poly$_{\deg \geq 2}$, 
we rarely have the situation that 
``$\G \setminus \G _{\min }$ is not compact.''
\begin{prop}
\label{nonminnoncpt}
Let $G$ be a polynomial semigroup generated by a compact 
subset $\G $ in {\em Poly }$_{\deg \geq 2}.$ 
Suppose that $G\in {\cal G}_{dis}$ and that 
$\G \setminus \G _{\min }$ is not compact. 
Then, all of the following statements \ref{nonminnoncpt1}, 
\ref{nonminnoncpt2}, \ref{nonminnoncpt3}, and 
\ref{nonminnoncpt4} hold.
\begin{enumerate}
\item \label{nonminnoncpt1}
Let $h\in \G _{\min }.$ 
Then, $J(h)=J_{\min }(G), K(h)=\hat{K}(G)$, and  
{\em int}$(K(h))$ is a non-empty connected set.
\item \label{nonminnoncpt2}
Either 

\begin{enumerate}
\item \label{nonminnoncpt2-1}
for each $h\in \G _{\min }$, $h$ is hyperbolic and $J(h)$ is a 
quasicircle; or 
\item \label{nonminnoncpt2-2}
for each $h\in \G _{\min }$, {\em int}$(K(h))$ is an immediate 
parabolic basin of a parabolic fixed point of $h.$ 
\end{enumerate} 

\item \label{nonminnoncpt3}
For each $\g \in \GN $, each limit function of 
$\{ f_{\g ,n}\} _{n\in \NN }$ in each connected component of 
$F_{\g }(f)$ is constant.
\item \label{nonminnoncpt4}
Suppose that (a) in statement \ref{nonminnoncpt2} holds. 
Then, $G_{\min ,\G }$ is hyperbolic and $G$ is semi-hyperbolic.
\end{enumerate}
\end{prop}
\begin{df}
Let $\G $ and $S$ be non-empty subsets of Poly$_{\deg \geq 2}$ 
with $S\subset \G .$ We set 
$$R(\G ,S):= 
\left\{ \g =(\g _{1},\g _{2}, \ldots )\in \GN \mid 
\sharp (\{ n\in \NN \mid \g _{n}\in S\} )=\infty \right\} .$$ 
\end{df}

\begin{df}
Let $f:X\times \CCI 
\rightarrow X\times \CCI $ be a 
rational skew product over 
$g:X\rightarrow X.$ 
We set
$$C(f):= \{ (x,y)\in X\times \CCI \mid y \mbox{ is a critical point of }
f_{x,1}\} .$$ 
Moreover, we set   
$P(f):=\overline{\cup _{n\in \NN  }f^{n}(C(f))}, $
where the closure is taken in the product space $X\times \CCI .$ 
This $P(f)$ is called the {\bf fiber-postcritical set} of 
$f.$ 

 We say that $f$ is hyperbolic 
(along fibers) if 
$P(f)\subset F(f).$

\end{df}
We present a result which describes 
the details of the fiberwise dynamics along 
$\g $ in $R(\G ,\G \setminus \G _{\min }).$  
\begin{thm}
\label{mainth3}
Let $G$ be a polynomial 
semigroup generated by a compact subset $\G $ of 
{\em Poly}$_{\deg \geq 2}.$  Suppose $G\in {\cal G}_{dis}.$ 
Let $f:\GN \times \CCI 
\rightarrow \GN \times \CCI $ be 
the  skew product associated with 
the family $\G $  
of polynomials. 
Then, 
 all of the following statements 
\ref{mainth3-0},\ref{mainth3-1}, and \ref{mainth3-3} hold. 
\begin{enumerate}
\item \label{mainth3-0}
Let $\g \in R(\G ,\G \setminus \G _{\min }).$ 
Then, each limit function of 
$\{ f_{\g ,n}\} _{n\in \NN }$ in each connected component 
of $F_{\g }(f)$ is constant. 
\item \label{mainth3-1}
Let $S$ be a non-empty compact subset of $\G \setminus \G _{\min }.$ 
Then, for each $\g \in R(\G ,S)$, we have the following.

 \begin{enumerate}
 
 \item \label{mainth3-1-1}
 There exists exactly one bounded component $U_{\g }$ 
 of $F_{\g }(f).$
  Furthermore,\ 
 $\partial U_{\g }=\partial A_{\gamma }(f)=J_{\g }(f).$ 
 
 \item \label{mainth3-1-2}
For each $y\in U_{\g }$,  there exists a number 
 $n\in \NN $ such that $f_{\g ,n}(y)
 \in $ {\em int}$(\hat{K}(G)).$
 
 \item \label{mainth3-1-3}
 $\hat{J}_{\g }(f)=J_{\g }(f).$ 
 Moreover, the map $\omega \mapsto J_{\omega }(f)$ defined on 
 $\GN $ 
 is continuous at $\g $, with respect to the Hausdorff topology 
 in the space of non-empty compact subsets of $\CCI .$ 
 \item \label{mainth3-1-4}
 The 2-dimensional Lebesgue measure of 
 $\hat{J}_{\g }(f)=J_{\g }(f)$
 is equal to zero. 
 
 \end{enumerate}

\item 
\label{mainth3-3}
Let $S$ be a non-empty compact subset of 
$\G \setminus \G _{\min }.$ For each $p \in \NN ,\ $ 
we denote by $W_{S,p}$ the set of 
elements $\g =(\g _{1},\g _{2},\ldots )\in \GN $ such that 
for each $l\in \NN $, at least one of 
$\g _{l+1},\ldots ,\g _{l+p}$ belongs to $S.$
Let $\overline{f}:=
f| _{W_{S,p}\times \CCI }:
W_{S,p}\times \CCI \rightarrow 
W_{S,p}\times \CCI .$ Then,  
$\overline{f} $ is a hyperbolic skew product 
over the shift map $\sigma :W_{S,p}\rightarrow W_{S,p}$,  
and there exists a constant 
$K_{S,p}\geq 1$ such that
for each $\g \in W_{S,p},\ $ 
$\hat{J}_{\g }(f)=J_{\g }(f)
=J_{\g }(\overline{f})$ is a $K_{S,p}$-quasicircle.  
Here, a Jordan curve $\xi $ in $\CCI $ is said to be 
a $K$-quasicircle, if $\xi $ is the image of $S^{1}(\subset \CC )$  
under a $K$-quasiconformal 
homeomorphism $\varphi :\CCI \rightarrow \CCI .$  
(For the definition of a quasicircle and a quasiconformal homeomorphism, see 
\cite{LV}.) 
\end{enumerate}
\end{thm}
We now present some results on semi-hyperbolic polynomial semigroups 
in ${\cal G}_{dis}.$ 
\begin{thm}
\label{shshprop}
Let $G$ be a polynomial semigroup generated by a non-empty compact subset 
$\G $ of {\em Poly}$_{\deg \geq 2}.$ Suppose that $G\in {\cal G}_{dis}.$ 
If $G_{\min ,\G }$ is semi-hyperbolic, then $G$ is semi-hyperbolic.
\end{thm}
\begin{thm}
\label{hhprop}
Let $G$ be a polynomial semigroup generated by a 
non-empty compact subset $\G $ of 
{\em Poly}$_{\deg \geq 2}.$ Suppose that 
$G\in {\cal G}_{dis}.$ If 
$G_{\min ,\G }$ is hyperbolic and 
$(\cup _{h\in \G \setminus \G _{\min }}CV^{\ast }(h))\cap 
J_{\min }(G)=\emptyset $, 
then $G$ is hyperbolic.
\end{thm}
\begin{rem}
\label{hhrem}
In \cite{SS}, it will be shown that 
in Theorem~\ref{hhprop}, the condition 
$(\cup _{h\in \G \setminus \G _{\min }}CV^{\ast }(h))\cap 
J_{\min }(G)=\emptyset $ is necessary. 
\end{rem}
\begin{thm}
\label{mainth3-2}
Let $G$ be a polynomial semigroup generated by a compact 
subset $\G $ of {\em Poly}$_{\deg \geq 2}.$ 
Let $f:\GNCR $ be the skew product associated with 
the family $\G .$ Suppose that $G\in {\cal G}_{dis}$ and that 
$G$ is semi-hyperbolic. Let $\g \in R(\G , \G \setminus \G _{\min })$ be 
any element. Then, $\hat{J}_{\g }(f)=J_{\g }(f)$ and 
$J_{\g }(f)$ is a Jordan curve. Moreover, 
for each point $y_{0}\in $ {\em int}$(K_{\g }(f))$, there exists an $n\in \NN $ 
such that $f_{\g ,n}(y_{0})\in $ {\em int}$(\hat{K}(G)).$

\end{thm}
We next present a result that there exist    
 families of uncountably many mutually disjoint quasicircles 
with uniform distortion, densely inside the Julia set of 
a semigroup in ${\cal G}_{dis}.$  
\begin{thm}
\label{cantorqc}
{\bf (Existence of a Cantor family of quasicircles.)} 
Let $G\in {\cal G}_{dis}$ (possibly generated by a non-compact 
family) 
and let $V$ be an open subset of $\CCI $ 
with $V\cap J(G)\neq \emptyset. $  
Then, there exist elements $g_{1}$ and 
 $g_{2}$ in $G$ such that all of the following hold.
 \begin{enumerate}
 \item 
 \label{cantorqc1}
 $H=\langle g_{1},g_{2}\rangle $ satisfies 
 that $J(H)\subset J(G).$ 
 \item 
 \label{cantorqc2}
 There exists a non-empty open set $U$ in $\CCI $ such that    
 $g_{1}^{-1}(\overline{U})\cup g_{2}^{-1}(\overline{U}) 
 \subset U$, and such that $g_{1}^{-1}(\overline{U})
 \cap g_{2}^{-1}(\overline{U})=\emptyset .$  
 
 \item 
 \label{cantorqc3}
 $H=\langle g_{1},g_{2}\rangle $ is a 
 hyperbolic polynomial semigroup.
 
 \item  
 \label{cantorqc4}
 Let 
 $f:\GN \times \CCI \rightarrow 
 \GN \times \CCI $ be the skew product  
 associated with the family  
 $\G= \{ g_{1},g_{2}\} $ of polynomials. Then, 
 we have the following. 
 \begin{enumerate}
 \item \label{cantorqc4a}
 $J(H)=\bigcup _{\g \in \GN }
 J_{\g }(f)$ (disjoint union). 
 \item \label{cantorqc4b}
 For each connected component $J$ of $J(H)$, there exists an 
 element $\g \in \GN $ such that   
 $J=J_{\g }(f).$  
  
 \item \label{cantorqc4c}
 There exists a constant $K\geq 1$
 independent of $J$ 
  such that each connected component $J$ of $J(H)$ 
  is a $K$-quasicircle. 
  \item 
 \label{cantorqc5}
 The map $\g \mapsto J_{\g }(f)$, defined for all 
 $\g \in \GN $, is continuous with respect to 
 the Hausdorff topology in the space of non-empty compact subsets of 
 $\CCI $, and 
 injective.   
 \item 
 \label{cantorqc6}
For each element $\g \in \GN ,$  
 $J_{\g }(f)\cap V\neq \emptyset .$ 
\item 
\label{cantorqc7}
Let $\omega \in \GN $ be an element such that 
$\sharp (\{ j\in \NN \mid \omega _{j}=g_{1}\} )=\infty $ 
and such that $\sharp (\{ j\in \NN \mid \omega _{j}=g_{2}\} )=\infty .$ 
Then, 
 $J_{\omega }(f)$ does not meet the boundary of any connected component of 
 $F(G).$  
 \end{enumerate} 
\end{enumerate}
\end{thm}
\subsection{Fiberwise Julia sets that are Jordan curves but not quasicircles}
\label{fjjq}
We present a result on a sufficient condition for a fiberwise Julia set 
$J_{x}(f)$ to be a Jordan curve but not a quasicircle. 
The proofs are given in Section~\ref{Proofs of fjjq}.
\begin{df}
Let $V$ be a subdomain of $\CCI $ such that 
$\partial V\subset \CC .$  
We say that $V$ is a John domain if there exists a 
constant $c>0$ and a point $z_{0}\in V$ ($z_{0}=\infty $ when 
$\infty \in V$) satisfying the following:  
for all $z_{1}\in V$ there exists an arc $\xi \subset V$ connecting 
$z_{1}$ to $z_{0}$ such that 
for any $z\in \xi $, we have 
$\min \{ |z-a|\mid a\in \partial V\} \geq c|z-z_{1}|.$
\end{df}
\begin{rem}
Let $V$ be a simply connected domain in $\CCI $ such that 
$\partial V\subset \CC .$ 
It is well-known that  
if $V$ is a John domain, then 
$\partial V$ is locally connected (\cite[page 26]{NV}). 
Moreover, a Jordan curve $\xi \subset \CC $ is a 
quasicircle if and only if both components of $\CCI \setminus \xi $ are 
John domains (\cite[Theorem 9.3]{NV}).  
\end{rem}
\begin{thm}
\label{mainthjbnq}
Let $G$ be a polynomial 
semigroup generated by a compact subset $\G $ of 
{\em Poly}$_{\deg \geq 2}.$  Suppose that $G\in {\cal G}_{dis}.$ 
Let $f:\GN \times \CCI 
\rightarrow \GN \times \CCI $ be 
the  skew product associated with 
the family $\G $  
of polynomials. 
Let $m\in \NN $ and suppose that there exists an element 
$(h_{1},h_{2},\ldots ,h_{m})\in \G ^{m}$  
such that $J(h_{m}\circ \cdots \circ h_{1})$ is not a quasicircle.
Let $\alpha =(\alpha _{1},\alpha _{2},\ldots )\in \GN $ be the  
element such that for each $k,l\in \NN \cup \{ 0\} $ with $1\leq l\leq m$, 
$\alpha _{km+l}=h_{l}.$ 
Then, the following statements \ref{mainthjbnq1} and 
\ref{mainthjbnq2} hold.
\begin{enumerate}
\item \label{mainthjbnq1}
Suppose that $G$ is hyperbolic. 
Let $\g \in R(\G ,\G \setminus \G _{\min })$ be an element 
such that there exists a sequence $\{ n_{k}\} _{k\in \NN }$ 
of positive integers satisfying that 
$\sigma ^{n_{k}}(\g )\rightarrow \alpha $ as $k\rightarrow \infty .$ 
Then, $J_{\g }(f)$ is a Jordan curve but not a quasicircle. 
Moreover, the unbounded component $A_{\g }(f)$ of 
$F_{\g }(f)$ is a John domain, but the unique bounded component 
$U_{\g }$ of $F_{\g }(f)$ is not a John domain.
\item \label{mainthjbnq2}
Suppose that $G$ is semi-hyperbolic. 
Let $\rho _{0}\in 
\G \setminus \G _{\min }$ be any element and 
let $\beta := (\rho _{0},\alpha _{1},\alpha _{2},\ldots )\in \GN .$ 
Let $\g \in R(\G ,\G \setminus \G _{\min })$ be an element 
such that there exists a sequence $\{ n_{k}\} _{k\in \NN }$ 
of positive integers satisfying that 
$\sigma ^{n_{k}}(\g )\rightarrow \beta $ as $k\rightarrow \infty .$ 
Then, $J_{\g }(f)$ is a Jordan curve but not a quasicircle. 
Moreover, the unbounded component $A_{\g }(f)$ of 
$F_{\g }(f)$ is a John domain, but the unique bounded component 
$U_{\g }$ of $F_{\g }(f)$ is not a John domain.

\end{enumerate}
\end{thm}
\begin{ex}
\label{jbnqexfirst}
Let $g_{1}(z):=z^{2}-1$ and $g_{2}(z):= \frac{z^{2}}{4}.$ 
Let $\G :=\{ g_{1}^{2}, g_{2}^{2}\} .$ 
Moreover, let $G$ be the polynomial semigroup generated by 
$\G .$ Let $D:= \{ z\in \CC \mid |z|<0.4\} .$ Then, 
it is easy to see $g_{1}^{2}(D)\cup g_{2}^{2}(D)\subset D.$ Hence, 
$D\subset F(G).$ Moreover, by Remark~\ref{pcbrem}, we have that 
$P^{\ast }(G)=
\overline{\cup _{g\in G\cup \{ Id\} }g(\{ 0,-1\} )} 
\subset D\subset F(G).$ Hence, $G\in {\cal G}$ and $G$ is hyperbolic.  
Furthermore, let $K:=\{ z\in \CC \mid 0.4\leq |z|\leq 4\} .$ Then, 
it is easy to see that $(g_{1}^{2})^{-1}(K)\cup (g_{2}^{2})^{-1}(K)\subset K$ and 
$(g_{1}^{2})^{-1}(K)\cap (g_{2}^{2})^{-1}(K)=\emptyset .$ 
Combining it with Lemma~\ref{hmslem}-\ref{backmin} and Lemma~\ref{hmslem}-\ref{bss}, 
we obtain that $J(G)$ is disconnected. Therefore, 
$G\in {\cal G}_{dis}.$ 
Moreover, it is easy to see that 
$\G _{\min }=\{ g_{1}^{2}\} .$ 
Since $J(g_{1}^{2})$ is not a Jordan curve, 
we can apply Theorem~\ref{mainthjbnq}. Setting 
$\alpha := (g_{1}^{2},g_{1}^{2},g_{1}^{2},\ldots )\in \GN $, 
it follows that 
for any $$\g \in {\cal I}:= \left\{ \omega \in R(\G ,\G \setminus \G _{\min }) 
\mid 
\exists (n_{k}) \mbox { with } \sigma ^{n_{k}}(\omega )\rightarrow 
\alpha \right\} ,$$ 
$J_{\g }(f)$ is a Jordan curve but not a quasicircle, and  
$A_{\g }(f)$ is a John domain but the bounded component 
of $F_{\g }(f)$ is not a John domain. 
(See figure~\ref{fig:dcgraph}: the Julia set of $G$. In this example,  
$\hat{{\cal J}}_{G}=\{ J_{\g }(f)\mid \g \in \GN \} $ 
and if $\g \neq \omega , J_{\g }(f)\cap J_{\omega }(f)=\emptyset .$) 
Note that by Theorem~\ref{mainth3}-\ref{mainth3-3}, if $\gamma \not\in {\cal I}$, 
then either $J_{\gamma }(f)$ is not a Jordan curve or $J_{\gamma }(f)$ is a quasicircle.
\begin{figure}[htbp]
\caption{The Julia set of $G=\langle g_{1}^{2},g_{2}^{2}\rangle .$}
\ \ \ \ \ \ \ \ \ \ \ \ \ \ \ \ \ \ \ \ \ \ \ \ \ \ \ \ \ \ \ 
\includegraphics[width=5cm,width=5cm]{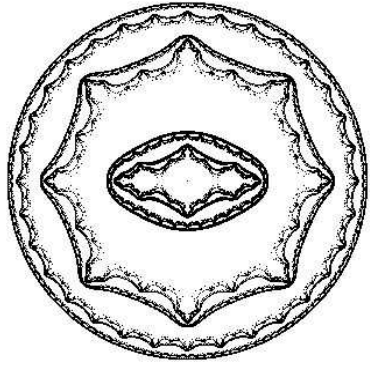}
\label{fig:dcgraph}
\end{figure}
\end{ex}

\subsection{Random dynamics of polynomials and 
 classification of compactly generated, 
 (semi-)hyperbolic, polynomial semigroups $G$ in 
 ${\cal G}$ }
\label{random}
In this section, we present some results on 
the random dynamics of polynomials. 
Moreover, we present some results on classification 
of compactly generated, (semi-) hyperbolic, polynomial semigroups 
$G$ in ${\cal G}.$ The proofs are given in Section~\ref{Proofs of random}.
 
 Let $\tau $ be a Borel probability measure on Poly$_{\deg \geq 2}.$  
We consider the i.i.d. random dynamics on $\CCI $ such that 
at every step we choose a polynomial map 
$h:\CCI \rightarrow \CCI $ according to the distribution $\tau .$ 
(Hence, this is a kind of Markov process on $\CCI .$ ) \\ 
{\bf Notation:} 
For a Borel probability measure $\tau $ on  Poly$_{\deg \geq 2}$, 
we denote by $\G_{\tau }$ the support of $\tau $ on 
Poly$_{\deg \geq 2}.$ (Hence, $\G_{\tau }$ is a closed set in 
Poly$_{\deg \geq 2}.$) 
Moreover, we denote by $\tilde{\tau } $ the infinite product measure $\otimes _{j=1}^{\infty } \tau .$  
This is a Borel probability measure on 
$\G_{\tau }^{\NN }.$ 
 Furthermore, we denote by $G_{\tau }$ the polynomial semigroup 
generated by $\G _{\tau }.$  
\begin{df}
Let $X$ be a complete metric space. 
A subset $A$ of $X$ is said to be residual if 
$X\setminus A$ is a countable union of nowhere dense subsets of 
$X.$ Note that by Baire Category Theorem, a residual set $A$ is dense in $X.$ 
\end{df}
\begin{cor}\label{rancor1}
(Corollary of Theorem~\ref{mainth3}-\ref{mainth3-1})
Let $\G $ be a non-empty compact subset of 
{\em Poly}$_{\deg \geq 2}.$ 
Let 
$f:\G ^{\NN }\times \CCI \rightarrow 
\G ^{\NN }\times \CCI $ be the skew product 
associated with the family $\G $ of polynomials.
Let $G$ be the polynomial semigroup generated by 
$\G .$ 
Suppose $G\in {\cal G}_{dis}.$ 
Then, there exists a residual subset ${\cal U}$ of $\G ^{\NN }$  
such that for each Borel probability measure 
$\tau $ on {\em Poly}$_{\deg \geq 2}$ with $\G _{\tau }=\G $, 
we have  
 $\tilde{\tau }({\cal U})=1$, and such that 
each $\g \in {\cal U}$ satisfies   
all of the following.
\begin{enumerate}
 \item \label{rancor1-1}
 There exists exactly one bounded component $U_{\g }$ 
 of $F_{\g }(f).$
  Furthermore,
 $\partial U_{\g }=\partial A_{\gamma }(f)=J_{\g }(f).$ 
 
 \item \label{rancor1-2}
 Each limit function of 
 $\{ f_{\g ,n} \} _{n}$ in 
 $U_{\g }$ is constant. Moreover, 
 for each $y\in U_{\g }$,  there exists a number 
 $n\in \NN $ such that $f_{\g ,n}(y)
 \in $ {\em int}$(\hat{K}(G)).$
 
 \item \label{rancor1-3}
 $\hat{J}_{\g }(f)=J_{\g }(f).$ Moreover, 
 the map $\omega \mapsto J_{\omega }(f)$ defined on 
 $\GN $ is continuous at 
 $\g $, with respect to the Hausdorff topology in the space of 
 non-empty compact subsets of $\CCI .$ 
 \item \label{rancor1-4}
The 2-dimensional Lebesgue measure of 
 $\hat{J}_{\g }(f)=J_{\g }(f)$
 is equal to zero.

\end{enumerate}

\end{cor}
Next we present a result on compactly generated, semi-hyperbolic, 
polynomial semigroups in ${\cal G}.$ 
\begin{thm}
\label{mainthran1}
Let $\G $ be a non-empty compact subset of 
{\em Poly}$_{\deg \geq 2}.$ 
Let 
$f:\G ^{\NN }\times \CCI \rightarrow 
\G ^{\NN }\times \CCI $ be the skew product 
associated with the family $\G $ of polynomials.
Let $G$ be the polynomial semigroup generated by $\G .$ 
Suppose that $G \in {\cal G}$ and 
that $G$ is semi-hyperbolic. 
Then, exactly one of the following three statements 
\ref{mainthran1-1}, \ref{mainthran1-2}, and 
\ref{mainthran1-3} holds.
\begin{enumerate}
\item 
\label{mainthran1-1}
$G$ is hyperbolic. Moreover, there exists a constant 
$K\geq 1$ such that for each $\g \in \G ^{\NN }$, 
$J_{\g }(f)$ is a $K$-quasicircle.
\item \label{mainthran1-2}
There exists a residual subset ${\cal U}$ of 
$\G^{\NN }$ such that 
for each Borel probability measure 
$\tau $ on {\em Poly}$_{\deg \geq 2}$ with $\G _{\tau }=\G $, we have 
$\tilde{\tau }({\cal U})=1$, and such that 
for each $\g \in {\cal U}$, 
$J_{\g }(f)$ is a Jordan curve but not a quasicircle,  
$A_{\g }(f)$ is a John domain, 
and the bounded component of $F_{\g }(f)$ is not a 
John domain.
Moreover, there exists a dense subset ${\cal V}$ of 
$\G ^{\NN }$ such that 
for each $\g \in {\cal V}$, $J_{\g }(f)$ is not a 
Jordan curve.  
Furthermore, there exist two elements $\alpha ,\beta \in \G ^{\NN }$ such that 
$J_{\beta }(f)<J_{\alpha }(f).$  
\item \label{mainthran1-3}
There exists a dense subset ${\cal V}$ of $\G ^{\NN }$ such that 
for each $\g \in {\cal V}$, $J_{\g }(f)$ is not a Jordan curve. 
Moreover, for each  $\alpha ,\beta \in \G ^{\NN }$, 
$J_{\alpha }(f)\cap J_{\beta }(f)
\neq \emptyset .$ Furthermore, 
$J(G)$ is arcwise connected. 
\end{enumerate}
\end{thm} 
\begin{cor}
\label{rancor2}
Let 
$\G $ be a non-empty compact subset of 
{\em Poly}$_{\deg \geq 2}.$  
Let 
$f:\G ^{\NN }\times \CCI \rightarrow 
\G ^{\NN }\times \CCI $ be the skew product 
 associated with the family $\G $ of polynomials. 
Let $G$ be the polynomial semigroup generated by 
$\G .$ 
Suppose that $G \in {\cal G}_{dis}$ and that $G$ is 
semi-hyperbolic. Then, either statement \ref{mainthran1-1} or 
statement \ref{mainthran1-2} in Theorem~\ref{mainthran1} holds. 
In particular, for any Borel Probability measure 
$\tau $ on {\em Poly}$_{\deg \geq 2}$ with $\G _{\tau }=\G $, 
for almost every $\g \in \G _{\tau }^{\NN }$ with respect to 
$\tilde{\tau }$, $J_{\g }(f)$ is a Jordan curve.  
\end{cor} 
We now classify compactly generated, hyperbolic, polynomial semigroups 
in ${\cal G}.$ 
\begin{thm}
\label{mainthran2}
Let $\G $ be a non-empty compact subset of {\em Poly}$_{\deg \geq 2}.$ 
Let $f:\GN \times \CCI \rightarrow \GN \times \CCI $ be the 
skew product associated with the family $\G .$ 
Let $G$ be the polynomial semigroup generated by $\G .$ 
Suppose that $G\in {\cal G}$ and that $G$ is hyperbolic. 
Then, exactly one of the following three statements \ref{mainthran2-1}, 
\ref{mainthran2-2}, \ref{mainthran2-3}
 holds. 
\begin{enumerate}
\item 
\label{mainthran2-1}
There exists a constant $K\geq 1$ such that for each 
$\g \in \GN $, $J_{\g }(f)$ is a $K$-quasicircle.
\item 
\label{mainthran2-2}
There exists a residual subset ${\cal U}$ of $\GN $ such that 
for each Borel probability measure $\tau $ on {\em Poly}$_{\deg \geq 2}$ 
with $\G _{\tau }=\G $, we have $\tilde{\tau }({\cal U})=1$, 
and such that for each $\g \in {\cal U}$, $J_{\g }(f)$ is a 
Jordan curve but not a quasicircle, 
$A_{\g }(f)$ is a John domain, 
and the bounded component 
of $F_{\g }(f)$ is not a John domain.   
Moreover, there exists a dense subset ${\cal V}$ of 
$\GN $
 such that for each $\g \in {\cal V}$, 
$J_{\g }(f)$ is a quasicircle.
Furthermore, there exists a dense subset 
${\cal W}$ of $\GN $ such that for 
each $\g \in {\cal W}$, there are infinitely many bounded connected 
components of $F_{\g }(f).$  
\item 
\label{mainthran2-3}
For each $\g \in \GN $, there are infinitely many bounded 
connected components 
of $F_{\g }(f)$. Moreover, 
for each $\alpha ,\beta \in \GN $, 
$J_{\alpha }(f)\cap J_{\beta }(f)\neq 
\emptyset .$ Furthermore, 
$J(G)$ is arcwise connected.  
\end{enumerate}

\end{thm}
\begin{ex}
\label{jbnqex}
Let $h_{1}(z):= z^{2}-1$ and $h_{2}(z):=az^{2}$, 
where $a\in \CC $ with $0<|a|<0.1.$ 
Let $\G := \{ h_{1},h_{2}\} .$ Moreover, 
let $G:=\langle h_{1},h_{2}\rangle .$ 
Let $U:= \{ |z|<0.2\} .$ Then, it is easy to see that 
$h_{2}(U)\subset U,\ h_{2}(h_{1}(U))\subset U,$ and 
$h_{1}^{2}(U)\subset U.$ Hence, $U\subset 
F(G).$ It follows that $P^{\ast }(G)\subset $ 
int$(\hat{K}(G))\subset F(G).$ Therefore, $G\in {\cal G}$ and 
$G$ is hyperbolic. Since $J(h_{1})$ is not a Jordan curve and 
$J(h_{2})$ is a Jordan curve, Theorem~\ref{mainthran2} implies that 
there exists a residual subset ${\cal U} $ of $ \GN $ such that 
for each Borel probability measure $\tau $ 
on Poly$_{\deg \geq 2}$ with $\G _{\tau }=\G $, we have 
$\tilde{\tau }({\cal U})=1$, and such that for each 
$\g \in {\cal U}$, 
$J_{\g }(f)$ is a Jordan curve but not a quasicircle. Moreover, 
for each $\g \in {\cal U}$, $A_{\g }(f)$ is a John domain, but 
the bounded component of $F_{\g }(f)$ is not a John domain. 
Furthermore, by Theorem~\ref{polyandrathm3}, $J(G)$ is connected.
\end{ex} 
\begin{rem}
Let $h\in $ Poly$_{\deg \geq 2}$ be a polynomial. 
Suppose that $J(h)$ is a Jordan curve but not a quasicircle. 
Then, it is easy to see that there exists a parabolic fixed point 
of $h$ in $\CC $ and the bounded connected component of $F(h)$ is the 
immediate parabolic basin. Hence, $\langle h\rangle $ is not semi-hyperbolic.
Moreover, by \cite{CJY}, 
$F_{\infty }(h)$ is not a John domain. 

 Thus what we see in statement \ref{mainthran1-2} in Theorem~\ref{mainthran1} and statement \ref{mainthran2-2} in Theorem~\ref{mainthran2}, 
 as illustrated in Example ~\ref{jbnqexfirst} and Example~\ref{jbnqex}, is a special and new 
 phenomenon which can hold in the {\em random} dynamics 
 of a family of polynomials, but cannot hold in the usual 
 iteration dynamics of a single polynomial. 
 Namely, it can hold that for almost every $\g \in \GN $, 
$J_{\g }(f)$ is a Jordan curve and fails to be a quasicircle 
all while the basin of infinity $A_{\g }(f)$ is still a John domain. 
Whereas, if $J(h)$, for some polynomial $h$, is a Jordan curve which 
fails to be a quasicircle, then the basin of infinity $F_{\infty }(h)$ 
is necessarily {\bf not} a John domain.   

 Pilgrim and Tan Lei (\cite{PT}) showed that there exists a  
hyperbolic rational map $h$ with {\em disconnected} Julia set such that  
``almost every'' connected component of $J(h)$ 
is a Jordan curve but not a quasicircle.  
\end{rem}
We give a sufficient condition so that statement~\ref{mainthran2-1} in Theorem~\ref{mainthran2} holds. 
\begin{prop}
\label{ranprop1}
Let $\G $ be a non-empty compact subset of {\em Poly}$_{\deg \geq 2}.$ 
Let $f:\GN \times \CCI \rightarrow \GN \times \CCI $ be the 
skew product associated with the family $\G .$ 
Let $G$ be the polynomial semigroup generated by $\G .$ 
Suppose that $P^{\ast }(G)$ is included in a connected component 
of {\em int}$(\hat{K}(G)).$ Then, there exists a constant $K\geq 1$ such that 
for each $\g \in \GN $, $J_{\g }(f)$ is a $K$-quasicircle.
\end{prop}
\subsection{Construction of examples}
\label{Const}
We present a way 
to construct examples of semigroups $G$ in ${\cal G}_{dis}.$ 
\begin{prop}
\label{Constprop}
Let $G$ be a 
polynomial semigroup generated by 
a compact subset $\G $ of {\em Poly}$_{\deg \geq 2}.$ 
Suppose that $G\in {\cal G}$ and  
{\em int}$(\hat{K}(G))\neq \emptyset .$ 
Let $b\in $ {\em int}$(\hat{K}(G)).$ 
Moreover, let $d\in \NN $ be any positive integer such that 
$d\geq 2$, and such that 
$(d, \deg (h))\neq (2,2)$ for each $h\in \G .$ 
Then, there exists a number $c>0$ such that 
for each $a\in \CC $ with $0<|a|<c$, 
there exists a compact neighborhood $V$ of 
$g_{a}(z)=a(z-b)^{d}+b$ in {\em Poly}$_{\deg \geq 2}$ satisfying 
that for any non-empty subset $V'$ of $V$,  
the polynomial semigroup 
$H_{\G, V'} $ generated by the family $\G \cup V'$ 
belongs to ${\cal G}_{dis}$, $\hat{K}(H_{\G, V'})=\hat{K}(G)$ 
 and $(\G \cup V')_{\min }\subset \G .$ 
Moreover, in addition to the assumption above, 
if $G$ is semi-hyperbolic (resp. hyperbolic), 
then the above $H_{\G, V'}$ is semi-hyperbolic (resp. hyperbolic).   
\end{prop}
\begin{rem}
By Proposition~\ref{Constprop}, 
there exists a $2$-generator polynomial semigroup 
$G=\langle h_{1},h_{2}\rangle $ in ${\cal G}_{dis}$ such that 
$h_{1}$ has a Siegel disk.  
\end{rem}
\begin{df}
Let $d\in \NN $ with $d\geq 2.$ 
We set 
${\cal Y}_{d}:=\{ h\in \mbox{Poly} \mid \deg (h)=d\} $ 
endowed with the relative topology from Poly.
\end{df}
\begin{thm}
\label{shshfinprop}
Let $m\geq 2$ and let $d_{2},\ldots ,d_{m}\in \NN $ be such that  
$d_{j}\geq 2$ for each $j=2,\ldots ,m.$ Let 
$h_{1}\in {\cal Y}_{d_{1}}$ with {\em int}$(K(h_{1}))\neq \emptyset $ be  
such that $\langle h_{1}\rangle \in {\cal G}.$ 
Let $b_{2},b_{3},\ldots ,b_{m}\in $ {\em int}$(K(h_{1})).$ 
Then, all of the following statements hold. 
\begin{enumerate}
\item \label{shshfinprop1}
Suppose that $\langle h_{1}\rangle $ is 
semi-hyperbolic (resp. hyperbolic). 
Then, there exists a number $c>0$ such that 
for each $(a_{2},a_{3},\ldots ,a_{m})\in \CC ^{m-1}$ with 
$0<|a_{j}|<c$ ($j=2,\ldots ,m$), 
setting $h_{j}(z)=a_{j}(z-b_{j})^{d_{j}}+b_{j}$ ($j=2,\ldots ,m$), 
the polynomial semigroup 
$G=\langle h_{1},\ldots ,h_{m}\rangle $ satisfies that 
$G\in {\cal G}$, $\hat{K}(G)=K(h_{1})$ and $G$ is semi-hyperbolic (resp. hyperbolic). 
\item \label{shshfinprop2}
Suppose that $\langle h_{1}\rangle $ is 
semi-hyperbolic (resp. hyperbolic). Suppose also that 
either (i) there exists a $j\geq 2$ with $d_{j}\geq 3$, or 
(ii) $\deg(h_{1})=3$, $b_{2}=\cdots =b_{m}.$ Then, there exist 
$a_{2},a_{3},\ldots ,a_{m}>0$ such that setting 
$h_{j}(z)=a_{j}(z-b_{j})^{d_{j}}+b_{j}$ ($j=2,\ldots ,m$), 
the polynomial semigroup $G=\langle h_{1},h_{2},\ldots ,h_{m}\rangle $ 
satisfies that $G\in {\cal G}_{dis}$, $\hat{K}(G)=K(h_{1})$ and 
$G$ is semi-hyperbolic (resp. hyperbolic).  

\end{enumerate} 
\end{thm}
\begin{df}
Let $m\in \NN .$ We set 
\begin{itemize}
\item ${\cal H}_{m}:= 
\{ (h_{1},\ldots ,h_{m})\in 
(\mbox{Poly}_{\deg \geq 2 })^{m}\mid 
\langle h_{1},\ldots ,h_{m}\rangle  \mbox{ is hyperbolic}\}$, 
\item 
${\cal B}_{m}:= 
 \{ (h_{1},\ldots ,h_{m})\in 
(\mbox{Poly}_{\deg \geq 2 })^{m}\mid 
\langle h_{1},\ldots ,h_{m}\rangle \in {\cal G} \}$, and 
\item 
${\cal D}_{m}:= 
\{ (h_{1},\ldots ,h_{m})\in 
(\mbox{Poly}_{\deg \geq 2 })^{m}\mid 
J(\langle h_{1},\ldots ,h_{m}\rangle ) \mbox{ is disconnected}\} $.  
\end{itemize}
Moreover, let $\pi _{1}: (\mbox{Poly}_{\deg \geq 2 })^{m}\rightarrow 
\mbox{Poly}_{\deg \geq 2 }$ be the projection defined by 
$\pi (h_{1},\ldots ,h_{m})=h_{1}.$ 
\end{df}
\begin{thm}
\label{sphypopen}
Under the above notation, all of the following statements hold.
\begin{enumerate}
\item \label{sphypopen1}
${\cal H}_{m}, {\cal H}_{m}\cap {\cal B}_{m},$ 
${\cal H}_{m}\cap {\cal D}_{m}$, 
and ${\cal H}_{m}\cap {\cal B}_{m}\cap {\cal D}_{m}$ 
are open in $(\mbox{{\em Poly}}_{\deg \geq 2})^{m}.$ 
\item \label{sphypopen2}
Let $d_{1},\ldots ,d_{m}\in \NN $ be such that  
$d_{j}\geq 2$ for each $j=1,\ldots ,m.$\\  
Then, 
$\pi _{1}: {\cal H}_{m}\cap {\cal B}_{m}\cap 
({\cal Y}_{d_{1}}\times \cdots \times 
{\cal Y}_{d_{m}})\rightarrow 
{\cal H}_{1}\cap {\cal B}_{1}\cap {\cal Y}_{d_{1}}$ is 
surjective. 
\item \label{sphypopen3} 
Let $d_{1},\ldots ,d_{m}\in \NN $ be such that  
$d_{j}\geq 2$ for each $j=1,\ldots ,m$ and such that  
$(d_{1},\ldots ,d_{m})\neq (2,2,\ldots ,2).$ Then, 
 $\pi _{1}: {\cal H}_{m}\cap {\cal B}_{m}\cap {\cal D}_{m}\cap 
({\cal Y}_{d_{1}}\times \cdots \times 
{\cal Y}_{d_{m}})\rightarrow 
{\cal H}_{1}\cap {\cal B}_{1}\cap {\cal Y}_{d_{1}}$ is 
surjective. 
\end{enumerate}
\end{thm}
\begin{rem}
Combining Proposition~\ref{Constprop}, 
Theorem~\ref{shshfinprop}, and Theorem~\ref{sphypopen}, 
we can construct many examples of semigroups $G$ 
in ${\cal G} $ (or ${\cal G}_{dis}$) with some 
additional properties (semi-hyperbolicity, hyperbolicity, etc.). 
\end{rem}
\section{Tools}
\label{Tools}
To show the main results, we need 
some tools in this section.
\subsection{Fundamental properties of rational semigroups}
{\bf Notation:} 
For a rational semigroup $G$, we set 
$E(G):=\{ z\in \CCI \mid \sharp (\cup _{g\in G}g^{-1}(\{ z\} ))<\infty \} .$ 
This is called the exceptional set of $G.$ 
\begin{lem}[\cite{HM1,GR,S1}]
\label{hmslem}
Let $G$ be a rational semigroup.\
\begin{enumerate}
\item 
\label{invariant}
For each $h\in G,\ $ we have 
$h(F(G))\subset F(G)$ and $h^{-1}(J(G))\subset J(G).$ Note that we do not 
have that the equality holds in general.
\item
\label{bss}
If $G=\langle h_{1},\ldots ,h_{m}\rangle $, then 
$J(G)=h_{1}^{-1}(J(G))\cup \cdots \cup h_{m}^{-1}(J(G)).$ 
More generally, if $G$ is generated by a compact subset 
$\G $ of {\em Rat}, then 
$J(G)=\cup _{h\in \G}h^{-1}(J(G)).$ 
(We call this property of the Julia set of a compactly generated rational semigroup ``backward self-similarity." )
\item
\label{Jperfect}
If \( \sharp (J(G)) \geq 3\) ,\ then \( J(G) \) is a 
perfect set.\ 
\item
\label{egset}
If $\sharp (J(G))\geq 3$ ,\ then 
$ \sharp (E(G)) \leq 2. $
\item
\label{o-set}
If a point \( z\) is not in \( E(G),\ \) then 
 \( J(G)\subset \overline{\cup _{g\in G}g^{-1}(\{ z\} )} .\) In particular 
if a point\ \( z\)  belongs to \ \( J(G)\setminus E(G), \) \ then
$ \overline{\cup _{g\in G}g^{-1}(\{ z\})}=J(G). $
\item
\label{backmin}
If \( \sharp (J(G))\geq 3 \) ,\ 
then
\( J(G) \) is the smallest closed backward invariant set containing at least
three points. Here we say that a set $A$ is backward invariant under $G$ if
for each \( g\in G,\ g^{-1}(A)\subset A.\ \)
\end{enumerate}
\end{lem}
\begin{thm}[\cite{HM1,GR,S1}]
\label{repdense}
Let $G$ be a rational semigroup.
If $\sharp (J(G))\geq 3$, then 
$ J(G)=\overline{\{ z\in \CCI \mid 
\exists g \in G,\ g(z)=z,\ |g'(z)|>1 \} } .$ 
In particular,\ $J(G)=\overline{\bigcup _{g\in G}J(g)}.$
\end{thm}
\begin{rem}
If a rational semigroup $G$ contains an element $g$ with $\deg (g)\geq 2$, 
then $\sharp (J(g))\geq 3$, which implies that $\sharp (J(G))\geq 3.$  
\end{rem}
\subsection{Fundamental properties of fibered rational maps}
\begin{lem}
\label{fibfundlem}
Let $f:X\times \CCI \rightarrow 
X\times \CCI $ be a rational skew product 
over $g:X\rightarrow X.$ Then, 
we have the following.
\begin{enumerate}
\item \label{fibfundlem1}
{\bf (\cite[Lemma 2.4]{S1})}
For each $x\in X$, 
$(f_{x,1})^{-1}(J_{g (x)}(f))
=J_{x}(f).$ Furthermore, we have 
$\hat{J}_{x}(f)\supset 
J_{x}(f).$ Note that 
{\bf equality $\hat{J_{x}}(f)=J_{x}(f)$ does not 
hold in general.}

 If $g:X\rightarrow X$ is 
a surjective and open map, then
$f^{-1}(\tilde{J}(f))=\tilde{J}(f)=f(\tilde{J}(f))$, and 
 for each $x\in X$, 
$(f_{x,1})^{-1}(\hat{J}_{g (x)}(f))=
\hat{J}_{x}(f).$
\item \label{fibfundlem2}
({\bf \cite{J, S1}})
If $d(x)\geq 2$ 
for each $x\in X $,  
then 
for each $x\in X$, $J_{x}(f)$ is a non-empty perfect set with $\sharp (J_{x}(f))\geq 3.$
Furthermore, the map $x\mapsto J_{x}(f)$ is 
lower semicontinuous; i.e.,  
for any point $(x,y)\in X\times \CCI $ 
with $y\in J_{x}(f)$ and 
any sequence $\{ x^{n}\} _{n\in \NN }$ in $X$ 
with $x^{n}\rightarrow x,\ $  
there exists a sequence $\{ y^{n}\} _{n\in \NN }$ 
in $\CCI $ with 
$y^{n}\in J_{x^{n}}(f)$ for each $n\in \NN $ such that 
$y^{n}\rightarrow y.$ However,  
$x\mapsto J_{x}(f)$ is {\bf NOT}
 continuous with respect to the Hausdorff topology in 
general.
\item \label{fibfundlemast}
If $d(x)\geq 2$ for each $x\in X$, 
then $\inf _{x\in X}$diam$_{S}J_{x}(f)>0$, 
where diam$_{S}$ denotes the diameter with respect to the 
spherical distance.
\item \label{fibfundlem4}
If $f:X\times \CCI \rightarrow 
X\times \CCI $ is a polynomial skew product 
and $d(x)\geq 2$     
for each $x\in X $,  
then  we have that there exists a ball $B$ 
around $\infty $ such that for each $x\in 
X$, $B\subset A_{x}(f)\subset F_{x}(f)$, and that 
for each $x\in X$, 
$J_{x}(f)=\partial K_{x}(f)=\partial A_{x}(f) .$
Moreover, for each $x\in X$, 
$A_{x}(f)$ is connected.  
\item \label{fibfundlema}
If $f:X\times \CCI \rightarrow 
X\times \CCI $ is a polynomial skew product 
and $d(x)\geq 2$     
for each $x\in X $, and if 
$\omega \in X$ is a point such that 
{\em int}$(K_{w}(f))$ is a non-empty 
set, then 
$\overline{\mbox{{\em int}}(K_{\omega }(f))}=K_{\omega }(f)$ 
and $\partial (\mbox{{\em int}}(K_{\omega }(f)))=J_{\omega }(f).$ 
\end{enumerate}
\end{lem}
\begin{proof}
For the proof of statement \ref{fibfundlem1}, 
see \cite[Lemma 2.4]{S1}.
For the proof of statement \ref{fibfundlem2}, 
see \cite{J} and \cite{S1}. 

 By statement \ref{fibfundlem2}, 
it is easy to see that statement \ref{fibfundlemast} holds. 
Moreover, 
it is easy to see that 
statement \ref{fibfundlem4} holds. 

 To show statement \ref{fibfundlema}, 
let $y\in J_{\omega }(f)$ be a point. 
Let $V$ be an arbitrary 
neighborhood of $y$ in $\CCI .$ 
Then, by the self-similarity of Julia sets (see \cite{Bu1}), 
there exists an $n\in \NN $ such that 
$f_{\omega ,n}(V\cap J_{\omega }(f))=J_{g^{n}(\omega )}(f).$ 
Since $\partial (\mbox{int}(K_{g^{n}(\omega )}(f)))\subset 
J_{g^{n}(\omega )}(f)$ and $(f_{\omega ,n})^{-1}
(K_{g^{n}(\omega )}(f))=K_{\omega }(f)$, 
it follows that 
$V\cap \partial (\mbox{int}(K_{\omega }(f)))\neq \emptyset .$ 
Hence, we obtain $J_{\omega }(f)=\partial 
(\mbox{int}(K_{\omega }(f)).$ Therefore, we have proved 
statement \ref{fibfundlema}.  
\end{proof}
\begin{lem}
\label{fiblem}
Let $f:\GN \times 
\CCI \rightarrow \GN \times 
\CCI $ be a skew product  
associated with a compact subset 
$\G $ of {\em Rat}. 
Let $G$ be the rational semigroup generated by $\G .$ 
Suppose that $\sharp (J(G))\geq 3.$ 
Then, we have the following. 
\begin{enumerate}
\item  \label{pic}
$\pi _{\CCI }(\tilde{J}(f))=J(G).$
\item \label{fibfundlem5}
For each $\g =(\g _{1},\g _{2},\ldots ,)\in \GN $, 
$\hat{J}_{\g }(f)=\cap _{j=1}^{\infty }\g_{1}^{-1}
\cdots \g_{j}^{-1}(J(G)).$
\end{enumerate}
\end{lem}
\begin{proof}
First, we show statement \ref{pic}. 
Since $J_{\g }(f)\subset J(G)$ for each $\g \in \G $, 
we have $\pi _{\CCI }(\tilde{J}(f))\subset J(G).$ 
By Theorem~\ref{repdense}, we have 
$J(G)=\overline{\cup _{g\in G}J(g)}.$ 
Since $\cup _{g\in G}J(g)\subset \pi _{\CCI }(\tilde{J}(f))$, 
we obtain $J(G)\subset \pi _{\CCI }(\tilde{J}(f)).$ 
Therefore, we obtain $\pi _{\CCI }(\tilde{J}(f))=J(G).$ 

 We now show statement \ref{fibfundlem5}. 
Let $\g =(\g _{1},\g _{2},\ldots )\in \GN .$ 
By statement \ref{fibfundlem1} in Lemma~\ref{fibfundlem}, 
we see that for each $j\in \NN $, 
$\g_{j}\cdots \g_{1}(
\hat{J}_{\g }(f))=
\hat{J}_{\sigma ^{j}(\g )}(f)
\subset J(G).$ 
Hence, 
$\hat{J}_{\g }(f)\subset 
\cap _{j=1}^{\infty }
\g _{1}^{-1}\cdots 
\g _{j}^{-1}(J(G)).$
 Suppose that there exists a point 
 $(\g ,y)\in \GN \times \CCI $ such that  
 $y\in $ $\left( \cap _{j=1}^{\infty }
\g _{1}^{-1}\cdots 
\g _{j}^{-1}(J(G))\right) \setminus 
\hat{J}_{\g }(f).$ 
Then, we have 
$(\g ,y)\in (\GN \times \CCI ) 
\setminus \tilde{J}(f).$ 
Hence, there exists a 
neighborhood $U$
of $\g $ in $\GN $ and a 
neighborhood $V$ of $y$ in $\CCI $ such that 
$U\times V\subset \tilde{F}(f).$ 
Then, there exists an $n\in \NN $ such that 
$\sigma ^{n}(U)=\GN .$ 
Combining it with Lemma~\ref{fibfundlem}-\ref{fibfundlem1}, 
we obtain $\tilde{F}(f)\supset 
f^{n}(U\times V)\supset \GN \times 
\{ f_{\g ,n}(y)\} .$ 
Moreover, 
since we have 
$f_{\g ,n}(y)\in J(G)=
\pi _{\CCI }(\tilde{J}(f))$, 
where the last equality holds 
by statement \ref{pic}, 
we get that there exists an element 
$\g '\in \GN $ such that 
$(\g ', f_{\g ,n}(y))\in 
\tilde{J}(f).$ However, it contradicts 
$(\g ',f_{\g ,n}(y))\in 
\GN \times \{ f_{\g ,n}(y)\} 
\subset \tilde{F}(f).$ 
Hence, we obtain 
$\hat{J}_{\g }(f)= 
\cap _{j=1}^{\infty }
\g_{1}^{-1}\cdots 
\g_{j}^{-1}(J(G)).$
\end{proof}

\begin{lem}
\label{fibconnlem}
Let $f:X\times \CCI \rightarrow 
X\times \CCI $ be a polynomial skew product 
over $g:X\rightarrow X$ such that  
for each $x\in X$, 
$d(x)\geq 2.$ Then, 
the following are equivalent.
\begin{enumerate}
\item \label{fcc1}
$\pi _{\CCI }(P(f))\setminus \{ \infty \} $ 
is bounded in $\CC .$ 
\item \label{fcc2}
For each $x\in X$, $J_{x}(f)$ is 
connected.
\item \label{fcc3}
For each $x\in X$, $\hat{J}_{x}(f)$ is 
connected.
\end{enumerate} 
\end{lem}
\begin{proof} 
First, we show \ref{fcc1} $\Rightarrow $\ref{fcc2}. 
Suppose that \ref{fcc1} holds.
Let 
$R>0$ be a number such that 
for each $x\in X$, 
$B:=\{ y\in \CCI 
\mid |y|>R\} \subset A_{x}(f)$ 
and $\overline{f_{x,1}(B)}\subset B.$ 
Then, for each $x\in X$, we have 
$A_{x}(f)=\cup _{n\in \NN }(f_{x,n})^{-1}(B)$ 
and $(f_{x,n})^{-1}(B)\subset 
(f_{x,n+1})^{-1}(B)$, 
for each $n\in \NN .$ Furthermore, since 
we assume \ref{fcc1}, we see that 
for each $n\in \NN $, $(f_{x,n})^{-1}(B)$ 
is a simply connected domain, by the Riemann-Hurwitz formula. 
Hence, for each $x\in X$, 
$A_{x}(f)$ is a simply connected domain.   
Since $\partial A_{x}(f)=J_{x}(f)$ for each $x\in X,$ 
we conclude that for each $x\in X$, 
$J_{x}(f)$ is connected. 
Hence, we have shown \ref{fcc1} $\Rightarrow $ \ref{fcc2}. 

 Next, we show 
 \ref{fcc2} $\Rightarrow $ \ref{fcc3}. 
 Suppose that \ref{fcc2} holds. 
 Let $z_{1}\in \hat{J}_{x}(f)$ and $z_{2}\in J_{x}(f)$ be two points. 
Let 
 $\{ x^{n}\} _{n\in \NN }$ be a sequence 
 such that $x^{n}\rightarrow x$ as $n\rightarrow \infty $, and 
 such that 
 $d(z_{1},J_{x^{n}}(f))\rightarrow 0$ as 
 $n\rightarrow \infty .$ We may assume that 
 there exists a non-empty compact set $K$ in 
 $\CCI $ such that 
 $J_{x^{n}}(f)\rightarrow K$ as $n\rightarrow \infty $, 
 with respect to the Hausdorff topology in the 
 space of non-empty compact sets in $\CCI .$ 
 Since we assume \ref{fcc2}, 
 $K$ is connected. 
 By Lemma~\ref{fibfundlem}-\ref{fibfundlem2}, 
 we have $d(z_{2}, J_{x^{n}}(f))\rightarrow 0$ as 
 $n\rightarrow \infty .$ Hence, 
 $z_{i}\in K$ for each $i=1,2.$ 
 Therefore, $z_{1}$ and $z_{2}$ belong to the same connected 
 component of $\hat{J}_{x}(f).$ Thus, 
 we have shown \ref{fcc2} $\Rightarrow $ \ref{fcc3}. 

 Next, we show \ref{fcc3} $\Rightarrow $ \ref{fcc1}. 
 Suppose that \ref{fcc3} holds.
 It is easy to see that 
 $A_{x}(f)\cap \hat{J}_{x}(f)=\emptyset $ for each 
 $x\in X.$ 
 Hence, $A_{x}(f)$ is a connected component of 
 $\CCI \setminus \hat{J}_{x}(f).$ 
 Since we assume \ref{fcc3}, we have that 
 for each $x\in X$, 
 $A_{x}(f)$ is a simply connected domain. 
 Since $(f_{x,1})^{-1}(A_{g(x)}(f))=
 A_{x}(f)$, the Riemann-Hurwitz formula implies that 
 for each $x\in X$, 
 there exists no critical point of 
 $f_{x,1}$ in $A_{x}(f)\cap \CC .$ 
 Therefore, we obtain \ref{fcc1}. 
 Thus, we have shown \ref{fcc3} $\Rightarrow $ \ref{fcc1}.
\end{proof}   
\begin{cor}
\label{fibconncor}
Let $G=\langle h_{1},h_{2}\rangle \in {\cal G}.$ Then, 
$h_{1}^{-1}(J(h_{2}))$ is connected.
\end{cor} 
\begin{proof}
Let $f:\GN \times \CCI \rightarrow \GN \times \CCI $ 
be the skew product associated with the family 
$\G =\{ h_{1},h_{2}\} .$ 
Let $\g =(h_{1},h_{2},h_{2},h_{2},h_{2},\ldots )\in \GN .$ 
Then, by Lemma~\ref{fibfundlem}-\ref{fibfundlem1},
 we have $J_{\g }(f)=h_{1}^{-1}(J(h_{2})).$ 
From Lemma~\ref{fibconnlem}, it follows that 
$h_{1}^{-1}(J(h_{2}))$ is connected.
\end{proof} 
 
\begin{lem}
\label{fibconnlem2}
Let $G$ be a polynomial semigroup
generated by a compact subset $\G $ of 
{\em Poly}$_{\deg \geq 2}.$ 
Let $f:\GN \times \CCI 
\rightarrow \GN \times \CCI $ be 
the skew  product associated with the 
family $\G .$
Suppose that  
$G\in {\cal G}.$ Then 
for each $\g =(\g _{1},\g _{2},\ldots ,)\in 
\GN $, 
the sets $J_{\g }(f),\ \hat{J} _{\g }(f)$, and 
$\cap _{j=1}^{\infty }
 \g _{1}^{-1}\cdots \g _{j}^{-1}(J(G))$ 
 are connected.
\end{lem}
\begin{proof}
From Lemma~\ref{fiblem}-\ref{fibfundlem5} and 
Lemma~\ref{fibconnlem}, the lemma follows.
\end{proof}
\begin{lem}
\label{fiborder}
Under the same assumption as that in Lemma~\ref{fibconnlem2}, 
let $\g ,\rho \in \GN $ be two elements with 
$J_{\g }(f)\cap J_{\rho }(f)=\emptyset .$ 
Then, either $J_{\g }(f)< J_{\rho }(f)$ 
or $J_{\rho }(f)< J_{\g }(f).$
\end{lem}
\begin{proof}
Let $\g ,\rho \in \GN $ with 
$J_{\g }(f)\cap J_{\rho }(f)=\emptyset .$ 
Suppose that the statement 
``either $J_{\g }(f)<J_{\rho }(f)$ or $J_{\rho }(f)<J_{\g }(f)$" 
is not true. 
Then, Lemma~\ref{fibconnlem} implies that 
$J_{\g }(f)$ is included in the unbounded 
component of $\CC \setminus J_{\rho }(f)$, and that 
$J_{\rho }(f)$ is included in the unbounded 
component of $\CC \setminus J_{\g}(f).$ 
From Lemma~\ref{fibfundlem}-\ref{fibfundlem4}, it follows that 
$K_{\rho }(f)$ is included in the unbounded 
component $A_{\g }(f)\setminus \{ \infty \} $ 
of $\CC \setminus J_{\g }(f).$ 
However, it causes a contradiction, since 
$\emptyset \neq P^{\ast }(G)\subset \hat{K}(G)\subset 
K_{\rho }(f)\cap K_{\g }(f).$  
\end{proof}
\begin{df}
Let 
$f:\GNCR $ be a polynomial skew product 
over $g:X\rightarrow X.$ 
Let $p\in \CC $ and 
$\epsilon >0.$ 
We set \\ 
${\cal F}_{f,p,\epsilon }:= 
\{\alpha :D(p,\epsilon )\rightarrow \CC \mid 
\alpha \mbox{ is a well-defined inverse branch of }
(f_{x,n})^{-1}, x\in X, n\in \NN \} .$  
\end{df}
\begin{lem}
\label{invnormal}
Let $f:\GNCR $ be a polynomial skew product over 
$g:X\rightarrow X$ such that for each $x\in X$, 
$d(x)\geq 2.$ 
Let $R>0,\epsilon >0$, and \\ 
${\cal F}:=
\{ \alpha \circ \beta :D(0,1)\rightarrow \CC 
\mid 
\beta :D(0,1)\cong D(p,\epsilon ),\ 
\alpha :D(p,\epsilon )\rightarrow \CC ,\ 
\alpha \in {\cal F}_{f,p,\epsilon },\ p\in D(0,R)\} .$ 
Then, ${\cal F} $ is normal on $D(0,1).$ 
\end{lem}
\begin{proof}
Since $d(x)\geq 2$ for each $x\in X$, 
there exists a ball $B$ around $\infty $ 
with $B\subset \CCI \setminus D(0,R+\epsilon )$ such that 
for each $x\in X$, $f_{x,1}(B)\subset B.$ 
Let $p\in D(0,R).$ Then, for each 
$\alpha \in {\cal F}_{f,p,\epsilon }$, 
$\alpha (D(p,\epsilon ))\subset \CCI \setminus B.$ 
Hence, ${\cal F} $ is normal in $D(0,1).$ 
\end{proof}
\subsection{A lemma from general topology}
\begin{lem}[\cite{N}]
\label{nadlem}
Let $X$ be a compact metric space 
and let $f:X\rightarrow X$ be a continuous 
open map. 
Let $A$ be a compact connected subset of 
$X.$ Then for each connected component 
$B$ of $f^{-1}(A)$, we have 
$f(B)=A.$ 

\end{lem}
\section{Proofs of the main results}
\label{Proofs} 
In this section, we demonstrate the main results.
\subsection{Proofs of results in \ref{concompsec}}
\label{pfconcompsec}
{\bf Proof of Theorem~\ref{mainth0}:}
First, we show the following:\\ 
Claim: For any $\lambda \in \Lambda $, 
$h_{\lambda }^{-1}(A)\subset A.$ 

 To show  the claim, 
let $\lambda \in \Lambda $ with $J(h_{\lambda })\neq \emptyset $
 and 
let $B$ be a connected component of 
$h_{\lambda }^{-1}(A).$ Then 
by Lemma~\ref{nadlem}, 
$h_{\lambda }(B)=A.$
Combining this with 
$h_{\lambda }^{-1}(J(h_{\lambda }))=
J(h_{\lambda })$, we obtain 
$B\cap J(h_{\lambda })\neq \emptyset .$ 
Hence $B\subset A.$ 
This means that 
$h_{\lambda }^{-1}(A)\subset A$ for each 
$\lambda \in \Lambda $ with $J(h_{\lambda })\neq \emptyset .$ 
Next, let $\lambda \in \Lambda $ with 
$J(h_{\lambda })=\emptyset .$ 
Then $h_{\lambda }$ is either identity or an 
elliptic M\"{o}bius transformation. 
By hypothesis and Lemma~\ref{hmslem}-\ref{invariant},  
 we obtain $h_{\lambda }^{-1}(A)\subset A.$ 
Hence, we have shown the claim.

 Combining the above claim with 
$\sharp A\geq 3$, 
by 
Lemma~\ref{hmslem}-\ref{backmin} 
we obtain $J(G)\subset A.$ Hence $J(G)=A$ and 
$J(G)$ is connected.
\qed  

\

\noindent {\bf Notation:}
We denote by $d$ the spherical distance on $\CCI .$ 
Given $A\subset \CCI $ and $z\in \CCI $, 
we set $d(z,A):=\inf \{ d(z,w)\mid w\in A\} .$ 
Given $A\subset \CCI $ and $\epsilon >0$,\ we set 
$B(A,\epsilon ):= \{ a\in \CCI \mid d(a,A)<\epsilon \} .$
Furthermore, given $A\subset \CC $, $z\in \CC $, 
and $\epsilon >0$, 
we set 
$d_{e}(z,A):=\inf \{ |z-w| \mid w\in A\} $ and 
$D(A,\epsilon ):= \{ a\in \CC \mid d_{e}(a,A)<\epsilon \} .$  

\ 

 We need the following lemmas to prove the main results.
\begin{lem}
\label{appjlem}
Let $G\in {\cal G}$ and let   
$J$ be a connected component of 
$J(G)$, $z_{0}\in J$ a point, and 
$\{ g_{n}\} _{n\in \NN }$ a sequence in  
$G$ such that 
$d(z_{0},J(g_{n}))\rightarrow 
0$ as $n\rightarrow \infty .$ 
Then 
$\sup\limits _{z\in J(g_{n})}d(z,J)\rightarrow 
0$ as $n\rightarrow \infty .$ 
\end{lem}
\begin{proof}
Suppose 
 there exists a connected component 
 $J'$ of $J(G)$ with $J'\neq J$ 
 and a subsequence $\{ g_{n_{j}}\} _{j\in \NN }$ 
 of $\{ g_{n}\} _{n\in \NN }$ such that 
 $\min\limits _{z\in J(g_{n_{j}})}
 d(z,J')\rightarrow 0$ as 
 $j\rightarrow \infty .$ 
 Since $J(g_{n_{j}})$ is compact and 
 connected for each $j$, 
 we may assume, passing to a subsequence, that there exists a non-empty 
 compact connected subset $K$ of $\CCI $ such that 
 $J(g_{n_{j}})\rightarrow K$ as 
 $j\rightarrow \infty $, with respect to 
 the Hausdorff topology. 
 Then $K\cap J \neq \emptyset $ and 
 $K\cap J'\neq \emptyset .$ Since 
 $K\subset J(G)$ and $K$ is connected, 
 it contradicts $J'\neq J.$
\end{proof}
\begin{lem}
\label{appj2lem}
Let $G\in {\cal G}.$ Then given $J\in {\cal J}$ and 
$\epsilon >0$,\ there exists an element $g\in G$ such that 
$J(g)\subset B(J,\epsilon ).$  

\end{lem}
\begin{proof}
We take a point $z\in J.$ Then, by Theorem~\ref{repdense}, 
there exists a sequence $\{ g_{n}\} _{n\in \NN }$ in $G$ such that 
$d(z,J(g_{n}))\rightarrow 0$ as $n\rightarrow \infty .$ 
By Lemma~\ref{appjlem}, we conclude that there exists 
an $n\in \NN $ such that 
$J(g_{n})\subset B(J,\epsilon ).$  
\end{proof}
\begin{lem}
\label{inftyj1}
Let $G$ be a polynomial semigroup. 
Suppose that $J(G)$ is disconnected, and  
 $\infty \in J(G).$ Then,\ 
 the connected component $A$ of 
 $J(G)$ containing $\infty $ is equal to 
 $\{ \infty \} .$ 
\end{lem}
\begin{proof}
 By Lemma~\ref{nadlem}, we obtain  
 $g^{-1}(A)\subset A$ for each 
 $g\in G.$ Hence, if $\sharp A\geq 3$, then 
 $J(G)\subset A$,\ by Lemma~\ref{hmslem}-\ref{backmin}. 
 Then $J(G)=A$ and it causes a contradiction, since 
 $J(G)$ is disconnected.  
\end{proof}

We now demonstrate Theorem~\ref{mainth1}.\\ 
\noindent {\bf Proof of Theorem~\ref{mainth1}:}
First, we show statement \ref{mainth1-1}.
Suppose the statement is false. Then, 
there exist elements $J_{1},J_{2}\in {\cal J}$ 
such that $J_{2}$ is included in the 
unbounded component $A_{1}$  
of $\CC \setminus J_{1}$, and such that 
$J_{1}$ is included in the unbounded 
component $A_{2}$  of 
$\CC \setminus J_{2}.$ 
Then we can find an $\epsilon >0$ such that 
$\overline{B(J_{2},\epsilon )}$ is 
included in the unbounded component 
of $\CC \setminus \overline{B(J_{1},\epsilon )}$, 
and such that 
$\overline{B(J_{1},\epsilon )}$ is 
included in the unbounded component of 
$\CC \setminus \overline{B(J_{2},\epsilon )}.$  
By Lemma~\ref{appj2lem}, for each $i=1,2$, 
 there exists an element $g_{i}\in G$ such that  
 $J(g_{i})\subset B(J_{i},\epsilon ).$ 
 This implies that 
 $J(g_{1})\subset A_{2}'$ and 
 $J(g_{2})\subset A_{1}'$, where 
 $A_{i}'$ denotes the unbounded 
 component of $\CC \setminus J(g_{i}).$  
 Hence we obtain 
 $K(g_{2})\subset A_{1}'.$ 
 Let  $v$ be a critical value 
 of $g_{2}$ in $\CC .$ Since 
 $P^{\ast }(G) $ is 
 bounded in $\CC $, we have 
 $v\in K(g_{2}).$  It implies  
 $v\in A_{1}'.$ 
  Hence $g_{1}^{l}(v)\rightarrow \infty $ as $l\rightarrow \infty .$ 
  However, this implies a 
 contradiction since 
 $P^{\ast }(G) $ is 
 bounded in $\CC .$ 
 Hence we have shown statement \ref{mainth1-1}.

 Next, we show statement \ref{mainth1-2}. 
 Let $F_{1}$ be a connected component of 
 $F(G).$ 
 Suppose that there exist three connected 
 components $J_{1},J_{2}$ and $J_{3}$ of 
 $J(G)$ such that they are mutually disjoint and 
 such that 
 $\partial F_{1}\cap J_{i}\neq \emptyset $ 
 for each $i=1,2,3.$ 
 Then, by statement \ref{mainth1-1}, 
 we may assume that we have 
 either (1): $J_{i}\in {\cal J}$ for each 
  $i=1,2,3$ and $J_{1}<J_{2}<J_{3}$, or 
  (2): $J_{1},J_{2}\in {\cal J},\ J_{1}<J_{2}$, and 
  $\infty \in J_{3}.$ 
  Each of these cases implies that 
  $J_{1}$ is included in a bounded component 
  of $\CC \setminus J_{2}$ and 
  $J_{3}$ is included in the unbounded component of 
  $\CCI \setminus J_{2}.$ However, it 
  causes a contradiction, 
  since $\partial F_{1}\cap J_{i}\neq \emptyset $ 
 for each $i=1,2,3.$ Hence, we have shown that 
 we have either \\  
 Case I: $\sharp \{ J:\mbox{component of }J(G)\mid 
 \partial F_{1}\cap J\neq \emptyset \} =1$ or \\ 
 Case II: $\sharp \{ J:\mbox{component of }J(G)\mid 
 \partial F_{1}\cap J\neq \emptyset \} =2.$

  Suppose that we have Case I. Let 
 $J_{1}$ be the connected component of $J(G)$ such that 
 $\partial F_{1}\subset J_{1}.$ 
 Let $D_{1}$ be the connected component of 
 $\CCI \setminus J_{1}$ containing $F_{1}.$ 
 Since $\partial F_{1}\subset J_{1},$
  we have $\partial F_{1}\cap D_{1}=\emptyset .$ 
  Hence, we have $F_{1}=D_{1}.$ Therefore, 
  $F_{1}$ is simply connected. 

  Suppose that we have Case II. 
 Let $J_{1}$ and $J_{2}$ be the two connected components 
 of $J(G)$ such that 
 $J_{1}\neq J_{2}$ and 
 $\partial F_{1}\subset J_{1}\cup J_{2}.$ 
 Let $D$ be the connected component of 
 $\CCI \setminus ( J_{1}\cup J_{2})$ 
 containing $F_{1}.$ 
 Since $\partial F_{1}\subset J_{1}\cup J_{2},$
 we have $\partial F_{1}\cap D=\emptyset .$ Hence, 
 we have $F_{1}=D.$ Therefore, $F_{1}$ is 
 doubly connected. Thus, we have shown 
statement \ref{mainth1-2}.

 We now show statement \ref{mainth1-3}.
 Let $g\in G$ be an element and 
$J$ a connected component 
of $J(G).$ 
 Suppose that $g^{-1}(J)$ is disconnected. 
 Then, by Lemma~\ref{nadlem}, 
 there exist at most finitely 
 many connected components $C_{1},\ldots ,C_{r}$
 of $g^{-1}(J)$ with $r\geq 2.$  
 Then there exists a positive number 
 $\epsilon $ such that 
 denoting by $B_{j}$ the connected 
 component of $g^{-1}(B(J,\epsilon ))$ 
 containing $C_{j}$ for each 
 $j=1,\ldots ,r$, 
 $\{ B_{j}\} $ are mutually 
 disjoint. By 
 Lemma~\ref{nadlem}, 
 we see that, for each 
 connected component $B$ of 
 $g^{-1}(B(J,\epsilon )), $
  $g(B)=B(J,\epsilon )$ 
 and $B\cap C_{j}\neq \emptyset $ for 
 some $j.$ Hence we get that 
 $g^{-1}(B(J,\epsilon ))=
 \bigcup _{j=1}^{r}B_{j}$ 
  (disjoint union) and 
  $g(B_{j})=B(J,\epsilon )$ for each 
  $j.$  
By Lemma~\ref{appj2lem}, there exists an element 
$h\in G$ such that $J(h)\subset B(J,\epsilon ).$  
Then it follows that 
$g^{-1}(J(h))\cap B_{j}\neq \emptyset 
$ for each $j=1,\ldots ,r.$ Moreover, 
we have $g^{-1}(J(h))\subset 
g^{-1}(B(J,\epsilon ))=\bigcup _{j=1}^{r}B_{j}.$
On the other hand, by 
Corollary~\ref{fibconncor}, 
we have that $g^{-1}(J(h))$ is connected. 
This is a contradiction. 
Hence, we have shown that, for each $g\in G$ and 
each connected component $J$ of $J(G)$,\   
$g^{-1}(J)$ is connected. 

 By Lemma~\ref{inftyj1}, we get that 
if $J\in {\cal J}$, then $g^{\ast }(J)\in {\cal J}.$ 
 Let $J_{1}$ and $J_{2}$ be two elements 
 of ${\cal J}$ such that 
 $J_{1}\leq J_{2}.$ 
 Let $U_{i}$ be the unbounded component 
 of $\CC \setminus J_{i}$, for each 
 $i=1,2.$ 
 Then 
 $U_{2}\subset U_{1}.$ 
  Let $g\in G$ be an element. 
 Then $g^{-1}(U_{2})\subset g^{-1}(U_{1}).$ 
 Since $g^{-1}(U_{i})$ is the unbounded 
 connected component of $\CC \setminus 
 g^{-1}(J_{i})$ for each $i=1,2$, 
 it follows that 
 $g^{-1}(J_{1})\leq g^{-1}(J_{2}).$ 
 Hence $g^{\ast }(J_{1})\leq g^{\ast }(J_{2})$, 
 otherwise 
$g^{\ast }(J_{2})<g^{\ast }(J_{1})$, 
and 
 it contradicts $g^{-1}(J_{1})\leq g^{-1}(J_{2}).$    
\qed \\ 
\subsection{Proofs of results in \ref{Upper}}
\label{Proof of Upper}
In this section, we prove results in Section \ref{Upper}, 
except Theorem~\ref{polyandrathm1}-\ref{polyandrathm1-2} and 
Theorem~\ref{polyandrathm1}-\ref{polyandrathm1-3}, which will be 
proved in Section~\ref{Proof of Properties}.

 To demonstrate Theorem~\ref{polyandrathm1}, we need the 
following lemmas.
\begin{lem}
\label{j1aj2lem}
Let $G$ be a polynomial semigroup in ${\cal G}_{dis}.$ 
Let $J_{1},J_{2}\in \hat{{\cal J}}$ be two elements with 
$J_{1}\neq J_{2}.$ 
Then, we have the following.
\begin{enumerate}
\item \label{j1aj2lem1}
If $J_{1},J_{2}\in {\cal J}$ and 
$J_{1}<J_{2}$, then 
there exists a doubly connected component $A$ of $F(G)$ such that 
$J_{1}<A<J_{2}.$ 
\item \label{j1aj2lem2}
If $\infty \in J_{2}$, then 
there exists a doubly connected component $A$ of $F(G)$ such that 
$J_{1}<A.$ 
\end{enumerate}
\end{lem}
\begin{proof}
First, we show statement \ref{j1aj2lem1}. 
Suppose that $J_{1},J_{2}\in {\cal J}$ and $J_{1}<J_{2}.$ 
We set $B=\cup _{J\in {\cal J}, J_{1}\leq J\leq J_{2}}J.$ 
Then, $B$ is a closed disconnected set. 
Hence, there exists a multiply connected component 
$A'$ of $\CCI \setminus B.$ Since $A'$ is multiply connected, 
we have that 
$A'$ is included in the unbounded component of $\CCI \setminus J_{1}$,  
and that 
$A'$ is included in a bounded component of $
\CCI \setminus J_{2}.$ 
This 
implies that $A'\cap J(G)=\emptyset .$ 
Let $A$ be the connected component of $F(G)$ such that  
$A'\subset A.$ 
Since $B\subset J(G)$, we have 
$F(G)\subset \CCI \setminus B.$  
Hence, $A$ must be equal to $A'.$ 
Since $A'$ is multiply connected, 
Theorem~\ref{mainth1}-\ref{mainth1-2} implies 
that $A=A'$ is doubly connected. 
Let $J$ be the connected component $J(G)$ such that 
$J<A$ and $J\cap \partial A\neq \emptyset .$  Then, since 
$A'=A$ is included in the unbounded component of 
$\CCI \setminus J_{1},$ 
we have that $J$ does not meet any bounded component 
of $\CC \setminus J_{1}.$ Hence, we obtain $J_{1}\leq J$, 
which implies that $J_{1}\leq J<A.$ 
Therefore, 
$A$ is a doubly connected component of $F(G)$ such that 
$J_{1}<A<J_{2}.$ Hence, we have shown statement \ref{j1aj2lem1}.

 Next, we show statement \ref{j1aj2lem2}. 
 Suppose that $\infty \in J_{2}.$ 
We set $B=(\cup _{J\in {\cal J}, J_{1}\leq J}J )\cup J_{2}.$ 
Then, $B$ is a disconnected closed set. 
Hence, there exists a multiply connected component 
$A'$ of $\CCI \setminus B.$ By the same method as that of proof 
of statement \ref{j1aj2lem1}, we see that 
$A'$ is equal to a doubly connected component $A$ of 
$F(G)$ such that $J_{1}<A.$ 
 Hence, we have shown statement \ref{j1aj2lem2}. 
\end{proof}
\begin{lem}
\label{cptmsetlem}
Let $H_{0}$ be a real affine semigroup 
generated by a compact set $C$ in {\em RA}. 
Suppose that each element $h\in C$ is of the form 
$h(x)=b_{1}(h)x+b_{2}(h)$, where 
$b_{1}(h), b_{2}(h)\in \RR $, 
$|b_{1}(h)|>1.$   
Then, for any subsemigroup $H$ of $H_{0}$, we have 
$M(H)=J(\eta (H))\subset \RR .$ 

\end{lem}
\begin{proof}
From the assumption, there exists a number $R>0$ such that 
for each $h\in C$, 
$\eta (h)(B(\infty ,R))\subset B(\infty ,R).$  
Hence, we have $B(\infty ,R)\subset F(\eta (H))$, which 
implies that $J(\eta (H))$ is a bounded subset of $\CC .$  
We consider the following cases:\\ 
Case 1: $\sharp (J(\eta (H)))\geq 3.$ \\ 
Case 2: $\sharp (J(\eta (H)))\leq 2.$ 

 Suppose that we have case 1. Then, 
from Theorem~\ref{repdense}, it follows that 
$M(H)=J(\eta (H))\subset \RR .$ 

 Suppose that we have case 2. 
Let $b(h)$ be the unique fixed point of $h\in H$ 
in $\RR .$ From the hypothesis, 
we have that for each $h\in H$, 
$b(h)\in J(\eta (H)).$ Since we assume $\sharp (J(\eta (H)))\leq 2$, 
Lemma~\ref{hmslem}-\ref{invariant} implies that 
 there exists a point $b\in \RR $ such that 
for each $h\in H$, we have $b(h)=b.$ Then 
any element $h\in H$ is of the form 
$h(x)=c_{1}(h)(x-b)+c_{2}(h),$ where $c_{1}(h), c_{2}(h)\in \RR , 
|c_{1}(h)|>1.$ Hence, $M(H)=\{ b\} \subset  J(\eta (H)).$ 
Suppose that there exists a point $c$ in $J(\eta (H))\setminus 
\{ b\} .$ Since $J(\eta (H))$ is a bounded set of $\CC $, 
and since we have $h^{-1}(J(\eta (H)))\subset J(\eta (H))$ for each $h\in H$ 
(Lemma~\ref{hmslem}-\ref{invariant}), 
we get that $h^{-1}(c)\in J(\eta (H))\setminus (\{ b\} \cup \{ c\} ) $, for each element $h\in H.$ 
This implies that $\sharp (J(\eta (H)))\geq 3$, which is a contradiction.
Hence, we must have that 
$J(\eta (H))=\{ b\} =M(H).$     
\end{proof}
We need the notion of Green's functions, in order to demonstrate 
Theorem~\ref{polyandrathm1}.
\begin{df}
Let $D$ be a domain in $\CCI $ with 
$\infty \in D.$ 
We denote by $\varphi (D,z)$ Green's function on 
$D$ with pole at $\infty .$ 
By definition, this is the unique function on $D\cap \CC $ with the 
properties:
\begin{enumerate}
\item $\varphi (D,z)$ is harmonic and positive in $D\cap \CC $;
\item $\varphi (D,z)-\log |z|$ is bounded in a neighborhood of 
$\infty $; and 
\item $\varphi (D,z)\rightarrow 0$ as $z\rightarrow \partial D.$ 
\end{enumerate} 
\end{df}
\begin{rem}\ 
\begin{enumerate}
\item The limit $\lim\limits _{z\rightarrow \infty }(\varphi (D,z)-\log |z|)$ 
exists and this is called {\em Robin's constant} of $D.$ 
\item If $D$ is a simply connected domain with $\infty \in D$, 
then we have $\varphi (D,z)=-\log |\psi (z)|$, where 
$\psi :D\rightarrow \{ z\in \CC \mid |z|<1\} $ denotes 
a biholomorphic map with $\psi (\infty )=0.$ 
\item 
It is well-known that for any $g\in $ Poly$_{\deg \geq 2},$  
\begin{equation}
\label{polygreeneq}
\varphi (F_{\infty }(g),z)=\log |z|+
\frac{1}{\deg (g)-1}\log |a(g)|+o(1)\ \mbox{ as }z\rightarrow \infty .
\end{equation}
(See \cite[p147]{Ste}.) Note that the point 
$-\frac{1}{\deg (g)-1}\log |a(g)|\in \RR $ is the unique fixed point 
of $\Psi (g)$ in $\RR .$  
\end{enumerate}
\end{rem}
\begin{lem}
\label{greenk1k2}
Let $K_{1}$ and $K_{2}$ be two non-empty connected compact sets in 
$\CC $ such that $K_{1}<K_{2}.$ 
Let $A_{i}$ denote the unbounded component of 
$\CCI \setminus K_{i}$, for each $i=1,2.$ 
Then, we have 
$\lim _{z\rightarrow \infty }(\log |z|-\varphi (A_{1},z)) 
<\lim _{z\rightarrow \infty }(\log |z|-\varphi (A_{2},z)).$
\end{lem}
\begin{proof}
The function 
$\phi (z):= 
\varphi (A_{2},z)-\varphi (A_{1},z)=
(\log |z|-\varphi (A_{1},z))-(\log |z|-\varphi (A_{2},z))$ 
is harmonic on $A_{2}\cap \CC .$ This $\phi $ is bounded 
around $\infty .$ Hence $\phi $ extends to a 
harmonic function on $A_{2}.$ 
Moreover, since $K_{1}<K_{2}$, 
we have $\limsup _{z\rightarrow \partial A_{2}}\phi (z)$ $<0.$ 
From the maximum principle, 
it follows that 
$\phi (\infty )<0.$ Therefore, the statement of our lemma holds. 
\end{proof}
In order to demonstrate Theorem~\ref{polyandrathm1}-\ref{polyandrathm1-1}, 
we will prove the following lemma.  
 (Theorem~\ref{polyandrathm1}-\ref{polyandrathm1-2} and 
 Theorem~\ref{polyandrathm1}-\ref{polyandrathm1-3} will be proved 
 in Section~\ref{Proof of Properties}.) 
\begin{lem}
\label{part1-1lem}
Let $G$ be a polynomial semigroup in ${\cal G}.$ Then, there exists an 
injective map $\tilde{\Psi }:
\hat{{\cal J}}_{G}\rightarrow 
{\cal M}_{\Psi (G)}$ such that: 
\begin{enumerate}
\item  
if $J_{1},J_{2}\in {\cal J}_{G}$ and 
$J_{1}<J_{2}$, 
then $\tilde{\Psi }(J_{1})<_{r}\tilde{\Psi }(J_{2})$;  
\item if $J\in \hat{{\cal J}}_{G}$ and $\infty \in J$, then 
$+\infty \in \tilde{\Psi }(J)$; and 
\item if $J\in {\cal J}_{G}$, then 
$\tilde{\Psi }(J)\subset \hat{\RR }\setminus \{ +\infty \} .$   
\end{enumerate}
\end{lem}  
\begin{proof}
We first show the following claim.\\  
\noindent Claim 1: 
In addition to the assumption of Lemma~\ref{part1-1lem}, 
if we have $\infty \in F(G)$, then 
$M(\Psi (G))\subset \hat{\RR } \setminus \{ +\infty \} .$ 
 
 To show this claim, let $R>0$ be a number such that 
 $J(G)\subset D(0,R).$ Then, 
 for any $g\in G$, we have 
 $K(g)<\partial D(0,R).$  
By Lemma~\ref{greenk1k2}, we get that 
there exists a constant $C>0$ such that for each 
$g\in G$, $\frac{-1}{\deg (g)-1}\log |a(g)|\leq C.$ 
Hence, it follows that $M(\Psi (G))\subset 
[-\infty ,C].$ Therefore, we have shown Claim 1.

We now prove the statement of the lemma in the case $G\in {\cal G}_{con} .$ 
If $\infty \in F(G)$, then claim 1 implies that 
$ M({\Phi (G)}) \subset \hat{\RR }\setminus \{ +\infty \} $ and the statement of the lemma holds. 
We now suppose $\infty \in J(G).$  
We put $L_{g}:= \max _{z\in J(g) } |z|$ for each $g\in G.$  
Moreover, for each non-empty compact subset $E$ of $\CC $, we denote by Cap $(E)$ the logarithmic capacity of $E.$ 
We remark that Cap$(E)=\exp (\lim _{z\rightarrow \infty }(\log |z|-\varphi (D_{E},z)))$, where $D_{E}$ denotes the 
connected component of $\CCI \setminus E$ containing $\infty .$   
We may assume that $0\in P^{\ast }(G).$ Then, by \cite{A}, we have 
Cap$(J(g))\geq $ Cap $([0,L_{g}])\geq L_{g}/4$ for each $g\in G.$  
Combining this with $\infty \in J(G)$ and Theorem~\ref{repdense}, 
we obtain $+\infty \in M_{\Psi (G)}$ and the statement of the lemma holds. 
  
We now prove the statement of the lemma in the case $G\in {\cal G}_{dis}.$  
Let $\{ J_{\lambda }\} _{\lambda \in \Lambda }$
be the set $\hat{{\cal J}}_{G}$ of all connected components 
of $J(G).$ 
By Lemma~\ref{appj2lem}, 
for each $\lambda \in \Lambda $ and each $n\in \NN $, 
there exists an element $g_{\lambda ,n}\in G$ such that 
\begin{equation}
\label{polyandrathm1pf0a}
J(g_{\lambda ,n})\subset B(J_{\lambda },\frac{1}{n}).
\end{equation} 
We have that the fixed point of 
$\Psi (g_{\lambda ,n})$ in $\RR $ is 
equal to $\frac{-1}{\deg (g_{\lambda ,n})-1}\log |a(g_{\lambda ,n})|.$ 
We may assume that 
$\frac{-1}{\deg (g_{\lambda ,n})-1}\log |a(g_{\lambda ,n})|\rightarrow 
\alpha _{\lambda }$ as $n\rightarrow \infty $, 
where $\alpha _{\lambda }$ is an element of $\hat{\RR }.$ 
For each $\lambda \in \Lambda $, let $B_{\lambda }\in 
{\cal M}_{\Psi (G)}$ be the element with 
$\alpha _{\lambda }\in B_{\lambda }.$ 
We will show the following claim.\\ 
Claim 2:  
If $\lambda ,\xi $ are two elements in $\Lambda $ with 
$\lambda \neq \xi $, then 
$B_{\lambda }\neq B_{\xi }.$ 
Moreover, if $J_{\lambda }, J_{\xi }\in {\cal J}_{G}$ and 
$J_{\lambda }<J_{\xi }$, then $B_{\lambda }<_{r}B_{\xi }. $
Furthermore, if $J_{\xi }\in \hat{{\cal J}}_{G}$ with 
$\infty \in J_{\xi }$, then $+\infty \in B_{\xi }.$ 
 
 To show this claim,  let $\lambda $ and $\xi $ be two elements in 
 $\Lambda $ with $\lambda \neq \xi .$ 
 We have the following two cases:\\ 
Case 1: $J_{\lambda },J_{\xi }\in {\cal J}_{G}$ and 
$J_{\lambda }<J_{\xi }.$ \\ 
Case 2:  $J_{\lambda }\in {\cal J}_{G}$ and $\infty \in J_{\xi } .$  

 Suppose that we have case 1. 
 By Lemma~\ref{j1aj2lem}, 
 there exists a doubly connected component $A$ of $F(G)$ such that 
\begin{equation}
 \label{polyandrathm1pf0b}
 J_{\lambda }<A<J_{\xi }.
\end{equation} 
 Let $\zeta _{1}$ and $\zeta _{2}$ be two Jordan curves in 
 $A$ such that they are not null-homotopic in $A$, and such that 
 $\zeta _{1}<\zeta _{2}.$   
For each $i=1,2,$ let $A_{i}$ be the unbounded component 
of 
$\CCI \setminus \zeta _{i}.$ 
Moreover, we set 
$\beta _{i}:=\lim _{z\rightarrow \infty }
(\log |z|-\varphi (A_{i},z))$, for each $i=1,2.$ 
By Lemma~\ref{greenk1k2}, we have $\beta _{1}<\beta _{2}.$ 
Let $g\in G$ be any element. 
By (\ref{polyandrathm1pf0a}) and (\ref{polyandrathm1pf0b}), 
there exists an $m\in \NN $ such that 
$J(g_{\lambda ,m})<\zeta _{1}.$ 
Since $P^{\ast }(G)\subset K(g_{\lambda ,m})$, 
it follows that $P^{\ast }(G)$ is included in the bounded component 
of $\CC \setminus \zeta _{1}.$ 
Hence, 
we see that  
\begin{equation}
\label{polyandrathm1pf0}
\mbox{either }J(g)<\zeta _{1}, \mbox{ or }
\zeta _{2}<J(g). 
\end{equation}
From Lemma~\ref{greenk1k2}, it follows that 
either 
$\frac{-1}{\deg (g)-1}\log |a(g)|<\beta _{1}$, or 
$\beta _{2}<\frac{-1}{\deg (g)-1}\log |a(g)|.$ 
This implies that 
\begin{equation}
\label{polyandrathm1pf1}
M(\Psi (G))\subset \hat{\RR }\setminus (\beta _{1},\beta _{2}),
\end{equation}
where $(\beta _{1},\beta _{2}):=\{ x\in \RR \mid \beta _{1}<x<\beta _{2}\}.$ 
Moreover, combining (\ref{polyandrathm1pf0a}), (\ref{polyandrathm1pf0b}), 
and (\ref{polyandrathm1pf0}), we get that  
there exists a number $n_{0}\in \NN $ such that 
for each $n\geq n_{0}$, 
$J(g_{\lambda ,n})<\zeta _{1}<\zeta _{2}<J(g_{\xi ,n}).$ 
From Lemma~\ref{greenk1k2}, it follows that 
\begin{equation}
\label{polyandrathm1pf2}
\frac{-1}{\deg (g_{\lambda ,n})-1}\log |a(g_{\lambda ,n})|
<\beta _{1}<\beta _{2}<
\frac{-1}{\deg (g_{\xi ,n})-1}\log |a(g_{\xi ,n})|,
\end{equation} 
for each $n\geq n_{0}.$ 
By (\ref{polyandrathm1pf1}) and (\ref{polyandrathm1pf2}), 
we obtain $B_{\lambda }<_{r} B_{\xi }.$ 

 We now suppose that we have case 2.
Then, by Lemma~\ref{j1aj2lem}, 
there exists a doubly connected component $A$ of $F(G)$ such that 
$J_{\lambda }<A.$ Continuing the same argument as that of case 1, 
we obtain $B_{\lambda }\neq B_{\xi }.$ 
In order to show $+\infty \in B_{\xi }$, 
let $R$ be any number such that $P^{\ast }(G)\subset D(0,R).$ 
Since $P^{\ast }(G)\subset K(g)$ for each $g\in G$, 
combining it with (\ref{polyandrathm1pf0a}) and 
Lemma~\ref{inftyj1}, we see that there exists an $n_{0}=n_{0}(R)$ 
such that for each $n\geq n_{0}$, $D(0,R)<J(g_{\xi ,n}).$ 
From Lemma~\ref{greenk1k2}, it follows that 
$\frac{-1}{\deg (g_{\xi ,n})-1}\log |a(g_{\xi ,n})|
\rightarrow +\infty .$ Hence, $+\infty \in B_{\xi }.$ 
Therefore, we have shown 
Claim 2. 
%

 Combining Claims 1 and 2, the statement of the lemma follows. 
 
 Therefore, we have proved Lemma~\ref{part1-1lem}. 
\end{proof}

\ 

We now demonstrate Theorem~\ref{polyandrathm1}-\ref{polyandrathm1-1}.\\ 
{\bf Proof of Theorem~\ref{polyandrathm1}-\ref{polyandrathm1-1}:} 
From Lemma~\ref{part1-1lem}, Theorem~\ref{polyandrathm1}-\ref{polyandrathm1-1}
follows.
\qed 

\ 

We now demonstrate Corollary~\ref{polyandracor}.\\ 
{\bf Proof of Corollary~\ref{polyandracor}:}
By Theorem~\ref{repdense}, 
we have 
$J(\Theta (G))=
\overline{ \cup _{h\in \Theta (G)}J(h) }=
\overline{\cup _{g\in G}J(\Theta (g))}$, 
where  the closure is taken in $\CCI .$  
Since $J(\Theta (g))=
\{ z\in \CC \mid |z|=|a(g)|^{-\frac{1}{\deg(g)-1}}\} $,  
we obtain 
\begin{equation}
\label{polyandracorpf1}
J(\Theta (G))=
\overline{\cup _{g\in G}\{ z\in \CC 
\mid |z|=|a(g)|^{-\frac{1}{\deg(g)-1}}\} },
\end{equation}
where the closure is taken in $\CCI .$ 
Hence, we see that 
$\sharp ({\hat{\cal J}}_{\Theta (G)})$ 
is equal to the cardinality of 
the set of all connected components of 
$J(\Theta (G))\cap [0,+\infty ].$
Moreover, let $\psi :[0,+\infty ]\rightarrow \hat{\RR }$ 
be the homeomorphism defined by 
$\psi (x):=\log (x)$ for $x\in (0,+\infty )$, 
$\psi (0):=-\infty $, and $\psi (+\infty )=+\infty .$ 
Then, 
(\ref{polyandracorpf1}) implies that, 
the map $\psi : 
[0,\infty ]\rightarrow \hat{\RR }$, 
maps $J(\Theta (G))\cap [0,+\infty ]$
onto $M(\Psi (\Theta (G))).$ 
For any $J\in \hat{{\cal J}}_{\Theta (G)}$, 
let $\tilde{\psi }(J)\in {\cal M}_{\Psi (\Theta (G))}=
{\cal M}_{\Psi (G)}$ be the element such that 
$\psi (J\cap [0,+\infty ])=\tilde{\psi }(J).$ 
Then,  
the map $\tilde{\psi }:\hat{{\cal J}}_{\Theta (G)}\rightarrow 
{\cal M}_{\Psi (\Theta (G))}={\cal M}_{\Psi (G)}$ is a 
bijection, and moreover, 
for any $J_{1},J_{2}\in {\cal J}_{\Theta (G)}$, 
we have that $J_{1}<J_{2}$ if and only if 
$\tilde{\psi }(J_{1})<_{r}\tilde{\psi }(J_{2}).$ 
Furthermore, for any $J\in \hat{{\cal J}}_{\Theta (G)}$, 
$\infty \in J$ if and only if $+\infty \in \tilde{\psi}(J).$ 
Let $\tilde{\Theta }:\hat{{\cal J}}_{G}\rightarrow 
\hat{{\cal J}}_{\Theta (G)}$ be the map defined by 
$\tilde{\Theta }=\tilde{\psi }^{-1}\circ \tilde{\Psi }$, 
where $\tilde{\Psi }: \hat{{\cal J}}_{G}\rightarrow {\cal M}_{\Psi (G)}$ is the 
map in Lemma~\ref{part1-1lem}.  
Then, by Lemma~\ref{part1-1lem}, 
$\tilde{\Theta }:\hat{{\cal J}}_{G}\rightarrow 
\hat{{\cal J}}_{\Theta (G)}$ is injective, and moreover, 
if $J_{1},J_{2}\in {\cal J}_{G}$ and $J_{1}<J_{2}$, then 
$\tilde{\Theta }(J_{1})\in {\cal J}_{\Theta (G)}$, 
$\tilde{\Theta }(J_{2})\in {\cal J}_{\Theta (G)}$, 
and $\tilde{\Theta }(J_{1})<\tilde{\Theta }(J_{2}).$ 

 Thus, we have proved Corollary~\ref{polyandracor}.
\qed  

\ 

 We now demonstrate Theorem~\ref{polyandrathm2}.

\noindent {\bf Proof of Theorem~\ref{polyandrathm2}:}
We have that for any $j=1,\ldots ,m$, 
$(\Psi (h_{j}))^{-1}(x)=\frac{1}{\deg (h_{j})}(x-\log |a_{j}|)
=\frac{1}{\deg (h_{j})}(x-\frac{-1}{\deg (h_{j})-1}\log |a_{j}|)+
\frac{-1}{\deg (h_{j})-1}\log |a_{j}|$, 
where $x\in \RR .$ 
Hence, it is easy to see that 
$\cup _{j=1}^{m}(\Psi (h_{j}))^{-1}([\alpha ,\beta ])
\subset [\alpha ,\beta ].$ 
From the assumption, it follows that 
\begin{equation}
\label{polyandrathm2pf1}
\cup _{j=1}^{m}(\Psi (h_{j}))^{-1}([\alpha ,\beta ])=
[\alpha ,\beta ].
\end{equation} 
Moreover, by Lemma~\ref{hmslem}-\ref{bss}, 
we have 
\begin{equation}
\label{polyandrathm2pf2}
\cup _{j=1}^{m}(\eta (\Psi (h_{j})))^{-1}(J(\eta (\Psi (G))))=J(\eta (\Psi (G))).
\end{equation}
Furthermore, by Lemma~\ref{cptmsetlem}, $J(\eta (\Psi (G)))$ is a compact 
subset of $\RR .$ 
Applying \cite[Theorem 2.6]{F}, it follows that 
$J(\eta (\Psi (G)))=[\alpha ,\beta ].$ 
Combined with Lemma~\ref{cptmsetlem}, we obtain 
$M(\Psi (G))=[\alpha ,\beta ].$ Hence, $M(\Psi (G))$ is connected. 
Therefore, from Theorem~\ref{polyandrathm1}-\ref{polyandrathm1-1}, 
it follows that 
$J(G)$ is connected.
\qed 

\ 

 We now demonstrate Theorem~\ref{polyandrathm3}.\\ 
{\bf Proof of Theorem~\ref{polyandrathm3}:}
Let $C$ be a set of polynomials of degree two such that $C$ generates $G.$ 
 Suppose that $J(G)$ is disconnected. Then, 
by Theorem~\ref{mainth0}, 
there exist two elements $h_{1},h_{2}\in C$ 
such that the semigroup 
$H=\langle h_{1},h_{2}\rangle $ satisfies that 
$J(H)$ is disconnected. 
For each $j=1,2$, let $a_{j}$ be the coefficient of 
the highest degree term of polynomial $h_{j}.$ 
Let $\alpha := 
\min _{j=1,2}\{ \frac{-1}{\deg (h_{j})-1}\log |a_{j}|\} $ 
and 
$\beta :=\max _{j=1,2}\{ \frac{-1}{\deg (h_{j})-1}\log |a_{j}|\}.$ 
Then we have that $\alpha =\min _{j=1,2}\{- \log |a_{j}|\} $ 
and $\beta =\max _{j=1,2}\{- \log |a_{j}|\} .$ 
Since 
$\Psi (h_{j})^{-1}(x)=\frac{1}{2}(x-\log |a_{j}|)
=\frac{1}{2}(x-(-\log |a_{j}|))+(-\log |a_{j}|)$ for each $j=1,2$, 
we obtain  
$[\alpha ,\beta ]= 
\cup _{j=1}^{2}(\Psi (h_{j}))^{-1}([\alpha ,\beta ]).$ 
Hence, by Theorem~\ref{polyandrathm2}, 
it must be true that 
$J(H)$ is connected. However, this is a contradiction. 
Therefore, $J(G)$ must be connected.
\qed 

\ 

We now demonstrate Theorem~\ref{polyandrathm4}. \\ 
{\bf Proof of Theorem~\ref{polyandrathm4}:}
For each $\lambda \in \Lambda $, let 
$b_{\lambda }$ be the fixed point of 
$\Psi (h_{\lambda })$ in $\RR .$ 
It is easy to see that 
$b_{\lambda }=\frac{-1}{\deg (h_{\lambda })-1}\log |a_{\lambda }|$, 
for each $\lambda \in \Lambda .$ 
From the assumption, 
it follows that 
there exists a point $b\in \RR $ such that 
for each $\lambda \in \Lambda $, 
$b_{\lambda }=b.$ 
This implies that 
for any element $g\in G$, 
the fixed point $b(g)\in \RR $ of 
$\Psi (g)$ in $\RR $ is equal to $b.$ 
Hence, we obtain  
$M(\Psi (G))=\{ b\} .$ 
Therefore, $M(\Psi (G))$ is connected. 
From Theorem~\ref{polyandrathm1}-\ref{polyandrathm1-1}, 
it follows that 
$J(G)$ is connected.
\qed 

\ 
\subsection{Proofs of results in \ref{Properties}}
\label{Proof of Properties}
In this section, we prove results in \ref{Properties}, 
Theorem~\ref{polyandrathm1}-\ref{polyandrathm1-2} and 
Theorem~\ref{polyandrathm1}-\ref{polyandrathm1-3}. 

 In order to demonstrate Theorem~\ref{mainth2},  
Theorem~\ref{polyandrathm1}-\ref{polyandrathm1-2}, and 
Theorem~\ref{polyandrathm1}-\ref{polyandrathm1-3}, we need the following lemma.
\begin{lem}
\label{gdisinflem}
If $G\in {\cal G}_{dis}$, then 
$\infty \in F(G).$ 
\end{lem}
\begin{proof}
Suppose that $G\in {\cal G}_{dis}$ and  
$\infty \in J(G).$ We will deduce a contradiction. 
By Lemma~\ref{inftyj1}, 
the element $J\in \hat{{\cal J}}_{G}$ with 
$\infty \in J$ satisfies that 
$J=\{ \infty \} .$ 
Hence, by Lemma~\ref{appj2lem}, 
for each $n\in \NN $, 
there exists an element $g_{n}\in G$ such that 
$J(g_{n})\subset B(\infty ,\frac{1}{n}).$ 
Let $R>0$ be any number which is sufficiently large so that 
$P^{\ast }(G)\subset B(0,R).$ Since 
we have that $P^{\ast }(G)\subset K(g)$ for each $g\in G$, 
it must hold that there exists an number $n_{0}=n_{0}(R)\in \NN $ 
such that for each $n\geq n_{0}$, 
$B(0,R)<J(g_{n}).$ From Lemma~\ref{greenk1k2}, 
it follows that 
$\lim _{z\rightarrow \infty }
(\log|z|-\varphi (F_{\infty }(g_{n}),z))\rightarrow +\infty 
$ as $n\rightarrow \infty .$
Hence, we see that 
$\frac{-1}{\deg (g_{n})-1}\log |a(g_{n})|\rightarrow 
+\infty $, as $n\rightarrow \infty .$ 
This implies that 
\begin{equation}
\label{typeproppf0}
|a(g_{n})|^{-\frac{1}{\deg(g_{n})-1}}\rightarrow \infty , \mbox{ as }
n\rightarrow \infty .
\end{equation}
Furthermore, by Theorem~\ref{polyandrathm1}-\ref{polyandrathm1-1}, 
we must have that $M(\Psi (G))$ is disconnected. 

 We now consider the polynomial semigroup 
 $H=\{ z\mapsto |a(g)|z^{\deg(g)}\mid g\in G\} \in {\cal G}.$
By Theorem~\ref{repdense}, 
we have 
$J(H)=\overline{\cup _{h\in H}J(h)}.$ 
Since the Julia set of polynomial $|a(g)|z^{\deg(g)}$ 
is equal to $\{ z\in \CC \mid |z|=|a(g)|^{-\frac{1}{\deg(g)-1}}\} $, 
it follows that 
\begin{equation}
\label{typeproppf1}
J(H)=\overline{\cup _{g\in G}\{ z\in \CC \mid 
|z|=|a(g)|^{-\frac{1}{\deg(g)-1}}\} }, 
\end{equation}
where the closure is taken in $\CCI .$ 
Moreover, $J(\Theta (G))=J(H).$ 
Combining it with (\ref{typeproppf0}), 
(\ref{typeproppf1}), and 
Corollary~\ref{polyandracor}, 
we see that 
\begin{equation}
\label{gdisinflempf1}
\infty \in J(H),  
 \mbox{ and } J(H) \mbox{ is disconnected.}
\end{equation} 
Let $\psi :[0,+\infty ]\rightarrow \hat{\RR }$ be 
the homeomorphism as in the proof of 
Corollary~\ref{polyandracor}. 
By (\ref{typeproppf1}), we have 
\begin{equation}
\label{gdisinflempf2}
\psi (J(H)\cap [0,+\infty ])= M(\Psi (H))=M(\Psi (G)). 
\end{equation}
Moreover, by Lemma~\ref{hmslem}-\ref{invariant}, we have 
\begin{equation}
\label{gdisinflempf3}
h(F(H)\cap [0,+\infty ])\subset F(H)\cap [0,+\infty ], 
\mbox{ for each } h\in H.
\end{equation}
Furthermore, we have that  
\begin{equation}
\label{gdisinflempf4}
\psi \circ h=\Psi (h)\circ \psi \mbox{ on } [0,+\infty ], 
\mbox{ for each } h\in H. 
\end{equation}
Combining (\ref{gdisinflempf2}), (\ref{gdisinflempf3}), 
and (\ref{gdisinflempf4}), we see that 
\begin{equation}
\label{gdisinflempf5}
\Psi (h)(\hat{\RR }\setminus M(\Psi (H)))\subset 
(\hat{\RR }\setminus M(\Psi (H))), \mbox{ for each } 
h\in H.
\end{equation}
By Lemma~\ref{inftyj1} and (\ref{gdisinflempf1}), 
we get that the connected component $J$ of 
$J(H)$ containing $\infty $ satisfies that 
\begin{equation}
\label{gdisinflempf6}
J=\{ \infty \} .
\end{equation} 
Combined with Lemma~\ref{appj2lem}, we see that 
for each $n\in \NN $, there exists an element $h_{n}\in H $ such that 
\begin{equation}
\label{gdisinflempf7}
J(h_{n})\subset B(\infty ,\frac{1}{n}).
\end{equation}
 Combining (\ref{typeproppf1}), (\ref{gdisinflempf2}), 
 (\ref{gdisinflempf6}), and (\ref{gdisinflempf7}),  
we obtain the following claim.\\ 
Claim 1:   
$+\infty $ is a non-isolated point of $M(\Psi (H))$ and 
the connected component of $M(\Psi (H))$ containing $+\infty $ is 
equal to $\{ +\infty \} .$ 

 Let $h\in H$ be an element. Conjugating $G$ by some 
 linear transformation, we may assume that 
$h$ is of the form $h(z)=z^{s}, s\in \NN , s>1.$ 
Hence $\Psi (h)(x)=sx, s>1.$  
Since $0$ is a fixed point of $\Psi (h)$, we have that 
$0\in M(\Psi (H)).$ By Claim 1, 
there exists $c_{1},c_{2}\in [0,+\infty )$ with  
$c_{1}<c_{2}$ such that the open interval $I=(c_{1},c_{2})$ is a 
connected component of $\hat{\RR }\setminus M(\Psi (H)).$ 
 We now show the following claim. \\ 
Claim 2:  
Let $Q=(r_{1},r_{2})\subset (0, + \infty )$ be any connected open interval in  
$\hat{\RR } \setminus M(\Psi (H))$, where 
$0\leq r_{1}<r_{2}<+\infty .$ Then, we have $r_{2}\leq sr_{1}.$ 

 To show this claim, suppose that $sr_{1}< r_{2}.$ Then, 
it implies that\\ $\cup _{n\in \NN \cup \{ 0\} }\Psi (h)^{n}(Q)=(r_{1},+\infty ).$  
However, by (\ref{gdisinflempf5}), 
we have $\cup _{n\in \NN \cup \{ 0\} }\Psi (h)^{n}(Q)\subset 
\hat{\RR }\setminus M(\Psi (H))$, which implies that 
the connected component $Q'$ of $\hat{\RR }\setminus M(\Psi (H))$ containing $Q$ satisfies that  
$Q'\supset (r_{1},+\infty ).$ 
This contradicts Claim 1. 
Hence, we obtain Claim 2.

 By Claim 2, we obtain $c_{1}>0.$ Let
 $c_{3}\in (0,c_{1})$ be a number so that 
 $c_{2}-c_{3}>s(c_{1}-c_{3}).$ 
Since $c_{1}\in M(\Psi (H))$, there exists an element 
$c\in (c_{3}, c_{1}]$ and an element $h_{1}\in H$ such that 
$\Psi (h_{1})(c)=c$ and $(\Psi (h_{1}))'(c)>1.$ 
Since $c_{2}-c_{3} >s(c_{1}-c_{3})$, we obtain 
\begin{equation}
\label{gdisinflempf8}
c_{2}-c>s(c_{1}-c).
\end{equation}
Let $t:=(\Psi (h_{1}) )'(c)>1.$      
Then, for each $n\in \NN $, 
we have 
$(\Psi (h_{1}) )^{n}(I)=
(t^{n}(c_{1}-c)+c, t^{n}(c_{2}-c)+c).$ 
From Claim 2, it follows that 
$t^{n}(c_{2}-c)+c\leq s(t^{n}(c_{1}-c)+c)$, for each $n\in \NN .$  
Dividing both sides by $t^{n}$ and then letting $n\rightarrow \infty $, 
we obtain $c_{2}-c\leq s(c_{1}-c).$ However, this contradicts
(\ref{gdisinflempf8}). 
Hence, we must have that 
$\infty \in F(G).$ Thus, we have proved Lemma~\ref{gdisinflem}.
\end{proof} 

\ 

We now demonstrate Proposition~\ref{fcprop}.\\ 
\noindent {\bf Proof of Proposition~\ref{fcprop}:}
Let $U$ be a connected component of $F(G)$ with $U\cap \hat{K}(G)\neq 
\emptyset .$ 
Let $g\in G$ be an element. 
Then we have $\hat{K}(G)\cap F(G)\subset $ int$(K(g)).$ 
Since $h(F(G))\subset F(G)$ and $h(\hat{K}(G)\cap F(G))\subset 
\hat{K}(G)\cap F(G)$ for each $h\in G,$ 
it follows that $h(U)\subset 
$ int$(K(g))$ for each $h\in G.$ Hence 
$U\subset $ int$(\hat{K}(G)).$ From this, 
it is easy to see 
$\hat{K}(G)\cap F(G)=$ int$(\hat{K}(G)).$ 
By the maximum principle, 
we see that $U$ is simply connected.
\qed \\ 

We now demonstrate Theorem~\ref{mainth2}. 

\noindent {\bf Proof of Theorem~\ref{mainth2}:}
%

 First, we show statement \ref{mainth2-2}. 
By Lemma~\ref{gdisinflem}, we have that 
$\infty \in F(G).$ 
 Let $F_{\infty }(G)$ be the connected component of 
 $F(G)$ containing $\infty .$ 
 Let $J\in {\cal J}$ be an element 
 such that $\partial F_{\infty }(G)\cap J\neq \emptyset .$ 
 Let $D$ be the unbounded component of 
 $\CCI \setminus J.$ Then 
 $F_{\infty }(G)\subset D$ and $D$ is 
 simply connected. We show 
 $F_{\infty }(G)=D.$ Otherwise, 
 there exists an element $J_{1}\in {\cal J}$ such that  
 $J_{1} \neq J$ and $J_{1}\subset D.$ 
 By Theorem~\ref{mainth1}-\ref{mainth1-1}, we have 
 either $J_{1}< J$ or $J< J_{1}.$ 
 Hence, it follows that $J< J_{1}$ and we have 
 that $J$ is included in a bounded component 
 $D_{0}$ of $\CC \setminus J_{1}.$ Since 
 $F_{\infty }(G)$ is included in the unbounded 
 component $D_{1}$ of $\CCI \setminus J_{1},$ 
 it contradicts $\partial F_{\infty }(G)\cap J\neq 
 \emptyset .$ Hence, $F_{\infty }(G)=D$ and 
 $F_{\infty }(G)$ is simply connected.  

 Next, let $J_{\max }$ be the element of 
 ${\cal J}$ with $\partial F_{\infty }(G)\subset J_{\max }$, and 
 suppose that there exists an element 
 $J\in {\cal J}$ such that 
 $J_{\max }<J.$ Then $J_{\max }$ is 
 included in a bounded component of 
 $\CC \setminus J.$ On the other hand, 
 $F_{\infty }(G)$ is included in the 
 unbounded component of $\CCI \setminus J.$ 
 Since $\partial F_{\infty }(G)\subset J_{\max },$ 
 we have a contradiction. Hence, 
 we have shown that $J\leq J_{\max }$ for each 
 $J\in {\cal J}.$ 

 Therefore, we have shown statement \ref{mainth2-2}.

 Next, we show statement \ref{mainth2-3}. 
 Since $\emptyset \neq P^{\ast }(G) \subset 
 \hat{K}(G),$ we have $\hat{K}(G)\neq \emptyset .$ 
 By Proposition~\ref{fcprop}, we have 
 $\partial \hat{K}(G)\subset J(G).$ 
 Let $J_{1}$ be a connected component 
 of $J(G)$  
 with $J_{1}\cap 
 \partial \hat{K}(G)\neq \emptyset .$ By Lemma~\ref{inftyj1}, $J_{1}\in {\cal J}.$  
 Suppose that there exists an element 
 $J\in {\cal J}$ such that 
 $J<J_{1}.$ Let $z_{0}\in J$ be a point. 
 By Theorem~\ref{repdense}, there exists a 
 sequence $\{ g_{n}\} _{n\in \NN }$ in $G$ such that 
 $d(z_{0},J(g_{n}))\rightarrow 0$ as $n\rightarrow \infty .$ 
 Then by Lemma~\ref{appjlem}, 
 $\sup\limits _{z\in J(g_{n})}d(z,J)\rightarrow 
0$ as $n\rightarrow \infty .$ 
 Since $J_{1}$ is included in the unbounded component 
 of $\CC \setminus J$, it follows that 
 for a large $n\in \NN ,$ 
 $J_{1}$ is included in the unbounded component 
 of $\CC \setminus J(g_{n}).$ 
 However, this causes a contradiction, since 
 $J_{1}\cap \hat{K}(G)\neq \emptyset .$ 
 Hence, by Theorem~\ref{mainth1}-\ref{mainth1-1}, it must hold that 
 $J_{1}\leq J$ for each $J\in {\cal J}.$ 
 This argument shows that 
 if $J_{1}$ and $J_{2}$ are two connected 
 components of $J(G)$ such that 
 $J_{i}\cap \partial \hat{K}(G)\neq \emptyset $ for each 
 $i=1,2$, then $J_{1}=J_{2}.$ Hence,\ 
 we conclude that there exists a unique minimal element 
 $J_{\min }$ in $({\cal J},\leq )$ and   
 $\partial \hat{K}(G)\subset J_{\min }.$ 

 Next, let $D$ be the unbounded component of 
 $\CC \setminus J_{\min }.$ 
 Suppose $D\cap \hat{K}(G)\neq \emptyset .$ 
 Let $x\in D\cap \hat{K}(G)$ be a point. 
 By Theorem~\ref{repdense} and Lemma~\ref{appjlem}, 
 there exists a sequence $\{ g_{n}\} _{n\in \NN }$ in $G$ 
 such that 
 $\sup\limits _{z\in J(g_{n})}d(z,J_{\min})\rightarrow 
0$ as $n\rightarrow \infty .$ 
Then, for a large $n\in \NN $, 
$x$ is in the unbounded component of 
$\CC \setminus J(g_{n}).$ 
However, this is a contradiction, 
since $g_{n}^{l}(x)\rightarrow \infty 
$ as $l\rightarrow \infty $,  
and $x\in \hat{K}(G).$
 Hence, we have shown statement \ref{mainth2-3}.

 Next, we show statement \ref{mainth2-3b}. 
 By Theorem~\ref{mainth0}, there exist 
 $\lambda _{1},\lambda _{2}\in \Lambda $ and 
 connected components $J_{1},J_{2}$ of $J(G)$ 
 such that 
 $J_{1}\neq J_{2}$ and 
 $J(h_{\lambda _{i}})\subset J_{i}$ for each $i=1,2.$ 
 By Lemma~\ref{inftyj1}, we have $J_{i}\in {\cal J}$ for each 
 $i=1,2.$ Then $J(h_{\lambda _{1}})\cap 
 J(h_{\lambda _{2}})=\emptyset .$ 
 Since $P^{\ast }(G) $ is bounded 
 in $\CC $, we may assume 
 $J(h_{\lambda _{2}})< J(h_{\lambda _{1}}).$ 
 Then we have $ K(h_{\lambda _{2}})\subset 
 $ int$(K(h_{\lambda _{1}})) $ and 
 $J_{2}<J_{1}.$ By statement \ref{mainth2-3}, 
 $J_{1}\neq J_{\min }.$ Hence 
 $J(h_{\lambda _{1}})\cap J_{\min }=\emptyset .$  
 Since $P^{\ast }(G)$ is bounded 
 in $\CC $, we have that 
 $ K(h_{\lambda _{2}})$ is connected. 
 Let $U$ be the connected component of 
 int$(K(h_{\lambda _{1}}))$ containing 
 $K(h_{\lambda _{2}}).$ 
 Since $P^{\ast }(G) \subset 
 K(h_{\lambda _{2}}),$ it follows that 
 there exists 
 an attracting fixed point $z_{1}$ of 
 $h_{\lambda _{1}}$ in $K(h_{\lambda _{2}})$ and  
 $U$ is the immediate attracting basin 
 for $z_{1}$ with respect to the 
 dynamics of $h_{\lambda _{1}}.$  
  Furthermore, 
 by Corollary~\ref{fibconncor}, 
$ h_{\lambda _{1}}^{-1}(J(h_{\lambda _{2}}))$ 
is connected. Therefore,  
$h_{\lambda _{1}}^{-1}(U)=U.$ Hence,  
int$(K(h_{\lambda _{1}}))=U.$ 

 Suppose that there exists an $n\in \NN $ such that 
 $h_{\lambda _{1}}^{-n}(J(h_{\lambda _{2}}))
  \cap J(h_{\lambda _{2}})\neq \emptyset .$ 
  Then, by Corollary~\ref{fibconncor}, $A:=\cup _{s\geq 0 } h_{\lambda _{1}}^{-ns}
  (J(h_{\lambda _{2}})) $ is connected and 
  its closure $\overline{A}$ contains 
  $J(h_{\lambda _{1}}).$ Hence 
  $J(h_{\lambda _{1}})$ and $J(h_{\lambda _{2}})$ 
  are included in the same connected component of $J(G).$ 
  This is a contradiction. 
  Therefore, for each $n\in \NN $, 
  we have $h_{\lambda _{1}}^{-n}(J(h_{\lambda _{2}}))
  \cap J(h_{\lambda _{2}})=\emptyset .$ Similarly, 
  for each $n\in \NN $, 
  we have $h_{\lambda _{2}}^{-n}(J(h_{\lambda _{1}}))
  \cap J(h_{\lambda _{1}})=\emptyset .$ 
  Combining $h_{\lambda _{1}}^{-1}(J(h_{\lambda _{2}}))
  \cap J(h_{\lambda _{2}})=\emptyset $ 
   with $z_{1}\in K(h_{\lambda _{2}})$, 
  we obtain $z_{1}\in $ int$(K(h_{\lambda _{2}})).$ 
Hence, we have proved statement \ref{mainth2-3b}. 
   
 We now prove statement \ref{mainth2ast1}. Let $g\in G$ be an element with 
 $J(g)\cap J_{\min }=\emptyset .$ 
We show the following:\\ 
Claim 2: $J_{\min }<J(g).$ 

 To show the claim, suppose 
 that $J_{\min }$ is included in 
 the unbounded component $U$ of 
 $\CC \setminus J(g).$ 
Since $\emptyset \neq \partial \hat{K}(G)\subset  J_{\min }$, 
it follows that $\hat{K}(G)\cap U\neq \emptyset .$ However, 
this is a contradiction.
Hence, 
 we have shown Claim 2. 

 Combining Claim 2, Theorem~\ref{repdense} and 
 Lemma~\ref{appjlem}, we get that 
 there exists an element $h_{1}\in G$ 
 such that $J(h_{1})<J(g).$  
 From an argument which we have used in the proof of 
 statement \ref{mainth2-3b},  
 it follows that $g$ has an attracting fixed 
 point $z_{g}$ in $\CC $ and 
 int$(K(g))$ consists of only one immediate  
 attracting basin 
 for $z_{g}.$ 
  Hence, we have shown 
  statement \ref{mainth2ast1}.

Next, we show statement \ref{mainth2-4}.
Suppose that int$(\hat{K}(G))=\emptyset .$ We will deduce a 
contradiction.  If int$(\hat{K}(G))=\emptyset $, then 
by Proposition~\ref{fcprop}, we obtain 
$F(G)\cap \hat{K}(G)=\emptyset .$ 
By statement \ref{mainth2-3b}, 
there exist two elements $g_{1}$ and 
$g_{2}$ of $G$ and two elements 
$J_{1}$ and $J_{2}$ of ${\cal J}$ 
such that $J_{1}\neq J_{2}$, 
such that $J(g_{i})\subset J_{i}$ for each 
$i=1,2$, such that $g_{1}$ has an attracting fixed 
point $z_{0}$ in int$(K(g_{2}))$, and such that 
$K(g_{2})\subset $ int$(K(g_{1})).$ 
Since we assume $F(G)\cap \hat{K}(G)=\emptyset ,$ 
we have $z_{0}\in P^{\ast }(G)\subset \hat{K}(G)\subset J(G).$ 
Let $J$ be the connected component of 
$J(G)$ containing $z_{0}.$ 
We now show $J=\{ z_{0}\} .$ 
Suppose $\sharp J\geq 2.$  
Then $J(g_{1})\subset 
\overline{\cup _{n\geq 0}g_{1}^{-n}(J)}.$
Moreover,  
by Theorem~\ref{mainth1}-\ref{mainth1-3},  
$g_{1}^{-n}J $ is connected 
for each $n\in \NN .$ 
Since $g_{1}^{-n}(J)\cap J\neq \emptyset $ 
for each $n\in \NN $, we see that 
$\overline{\cup _{n\geq 0}g_{1}^{-n}(J)}$ 
is connected. Combining this with 
$z_{0}\in $ int$(K(g_{2}))$,  
$K(g_{2})\subset $ int$(K(g_{1}))$, 
$z_{0}\in J$ and 
$J(g_{1})\subset 
 \overline{\cup _{n\geq 0}g_{1}^{-n}(J)}$, 
 we obtain 
 $ \overline{\cup _{n\geq 0}g_{1}^{-n}(J)}\cap J(g_{2})
 \neq \emptyset .$ 
 Then it follows that 
 $J(g_{1})$ and $J(g_{2})$ are included 
 in the same connected component of $J(G).$ 
 This is a contradiction. Hence, we have shown 
 $J=\{ z_{0}\} .$ By statement \ref{mainth2-3}, 
 we obtain $\{ z_{0}\} =J_{\min }
 =P^{\ast }(G).$ 
Let 
$\varphi (z):= \frac{1}{z-z_{0}}$ and 
let $\tilde{G}:= \{ \varphi g\varphi ^{-1}\mid g\in G\} .$ 
Then $\tilde{G}\in {\cal G}_{dis}.$ Moreover, 
since $z_{0}\in J(G)$, we have that 
$\infty \in J(\tilde{G}).$ This contradicts Lemma~\ref{gdisinflem}. 
Therefore, we must have that int$(\hat{K}(G))\neq \emptyset .$ 
%

 Since $\partial \hat{K}(G)\subset 
 J_{\min }$ (statement \ref{mainth2-3}) 
 and $\hat{K}(G)$ is bounded, 
 it follows that $\CC \setminus J_{\min } $ 
 is disconnected and 
 $\sharp J_{\min }\geq 2.$ Hence, $\sharp J\geq 2$ 
 for each $J\in {\cal J}=\hat{{\cal J}}.$ 
 Now, let $g\in G$ be an element with 
 $J(g)\cap J_{\min }=\emptyset .$ 
 we show $J_{\min }\neq g^{\ast }(J_{\min }).$ 
 If $J_{\min }=g^{\ast }(J_{\min })$, then 
 $g^{-1}(J_{\min })\subset 
 J_{\min }.$ Since $\sharp J_{\min }\geq 3,$ 
 it follows that $J(g)\subset J_{\min }$, which  
 is a contradiction. 
 Hence, $J_{\min }\neq g^{\ast }(J_{\min })$, and so 
 $J_{\min }<g^{\ast }(J_{\min }).$ Combined with 
 Theorem~\ref{mainth1}-\ref{mainth1-3}, 
 we obtain $g^{-1}(J(G))\cap J_{\min }=\emptyset .$ 
 Since $g(\hat{K}(G))\subset \hat{K}(G)$, 
 we have $g($int$(\hat{K}(G)))\subset $ int$(\hat{K}(G)).$ 
 Suppose $g(\partial \hat{K}(G))\cap 
 \partial \hat{K}(G)\neq \emptyset .$ 
 Then, since $\partial \hat{K}(G)\subset J_{\min }$ 
 (statement \ref{mainth2-3}), 
 we obtain $g(J_{\min })\cap J_{\min }\neq 
 \emptyset .$ This implies 
 $g^{-1}(J_{\min })\cap J_{\min }\neq 
 \emptyset $, 
 which contradicts $g^{-1}(J(G))\cap J_{\min }=\emptyset .$ 
Hence, 
 it must hold 
 $g(\partial \hat{K}(G))\subset $ int$(\hat{K}(G))$, and so 
 $g(\hat{K}(G))\subset $ int$(\hat{K}(G)).$  
Moreover, since $g^{-1}(J(G))\cap J_{\min }=\emptyset $, 
we have that $g(J_{\min })$ is a connected subset of $F(G).$ 
Since $\partial \hat{K}(G)\subset J_{\min }$ 
and $g(\partial \hat{K}(G))\subset $ int$(\hat{K}(G))$, 
Proposition~\ref{fcprop} implies that $g(J_{\min })$ must be contained in int$(\hat{K}(G)).$ 

 By statement \ref{mainth2ast1}, 
 $g$ has a unique attracting fixed point 
 $z_{g}$ in $\CC .$ 
Then, $z_{g}\in P^{\ast }(G) 
\subset \hat{K}(G).$ Hence, 
$z_{g}=g(z_{g})\in g(\hat{K}(G))
\subset $ int$(\hat{K}(G)).$ 
 Hence, we have shown statement \ref{mainth2-4}.

 We now show statement \ref{mainth2ast2}. 
 Since $F_{\infty }(G)$ is simply connected 
 (statement \ref{mainth2-2}), 
we have $\cup _{A\in {\cal A}}A\subset \CC .$ 
 Suppose that there exist two distinct elements 
 $A_{1}$ and $A_{2}$ in ${\cal A}$ such that 
 $A_{1}$ is included in the unbounded component of 
 $\CC \setminus A_{2}$, and such that $A_{2}$ is included in the 
 unbounded component of $\CC \setminus A_{1}.$ 
For each $i=1,2$, Let $J_{i}\in {\cal J}$ be the 
element that intersects the bounded component of $\CC \setminus A_{i}.$ 
Then, $J_{1}\neq J_{2}.$ 
Since $({\cal J},\leq )$ is 
totally ordered (Theorem~\ref{mainth1}-\ref{mainth1-1}), 
we may assume that $J_{1}<J_{2}.$ 
Then, it implies that $A_{1}<J_{2}<A_{2}$, which is a contradiction. 
Hence, $({\cal A},\leq )$ is totally ordered. 
Therefore, we have proved statement \ref{mainth2ast2}. 

Thus, we have proved Theorem~\ref{mainth2}.
\qed 

\ 

 We now demonstrate Theorem~\ref{mainupthm}.\\ 
{\bf Proof of Theorem~\ref{mainupthm}:} 
First, we show Theorem~\ref{mainupthm}-\ref{mainupthm1}.
If 
$G\in {\cal G}_{con}$, then $J(G)$ is uniformly perfect. 

 We now suppose that $G\in {\cal G}_{dis}.$ 
Let $A$ be an annulus separating $J(G).$
Then $A$ separates $J_{\min }$ and $J_{\max }.$ 
 Let $D$ be the unbounded component 
 of $\CC \setminus J_{\min }$ and 
 let $U$ be the connected component 
 of $\CC \setminus J_{\max } $ containing 
 $J_{\min}.$ Then 
 it follows that $A\subset U\cap D.$ 
 Since $\sharp J_{\min }>1$ and $\infty \in F(G)$ 
 (Theorem~\ref{mainth2}), 
 we get that 
 the doubly connected domain 
 $U\cap D$ satisfies mod $(U\cap D) <\infty .$ 
 Hence, we obtain 
 mod $A\leq $ mod $(U\cap D )<\infty .$ 
 Therefore, $J(G)$ is uniformly perfect.  

 If a point $z_{0}\in J(G)$ is a superattracting fixed point of 
 an element $g\in G$, then, combining uniform perfectness of $J(G)$ and 
 \cite[Theorem 4.1]{HM2}, it follows that $z_{0}\in $ int$(J(G)).$  
Thus, we have shown Theorem~\ref{mainupthm}-\ref{mainupthm1}.
%

Next, we show Theorem~\ref{mainupthm}-\ref{mainupthm2}.
If $G\in {\cal G}$ and $\infty \in J(G)$, then 
by Lemma~\ref{gdisinflem}, we obtain  
$G\in {\cal G}_{con}.$ Moreover, 
Theorem~\ref{mainupthm}-\ref{mainupthm1} 
implies that 
$\infty \in $ int$(J(G)).$ Therefore, we have shown 
Theorem~\ref{mainupthm}-\ref{mainupthm2}.

 We now show Theorem~\ref{mainupthm}-\ref{mainupthm3}.
Suppose that $G\in {\cal G}_{dis}.$ 
 Let $g\in G$ and let $z_{1}\in J(G)\cap \CC $ with 
 $g(z_{1})=z_{1}$ and $g'(z_{1})=0.$ Then, 
$z_{1}\in P^{\ast }(G)\subset \hat{K}(G).$  
By Theorem~\ref{mainth2}-\ref{mainth2-3},  
  we obtain 
 $z_{1}\in J_{\min }.$
%
Moreover, 
Theorem~\ref{mainupthm}-\ref{mainupthm1} 
implies that 
 $z_{1}\in $ int$(J(G)).$ Combining this and 
 $z_{1}\in J_{\min }$, we obtain $z_{1}\in $ int$(J_{\min }).$ 
 By Theorem~\ref{mainth2}- \ref{mainth2-4-2}, we 
 obtain  
 $J(g)\subset J_{\min }.$ 

 Hence, we have shown Theorem~\ref{mainupthm}.
\qed 

\ 

We now demonstrate Theorem~\ref{polyandrathm1}-\ref{polyandrathm1-2}.\\ 
{\bf Proof of Theorem~\ref{polyandrathm1}-\ref{polyandrathm1-2}:}
Suppose $G\in {\cal G}_{dis}. $ Then, by Lemma~\ref{gdisinflem}, 
we obtain $\infty \in F(G).$ Hence, 
there exists a number $R>0$ such that 
for each $g\in G$, $J(g)<\partial B(0,R).$ 
From Lemma~\ref{greenk1k2}, it follows that 
there exists a constant $C_{1}>0$ such that 
for each $g\in G$, 
$\frac{-1}{\deg(g)-1}\log |a(g)|< C_{1}.$ 
This implies that there exists a constant $C_{2}\in \RR $ such that 
\begin{equation}
\label{polyandrathm1-2pf1}
M(\Psi (G))\subset [-\infty ,C_{1}].
\end{equation} 
Moreover, by Theorem~\ref{mainth2}-\ref{mainth2-4}, 
we have that int$(\hat{K}(G))\neq \emptyset .$ 
Let $B$ be a closed disc in int$(\hat{K}(G)).$ 
Then it must hold that 
for each $g\in G$, $B<J(g).$ 
Hence, by Lemma~\ref{greenk1k2}, there exists a constant 
$C_{3}\in \RR $ such that 
for each $g\in G$, 
$C_{3}\leq \frac{-1}{\deg(g)-1}\log |a(g)|.$ 
Therefore, we obtain  
\begin{equation}
\label{polyandrathm1-2pf2}
M(\Psi (G))\subset [C_{3},+\infty ].
\end{equation}
Combining (\ref{polyandrathm1-2pf1}) and (\ref{polyandrathm1-2pf2}), 
we obtain $M(\Psi (G))\subset \RR .$ 
Let $C_{4}$ be a large number so that 
$M(\Psi (G))\subset D(0,C_{4}).$ 
Since for each $g\in G$, 
the repelling fixed point 
$-\frac{1}{\deg(g)-1}\log |a(g)|$ of 
$\eta (\Psi (g))$ belongs to $D(0,C_{4})\cap \RR $, 
we see that for each $z\in \CC \setminus D(0,C_{4})$, 
$|\eta (\Psi (g))(z)|=
|\deg (g)(z-\frac{-1}{\deg(g)-1}\log |a(g)|)+
\frac{-1}{\deg(g)-1}\log |a(g)||
\geq \deg (g)C_{4}-(\deg (g)-1)C_{4}=
C_{4}.$ It follows that 
$\infty \in F(\eta (\Psi (G))).$ 
Combining this and Theorem~\ref{repdense}, 
we obtain  
$M(\Psi (G))=J(\eta (\Psi (G)))$, 
if $\sharp (J(\eta (\Psi (G))))\geq 3.$ 

 Suppose that $\sharp (J(\eta (\Psi (G))))= 2.$ 
Let $g\in G$ be an element and let $b\in \RR $ be 
the unique fixed point of $\Psi (g)$ in $\RR .$ 
Then, since $\infty \in F(\eta (\Psi (G)))$, 
there exists a point $c\in (J(\eta (\Psi (G)))\cap \CC )\setminus 
\{ b\} .$ By Lemma~\ref{hmslem}-\ref{invariant}, 
$(\eta (\Psi (g)))^{-1}(c)\in J(\eta (\Psi (G)))\setminus \{ b,c\} .$ 
This contradicts $\sharp (J(\eta (\Psi (G))))= 2.$ Hence it must hold that  
 $\sharp (J(\eta (\Psi (G))))\neq  2.$

 Suppose that $\sharp (J(\eta (\Psi (G))))=1 .$   
Since $M(\Psi (G))\subset \RR $ 
and $M(\Psi (G))\cap \RR \subset J(\eta (\Psi (G)))$, 
it follows that $M(\Psi (G))=J(\eta (\Psi (G))).$ 

 Therefore, we always have that $M(\Psi (G))=J(\eta (\Psi (G))).$ 
Thus, we have proved Theorem~\ref{polyandrathm1}-\ref{polyandrathm1-2}.
\qed 

\ 

We now demonstrate Theorem~\ref{polyandrathm1}-\ref{polyandrathm1-3}.\\ 
{\bf Proof of Theorem~\ref{polyandrathm1}-\ref{polyandrathm1-3}:} 
By Theorem~\ref{polyandrathm1}-\ref{polyandrathm1-1} and 
Theorem~\ref{polyandrathm1}-\ref{polyandrathm1-2}, 
the statement holds.
\qed 

\ 

We now demonstrate Proposition~\ref{orderjprop}.\\ 
\noindent {\bf Proof of Proposition~\ref{orderjprop}:}
First, we show statement \ref{orderjprop1}. Let $g\in Q_{1}.$ 
 We show the following:\\ 
\noindent Claim 1: For any element $J_{3}\in {\cal J}$ with 
$J_{1}\leq J_{3}$, we have 
$J_{1}\leq g^{\ast }(J_{3}).$ \\ 
To show this claim, let 
$J\in {\cal J}$ be an element with $J(g)\subset J.$
We consider the following two cases;\\   
Case 1: $J\leq J_{3}$, and\\ 
Case 2: $J_{1}\leq J_{3}\leq J.$ 

 Suppose that we have Case 1. Then, $J_{1}\leq J=g^{\ast }(J)\leq 
 g^{\ast }(J_{3}).$ Hence, the statement of Claim 1 is true. 

 Suppose that we have Case 2. 
If we have  $g^{\ast }(J_{3})< J_{3}$,  
then, we have $(g^{n})^{\ast }(J_{3})\leq g^{\ast }(J_{3})<J_{3}\leq J$ 
for each $n\in \NN .$ Hence, $\inf \{ d(z,J)\mid z\in g^{-n}(J_{3}),
n\in \NN \} >0.$ However, since $J(g)\subset J$ and 
$\sharp J_{3}\geq 3$, we obtain a contradiction. 
Hence, we must have $J_{3}\leq g^{\ast }(J_{3})$, which implies 
$J_{1}\leq J_{3}\leq g^{\ast }(J_{3}).$ 
Hence, we conclude that Claim 1 holds. 

 Now, let $K_{1}:= J(G)\cap (J_{1}\cup A_{1}).$ Then, 
 by Claim 1,\ we obtain $g^{-1}(K_{1})\subset K_{1}$, 
 for each $g\in Q_{1}.$ From Lemma~\ref{hmslem}-\ref{backmin}, 
 it follows that $J(H_{1})\subset K_{1}.$  
Hence, we have shown statement \ref{orderjprop1}.

 Next, we show statement \ref{orderjprop2}. Let $g\in Q_{2}.$ 
Then, by the same method as that of the proof 
 of Claim 1, we obtain the following.\\ 
 \noindent Claim 2: For any element $J_{4}\in {\cal J}$ with 
 $J_{4}\leq J_{2}$, we have 
 $g^{\ast }(J_{4})\leq J_{2}.$ 

 Now, let $K_{2}:=J(G)\cap (\CC \setminus A_{2}).$ 
 Then, by Claim 2, we obtain $g^{-1}(K_{2})\subset 
 K_{2}$, for each $g\in Q_{2}.$ From Lemma~\ref{hmslem}-\ref{backmin}, 
 it follows that $J(H_{2})\subset K_{2}.$ Hence, we have shown 
 statement \ref{orderjprop2}.

 Next, we show statement \ref{orderjprop3}. 
 By statements \ref{orderjprop1} and \ref{orderjprop2}, 
 we obtain 
 $J(H)\subset J(H_{1})\cap J(H_{2})\subset K_{1}\cap K_{2}
 \subset (\CC \setminus A_{2})\cap (J_{1}\cup A_{1})
 \subset J_{1}\cup (A_{1}\setminus A_{2}).$ 

 Hence, we have proved Proposition~\ref{orderjprop}.   
\qed 

\ 

We now demonstrate Proposition~\ref{bminprop}.\\ 
{\bf Proof of Proposition~\ref{bminprop}:}
By Theorem~\ref{mainth1} and 
Theorem~\ref{mainth2}, 
$({\cal J},\leq )$ is totally ordered 
and there exists a maximal element 
$J_{\max }$ and a minimal element 
$J_{\min }.$ 
Suppose that for any $h\in \G $, $J(h)\cap J_{\max }=\emptyset .$ 
Then, since $\sharp J_{\max}\geq 3$ (Theorem~\ref{mainth2}-\ref{mainth2-4-1}), 
we get that for any $h\in \G $, 
$h^{-1}(J_{\max })\cap J_{\max }=\emptyset .$ 
Combining it with Theorem~\ref{mainth1}-\ref{mainth1-3}, 
it follows that for any $h\in \G $, 
$h^{-1}(J(G))\cap J_{\max }=\emptyset .$ 
However, since 
$J(G)=\cup _{h\in \G}h^{-1}(J(G))$ (Lemma~\ref{hmslem}-\ref{bss}), 
it causes a contradiction. Hence, there must be an element 
$h_{1}\in \G $ such that $J(h_{1})\subset J_{\max }.$

 By the same method as above, we can show that there exists an 
 element $h_{2}\in \G $ such that $J(h_{2})\subset J_{\min }.$  
\qed 

\subsection{Proofs of results in \ref{Poly}}
\label{Proof of Poly}
In this section, we prove results in \ref{Poly}.

We now prove Theorem~\ref{fcthm}.\\ 
\noindent {\bf Proof of Theorem~\ref{fcthm}:} 
Combining the assumption and Theorem~\ref{mainth1}-\ref{mainth1-3}, 
we get that for each $h\in \G $ and each $j\in \{ 1,\ldots ,n\} $, 
there exists a $k\in \{ 1,\ldots ,n\} $ with 
$h^{-1}(J_{j})\subset J_{k}.$ 
Hence, 
\begin{equation}
\label{fcthmpfeq1}
h^{-1}(\cup _{j=1}^{n}J_{j})\subset \cup _{j=1}^{n}J_{j}, 
\mbox{ for each } h\in \G .
\end{equation}

Moreover, by Theorem~\ref{mainth2}-\ref{mainth2-4-1}, 
we obtain 
\begin{equation}
\label{fcthmpfeq2}
\sharp (\cup _{j=1}^{n}J_{j})\geq 3.
\end{equation}
Combining (\ref{fcthmpfeq1}), (\ref{fcthmpfeq2}), and 
Lemma~\ref{hmslem}-\ref{backmin}, 
it follows that $J(G)\subset \bigcup _{j=1}^{n}J_{j}.$ 
Hence, $J(G)=\cup _{j=1}^{n}J_{j}.$ 
Therefore, we have proved Theorem~\ref{fcthm}.

\qed 

 We now prove Proposition~\ref{fincomp}.

\noindent {\bf Proof of Proposition~\ref{fincomp}:} 
Let $n\in \NN $ with $n>1$ and let 
$\epsilon $ be a number with $0<\epsilon <\frac{1}{2}.$ 
For each $j=1,\ldots ,n$, 
let $\alpha _{j}(z)=\frac{1}{j}z^{2}$ and  
let $\beta _{j}(z)=\frac{1}{j}(z-\epsilon )^{2}+\epsilon .$  

For any large $l\in\NN $, 
there exists an open neighborhood $U$  
 of $\{ 0,\epsilon \} $ with  
 $U\subset \{ z\mid |z|<1\} $ and 
 a open neighborhood $V$ of 
 $(\alpha _{1}^{l},\ldots ,\alpha _{n}^{l},
 \beta _{1}^{l},\ldots ,\beta _{n}^{l})$ in 
 (Poly)$^{2n}$ such that 
 for each $(h_{1},\ldots ,h_{2n})\in V$, 
 we have $\cup _{j=1}^{2n}h_{j}(U)\subset U$ and 
 $\cup _{j=1}^{m}C(h_{j})\cap \CC \subset U$, 
 where $C(h_{j})$ denotes the set of all critical points of 
 $h_{j}.$ Then, by Remark~\ref{pcbrem}, for each $(h_{1},\ldots ,h_{m})\in V $,  
 $\langle h_{1},\ldots ,h_{2n}\rangle \in {\cal G}.$  
If $l$ is large enough and $V$ is so small, then, 
for each $(h_{1},\ldots ,h_{2n})\in V$, the set  
 $I_{j}:=
J(h_{j})\cup 
J(h_{j+n})$ is connected, for each $j=1,\ldots ,n$, and   
we have:
\begin{equation}
\label{finpropeq1}
(h_{i})^{-1}(I_{j})\cap I_{i}\neq \emptyset , 
(h_{i+n})^{-1}(I_{j})\cap I_{i}\neq \emptyset ,
\end{equation}
for each $(i,j).$ 
Furthermore, for a closed annulus $A=\{z\mid \frac{1}{2}
\leq |z|\leq n+1\} $,   
if $l\in \NN $ is large enough and $V$ is so small, then 
for each $(h_{1},\ldots ,h_{m})\in V$, 
$\cup _{j=1}^{2n}(h_{j})^{-1}(A)
\subset \mbox{int}(A)$  
 and 
$\{ (h_{j})^{-1}(A)\cup (h_{j+n})^{-1}(A)\} _{j=1}^{n}$ 
are mutually disjoint. 
Combining it with Lemma~\ref{hmslem}-\ref{backmin} and 
Lemma~\ref{hmslem}-\ref{bss}, we get 
that for each $(h_{1},\ldots ,h_{2n})\in V$,\ 
$J(\langle h_{1},\ldots ,h_{2n}\rangle)\subset A$ and 
$\{ J_{j}\} _{j=1}^{n}$ are 
mutually disjoint, where 
$J_{j}$ denotes the connected component of 
$J(\langle h_{1},\ldots ,h_{2n}\rangle )$ containing 
$I_{j}=J(h_{j})\cup J(h_{j+n}).$ 
Combining it with (\ref{finpropeq1}) and Theorem~\ref{fcthm}, 
it follows that for each $(h_{1},\ldots ,h_{2n})\in V$, 
 the polynomial semigroup 
$G=\langle h_{1},\ldots ,h_{2n}\rangle $ satisfies that  
$\sharp (\hat{{\cal J}}_{G})=n.$ 
\qed  

\ 

 To prove Theorem~\ref{countthm}, 
we need the following notation.
\begin{df} 
\ 
\begin{enumerate}
\item 
 Let $X$ be a metric space.
 Let $h_{j}:X\rightarrow X\ (j=1,\ldots ,m)$ be a 
 continuous map. 
 Let $G=\langle h_{1},\ldots ,h_{m}\rangle $ be 
 the semigroup generated by $\{ h_{j}\} .$
 A non-empty compact subset $L$ of $X$ is 
 said to be a {\bf backward self-similar set 
 with respect to $\{ h_{1},\ldots ,h_{m}\} $} 
 if  
 \begin{enumerate}
 \item $ L=\bigcup _{j=1}^{m}h_{j}^{-1}(L)$ 
 and 
 \item $g^{-1}(z)\neq \emptyset $ for each 
 $z\in L$ and $g\in G.$
 \end{enumerate}
 For example, if 
 $G=\langle h_{1},\ldots ,h_{m}\rangle $
 is a  
 finitely generated rational semigroup, then  
 the Julia set $J(G)$ is 
 a backward self-similar set with respect to 
 $\{ h_{1},\ldots ,h_{m}\} .$ 
 (See Lemma~\ref{hmslem}-\ref{bss}.)
\item We set $\Sigma _{m}:=
\{ 1,\ldots ,m\} ^{\NN }.$
For each $x=(x_{1},x_{2},\ldots ,)\in \Sigma _{m} $,  
we set 
$ L_{x}:=\bigcap _{j=1}^{\infty }
h_{x_{1}}^{-1}\cdots h_{x_{j}}^{-1}(L) \ (\neq \emptyset ).$ 

\item 
For a finite word $w=(w_{1},\ldots ,w_{k})\in 
\{ 1\ldots ,m\} ^{k}$, 
we set 
$h_{w}:=h_{w_{k}}\circ \cdots \circ h_{w_{1}}.$
 
\item Under the notation of \cite[page 110--page 115]{Sp},  
for any $k\in \NN $, 
let $\Omega _{k}=\Omega _{k}(L,\{ h_{1},\ldots ,h_{m}\} )$  
be the graph (one-dimensional simplicial complex) whose vertex set is  
$\{ 1,\ldots ,m\}^{k}$ and 
that satisfies that mutually different 
$w^{1},w^{2}\in \{ 1,\ldots ,m\}^{k}$ makes 
a $1$-simplex if and only if 
$\bigcap _{j=1}^{2} h_{w^{j}}^{-1}(L)\neq \emptyset .$  

 Let $\varphi _{k}:\Omega _{k+1}\rightarrow 
\Omega _{k}$ be the simplicial map 
defined by:\\ 
$(w_{1},\ldots ,w_{k+1})\mapsto (w_{1},\ldots ,w_{k})$   
for each 
$ (w_{1},\ldots ,w_{k+1})\in \{ 1,\ldots ,m\} ^{k+1}.$
Then $\{ \varphi _{k}:\Omega _{k+1}\rightarrow 
\Omega _{k}\} _{k\in \NN }$  makes an inverse system of 
simplicial maps. 
\item 
Let ${\cal C}(|\Omega _{k})|)$ be the 
set of all connected components of the realization 
$|\Omega _{k} |$ of $\Omega _{k}.$
Let $ \{ (\varphi _{k})_{\ast }: {\cal C}(|\Omega _{k+1}|)
\rightarrow {\cal C}(|\Omega _{k}|) \} _{k\in \NN }$ be 
the inverse system induced by 
$\{ \varphi _{k}\} _{k}.$      
\end{enumerate}
\end{df}
\noindent {\bf Notation:} 
We fix an $m\in \NN .$  
We set ${\cal W}^{\ast }:= 
\cup _{k=1}^{\infty } 
\{ 1,\ldots ,m\} ^{k}$ (disjoint union) and 
$\tilde{{\cal W}}:= {\cal W}^{\ast }\cup 
\Sigma _{m} $ (disjoint union). 
For an element $x\in \tilde{{\cal W}}$, 
we set $|x|= k$ if $x\in 
\{ 1,\ldots ,m\} ^{k} $, and 
$|x|=\infty $ if $x\in \Sigma _{m}.$ 
(This is called the word length of $x.$)  
For any $x\in \tilde{{\cal W}} $ and any 
$j\in \NN $ with $j\leq |x|$, we set 
$x|j:= (x_{1},\ldots x_{j})\in 
\{ 1,\ldots ,m\} ^{j} .$   
For any $x^{1}=(x_{1}^{1},\ldots ,
x_{p}^{1})\in {\cal W}^{\ast }$ 
and any $x^{2}=(x_{1}^{2},x_{2}^{2},\ldots )\in \tilde{{\cal W}}$, 
we set $
x^{1}x^{2}:=(x_{1}^{1},\ldots ,x_{p}^{1},
x_{1}^{2},x_{2}^{2},\ldots )\in \tilde{{\cal W}}.$

 To prove Theorem~\ref{countthm}, 
we need the following lemmas.
\begin{lem}
\label{phionto} 
Let $L$ be a backward self-similar set with respect to 
$\{ h_{1},\ldots ,h_{m}\} .$ 
Then, for each $k\in \NN $, the map 
$|\varphi _{k}|: |\Omega _{k+1}|\rightarrow 
|\Omega _{k}|$ induced from 
$\varphi _{k}:\Omega _{k+1}\rightarrow \Omega _{k}$ 
is surjective. 
In particular, 
$ (\varphi _{k})_{\ast }: 
{\cal C}(|\Omega _{k+1}|)\rightarrow 
{\cal C}(|\Omega _{k}|)$ is surjective.
\end{lem} 
\begin{proof}
 Let $x^{1},x^{2}\in 
 \{1 ,\ldots ,m\} ^{k}$ and 
 suppose that $\{ x^{1},x^{2}\} $ makes a $1$-simplex 
 in $\Omega _{k}.$ Then 
 $h_{x^{1}}^{-1}(L)\cap h_{x^{2}}^{-1}(L)
 \neq \emptyset .$ Since 
 $L=\cup _{j=1}^{m}h_{j}^{-1}(L)$, 
 there exist $x^{1}_{k+1}$ and 
 $x^{2}_{k+1}$ in $\{ 1,\ldots ,m\} $ 
 such that $
 h_{x^{1}}^{-1}h_{x^{1}_{k+1}}^{-1}(L)
 \cap h_{x^{2}}^{-1}h_{x^{2}_{k+1}}^{-1}(L)
 \neq \emptyset .$ 
 Hence, 
 $\{ x^{1}x^{1}_{k+1},x^{2}x^{2}_{k+1}\} $ 
 makes a $1$-simplex in 
 $\Omega _{k+1}.$ Hence 
 the lemma holds.
\end{proof}
\begin{lem}
\label{1disconlem}
Let $m\geq 2$ and 
let $L$ be a backward self-similar 
set with respect to 
$\{ h_{1},\ldots ,h_{m}\} .$ 
Suppose 
that for each $j$ with 
$j\neq 1$, $h_{1}^{-1}(L)\cap 
h_{j}^{-1}(L)=\emptyset .$
For each $k$, let $C_{k}\in {\cal C}(|\Omega _{k}|)$ 
 be the element containing 
 $(1,\ldots ,1)\in \{ 1,\ldots ,m\} ^{k}.$ 
 Then, we have the following.
 \begin{enumerate}
 \item \label{1disconlem1}
 For each $k\in \NN $, 
 $C_{k}=\{ (1,\ldots ,1)\} .$
  
\item \label{1disconlem2} For each $k\in \NN $, 
$\sharp ({\cal C}(|\Omega _{k}|))
<\sharp ({\cal C}(|\Omega _{k+1}|)).$
\item \label{1disconlem3}
$L$ has infinitely many connected components. 
\item \label{1disconlem4}
Let $x:=(1,1,1,\ldots )\in \Sigma _{m}$ and 
$x'\in \Sigma _{m}$ an element with $x\neq x'.$ Then, 
for any $y\in L_{x}$ and $y'\in L_{x'}$, there exists no 
connected component $A$ of $L$ such that 
$y\in A$ and $y'\in A.$  
\end{enumerate}
\end{lem}
\begin{proof}
 We show statement \ref{1disconlem1} by induction on $k.$ 
 We have  $C_{1}=\{ 1\} .$ 
 Suppose $C_{k}=\{ (1,\ldots ,1)\} .$ 
 Let $w\in \{ 1,\ldots ,m\} ^{k+1}\cap C_{k+1}$ be 
 any element. 
 Since 
 $(\varphi _{k})_{\ast }(C_{k+1})=C_{k}$, 
 we have $\varphi _{k}(w)=(1,\ldots ,1)\in \{ 1,\ldots ,m\} ^{k}.$ 
 Hence, $w|k=(1,\ldots ,1)\in 
 \{ 1,\ldots ,m\} ^{k}.$ 
 Since $
 h_{1}^{-1}(L)\cap h_{j}^{-1}(L)=\emptyset $ 
for each $j\neq 1$, 
we obtain $w=(1,\ldots ,1)\in \{ 1,\ldots ,m\} ^{k+1}.$ 
Hence, the induction is completed. Therefore, 
we have shown statement \ref{1disconlem1}. 

 Since both $(1,\ldots 1,1)\in \{ 1,\ldots ,m\} ^{k+1}$ 
 and $(1,\ldots ,1,2)\in \{ 1,\ldots ,m\} ^{k+1}$ are 
 mapped to $(1,\ldots ,1)\in \{ 1,\ldots ,m\} ^{k}$ under 
 $\varphi _{k}$, by statement \ref{1disconlem1} and 
 Lemma~\ref{phionto}, 
 we obtain statement \ref{1disconlem2}. 
 For each $k\in \NN $, 
 we have 
\begin{equation}
\label{1disconlemeq1}
 L=\coprod _{C\in {\cal C}(|\Omega _{k}|)}
 \ \bigcup _{w\in \{ 1,\ldots ,m\} ^{k}\cap C}
 h_{w}^{-1}(L).
\end{equation} 
 Hence, by statement \ref{1disconlem2}, we conclude  
 that $L$ has infinitely many connected components.

 We now show statement \ref{1disconlem4}. 
 Let $k_{0}:=\min \{ l\in \NN \mid x_{l}'\neq 1\} . $ 
 Then, by (\ref{1disconlemeq1}) and statement 
 \ref{1disconlem1}, we get that there exist compact sets 
 $B_{1}$ and $B_{2}$ in $L$  such that 
 $B_{1}\cap B_{2}=\emptyset ,\ B_{1}\cup B_{2}=L,\ 
L_{x}\subset (h_{1}^{k_{0}})^{-1}(L)\subset  B_{1},$ and 
$L_{x'}\subset h_{x_{1}'}^{-1}\cdots h_{x_{k_{0}}'}^{-1}(L)\subset B_{2}.$ 
Hence, statement \ref{1disconlem4} holds.
\end{proof}
We now demonstrate Theorem~\ref{countthm}.

\noindent {\bf Proof of Theorem~\ref{countthm}:} 
By Theorem~\ref{mainth2}-\ref{mainth2-2} or Remark~\ref{hatjcptrem}, we have 
$\hat{{\cal J}}={\cal J}.$
Let $J_{1}\in \hat{{\cal J}}$ be the element containing 
$J(h_{m}).$ By Theorem~\ref{mainth0}, we must have 
$J_{0}\neq J_{1}.$ Then, by Theorem~\ref{mainth1}-\ref{mainth1-1}, 
we have the following two possibilities.\\ 
Case 1. $J_{0}<J_{1}.$\\ 
Case 2. $J_{1}<J_{0}.$

 Suppose we have case 1. Then, by Proposition~\ref{bminprop}, 
 we have that $J_{0}=J_{\min }$ and $J_{1}=J_{\max }.$ 
 Combining it with the assumption and Theorem~\ref{mainth1}-\ref{mainth1-3}, 
 we obtain 
\begin{equation}
\label{cteq1}
\cup _{j=1}^{m-1}h_{j}^{-1}(J_{\max })\subset J_{\min }. 
\end{equation}
By (\ref{cteq1}) and Theorem~\ref{mainth1}-\ref{mainth1-3}, 
we get 
\begin{equation}
\label{cteq1-1}
\cup _{j=1}^{m-1}h_{j}^{-1}(J(G))\subset J_{\min }.
\end{equation}
Moreover, since $J(h_{m})\cap J_{\min }=\emptyset $, 
Theorem~\ref{mainth2}-\ref{mainth2-4-2} implies that 
\begin{equation}
\label{cteq2}
h_{m}^{-1}(J(G))\cap J_{\min }=\emptyset .
\end{equation}  
Then, by (\ref{cteq1-1}) and (\ref{cteq2}), we get 
\begin{equation}
\label{cteq3}
h_{m}^{-1}(J(G))\cap \left(\cup _{j=1}^{m-1}h_{j}^{-1}(J(G))\right)
=\emptyset .
\end{equation}
We now consider the backward self-similar set $J(G)$ 
with respect to $\{ h_{1},\ldots ,h_{m}\} .$ 
By Lemma~\ref{hmslem}-\ref{bss}, we have 
\begin{equation}
\label{cteqast1}
J(G)=\cup _{w\in \Sigma _{m}}(J(G))_{w}.
\end{equation}
Combining (\ref{cteq3}), Lemma~\ref{1disconlem}, 
Lemma~\ref{fibconnlem2}, and (\ref{cteqast1}), 
we obtain 
\begin{equation}
\label{cteq4}
J_{\max }=(J(G))_{m^{\infty }}\supset J(h_{m}),
\end{equation}
where we set $m^{\infty }:=(m,m,m,\ldots )\in \Sigma _{m}.$ 
Furthermore, by (\ref{cteq3}) and Lemma~\ref{1disconlem}, 
we get 
\begin{equation}
\label{cteq5}
\sharp (\hat{{\cal J}})\geq \aleph _{0}.
\end{equation}
Let $x=(x_{1},x_{2},\ldots )\in \Sigma _{m}$ be any element with $x\neq m^{\infty }$ 
and let $l:= \min \{ s\in \NN \mid x_{s}\neq m\} .$ 
Then, by (\ref{cteq1-1}), we have 
\begin{equation}
\label{cteq6}
(J(G))_{x}=\cap _{j=1}^{\infty }h_{x_{1}}^{-1}\cdots 
h_{x_{j}}^{-1}(J(G))\subset (h_{m}^{l-1})^{-1}(J_{\min }).
\end{equation}
Combining (\ref{cteqast1}) with (\ref{cteq4}) and (\ref{cteq6}),  
we obtain 
\begin{equation}
\label{cteq7}
J(G)=J_{\max }\cup \bigcup _{n\in \NN \cup \{ 0\} }
h_{m}^{-n}(J_{\min }).
\end{equation}
By (\ref{cteq5}) and (\ref{cteq7}), we get 
$\sharp (\hat{{\cal J}})=\aleph _{0}.$ 
Moreover, combining (\ref{cteq4}), (\ref{cteq7}), 
Theorem~\ref{mainth2}-\ref{mainth2ast1} and 
Theorem~\ref{mainth2}-\ref{mainth2-4-2}, we get that 
for each $J\in \hat{{\cal J}}$ with $J\neq J_{\max }$, 
there exists no sequence $\{ C_{j}\} _{j\in \NN }$ 
of mutually distinct elements of $\hat{{\cal J}}$ such that  
$\min _{z\in C_{j}}d(z,J)\rightarrow 0$ as $j\rightarrow \infty .$ 
 Furthermore, combining (\ref{cteq4}), Theorem~\ref{mainth2}-\ref{mainth2ast1} 
 and Theorem~\ref{mainth2}-\ref{mainth2-4-2}, 
we obtain $J_{\max }=(J(G))_{m^{\infty }}=J(h_{m}).$ 
Hence, all statements of Theorem~\ref{countthm} are true, provided that 
we have case 1.

 We now assume case 2: $J_{1}<J_{0}.$ 
Then, by Proposition~\ref{bminprop}, 
 we have that $J_{0}=J_{\max }$ and $J_{1}=J_{\min }.$ 
By the same method as that of case 1, we obtain  
\begin{equation}
\label{cteq8}
J_{\min }=(J(G))_{m^{\infty }}\supset J(h_{m}),
\end{equation}
\begin{equation}
\label{cteq9} 
J(G)=J_{\min }\cup \bigcup _{n\in \NN \cup \{ 0\} }
h_{m}^{-n}(J_{\max }),
\end{equation}
and 
\begin{equation}
\label{cteq10}
\sharp (\hat{{\cal J}})=\aleph _{0}. 
\end{equation}
Since $J(h_{j})\subset J_{0}$, for each $j=1,\ldots ,m-1$, 
and $J_{0}\neq J_{\min }$, 
Theorem~\ref{mainth2}-\ref{mainth2-4-2} 
implies that for each $j=1,\ldots ,m-1$, 
$h_{j}(J(h_{m}))\subset $ int$(K(h_{m})).$  
Hence, for each $j=1,\ldots ,m$, 
$h_{j}(K(h_{m}))\subset $ $K(h_{m}).$ 
Therefore, we have 
\begin{equation}
\label{cteq11}
 \mbox{int}(K(h_{m}))\subset F(G).
\end{equation}  
By (\ref{cteq11}) and (\ref{cteq8}), we obtain 
$J_{\min }=(J(G))_{m^{\infty }}=J(h_{m}).$ 
Moreover, by (\ref{cteq8}) and (\ref{cteq9}), 
we get that for each $J\in \hat{{\cal J}}$ with 
$J\neq J_{\min }$, there exists no sequence $\{ C_{j}\}_{j\in \NN }$ 
of mutually distinct elements of $\hat{{\cal J}}$ 
such that $\min _{z\in C_{j}}d(z,J)\rightarrow 0$ as $j\rightarrow \infty .$ 
Hence, we have shown Theorem~\ref{countthm}.     
\qed 

\ 

We now demonstrate Proposition~\ref{countprop}.

\noindent {\bf Proof of Proposition~\ref{countprop}:}
Let $0<\epsilon <\frac{1}{2}$ and let 
$\alpha _{1}(z):=z^{2},\alpha _{2}(z):=(z-\epsilon )^{2}+\epsilon , $ 
and $\alpha _{3}(z):=\frac{1}{2}z^{2}.$ 
If we take a large $l\in \NN $, then there exists an open 
neighborhood $U$ of $\{ 0,\epsilon \} $ with 
$U\subset \{ |z|<1\} $ and a neighborhood $V$ of 
$(\alpha _{1}^{l}, \alpha _{2}^{l},\alpha _{3}^{l})$ 
in (Poly)$^{3}$ such that for each $(h_{1},h_{2},h_{3})\in V$, 
we have 
$\cup _{j=1}^{3}h_{j}(U)\subset U$
and $\cup _{j=1}^{3}C(h_{j})\cap \CC \subset U$, where 
$C(h_{j})$ denotes the set of all critical points of $h_{j}.$ 
Then, by Remark~\ref{pcbrem}, for each $(h_{1},h_{2},h_{3})\in V$,  
$\langle h_{1},h_{2},h_{3}\rangle \in {\cal G}.$ 
Moreover, if we take an $l$ large enough and $V$ so small, 
then for each $(h_{1},h_{2},h_{3})\in V$, we have 
that: 
\begin{enumerate}
\item \label{cp0}
$J(h_{1})<J(h_{3})$;
\item \label{cp1}
$J(h_{1})\cup J(h_{2})$ is connected; 
\item \label{cp2}
$h_{i}^{-1}(J(h_{3}))\cap (J(h_{1})\cup J(h_{2}))\neq \emptyset ,$
for each $i=1,2$;   
\item \label{cp3}
$\cup _{j=1}^{3}h_{j}^{-1}(A)\subset A$, where 
$A=\{ z\in \CC \mid \frac{1}{2}\leq |z|\leq 3\}$; and 
\item \label{cp4}
$h_{3}^{-1}(A)\cap (\cup _{j=1}^{2}h_{i}^{-1}(A))=
\emptyset .$

\end{enumerate}
Combining statements \ref{cp3} and \ref{cp4} above,  
Lemma~\ref{hmslem}-\ref{backmin}, and Lemma~\ref{hmslem}-\ref{bss}, 
we get that for each $(h_{1},h_{2},h_{3})\in V$, 
$J(\langle h_{1},h_{2},h_{3}\rangle )\subset A$ and 
$J(\langle h_{1},h_{2},h_{3}\rangle )$ is disconnected.
Hence, for each $(h_{1},h_{2},h_{3})\in V$, we have 
$\langle h_{1},h_{2},h_{3}\rangle \in {\cal G}_{dis}.$ 
Combining it with statements \ref{cp1} and 
\ref{cp2} above and Theorem~\ref{countthm}, it follows that 
$J(h_{1})\cup J(h_{2})\subset J_{0}$ for some 
$J_{0}\in \hat{{\cal J}}_{\langle h_{1},h_{2},h_{3}\rangle }$, 
$h_{j}^{-1}(J(h_{3}))\cap J_{0}\neq \emptyset $ for each 
$j=1,2,$ 
and 
$\sharp (\hat{{\cal J}}_{\langle h_{1},h_{2},h_{3}\rangle })
=\aleph _{0}$, for each 
$(h_{1},h_{2},h_{3})\in V.$ 
Since $J(h_{1})<J(h_{3})$, 
Theorem~\ref{countthm} implies that 
the connected component $J_{0}$ should be equal to 
$J_{\min }(\langle h_{1},h_{2},h_{3}\rangle )$, 
and that $J_{\max }(\langle h_{1},h_{2},h_{3}\rangle )=J(h_{3}).$ 
 
 Thus, we have proved Proposition~\ref{countprop}.
\qed 

\ 

 We now show Proposition~\ref{countcomp}.

\noindent {\bf Proof of 
Proposition~\ref{countcomp}:} 
In fact, we show the following claim:\\ 
Claim: There exists a polynomial semigroup 
$G=\langle h_{1},h_{2},h_{3}\rangle $ 
in ${\cal G}$ such that 
all of the following hold.
\begin{enumerate}
\item $\sharp (\hat{{\cal J}})=\aleph _{0} .$
\item $J_{\min }\supset J(h_{1})\cup  
J(h_{2})$ and 
there exists a superattracting fixed point 
$z_{0}$ of $h_{1}$ with $z_{0}\in $ int$(J_{\min }).$ 
\item $J_{\max }=J(h_{3}).$
\item There exists a sequence $\{ n_{j}\} _{j\in \NN }$ of 
positive integers such that 
$\hat{{\cal J}}=
\{ J_{\min }\} \cup \{ J_{j}\mid j\in \NN \} $, 
where $J_{j}$ denotes the element of $\hat{{\cal J}}$ with 
$h_{3}^{-n_{j}}(J_{\min })\subset J_{j}.$ 
\item For any $J\in \hat{{\cal J}}$ with $J\neq J_{\max }$, 
there exists no sequence $\{ C_{j}\} _{j\in \NN }$ 
of mutually distinct elements of $\hat{{\cal J}}$ 
such that $\min _{z\in C_{j}}d(z,J)\rightarrow 0$ as 
$j\rightarrow \infty .$   
\item $G$ is sub-hyperbolic: i.e.,  
$\sharp (P(G)\cap J(G))<\infty  $ and 
$P(G)\cap F(G)$ is compact.
\end{enumerate} 
To show the claim, 
let $g_{1}(z)$ be the second iterate of $z\mapsto z^{2}-1.$ 
Let $g_{2}$ be a polynomial 
such that $J(g_{2})=\{ z\mid |z|=1\} $ and 
$g_{2}(-1)=-1.$ Then, we have 
$g_{1}(\sqrt{-1})=3\in \CCI \setminus K(g_{1}).$ 
Take a large, positive integer $m_{1}$, and 
let $a:=g_{1}^{m_{1}}(\sqrt{-1}).$ 
Then, 
\begin{equation}
\label{countcompeq1}
J(\langle g_{1}^{m_{1}},g_{2}\rangle )\subset 
\{ z\mid |z|<a\} .
\end{equation}  
Furthermore, 
since   
$a>\frac{1}{2}+\frac{\sqrt{5}}{2}$, we have   
\begin{equation}
\label{countcompeq1-5}
\overline{(g_{1}^{m_{1}})^{-1}(\{ z\mid |z|<a\} )}\subset 
\{ z\mid |z|<a\} .
\end{equation}
Let $g_{3}$ be a polynomial such that 
$J(g_{3})=\{ z\mid |z|=a\} .$ 
Since $-1$ is a superattracting fixed point 
of $g_{1}^{m_{1}}$ and it belongs to $J(g_{2})$, 
by \cite[Theorem 4.1]{HM2}, we see that 
for any $m\in \NN $, 
\begin{equation}
\label{countcompeq2}
-1\in \mbox{int}(J(\langle g_{1}^{m_{1}} , g_{2}^{m}\rangle )).
\end{equation}
Since $J(g_{2})\cap \mbox{int}(K(g_{1}^{m_{1}}))\neq \emptyset $ and 
$J(g_{2})\cap F_{\infty }(g_{1}^{m_{1}})\neq \emptyset $, 
we can take an $m_{2}\in \NN $ such that 
\begin{equation}
\label{countcompeq3}
(g_{2}^{m_{2}})^{-1}(\{ z\mid |z|=a\} )
\cap J(\langle g_{1}^{m_{1}},g_{2}^{m_{2}}\rangle )\neq \emptyset 
\end{equation} 
and 
\begin{equation}
\label{countcompeq4}
\overline{(g_{2}^{m_{2}})^{-1}(\{ z\mid |z|<a\} )}\subset 
\{ z\mid |z|<a\} .
\end{equation}
Take a small $r>0$ such that
\begin{equation}
\label{countcompeq5} 
\mbox{ for each }j=1,2,3,\ g_{j}(\{ z\mid |z|\leq r\})\subset 
\{ z\mid |z|<r\} . 
\end{equation}
Take an $m_{3}$ such that 
\begin{equation}
\label{countcompeq6}
(g_{3}^{m_{3}})^{-1}(\{ z\mid |z|=r\} )
\cap (\cup _{j=1}^{2}(g_{j}^{m_{j}})^{-1}(\{ z\mid |z|\leq a\} )) 
=\emptyset 
\end{equation}
and 
\begin{equation}
\label{countcompeq4-5}
g_{3}^{m_{3}}(-1)\in \{ z\mid |z|<r\} .
\end{equation}
Let $K:=\{ z\mid r\leq |z|\leq a\} .$ 
Then, by (\ref{countcompeq1-5}), (\ref{countcompeq4}), 
 (\ref{countcompeq5}) and (\ref{countcompeq6}), we have 
\begin{equation}
\label{countcompeq7}
(g_{j}^{m_{j}})^{-1}(K)\subset K, \mbox{ for } 
j=1,2,3, \mbox{and }
(g_{3}^{m_{3}})^{-1}(K)\cap 
(\cup _{j=1}^{2}(g_{j}^{m_{j}})^{-1}(K))= \emptyset .
\end{equation}     
Let $h_{j}:=g_{j}^{m_{j}}$, for each $j=1,2,3$, 
and let $G=\langle h_{1},h_{2},h_{3}\rangle .$ 
Then, by (\ref{countcompeq7}) and Lemma~\ref{hmslem}-\ref{backmin}, 
we obtain: 
\begin{equation}
\label{countcompeq8}
J(G)\subset K \mbox{ and }
h_{3}^{-1}(J(G))\cap (\cup _{j=1}^{2}
h_{j}^{-1}(J(G)))=\emptyset .
\end{equation}
Combining it with Lemma~\ref{hmslem}-\ref{bss}, 
it follows that 
$J(G)$ is disconnected. 
Furthermore, combining (\ref{countcompeq5}) and 
(\ref{countcompeq4-5}),  
we see $G\in {\cal G}$,
$P(G)\cap J(G)=\{ -1\} $, and 
that $P(G)\cap F(G)$ is compact. 
By Proposition~\ref{bminprop}, there exists 
a $j\in \{ 1,2,3\} $ with 
$J(h_{j})\subset J_{\min }.$ 
Since $ J(G)\subset K\subset \{ z\mid |z|\leq a\} $ and 
$J(h_{3})=\{ z\mid |z|=a\} $, 
we have 
\begin{equation}
\label{countcompeqmax}
J(h_{3})\subset J_{\max }.
\end{equation} 
Hence, either $J(h_{1})\subset J_{\min }$ or 
$J(h_{2})\subset J_{\min }.$ Since 
$J(h_{1})\cup J(h_{2})$ is connected, 
it follows that 
\begin{equation}
\label{countcompeqmin}
J(h_{1})\cup J(h_{2})\subset J_{\min }.
\end{equation} 
Combining this with Theorem~\ref{mainth1}-\ref{mainth1-3},
we have $h_{j}^{-1}(J_{\min })\subset J_{\min }$, for each 
$j=1,2.$ Hence, 
\begin{equation}
\label{countcompeq8-5}
J(\langle h_{1},h_{2}\rangle )\subset  
J_{\min }.
\end{equation}
Since $\sqrt{-1}\in J(h_{2})$ and $h_{1}(\sqrt{-1})=a\in J(h_{3})$,
 we obtain 
\begin{equation}
\label{countcompeqin}
h_{1}^{-1}(J(h_{3}))\cap J_{\min }\neq 
 \emptyset .
\end{equation} 
  
Similarly, by (\ref{countcompeq3}) and 
(\ref{countcompeq8-5}), 
we obtain  
\begin{equation}
\label{countcompeqin2}
h_{2}^{-1}(J(h_{3}))\cap J_{\min }\neq \emptyset .
\end{equation} 
Combining (\ref{countcompeqmax}), 
(\ref{countcompeqin}), (\ref{countcompeqin2}), and Theorem~\ref{countthm}, 
we obtain  
$\sharp (\hat{{\cal J}})=\aleph _{0}$, 
$J_{\max }=J(h_{3})$, 
$J(G)=J_{\max }\cup \bigcup _{n\in \NN \cup \{ 0\} }
h_{3}^{-n}(J_{\min })$, and 
that for any $J\in \hat{{\cal J}}$ with 
$J\neq J_{\max}$, 
there exists no sequence $\{C _{j}\} _{j\in \NN }$ 
of mutually distinct elements of 
$\hat{{\cal J}}$ such that 
$\min _{z\in C_{j}}d(z,J)\rightarrow 0$ as $j\rightarrow \infty .$ 
 
Moreover, by (\ref{countcompeq2}) and (\ref{countcompeq8-5})
(or by Theorem~\ref{mainupthm}-\ref{mainupthm3}), 
the superattracting fixed point $-1$ of $h_{1}$ belongs to 
int$(J_{\min }).$ 

 Hence, we have shown the claim.

Therefore, we have proved Proposition~\ref{countcomp}.
\qed

\subsection{Proofs of results in \ref{fibsec}}
\label{pffibsec}
In this section, we prove results in section \ref{fibsec}.

To prove results in \ref{fibsec}, we need the following notations and 
lemmas.
\begin{df}[\cite{S1}]
Let $f:X\times \CCI \rightarrow X\times \CCI $ be a rational skew 
product over $g:X\rightarrow X.$ Let $N\in \NN .$ 
We say that a point $(x_{0},y_{0})\in X\times \CCI $ belongs to 
$SH_{N}(f)$ if there exists a neighborhood $U$ of 
$x_{0}$ in $X$ and a positive number $\delta $ such that 
for any $x\in U$, any $n\in \NN $, any $x_{n}\in g^{-n}(x)$, 
and any connected component 
$V$ of $(f_{x_{n},n})^{-1}(B(y_{0},\delta ))$, 
$\deg (f_{x_{n},n}:V\rightarrow B(y_{0},\delta ))\leq N.$ 
Moreover, we set 
$UH(f):= (X\times \CCI )\setminus \cup _{N\in \NN }SH_{N}(f).$ 
We say that $f$ is semi-hyperbolic (along fibers) if 
$UH(f)\subset \tilde{F}(f).$ 
\end{df}
\begin{rem}
Under the above notation, we have $UH(f)\subset P(f).$ 
\end{rem}
\begin{rem}
\label{hypskewsemigrrem}
Let $\G $ be a compact subset of Rat and let 
$f:\GN \times \CCI \rightarrow \GN \times \CCI $ be the skew product 
associated with $\G .$  Let $G$ be the rational semigroup generated by 
$\G . $ Then, by Lemma~\ref{fiblem}-\ref{pic}, 
it is easy to see that $f$ is semi-hyperbolic if and only if 
$G$ is semi-hyperbolic. Similarly, it is easy to see that 
$f$ is hyperbolic if and only if $G$ is hyperbolic. 
\end{rem}

\begin{lem}
\label{nclimlem}
Let $f:X\times \CCI \rightarrow X\times \CCI $ be a 
polynomial skew product over $g:X\rightarrow X$ such that 
for each $\omega \in X$, $d(\omega )\geq 2.$ 
Let $x\in X$ be a point and  
 $y_{0}\in F_{x}(f)$ a point. 
 Suppose that there exists a strictly increasing sequence 
 $\{ n_{j}\} _{j\in \NN }$ of positive integers such that   
the sequence $\{ f_{x,n_{j}}\} _{j\in \NN }$  
converges to a non-constant map around $y_{0}$, and such that 
$\lim _{j\rightarrow \infty }f^{n_{j}}(x,y_{0})$ exists. 
We set $(x_{\infty }, y_{\infty }):=
\lim _{j\rightarrow \infty }f^{n_{j}}(x,y_{0}).$ 
Then, there exists a non-empty bounded open set $V$ in $\CC $
and 
a number $k\in \NN$ such that 
$\{ x_{\infty }\} \times \partial V\subset \tilde{J}(f)\cap UH(f)
\subset \tilde{J}(f)\cap P(f)$, 
and such that for each $j$ with $j\geq k$, 
$f_{x,n_{j}}(y_{0})\in V.$  
\end{lem}
\begin{proof}
We set 
$$V:= \{ y\in \CCI \mid 
\exists \epsilon >0, \lim_{i\rightarrow \infty }
\sup _{j>i}\sup _{d(\xi ,y)\leq \epsilon }
d(f_{g^{n_{i}}(x),n_{j}-n_{i}}(\xi ), \xi )=0\} .$$
Then, by \cite[Lemma 2.13]{S1}, 
 we have $\{ x_{\infty }\} \times \partial V
 \subset \tilde{J}(f)\cap UH(f)\subset 
 \tilde{J}(f)\cap P(f).$ 
 Moreover, since for each $x\in X$, 
 $f_{x,1}$ is a polynomial with $d(x)\geq 2$, 
Lemma~\ref{fibfundlem}-\ref{fibfundlem4} implies that 
there exists a ball $B$ around $\infty $ such that 
$B\subset \CCI \setminus V.$  

 From the assumption, 
 there exists a number $a>0$ and a non-constant map 
 $\varphi :D(y_{0},a)\rightarrow \CCI $ 
 such that $f_{x,n_{j}}\rightarrow \varphi $ 
 as $j\rightarrow \infty $,
 uniformly 
 on $D(y_{0},a).$   Hence,  
  $d(f_{x,n_{j}}(y),f_{x,n_{i}}(y))\rightarrow 0 $ as 
 $i,j\rightarrow \infty $, uniformly on $D(y_{0},a).$ 
 Moreover, since $\varphi $ is not constant, 
 there exists a positive number 
 $\epsilon $ such that 
 for each large $i$, $f_{x,n_{i}}(D(y_{0},a))\supset D(y_{\infty },\epsilon ).$ 
Therefore, it follows that 
$d(f_{g^{n_{i}}(x),n_{j}-n_{i}}(\xi ),\xi )\rightarrow 0$ as $i,j\rightarrow 
\infty $ uniformly on $D(y_{\infty },\epsilon ).$ 
Thus, $y_{\infty }\in V.$ Hence,  
there exists a number $k\in \NN $ such that 
for each $j\geq k$, $f_{x,n_{j}}(y_{0})\in V.$ 
Therefore, we have proved Lemma~\ref{nclimlem}.
\end{proof}
\begin{rem}
In \cite[Lemma 2.13]{S1} and \cite[Theorem 2.6]{S4}, 
the sequence $(n_{j})$ of positive integers should be strictly 
increasing.
\end{rem}
\begin{lem}
\label{ncintk}
Let $\G $ be a non-empty compact subset of 
{\em Poly}$_{\deg \geq 2}.$ 
Let $f:\GNCR $ be the skew product 
associated with $\G .$ 
Let $G$ be the polynomial semigroup generated by 
$\G .$ 
Let $\g \in \GN $ be a point. 
Let $y_{0}\in F_{\g }(f)$ and 
suppose that there exists a strictly increasing sequence 
$\{ n_{j}\} _{j\in \NN }$ of positive integers such that  
$\{ f_{\g ,n_{j}}\} _{j\in \NN }$ converges 
to a non-constant map around $y_{0}.$ 
Moreover, suppose that $G\in {\cal G}.$ 
Then, there exists a number $j\in \NN $ such that 
$f_{\g ,n_{j}}(y_{0})\in $ {\em int}$(\hat{K}(G)).$ 
\end{lem}
\begin{proof}
By Lemma~\ref{nclimlem}, 
there exists a bounded open set $V$ in $\CC $, 
a point $\g _{\infty }\in \GN $, and a number 
$j\in \NN $ such that 
$\{ \g _{\infty }\} \times \partial V\subset 
\tilde{J}(f)\cap P(f)$, and such that 
$f_{\g ,n_{j}}(y_{0})\in V.$  
Then, we have $\partial V\subset P^{\ast }(G).$ 
Since $g(P^{\ast }(G))\subset P^{\ast }(G)$ for each $g\in G$, 
the maximum principle implies that 
$V\subset $ int$(\hat{K}(G)).$ 
Hence, $f_{\g ,n_{j}}(y_{0})\in $ int$(\hat{K}(G)).$ 
Therefore, we have proved Lemma~\ref{ncintk}.
\end{proof}

\begin{lem}
\label{constlimlem}
Let $G$ be a  
polynomial semigroup generated by a compact subset 
$\G $ of {\em Poly}$_{\deg \geq 2}.$ 
If a sequence $\{ g_{n}\} _{n\in \NN }$ of elements of $G$ tends to a 
constant $w_{0}\in \CCI $ locally uniformly on a domain $V\subset \CCI $, 
then $w_{0}\in P(G).$   
\end{lem}
\begin{proof}
Since $\infty \in P(G)$, we may assume that 
$w_{0}\in \CC .$ 
Suppose $w_{0}\in \CC \setminus P(G).$ 
Then, there exists a $\delta >0$ such that 
$B(w_{0}, 2\delta )\subset \CC \setminus P(G).$ 
Let $z_{0}\in V$ be a point. Then, for each large $n\in \NN $, 
there exists a well-defined inverse branch 
$\alpha _{n}$
of $g_{n}^{-1}$ on $B(w_{0},2\delta )$ such that 
$\alpha _{n}(g_{n}(z_{0}))=z_{0}.$ 
Let $B:=B(w_{0},\delta ).$  
Since $\G $ is compact, there exists a connected component $F_{\infty }(G)$ of 
$F(G)$ containing $\infty .$ 
Let $C$ be a compact neighborhood of $\infty $ 
in $F_{\infty }(G).$   
Then, we must have that there exists a number $n_{0}$ such that  
$\alpha _{n}(B)\cap C=\emptyset $ for each $n\geq n_{0}$, 
since $g_{n}\rightarrow \infty 
$ uniformly on $C$ as $n\rightarrow \infty $, 
which follows from that 
$\deg (g_{n})\rightarrow \infty $ and 
local degree at $\infty $ of $g_{n}$ tends to 
$\infty $ as $n\rightarrow \infty .$ 
Hence, $\{ \alpha _{n} |_{B}\}_{n\geq n_{0}} $ is normal 
in $B.$ However, for a small 
$\epsilon $ so that $B(z_{0},2\epsilon )\subset V$, 
we have $g_{n}(B(z_{0},\epsilon ))\rightarrow w_{0}$ 
as $n\rightarrow \infty $, 
and this is a contradiction. Hence, we must have that 
$w_{0}\in P(G).$ 
\end{proof}
We now demonstrate Proposition~\ref{nonminnoncpt}-\ref{nonminnoncpt1},
\ref{nonminnoncpt}-\ref{nonminnoncpt2}, and 
\ref{nonminnoncpt}-\ref{nonminnoncpt3}. 
(Proposition~\ref{nonminnoncpt}-\ref{nonminnoncpt4} will be proved 
after Theorem~\ref{shshprop} is proved.)
\\ 
{\bf Proof of Proposition~\ref{nonminnoncpt}-\ref{nonminnoncpt1}, 
\ref{nonminnoncpt}-\ref{nonminnoncpt2}, and 
\ref{nonminnoncpt}-\ref{nonminnoncpt3} 
:} 
Since $\G \setminus \G _{\min }$ 
is not compact, 
there exists a sequence $\{ h_{j}\} _{j\in \NN }$ 
in $\G \setminus \G _{\min }$ and an element $h_{\infty }\in \G _{\min }$ 
such that $h_{j}\rightarrow h_{\infty }$ as $j\rightarrow \infty .$ 
By Theorem~\ref{mainth2}-\ref{mainth2-4-2}, 
for each $j\in \NN $, 
$h_{j}(K(h_{\infty }))$ is included in a connected component 
$U_{j}$ of int$(\hat{K}(G)).$ 
Let $z_{1}\in $ int$(\hat{K}(G))$ $(\subset $ int$(K(h_{\infty })))$ 
be a point. Then, 
$h_{\infty }(z_{1})\in $ int$(\hat{K}(G))$ and 
$h_{j}(z_{1})\rightarrow h_{\infty }(z_{1})$ as 
$j\rightarrow \infty .$  
 Hence, we may assume that there exists a connected 
 component $U$ of int$(\hat{K}(G))$ such that 
 for each $j\in \NN $, 
 $h_{j}(K(h_{\infty }))\subset U.$ 
 Therefore, $K(h_{\infty })=h_{\infty }(K(h_{\infty }))
 \subset \overline{U}.$ 
Since $\overline{U}\subset K(h_{\infty })$, 
we obtain $K(h_{\infty })=\overline{U}.$ 
Since $U\subset $ int$(K(h_{\infty }))\subset \overline{U}$ and 
$U$ is connected, it follows that 
int$(K(h_{\infty }))$ is connected. 
Moreover, we have 
$U\subset $ int$(K(h_{\infty }))\subset $ int$(\overline{U})\subset 
$ int$(\hat{K}(G)).$  
Thus, 
\begin{equation}
\label{nonminnoncptpfeq1}
\mbox{ int}(K(h_{\infty }))=U.
\end{equation} 
Furthermore, since
\begin{equation}
\label{nonminnoncptpfeq2}
J(h_{\infty })<J(h_{j}) \mbox{ for each } j\in \NN ,
\end{equation}
and $h_{j}\rightarrow h_{\infty }$ as $j\rightarrow \infty $, 
we obtain 
\begin{equation}
\label{nonminnoncptpfeq3}
J(h_{j})\rightarrow J(h_{\infty }) \mbox{ as } j\rightarrow \infty ,
\end{equation} 
with respect to the Hausdorff topology.
Combining that $h_{j}\in \G \setminus \G _{\min }$ for each $j\in \NN $ 
with Theorem~\ref{mainth2}-\ref{mainth2ast1}, 
(\ref{nonminnoncptpfeq1}), (\ref{nonminnoncptpfeq2}), 
and (\ref{nonminnoncptpfeq3}), we see that 
for each $h\in \G _{\min }$, $K(h)=K(h_{\infty }).$ 
Combining it with (\ref{nonminnoncptpfeq1}), 
(\ref{nonminnoncptpfeq2}) and (\ref{nonminnoncptpfeq3}),  
it follows that statements \ref{nonminnoncpt1} and 
\ref{nonminnoncpt2} in Proposition~\ref{nonminnoncpt} hold. 

 We now show that statement \ref{nonminnoncpt3} holds. 
Let $\g \in \GN $ and let $y_{0}\in $ int$(K_{\g }(f))$ be a point. 
Suppose that there exists a strictly increasing sequence $\{ n_{j}\} _{j\in \NN }$ 
of 
positive integers such that 
$f_{\g ,n_{j}}$ tends to a non-constant map as $j\rightarrow \infty $ 
around $y_{0}.$ 
Then, by Lemma~\ref{ncintk}, there exists a number $k\in \NN $ such that 
$f _{\g , n_{k}}(y_{0})\in $ int$(\hat{K}(G)).$ 
 Hence, the sequence 
 $\{ f_{\sigma ^{n_{k}}(\g ), n_{k+j}-n_{k}}\} _{j\in \NN }$ 
 converges to a non-constant map around $y_{1}:= 
f_{\g ,n_{k}}(y_{0})\in $ int$(\hat{K}(G)).$ 
However, combining Theorem~\ref{mainth2}-\ref{mainth2-4-2} and 
statements \ref{nonminnoncpt1} and \ref{nonminnoncpt2} in 
Proposition~\ref{nonminnoncpt}, 
we have that for each $h\in \G $, 
$h: $ int$(\hat{K}(G))\rightarrow $ int$(\hat{K}(G))$ is a contraction map 
with respect to the hyperbolic distance on int$(\hat{K}(G)).$ 
This causes a contradiction. 
 Therefore, statement \ref{nonminnoncpt3} in Proposition~\ref{nonminnoncpt} 
 holds. 

 Thus, we have proved Proposition~\ref{nonminnoncpt}-\ref{nonminnoncpt1}, 
 \ref{nonminnoncpt}-\ref{nonminnoncpt2}, and 
 \ref{nonminnoncpt}-\ref{nonminnoncpt3}. 
\qed 

\

We now demonstrate Theorem~\ref{mainth3}-\ref{mainth3-0} 
and Theorem~\ref{mainth3}-\ref{mainth3-1}.\\ 
{\bf Proof of Theorem~\ref{mainth3}-\ref{mainth3-0} and 
Theorem~\ref{mainth3}-\ref{mainth3-1}:}
First, we will show the following claim.\\ 
Claim 1. Let $\g \in R(\G ,\G \setminus \G _{\min }).$ Then, 
for any point $y_{0}\in F_{\g }(f)$, there exists 
no non-constant limit function of $\{ f_{\g ,n}\} _{n\in \NN }$ 
around $y_{0}.$ 

 To show this claim, by Proposition~\ref{nonminnoncpt}-\ref{nonminnoncpt3}, 
we may assume that $\G \setminus \G _{\min }$ is compact.  
 Suppose that there exists 
 a strictly increasing sequence $\{ n_{j}\} _{j\in \NN }$ of 
 positive integers
such that 
 $f_{\g ,n_{j}}$ tends to a non-constant map as $j\rightarrow \infty $ 
 around $y_{0}.$ 
By Lemma~\ref{ncintk}, there exists a 
number $k\in \NN $ such that 
$f_{\g ,n_{k}}(y_{0})\in $ int$(\hat{K}(G)).$ 
Hence, we get that the sequence 
$\{ f_{\sigma ^{n_{k}}(\g ), n_{k+j}-n_{k}}\} _{j\in \NN }$ converges to a non-constant map around 
the point $y_{1}:=f_{\g ,n_{k}}(y_{0})\in $ int$(\hat{K}(G)).$ 
However, since we are assuming that $\G \setminus \G _{\min }$ is 
compact, 
Theorem~\ref{mainth2}-\ref{mainth2-4-2} implies that  
$\cup _{h\in \G \setminus \G _{\min }}h(\hat{K}(G))$ is a compact 
subset of int$(\hat{K}(G))$, 
which implies 
that if we take the hyperbolic metric for each connected component 
of int$(\hat{K}(G))$, then there exists a constant $0<c<1$ such 
that for each $z\in $ int$(\hat{K}(G))$ and 
each $h\in \G \setminus \G _{\min }$, we have 
$\| h'(z)\| \leq c$, where $\| \cdot \| $ denotes 
the norm of the derivative measured from the 
hyperbolic metric on the connected component $W_{1}$ of int$(\hat{K}(G))$ 
containing $z$ to that of the connected component $W_{2}$ of int$(\hat{K}(G))$ 
containing $h(z).$ This causes a contradiction, since 
we have that 
$\g \in R(\G ,\G \setminus \G _{\min })$ 
and the sequence 
$\{ f_{\sigma ^{n_{k}}(\g ), n_{k+j}-n_{k}}\} _{j\in \NN }$ converges to a non-constant map around 
the point $y_{1}\in $ int$(\hat{K}(G)).$ Hence, we have shown Claim 1.

 Next, let $S$ be a non-empty 
 compact subset of $\G \setminus \G _{\min }$ and let 
 $\g \in R(\G ,S).$  
 We show the following claim.\\ 
Claim 2. For each point $y_{0}$ in each bounded component of $F_{\g }(f)$,  
there exists a number $n\in \NN $ such that 
$f_{\g ,n}(y_{0})\in $ int$(\hat{K}(G)).$ 

 To show this claim, suppose that there exists no $n\in \NN $ 
 such that $f_{\g ,n}(y_{0})\in $ int$(\hat{K}(G))$, and 
 we will deduce a contradiction. 
By Claim 1, $\{ f_{\g ,n}\} _{n\in \NN }$ has only constant limit functions 
around $y_{0}.$ Moreover, if a point $w_{0}\in \CC $ 
is a constant limit function of 
$\{ f_{\g ,n}\} _{n\in \NN }$, then by Lemma~\ref{constlimlem}, 
we must have $w_{0}\in P^{\ast }(G)\subset \hat{K}(G).$ Since 
we are assuming that there exists no $n\in \NN $ 
 such that $f_{\g ,n}(y_{0})\in $ int$(\hat{K}(G))$, 
 it follows that $w_{0}\in \partial \hat{K}(G).$ 
Combining it with Theorem~\ref{mainth2}-\ref{mainth2-3}, 
we obtain $w_{0}\in \partial \hat{K}(G)\subset J_{\min }.$  
From this argument, we get that 
\begin{equation}
\label{main3-1pfeq1}
d(f_{\g ,n}(y_{0}),J_{\min })\rightarrow 0, \mbox{ as }n\rightarrow \infty.
\end{equation}  
However, since  
$\g$  belongs to $R(\G ,S)$, 
the above (\ref{main3-1pfeq1}) implies that 
the sequence $\{ f_{\g ,n}(y_{0})\} _{n\in \NN }$ accumulates 
in the compact set 
$\cup _{h\in S}
h^{-1}(J_{\min})$, which is apart from 
$J_{\min }$, by Theorem~\ref{mainth2}-\ref{mainth2-4-2}. 
This contradicts (\ref{main3-1pfeq1}). Hence, we have shown that 
Claim 2 holds. 

  Next, we show the following claim.\\ 
Claim 3. There exists exactly one bounded 
component $U_{\g }$ of $F_{\g }(f).$ 

 To show this claim, we take an element $h\in \G_{\min }$ 
(Note that $\G_{\min }\neq \emptyset$, by Proposition~\ref{bminprop}).
We write the element $\g $ as $\g =(\g _{1},\g_{2},\ldots ).$ 
For any $l\in \NN $ with $l\geq 2$, 
let $s_{l}\in \NN $ be an integer 
with $s_{l}>l$ such that  
$\g_{s_{l}}\in S.$ 
We may assume that for each $l\in \NN $, 
$s_{l}<s_{l+1}.$ For each $l\in \NN $, 
let $\g ^{l}:= (\g _{1},\g _{2},\ldots ,\g _{s_{l}-1},h,h,h,\ldots )\in 
\GN $ and  
$\tilde{\g }^{l}:=\sigma ^{s_{l}-1}(\g )
=(\g _{s_{l}}, \g _{s_{l}+1},\ldots )
\in \GN .$ Moreover,  
let $\rho :=(h,h,h,\ldots )\in \GN .$ 
Since $h\in \G _{\min }$, 
we have  
\begin{equation}
\label{main3-1pfeq2}
J_{\rho }(f)=J(h)\subset 
J_{\min }.
\end{equation} 
Moreover, since $\g _{s_{l}}$ does not belong to 
$\G _{\min }$, combining it with Theorem~\ref{mainth2}-\ref{mainth2-4-2} 
we obtain  
$\g_{s_{l}}^{-1}(J(G))\cap J_{\min }=\emptyset .$ 
Hence, we have that for each $l\in \NN $, 
\begin{equation}
\label{main3-1pfeq3}
J_{\tilde{\g }^{l}}(f)=
\g _{s_{l}}^{-1}(J_{\sigma ^{s_{l}}(\g )}(f))
\subset  
\g_{s_{l}}^{-1}(J(G)) \subset \CCI \setminus J_{\min }.
\end{equation} 
Combining (\ref{main3-1pfeq2}),  (\ref{main3-1pfeq3}), 
and Lemma~\ref{fiborder}, 
we obtain 
\begin{equation}
\label{main3-1pfeq3a}
J_{\rho }(f)<J_{\tilde{\g }^{l}}(f), 
\end{equation}
which implies   
\begin{equation}
\label{main3-1pfeq4}
J_{\g ^{l}}(f)=
(f_{\g , s_{l}-1})^{-1}(J_{\rho }(f))
<(f_{\g ,s_{l}-1})^{-1}(J_{\tilde{\g }^{l}}(f))
=J_{\g }(f).
\end{equation} 
From Lemma~\ref{fiborder} and 
(\ref{main3-1pfeq4}), 
it follows that 
there exists a bounded component 
$U_{\g }$ of $F_{\g }(f)$ 
such that 
for each $l\in \NN $ with $l\geq 2$, 
\begin{equation}
\label{main3-1pfeq5}
J_{\g ^{l}}(f)\subset U_{\g }.
\end{equation}
We now suppose that there exists a bounded component $V$ of 
$F_{\g }(f)$ with $V\neq U_{\g }$, and we will deduce a contradiction.
Under the above assumption, we take a point $y\in V.$ 
Then, by Claim 2, we get that 
there exists a number $l\in \NN $ such that 
$f_{\g ,l}(y)\in $ int$(\hat{K}(G)).$ 
Since $s_{l}> l$, we obtain 
$f_{\g ,s_{l}-1}(y)\in $ int$(\hat{K}(G))\subset K(h)$, 
where, $h\in \G _{\min }$ is the element which we have taken before.  
By (\ref{main3-1pfeq3a}), 
we have that there exists a bounded component 
$B$ of $F_{\tilde{\g }^{l}}(f) $ 
containing $K(h).$ 
Hence, we have 
$f_{\g ,s_{l}-1}(y)\in B.$ 
Since the map $f_{\g ,s_{l}-1}:V\rightarrow B$ 
is surjective, it follows that 
$V\cap \left((f_{\g ,s_{l}-1})^{-1}(J(h))\right) 
\neq \emptyset .$ 
Combined with $ (f_{\g ,s_{l}-1})^{-1}(J(h))=
(f_{\g ^{l}, s_{l}-1})^{-1}(J(h))=
J_{\g ^{l}}(f)$, 
we obtain 
$V\cap J_{\g ^{l}}(f)   
 \neq \emptyset .$ 
 However, this causes a contradiction, 
since we have (\ref{main3-1pfeq5}) and $U_{\g }\cap V=\emptyset .$ 
Hence, we have shown Claim 3.

 Next, we show the following claim. \\ 
Claim 4. 
$\partial U_{\g }=\partial A_{\gamma }(f)= J_{\g }(f).$ 

 To show this claim, 
since $U_{\g }=\mbox{int}(K_{\g }(f))$, 
Lemma~\ref{fibfundlem}-\ref{fibfundlema} implies that 
$\partial U_{\g }=J_{\g}(f).$ Moreover, by Lemma~\ref{fibfundlem}-\ref{fibfundlem4}, we have 
$\partial A_{\gamma }(f)=J_{\gamma }(f).$  
Thus, we have shown Claim 4.

 We now show the following claim.\\ 
Claim 5. $\hat{J}_{\g }(f)=J_{\g }(f)$ and 
the map $\omega \mapsto J_{\omega }(f)$ is continuous at $\g $
 with respect to the Hausdorff topology in the space of non-empty 
 compact subsets of $\CCI .$  

 To show this claim, 
 suppose that there exists a point $z$ with 
 $z\in \hat{J}_{\g }(f)\setminus J_{\g }(f).$ 
Since $\hat{J}_{\g }(f)\setminus 
J_{\g }(f)$ is included in the union 
of bounded components of $F_{\g }(f)$, 
combining it with Claim 2,   
we get that there exists a number $n\in \NN $ such that 
$f_{\g ,n}(z)\in $ int$(\hat{K}(G))\subset F(G).$ 
However, since 
$z\in \hat{J}_{\g }(f)$, we must have 
that $f_{\g ,n}(z)=\pi _{\CCI }(f_{\g }^{n}(z))\in 
\pi _{\CCI }(\tilde{J}(f))=J(G).$ This is a contradiction. 
Hence, we obtain $\hat{J}_{\g }(f)=J_{\g }(f).$ 
Combining it with Lemma~\ref{fibfundlem}-\ref{fibfundlem2}, 
it follows that $\omega \mapsto J_{\omega }(f)$ is 
continuous at $\g .$   
Therefore, we have shown Claim 5. 

 Combining all Claims $1,\ldots ,5$, it follows   
 that statements \ref{mainth3-0}, \ref{mainth3-1-1}, \ref{mainth3-1-2}, 
 and \ref{mainth3-1-3} in Theorem~\ref{mainth3}
 hold.  

 We now show statement \ref{mainth3-1-4}. Let 
 $\g \in R(\G ,S)$ be an element. 
Suppose that $m_{2}(J_{\g }(f))>0$, where $m_{2}$ denotes the 
$2$-dimensional Lebesgue measure. Then, there exists a 
Lebesgue density point $b\in J_{\g }(f)$ so that 
\begin{equation}
\label{lebpteq}
\lim _{s\rightarrow 0}
\frac{m_{2}\left(D(b,s)\cap J_{\g }(f)\right)}
{m_{2}(D(b,s))}=1.
\end{equation} 
Since $\g $ belongs to $R(\G ,S)$, 
there exists an element $\g _{\infty }\in S$ and  
a sequence $\{ n_{j}\} _{j\in \NN }$ of positive integers  
such that $n_{j}\rightarrow \infty $ and 
$\g _{n_{j}}\rightarrow \g _{\infty }$ 
as $j\rightarrow \infty $,  
and such that for each $j\in \NN $, $\g _{n_{j}}\in S.$ 
We set $b_{j}:= f_{\g ,n_{j}-1}(b)$, for each $j\in \NN .$ 
We may assume that there exists a point $a\in \CC $ 
such that $b_{j}\rightarrow a$ as $j\rightarrow \infty .$ 
Since $\g _{n_{j}}(b_{j})=f_{\g ,n_{j}}(b)=
\pi _{\CCI }(f_{\g }^{n_{j}}(\g ,b))\in 
\pi _{\CCI }(\tilde{J}(f))=J(G)$, 
we obtain  
$a\in \g _{\infty }^{-1}(J(G)).$ Moreover, 
by Theorem~\ref{mainth2}-\ref{mainth2-4-2}, 
we obtain 
\begin{equation}
\label{mainth314pfeq0.5}
a\in \g _{\infty }^{-1}(J(G))\subset \CC \setminus J_{\min }.
\end{equation}
Combining it with Theorem ~\ref{mainth2}-\ref{mainth2-3}, 
it follows that 
\begin{equation}
\label{mainth314pfeq1}
r:=d_{e}(a,P(G))>0.
\end{equation}
Let $\epsilon $ be arbitrary number with $0<\epsilon <\frac{r}{10}.$
We may assume that for each $j\in \NN $, 
we have $b_{j}\in D(a,\frac{\epsilon }{2}).$  
 For each $j\in \NN ,$ let $\varphi _{j}$ be the well-defined inverse branch 
of $(f_{\g ,n_{j}-1})^{-1}$ on 
$D(a,r)$ such that $\varphi _{j}(b_{j})=b.$  
Let $V_{j}:= \varphi _{j}(D(b_{j},r-\epsilon ))$, for each 
$j\in \NN .$  
We now show the following claim.\\ 
Claim 6. diam $V_{j}\rightarrow 0$, 
as $\j\rightarrow \infty .$ 

 To show this claim, suppose that this is not true. Then, there exists a 
 strictly increasing sequence $\{ j_{k}\} _{k\in \NN }$ of  
 positive integers and a positive constant $\kappa $ 
 such that for each $k\in \NN $, diam $V_{j_{k}}\geq \kappa .$ 
 From Koebe distortion theorem, it follows that 
there exists a positive constant $c_{0}$ such that for each 
$k\in \NN $, 
$V_{j_{k}}\supset D(b,c_{0}).$ This implies that for each $k\in \NN $, 
$f_{\g ,v_{k}}(D(b,c_{0}))\subset 
D(b_{j_{k}},r-\epsilon )$, where $v_{k}:= n_{j_{k}}-1.$ 
Since $v_{k}\rightarrow \infty $ as $k\rightarrow \infty $ and 
$f_{\g ', n}|_{F_{\infty }(G)}\rightarrow \infty $ for any 
$\g '\in \GN $, 
it follows that 
for any $n\in \NN $, 
$f_{\g ,n}(D(b,c_{0}))\subset (\CCI \setminus F_{\infty }(G))$,  
which implies that $b\in F_{\g }(f).$ However,  
it contradicts $b\in J_{\g }(f).$ 
Hence, Claim 6 holds. 

 Combining Koebe distortion theorem and Claim 6, 
 we see that 
 there exist a constant $K>0$ and two sequences $\{ r_{j}\} _{j\in \NN }$ and 
 $\{ R_{j}\} _{j\in \NN }$ of positive numbers such that 
 $K\leq \frac{r_{j}}{R_{j}}<1$ and 
 $D(b,r_{j})\subset V_{j}\subset D(b,R_{j})$ for each $j\in \NN $, 
 and such that $R_{j}\rightarrow 0$ as $j\rightarrow \infty $. 
From (\ref{lebpteq}), it follows that 
\begin{equation}
\label{fatou0eq1}
\lim _{j\rightarrow \infty }
\frac{m_{2}\left(V_{j}\cap F_{\g }(f)\right)}
{m_{2}(V_{j})}=0.
\end{equation}   
For each $j\in \NN $, let 
$\psi _{j}: D(0,1)\rightarrow \varphi _{j}(D(a,r))$ be a biholomorphic 
map such that 
$\psi _{j}(0)=b.$ 
Then, there exists a constant $0<c_{1}<1$ such that 
for each $j\in \NN $, 
\begin{equation}
\label{psid0}
\psi _{j}^{-1}(V_{j})\subset 
D(0,c_{1}).
\end{equation}
Combining it with (\ref{fatou0eq1}) and Koebe distortion theorem, 
it follows that  
\begin{equation}
\label{fatou0eq2}
\lim _{j\rightarrow \infty }
\frac{m_{2}\left(\psi _{j}^{-1}(V_{j}\cap F_{\g }(f))\right)}
{m_{2}(\psi _{j}^{-1}(V_{j}))}=0.
\end{equation}   
Since $\varphi _{j}^{-1}(\psi _{j}(D(0,1))\subset 
D(a,r)$ for each $j\in \NN $, 
combining (\ref{psid0}) and 
Cauchy's formula yields that 
there exists a constant $c_{2}>0$ such that for any $j\in \NN $,  
\begin{equation}
\label{fatou0eq3}
|(f_{\g ,n_{j}-1}\circ \psi _{j})'(z)|\leq c_{2}\ \mbox{ on }
\psi _{j}^{-1}(V_{j}).
\end{equation}     
Combining (\ref{fatou0eq2}) and (\ref{fatou0eq3}), 
we obtain 
\begin{align*}
 \frac{m_{2}\left(D(b_{j},r-\epsilon )\cap 
F_{\sigma ^{n_{j}-1}(\g )}(f)\right)}
{m_{2}(D(b_{j},r-\epsilon ))} 
 =  \frac{m_{2}\left((f_{\g ,n_{j}-1}\circ \psi _{j})
(\psi _{j}^{-1}(V_{j}\cap 
F_{\g }(f)))\right)}
{m_{2}(D(b_{j},r-\epsilon ))}\\ 
 = 
\frac{\int _{\psi _{j}^{-1}(V_{j}\cap F_{\g }(f))}
|(f_{\g ,n_{j}-1}\circ \psi _{j})'(z)|^{2}\ dm_{2}(z)}
{m_{2}(\psi _{j}^{-1}(V_{j}))}
\cdot 
\frac{m_{2}(\psi _{j}^{-1}(V_{j}))}
{m_{2}(D(b_{j},r-\epsilon ))} 
 \rightarrow 0,
\end{align*}
as $j\rightarrow \infty .$ 
Hence, we obtain 
$$\lim _{j\rightarrow \infty }
\frac{m_{2}\left(D(b_{j},r-\epsilon )\cap 
J_{\sigma ^{n_{j}-1}(\g )}(f)\right)}
{m_{2}(D(b_{j},r-\epsilon ))}=1.$$
Since  $J_{\sigma ^{n_{j}-1}(\g )}(f)\subset J(G)$ 
for each $j\in \NN $, and $b_{j}\rightarrow a$ as $j\rightarrow \infty $, 
it follows that 
$$\frac{m_{2}(D(a,r-\epsilon )\cap J(G))}
{m_{2}(D(a,r-\epsilon ))}=1.$$ 
This implies that 
$D(a,r-\epsilon )\subset J(G).$ 
Since this is valid for any $\epsilon $, 
we must have that $D(a,r)\subset J(G).$ It follows that  
the point $a$ belongs to a connected component 
$J$ of $J(G)$ such that $J\cap P^{\ast }(G)\neq \emptyset .$ 
However, Theorem~\ref{mainth2}-\ref{mainth2-3} implies 
that the component $J$ is equal to $J_{\min }$, 
which causes a contradiction since we have  
(\ref{mainth314pfeq0.5}).
 Hence, we have shown statement \ref{mainth3-1-4} in 
 Theorem~\ref{mainth3}-\ref{mainth3-1}.

 Therefore, we have proved Theorem~\ref{mainth3}-\ref{mainth3-0} 
 and Theorem~\ref{mainth3}-\ref{mainth3-1}. 
\qed 

\ 
\ 

In order to demonstrate Theorem~\ref{mainth3}-\ref{mainth3-3}, 
we need the following result.
\begin{thm} 
\label{hypskewqc}
({\bf Uniform fiberwise quasiconformal surgery})
Let $f:X\times \CCI \rightarrow X\times \CCI $ be a 
polynomial skew product over $g:X\rightarrow X$ such that 
for each $x\in X$, $d(x)\geq 2.$ 
Suppose that $f$ is hyperbolic and that 
$\pi _{\CCI }(P(f))\setminus \{ \infty \} $ is bounded in 
$\CC .$ Moreover, suppose that 
for each $x\in X$,\ 
{\em int}$(K_{x}(f))$ is connected. 
Then, there exists a constant $K$ such that 
for each $x\in X,\ J_{x}(f)$ is a $K$-quasicircle.   
\end{thm}
\begin{proof}
Step 1: 
By \cite[Theorem 2.14-(4)]{S1}, the map $x\mapsto J_{x}(f)$ is 
continuous with respect to the Hausdorff topology. 
Hence, there exists a positive constant 
$C_{1}$ such that 
for each 
$x\in X,\ \inf \{ d(a,b)\mid a\in J^{x}(f),\ b\in \pi ^{-1}(\{ x\} )\cap P^{\ast }(f)\} >C_{1}$,  
where $P^{\ast }(f):= P(f)\setminus \pi _{\CCI }^{-1}(\{ \infty \} )$, 
and $d(\cdot ,\cdot )$ denotes the spherical distance, under the 
canonical identification $\pi ^{-1}(\{ x\} )\cong \CCI .$   
Moreover, from the assumption, we have that 
for each $x\in X$, int$(K_{x}(f))\neq \emptyset .$ 
Since $X$ is compact, 
it follows that for each $x\in X,$  
there exists an analytic 
Jordan curve 
$\zeta _{x}$ in 
$K^{x}(f)\cap F^{x}(f)$ 
such that: 
\begin{enumerate}
\item 
$\pi ^{-1}(\{ x\} )\cap P^{\ast }(f)$ is included in the 
bounded component $V_{x}$ of $\pi ^{-1}(\{ x\} )\setminus 
\zeta _{x}$;
\item $\inf _{z\in \zeta _{x}}
d(z,\ J^{x}(f)\cup (\pi ^{-1}(\{ x\} )\cap P^{\ast }(f)))
\geq C_{2}$, where $C_{2}$ is a positive constant 
independent of $x\in X$; and  
\item there exist finitely many Jordan curves 
$\xi _{1},\ldots ,\xi _{k}$ in $\CC $ such that 
for each $x\in X$, there exists a $j$ with 
$\pi _{\CCI }(\zeta _{x})=\xi _{j}.$ 
\end{enumerate}
Step 2: 
By 
\cite[Corollary 2.7]{S4}, 
there exists an $n\in \NN $ such that 
for each $x\in X,\ W_{x}:=(f_{x}^{n})^{-1}(V_{g ^{n}(x)})
\supset \overline{V_{x}}$, 
$\inf \{ d(a,b)\mid a\in \partial W_{x},b\in \partial V_{x}, x\in X\} >0$, 
 and 
mod $(W_{x}\setminus \overline{V_{x}})\geq C_{3},$ where 
$C_{3}$ is a positive constant independent of $x\in X.$ 
In order to prove the theorem, 
since $J_{x}(f^{n})=J_{x}(f)$ for each $x\in X$, 
replacing $f:X\times \CCI \rightarrow X\times \CCI $ by 
$f^{n}:X\times \CCI \rightarrow X\times \CCI $, 
we may assume $n=1.$ 

\noindent Step 3: For each $x\in X$, let 
$\varphi _{x}:\pi ^{-1}(\{ x\} )\setminus 
\overline{V_{x}}\rightarrow 
\pi ^{-1}(\{ x\} )\setminus 
\overline{D(0,\frac{1}{2})}$ be 
a biholomorphic map such that 
$\varphi _{x}(x,\infty )=(x,\infty ),$ 
under the canonical identification $\pi ^{-1}(\{ x\} )\cong \CCI .$ 
We see that $\varphi _{x}$ extends 
analytically over $\partial V_{x}=\zeta _{x}.$ 
For each $x\in X$, we define a quasi-regular map 
$h_{x}: \pi ^{-1}(\{ x\} )\cong \CCI \rightarrow 
\pi ^{-1}(\{ g (x)\} )\cong \CCI $ as 
follows:
$$h_{x}(z):=\begin{cases}
            \varphi _{g (x)}f_{x}\varphi _{x}^{-1}(z),\ 
           \mbox{ if } z\in \varphi _{x}(\pi ^{-1}(\{ x\} )\setminus W_{x}),\\ 
            z^{d(x)},\ \mbox{ if } z\in \overline{D(0,\frac{1}{2})},\\ 
            \tilde{h}_{x}(z),\ \mbox{\ if } z\in \varphi _{x}(W_{x}\setminus 
            \overline{V_{x}}),
            \end{cases}
$$
where 
$\tilde{h}_{x}:\varphi _{x}(W_{x}\setminus \overline{V_{x}}) 
\rightarrow 
D(0,\frac{1}{2})\setminus 
\overline{D(0,(\frac{1}{2})^{d(x)})}$ is a 
regular covering and 
a $K_{0}$-quasiregular map with dilatation constant 
$K_{0}$ independent of $x\in X.$ 

\noindent Step 4: 
For each $x\in X$, we define a 
Beltrami differential 
$\mu _{x}(z)\frac{d\overline{z}}{dz}$ on 
$\pi ^{-1}(\{ x\} )\cong \CCI $ as follows: 
$$\begin{cases}
   \frac{\partial _{\overline{z}}\tilde{h}_{x}}
        {\partial _{z}\tilde{h}_{x}}
  \frac{d\overline{z}}{dz},\   
  \mbox{ if }z\in \varphi _{x}(W_{x}\setminus 
  \overline{V_{x}}),\\
  (h_{g ^{m}(x)}\cdots h_{x})^{\ast }
  (\frac{\partial _{\overline{z}}\tilde{h}_{g ^{m}(x)}}
        {\partial _{z}\tilde{h}_{g ^{m}(x)}}
  \frac{d\overline{z}}{dz}),\   
  \mbox{ if } z\in (h_{g ^{m}(x)}\cdots h_{x})^{-1}
   (\varphi _{g^{m}(x)}(W_{g ^{m}(x)}\setminus 
   \overline{V_{g^{m}(x)}})),\\ 
   
  0,\ \ \ \mbox{otherwise.}     
  \end{cases}
$$
Then, there exists a constant $k$ with $0<k<1$ such that 
for each $x\in X,\ \| \mu _{x}\| _{\infty }\leq k.$ 
By the construction, we have 
$h_{x}^{\ast }(\mu _{g (x)}\frac{d\overline{z}}{dz})
=\mu _{x}\frac{d\overline{z}}{dz}$, for each $x\in X.$ 
By the measurable Riemann mapping theorem (\cite[page 194]{LV}), 
for each $x\in X$, 
there exists a quasiconformal map $\psi _{x}:\pi ^{-1}(\{ x\} )\rightarrow 
\pi ^{-1}(\{ x\} )$  
such that 
$\partial _{\overline{z}}\psi _{x}=\mu _{x}\partial _{z}\psi _{x},\ 
\psi _{x}(0)=0,\ \psi _{x}(1)=1, $ and $\psi _{x}(\infty )=\infty $, 
under the canonical identification $\pi ^{-1}(\{ x\} )\cong \CCI .$ 
For each $x\in X$, let 
$\hat{h}_{x}:=\psi _{g (x)}h_{x}\psi _{x}^{-1}
:\pi ^{-1}(\{ x\} )\rightarrow \pi ^{-1}(\{ g (x)\} ).$ 
Then, $\hat{h}_{x}$ is holomorphic on 
$\pi ^{-1}(\{ x\} ).$  By the construction, 
we see that 
$\hat{h}_{x}(z)=c(x)z^{d(x)}$, 
where $c(x)=\psi _{g (x)}h_{x}\psi _{x}^{-1}(1)
=\psi _{g(x)}h_{x}(1).$ 
Moreover, by the construction again, we see 
that there exists a positive constant $C_{4}$ such 
that for each $x\in X,\ \frac{1}{C_{4}}\leq |h_{x}(1)|\leq C_{4}.$ 
Furthermore, 
\cite[Theorem 5.1 in page 73]{LV} implies that 
under the canonical identification 
$\pi ^{-1}(\{ x\} )\cong \CCI $, 
the family $\{ \psi _{x}^{-1}\} _{x\in X}$ is normal in 
$\CCI .$ Therefore, it follows that 
there exists a positive constant $C_{5}$ such that 
for each $x\in X,\ \frac{1}{C_{5}}\leq |c(x)|\leq 
C_{5}.$ 
Let $\tilde{J}_{x}$ be the set of 
non-normality of the sequence 
$\{ \hat{h} _{g ^{m}(x)}\cdots \hat{h}_{x}\} _{m\in \NN }$ 
in $\pi ^{-1}(\{ x\} )\cong \CCI .$ 
Since $\hat{h}_{x}(z)=c(x)z^{d(x)}$ and 
$\frac{1}{C_{5}}\leq |c(x)|\leq C_{5}$ for each 
$x\in X$, we get that for each $x\in X$, 
$\tilde{J}_{x}$ is a round circle. 
Moreover, \cite[Theorem 5.1 in page 73]{LV} implies that 
 $\{ \psi _{x}\} _{x\in X}$ 
and $\{ \psi _{x}^{-1}\} _{x\in X}$ are normal in $\CCI $ 
(under the canonical identification $\pi ^{-1}(\{ x\} )\cong \CCI $). 
Combining it with \cite[Corollary 2.7]{S4},  
we see that for each $x\in X$,   
$J^{x}(f)=\varphi _{x}^{-1}(\psi _{x}^{-1}(\tilde{J}_{x}))$, 
and it follows that 
there exists a constant $K$ such that 
for each $x\in X,\ J_{x}(f)$ is a $K$-quasicircle. 

 Thus, we have proved Theorem~\ref{hypskewqc}.
\end{proof}
\begin{rem}
Theorem~\ref{hypskewqc} generalizes a result 
in \cite[TH\'{E}OR\`{E}ME 5.2]{Se}, 
where O. Sester investigated hyperbolic polynomial skew products 
$f:X\times \CCI \rightarrow X\times \CCI $ 
such that for each $x\in X$, $d(x)=2.$   
\end{rem}
We now demonstrate Theorem~\ref{mainth3}-\ref{mainth3-3}.\\ 
\noindent 
{\bf Proof of Theorem~\ref{mainth3}-\ref{mainth3-3}:}
First, we remark that the subset $W_{S,p}$ of 
$\GN $ is a $\sigma $-invariant compact set.
 Hence, $\overline{f}:
W_{S,p}\times \CCI \rightarrow W_{S,p}\times \CCI $ 
is a polynomial skew product over 
$\sigma :W_{S,p}\rightarrow W_{S,p}.$  
Suppose that $\tilde{J}(\overline{f})\cap 
P(\overline{f})\neq \emptyset $ 
and let $(\g ,y)\in \tilde{J}(\overline{f})\cap 
P(\overline{f})$ be a point. 
Then, since the point $\g =(\g _{1},\g _{2},\ldots )$ belongs to 
$W_{S,p}$, 
there exists a number $j\in \NN $ such that 
$\g _{j}\in S.$ 
Combining it with Theorem~\ref{mainth2}-\ref{mainth2-4-2} and 
Theorem~\ref{mainth2}-\ref{mainth2-3}, 
we have $\g_{j}^{-1}(J(G))\subset \CC \setminus 
\hat{K}(G) \subset \CC \setminus P(G).$ 
Moreover, we have that 
$\pi _{\CCI }(\overline{f}_{\g }^{j-1}(\g ,y))
= \pi _{\CCI }(f_{\g }^{j-1}(\g ,y))
\in J_{\sigma ^{j-1}(\g )}(f)
= \g _{j}^{-1}\left(J_{\sigma ^{j}(\g )}(f)\right)
\subset \g_{j}^{-1}(J(G)).$ 
Hence, we obtain 
\begin{equation}
\label{mainth3-3pfeq1}
\pi _{\CCI }(\overline{f}_{\g }^{j-1}(\g ,y))\in \CC \setminus 
P(G).
\end{equation} 
However, since $(\g ,y)\in P(\overline{f})$, 
we have that $\pi _{\CCI }(\overline{f}_{\g }^{j-1}(\g ,y))
\in \pi _{\CCI }(P(\overline{f}))
\subset P(G)$, which contradicts (\ref{mainth3-3pfeq1}). 
Hence, we must have that 
$\tilde{J}(\overline{f})\cap P(\overline{f})=
\emptyset .$ Therefore, 
$\overline{f}:W_{S,p}\times \CCI \rightarrow W_{S,p}\times \CCI $ 
is a hyperbolic polynomial skew product over 
the shift map $\sigma :W_{S,p}\rightarrow W_{S,p}.$

Combining this with Theorem~\ref{mainth3}-\ref{mainth3-1-1} and 
Theorem~\ref{hypskewqc}, we conclude that 
there exists a constant $K_{S,p}\geq 1$ such that 
for each $\g \in W_{S,p},\ J_{\g }(\overline{f})$ is 
a $K_{S,p}$-quasicircle.
Moreover, by Theorem~\ref{mainth3}-\ref{mainth3-1-3}, 
we have $J_{\g }(\overline{f})=J_{\g }(f)=\hat{J}_{\g }(f).$ 

Hence, we have shown Theorem~\ref{mainth3}-\ref{mainth3-3}. 
\qed \\ 

To demonstrate Theorem~\ref{shshprop}, 
we need the following.

\begin{lem}
\label{ksetmin}
Let $G$ be a  
polynomial semigroup generated by 
a non-empty compact set $\G $ in 
{\em Poly}$_{\deg \geq 2}.$ 
Suppose that $G\in {\cal G}_{dis}.$ 
Then, we have $\hat{K}(G_{\min ,\G })=\hat{K}(G).$ 
\end{lem}
\begin{proof}

Since $G_{\min ,\G }\subset G$, 
we have $\hat{K}(G)\subset \hat{K}(G_{\min ,\G }).$ 
Moreover, it is easy to see 
$\hat{K}(G_{\min ,\G })=\cap _{g\in G_{\min ,\G }}K(g).$ 
Let $g\in G_{\min ,\G }$ and 
$h\in \G \setminus \G _{\min }.$ 
For each $\alpha \in \G _{\min }$, 
we have $\alpha ^{-1}(J_{\min }(G))\subset J_{\min }(G).$ 
Since $\sharp (J_{\min }(G))\geq 3 $ (Theorem~\ref{mainth2}-\ref{mainth2-4-1}), 
Lemma~\ref{hmslem}-\ref{backmin} implies that 
$J(g)\subset J_{\min }(G).$ Hence, 
from Theorem~\ref{mainth2}-\ref{mainth2-4-2}, 
it follows that 
\begin{equation}
\label{ksetminpfeq1}
h(J(g))\subset \mbox{int}(\hat{K}(G))\subset \mbox{int}(\hat{K}(g)).
\end{equation}
Since $J(g)$ is connected and each connected component of int$(K(g))$ 
is simply connected, the above (\ref{ksetminpfeq1}) implies that 
$h(K(g))\subset K(g).$ 
Hence, we obtain 
$h(\hat{K}(G_{\min ,\G }))=
h(\cap _{g\in G_{\min ,\G }}K(g))
\subset \cap _{g\in G_{\min ,\G }}K(g)
=\hat{K}(G_{\min ,\G }).$
Combined with that $\alpha (\hat{K}(G_{\min ,\G }))\subset 
\hat{K}(G_{\min ,\G })$ for each $\alpha \in \G_{\min }$, 
it follows that for each $\beta \in G$, 
$\beta (\hat{K}(G_{\min ,\G }))\subset \hat{K}(G_{\min ,\G }).$ 
Therefore, we obtain $\hat{K}(G_{\min ,\G })\subset 
\hat{K}(G).$ Thus, it follows that  
 $\hat{K}(G_{\min ,\G })=\hat{K}(G).$
\end{proof}
\begin{df}
Let $G$ be a rational semigroup and $N$ a positive integer. 
We denote by $SH_{N}(G)$ the set of points $z\in \CCI $ 
satisfying that  
there exists a positive number $\delta $  
such that for each $g\in G$, 
$\deg (g:V\rightarrow B(z,\delta ))\leq N$, 
for each connected component $V$ of $g^{-1}(B(z,\delta )).$ 
Moreover, 
we set $UH(G):= \CCI \setminus \cup _{N\in \NN }SH_{N}(G).$ 

\end{df}
\begin{lem}
\label{060523lem1}
Let $G$ be a polynomial semigroup generated by a compact subset 
$\G $ of {\em Poly}$_{\deg \geq 2}.$ 
Suppose that $G\in {\cal G}_{dis }$ and that 
$\G \setminus \G _{\min }$ is not compact. 
Moreover, suppose that (a) in 
Proposition~\ref{nonminnoncpt}-\ref{nonminnoncpt2} holds.
Then, there exists an open neighborhood ${\cal U} $ of 
$\G _{\min }$ in $\G $ and an open set $U$ in {\em int}$(\hat{K}(G))$ 
with $\overline{U}\subset $ {\em int}$(\hat{K}(G))$ such that:
\begin{enumerate}
\item $\cup _{h\in {\cal U}}h(U)\subset U$; 
\item $\cup _{h\in {\cal U}}CV^{\ast }(h) \subset U$, 
and 
\item 
denoting by $H$ the polynomial semigroup generated by ${\cal U}$, 
we have that $P^{\ast }(H)\subset $ {\em int}$(\hat{K}(G))\subset F(H)$
and that $H$ is hyperbolic. 
\end{enumerate}  
\end{lem}
\begin{proof}
Let $h_{0}\in \G _{\min }$ be an element.  
Let 
${\cal E}:= 
\{ \psi (z)=az+b\mid a,b\in \CC, |a|=1, \psi (J(h_{0}))=J(h_{0})\} .$ 
Then, by \cite{Be}, ${\cal E}$ is compact in Poly. 
Moreover, by \cite{Be}, we have the following two claims:\\ 
Claim 1: If $J(h_{0})$ is a round circle with the center $b_{0}$ and 
radius $r$, then 
${\cal E}=\{ \psi (z)=a(z-b_{0})+b_{0}\mid |a|=r\} .$ \\ 
Claim 2: If $J(h_{0})$ is not a round circle, then 
$\sharp {\cal E}<\infty .$ 

 Let $z_{0}$ be the unique attracting fixed point of $h_{0}$ in $\CC .$ 
 Let $g\in G_{\min ,\G }.$ 
 By \cite{Be}, for each $n\in \NN $, 
 there exists an $\psi _{n}\in {\cal E}$ such that 
 $h_{0}^{n}g=\psi _{n}gh_{0}^{n}.$ Hence, 
 for each $n\in \NN $, 
 $h_{0}^{n}g(z_{0})=\psi _{n}gh_{0}^{n}(z_{0})=\psi _{n}g(z_{0}).$ 
 Combining it with Claim 1 and Claim 2, 
 it follows that there exists an $n\in \NN $ such that 
 $h_{0}^{n}(g(z_{0}))=z_{0}.$ For this $n$, 
 $g(z_{0})=\psi _{n}^{-1}(h_{0}^{n}(g(z_{0})))
 =\psi _{n}^{-1}(z_{0})\in 
 \cup _{\psi \in {\cal E}}\psi (z_{0}).$ 
 Combining it with Claim 1 and Claim 2 again, 
 we see that the set $C:= \overline{\cup _{g\in G_{\min ,\G}}\{ g(z_{0})\} } $ 
 is a compact subset of int$(\hat{K}(G)).$ 
Let $d_{H}$ be the hyperbolic distance on int$(\hat{K}(G)).$ 
Let $R>0$ be a large number such that 
setting $U:=\{ z\in \mbox{int}(\hat{K}(G))\mid 
\min _{a\in C}d_{H}(z,a)<R\} $, 
we have $\cup _{h\in \G _{\min }}CV^{\ast }(h)\subset U .$ 
 Then, for each $h\in \G _{\min }$, 
 $\overline{h(U)}\subset U.$ 
 Therefore, there exists an open neighborhood ${\cal U}$ of 
 $\G _{\min }$ in $\G $ such that 
 $\cup _{h\in {\cal U}}h(U)\subset U$, and such that 
$\cup _{h\in {\cal U}}CV^{\ast }(h) \subset U.$ 
Let $H$ be the polynomial semigroup generated by ${\cal U}.$ 
From the above argument, we obtain 
$P^{\ast }(H)=\overline{\cup _{g\in H}CV^{\ast }(g) }
\subset $ $\overline{\cup _{g\in H\cup \{ Id\} }g\left( \cup _{h\in {\cal U}}CV^{\ast }(h) \right) }   
\subset \overline{\cup _{g\in H\cup \{ Id\} }g(U)}\subset \overline{U}\subset 
$ int$(\hat{K}(G))\subset F(H).$ 
Hence, $H$ is hyperbolic. Thus, we have proved 
Lemma~\ref{060523lem1}. 
\end{proof}

We now demonstrate Theorem~\ref{shshprop}.\\ 
{\bf Proof of Theorem~\ref{shshprop}:}
Suppose that $G_{\min ,\G }$ is semi-hyperbolic.
We will consider the following two cases:\\ 
Case 1: $\G \setminus \G _{\min }$ is compact.\\ 
Case 2: $\G \setminus \G _{\min }$ is not compact.

 Suppose that we have Case 1.
Since $UH(G_{\min ,\G })\subset P(G_{\min ,\G })$, 
$G_{\min ,\G }\in {\cal G}$, and $G_{\min ,\G }$ is semi-hyperbolic, 
we obtain 
$UH(G_{\min ,\G })\cap \CC 
\subset F(G _{\min ,\G })\cap \hat{K}(G_{\min ,\G })$ 
=int$(\hat{K}(G_{\min ,\G })).$  
By Lemma~\ref{ksetmin}, we have 
$\hat{K}(G_{\min ,\G })=\hat{K}(G).$ 
Hence, we obtain 
\begin{equation}
\label{mainth3-2pfeq1}
UH(G_{\min ,\G })\cap \CC \subset  
\mbox{ int}(\hat{K}(G))\subset \CC \setminus J_{\min }(G).
\end{equation}
Therefore, there exists a positive integer $N$ and a 
positive number $\delta $ 
such that for each $z\in J_{\min }(G)$ and 
each $h\in G_{\min ,\G }$, 
we have 
\begin{equation}
\label{mainth3-2pfeq2}
\deg (h:V\rightarrow D(z,\delta ))\leq N,  
\end{equation}
for each connected component $V$ of $h^{-1}(D(z,\delta )).$ 
Moreover, combining Theorem~\ref{mainth2}-\ref{mainth2-4-2} and 
Theorem~\ref{mainth2}-\ref{mainth2-3}, we obtain  
$\cup _{\alpha \in \G \setminus \G_{\min }}
\alpha ^{-1}(J_{\min }(G))\cap P^{\ast }(G)=\emptyset .$ 
Hence, there exists a number $\delta _{1}$ such that 
for each\\ $z\in \cup _{\alpha \in \G \setminus \G_{\min }}
\alpha ^{-1}(J_{\min }(G))$ and each $\beta \in G\cup \{ Id \} $, 
\begin{equation}
\label{mainth3-2pfeq3}
\deg (\beta :W\rightarrow D(z,\delta _{1}))=1,
\end{equation}
for each connected component $W$ of $\beta ^{-1}(D(z,\delta _{1})).$ 
For this $\delta _{1}$, there exists a number $\delta _{2}>0$ 
such that for each $z\in J_{\min }(G)$ and  each $\alpha \in \G 
\setminus \G _{\min }$, 
\begin{equation}
\label{mainth3-2pfeq4}
\mbox{diam }B\leq \delta _{1},\  
\deg (\alpha : B\rightarrow D(z,\delta _{2}))
\leq \max\{ \deg (\alpha )\mid \alpha \in \G \setminus \G _{\min }\} 
\end{equation}
for each connected component $B$ of $\alpha ^{-1}(D(z,\delta _{2})).$
Furthermore, 
by \cite[Lemma 1.10]{S1} (or \cite{S2}), 
we have that there exists a constant $0<c<1$ such that 
for each $z\in J_{\min }(G)$, each $h\in G_{\min ,\G }\cup \{ Id \} $, 
and each connected component $V$ of $h^{-1}(D(z,c\delta ))$, 
\begin{equation}
\label{mainth3-2pfeq5}
 \mbox{diam }V\leq \delta _{2}.
\end{equation}
Let $g\in G$ be any element. 

 Suppose that $g\in G_{\min ,\G}.$ Then, 
 by (\ref{mainth3-2pfeq2}), for each $z\in J_{\min }(G)$,  we have 
$\deg (g:V\rightarrow D(z, c\delta  ))\leq N$, for each 
connected component $V$ of $g^{-1}(D(z,c\delta )).$  

 Suppose that 
 $g$ is of the form $g=h\circ \alpha \circ g_{0}$, where 
$h\in G_{\min ,\G }\cup \{ Id \}$, $ \alpha \in \G \setminus \G _{\min }$, 
and $g_{0}\in G\cup \{ Id \} .$ 
 Then, 
 combining (\ref{mainth3-2pfeq3}), (\ref{mainth3-2pfeq4}), 
 and (\ref{mainth3-2pfeq5}),  
 we get that for each $z\in J_{\min }(G)$, 
 $\deg (g: W\rightarrow D(z,c\delta ))\leq N\cdot \max 
 \{ \deg (\alpha )\mid \alpha \in \G \setminus \G _{\min }\} $, 
 for each connected component $W$ of $g^{-1}(D(z,c\delta )).$  

 From the above argument, we see that 
 $J_{\min }(G)\subset  SH_{N'}(G)$, 
 where $N':=N\cdot \max \{\deg( \alpha )\mid 
 \alpha \in \G \setminus \G _{\min }\} .$ 
Moreover, by Theorem~\ref{mainth2}-\ref{mainth2-3},
we see that for any point $z\in J(G)\setminus J_{\min }(G)$, 
$z\in SH_{1}(G).$ 
Hence, we have shown that 
$J(G)\subset \CCI \setminus UH(G).$ Therefore, 
$G$ is semi-hyperbolic, provided that we have Case 1.  

 We now suppose that we have Case 2. 
 Then, by Proposition~\ref{nonminnoncpt}, 
 we have that for each $h\in \G _{\min }$, 
 $K(h)=\hat{K}(G)$ and int$(K(h))$ is non-empty and connected. 
 Moreover, for each $h\in \G _{\min }$, 
 int$(K(h))$ is an immediate basin of an attracting fixed point 
 $z_{h}\in \CC .$ 
Let ${\cal U}$ be the open neighborhood of 
$\G _{\min }$ in $\G $ as in 
Lemma~\ref{060523lem1}. 
Denoting by $H$ the polynomial semigroup generated by 
${\cal U}$, we have 
$P^{\ast }(H)\subset $ int$(\hat{K}(G)).$ 
Therefore, there exists a number $\delta >0$ such that 
\begin{equation}
\label{shshproppfeq02}
D(J(G),\delta )\subset \CC \setminus P(H).
\end{equation} 
Moreover, combining Theorem~\ref{mainth2}-\ref{mainth2-4-2} 
and that $\G \setminus {\cal U}$ is compact, 
we see that there exists a number 
$\epsilon >0$ such that 
\begin{equation}
\label{shshproppfeq03}
\overline{ \bigcup_{\alpha \in \G \setminus {\cal U}}
\alpha ^{-1}(D(J_{\min }(G),\epsilon ))}\subset A_{0},  
\end{equation}
where $A_{0}$ denotes the unbounded component of $\CC \setminus J_{\min }(G).$ 
Combining it with Theorem~\ref{mainth2}-\ref{mainth2-3}, 
it follows that there exists a number 
$\delta _{1}>0$ such that 
\begin{equation}
\label{shshproppfeq04}
D\left( \bigcup _{\alpha \in 
\G \setminus {\cal U}}\alpha ^{-1}
(D(J_{\min }(G),\epsilon )),\ \delta _{1}\right) 
\subset \CC \setminus P(G).
\end{equation}
For this $\delta _{1}$, 
there exists a number $\delta _{2}>0$ such that 
for each $\alpha \in \G \setminus {\cal U}$ and each 
$x\in D(J_{\min }(G),\epsilon )$, 
\begin{equation}
\label{shshproppfeq05}
\mbox{diam }B\leq \delta _{1},\ 
\deg (\alpha :B\rightarrow D(x,\delta _{2}))\leq 
\max \{ \deg (\beta )\mid \beta \in \G \setminus {\cal U}\} 
\end{equation}
for each connected component $B$ of 
$\alpha ^{-1}(D(x,\delta _{2})).$ 
By Lemma~\ref{invnormal} and (\ref{shshproppfeq02}), 
there exists a constant $c>0$ such that 
for each $h\in H$ and each $z\in J_{\min }(G)$, 
\begin{equation}
\label{shshproppfeq06}
\mbox{diam }V\leq \min \{ \delta _{2},\epsilon \} , 
\end{equation}
for each connected component $V$ of $h^{-1}(D(z,c\delta )).$ 
Let $z\in J_{\min }(G)$ and $g\in G.$ We will show that 
$z\in \CC \setminus UH(G).$ 

 Suppose that $g\in H.$ Then, 
 (\ref{shshproppfeq02}) implies that 
 for each 
connected component $V$ of 
$g^{-1}(D(z,c\delta ))$, 
$\deg (g:V\rightarrow D(z,c\delta ))=1.$ 

 Suppose that $g$ is of the form $g=h\circ \alpha \circ 
 g_{0}$, where 
$h\in H\cup \{ Id\} , \alpha \in \G \setminus 
{\cal U}, g_{0}\in G\cup \{ Id \} .$ 
Let 
$W$ be a connected component of $g^{-1}(D(z,c\delta ))$ and 
let $W_{1}:=g_{0}(W)$ and $V:=\alpha (W_{1}).$ 
Let $z_{1}$ be the point such that 
$\{ z_{1}\} =V\cap h^{-1}(\{ z\} ).$ 
If $z_{1}\in \CC \setminus 
D(J_{\min }(G), \epsilon )$, then, 
by (\ref{shshproppfeq06}) and Theorem~\ref{mainth2}-\ref{mainth2-3}, 
$V\subset D(z_{1},\epsilon )\subset \CC \setminus P(G).$ 
 Hence, 
 $\deg (\alpha \circ g_{0}:W\rightarrow V)=1$, 
 which implies that 
 $\deg (g:W\rightarrow D(z,c\delta ))=1.$ 
If $z_{1}\in D(J_{\min }(G), \epsilon )$, 
then by (\ref{shshproppfeq06}), 
$V\subset D(z_{1},\delta _{2}).$ 
Combining it with 
(\ref{shshproppfeq04}) and (\ref{shshproppfeq05}), 
we obtain  
$\deg (\alpha \circ g_{0}:W\rightarrow V)=
\deg (\alpha :W_{1}\rightarrow V)\leq 
\max \{ \deg (\beta )\mid \beta \in \G \setminus {\cal U}\} .$ 
Therefore, 
$\deg (g:W\rightarrow D(z,c\delta ))\leq 
\max \{ \deg (\beta )\mid \beta \in \G \setminus {\cal U}\} .$ 
Thus, $J_{\min }(G)\subset \CC \setminus UH(G).$ 

 Moreover, Theorem~\ref{mainth2}-\ref{mainth2-3} 
 implies that 
 $J(G)\setminus J_{\min }(G)\subset 
 \CC \setminus P(G)\subset \CC \setminus UH(G).$ 
 Therefore, $J(G)\subset \CC \setminus UH(G)$, 
 which implies that $G$ is semi-hyperbolic.  

 Thus, we have proved Theorem~\ref{shshprop}.
\qed 

\ 

We now demonstrate Theorem~\ref{hhprop}.\\ 
{\bf Proof of Theorem~\ref{hhprop}:} 
We use the same argument as that in the proof of Theorem~\ref{shshprop}, 
but we modify it as follows: 
\begin{enumerate}
\item 
 In (\ref{mainth3-2pfeq1}), we replace 
$UH(G_{\min ,\G })\cap \CC $ by $P^{\ast }(G_{\min ,\G }).$  
\item  In (\ref{mainth3-2pfeq2}), 
we replace $N$ by $1.$  
\item  We replace (\ref{mainth3-2pfeq4}) by 
the following (\ref{mainth3-2pfeq4})' 
$\mbox{diam }B\leq \delta _{1},\  
\deg (\alpha : B\rightarrow D(z,\delta _{2}))=1.$  
\item  We replace (\ref{shshproppfeq05}) by the following 
(\ref{shshproppfeq05})' 
$\mbox{diam }B\leq \delta _{1},\ 
\deg (\alpha :B\rightarrow D(x,\delta _{2}))=1.$ 
(We take the number $\epsilon >0$ so small.)  
\end{enumerate} 
With these modification, it is easy to see that $G$ is hyperbolic. 

 Thus, we have proved Theorem~\ref{hhprop}. 
\qed 

\ 

We now demonstrate Proposition~\ref{nonminnoncpt}-\ref{nonminnoncpt4}.\\ 
{\bf Proof of Proposition~\ref{nonminnoncpt}-\ref{nonminnoncpt4}:}
Suppose that (a) in Proposition~\ref{nonminnoncpt}-\ref{nonminnoncpt2} 
holds. By Lemma~\ref{060523lem1}, 
$G_{\min ,\G }$ is hyperbolic. 
Combining it with Theorem~\ref{shshprop}, 
it follows that $G$ is semi-hyperbolic. 
Thus, we have proved Proposition~\ref{nonminnoncpt}-\ref{nonminnoncpt4}.
\qed

\ 
To demonstrate Theorem~\ref{mainth3-2}, 
we need the following  proposition.
\begin{prop}
\label{shonecomp}
Let $f:X\times \CCI \rightarrow X\times \CCI $ be a semi-hyperbolic 
polynomial skew product 
over $g:X\rightarrow X.$ 
 Suppose that for each $x\in X$, $d(x)\geq 2$, and that 
$\pc (P(f))\cap \CC $ is bounded in $\CC .$   
Let $\omega \in X$ be a point. 
If {\em int}$(K_{\omega }(f))$ is 
a non-empty connected set, then 
$J_{\omega }(f)$ is a Jordan curve.  
\end{prop}
\begin{proof}
By \cite[Theorem 1.12]{S4} and Lemma~\ref{fibconnlem}, we get that the 
unbounded component $A_{\omega }(f)$ of 
$F_{\omega }(f)$ is a John domain. Combining it,  
that $A_{\omega }(f)$ is simply connected (cf. Lemma~\ref{fibconnlem}), 
and \cite[page 26]{NV}, we see that 
$J_{\omega }(f)
=\partial (A_{\omega }(f))$ 
(cf. Lemma \ref{fibfundlem}) is locally connected. 
Moreover, by Lemma~\ref{fibfundlem}-\ref{fibfundlema}, 
we have 
$\partial (\mbox{int}(K_{\omega }(f)))=J_{\omega }(f).$ 
Hence, we see that $\CCI \setminus J_{\omega }(f)$ has 
exactly two connected components $A_{\g }(f)$ and int$(K_{\omega }(f))$, 
and that $J_{\omega }(f)$ is locally connected. From 
\cite[Lemma 5.1]{PT}, it follows that 
$J_{\g }(f)$ is a Jordan curve. Thus, we have proved 
Proposition~\ref{shonecomp}.
%
\end{proof}

 We now demonstrate Theorem~\ref{mainth3-2}.\\ 
\noindent {\bf Proof of Theorem~\ref{mainth3-2}:}
Let 
$\g \in R(\G ,\G \setminus \G _{\min })$ and 
$y\in $ int$(K_{\g }(f)).$ 
Combining Theorem~\ref{mainth3}-\ref{mainth3-0} 
and \cite[Lemma 1.10]{S1}, 
we obtain $\liminf _{n\rightarrow \infty }$ $
d(f_{\g,n}(y),$ $J(G))>0.$ 
Combining this with Lemma~\ref{constlimlem} and 
Theorem~\ref{mainth3}-\ref{mainth3-0}, 
we see that there exists a point 
$a\in P^{\ast }(G)\cap F(G)$ such that 
$\liminf _{n\rightarrow \infty }d(f_{\g ,n}(y),a)=0.$ 
Since $P^{\ast }(G)\cap F(G)\subset $ int$(\hat{K}(G))$, 
it follows that there exists a positive integer 
$l$ such that 
\begin{equation}
\label{mainth3-2pfeq0}
f_{\g ,l}(y)\in \mbox{int}(\hat{K}(G)).
\end{equation} 
Combining (\ref{mainth3-2pfeq0}) and the 
same method as that in the proof of Claim 3 in 
the proof of Theorem~\ref{mainth3}-\ref{mainth3-0} and 
Theorem~\ref{mainth3}-\ref{mainth3-1}, 
we get that there exists exactly one 
bounded component $U_{\g }$ of $F_{\g }(f).$ 
Combining it with Proposition~\ref{shonecomp}, 
it follows that $J_{\g }(f)$ is a Jordan curve. 
%
Moreover, by \cite[Theorem 2.14-(4)]{S1}, 
we have $\hat{J}_{\gamma }(f)=J_{\gamma }(f).$ 

 Thus, we have proved Theorem~\ref{mainth3-2}.
\qed

\  

 We now demonstrate Theorem~\ref{cantorqc}.

\noindent 
{\bf Proof of Theorem~\ref{cantorqc}:}
Let $V$ be an open set with $J(G)\cap V\neq \emptyset .$ 
We may assume that $V$ is connected.
Then, by Theorem~\ref{repdense}, 
there exists an element $\alpha _{1}\in G$ such that 
$J(\alpha _{1})\cap V\neq \emptyset .$ 
Since we have $G\in {\cal G}_{dis}$,  
Theorem~\ref{mainth0} implies that there exists an element 
$\alpha _{2}\in G$ such that 
no connected component $J$ of $J(G)$ satisfies  
$J(\alpha _{1})\cup J(\alpha _{2})\subset J.$ 
Hence, we have $\langle \alpha _{1}, \alpha _{2}\rangle \in {\cal G}_{dis}.$ 
Since $J(\alpha _{1})\cap V\neq \emptyset $, 
combined with Lemma~\ref{fibfundlem}-\ref{fibfundlem2}, 
we get that there exists an $l_{0}\in \NN $ such that 
for each $l$ with $l\geq l_{0}$, 
$J(\alpha _{2}\alpha _{1}^{l})\cap V\neq \emptyset .$ 
Moreover, since no connected component $J$ of $J(G)$ satisfies  
$J(\alpha _{1})\cup J(\alpha _{2})\subset J$, 
Lemma~\ref{fibfundlem}-\ref{fibfundlem2} implies that 
there exists an $l_{1}\in \NN $ such that 
for each $l$ with $l\geq l_{1}$, 
$J(\alpha _{2}\alpha _{1}^{l})\cap 
J(\alpha _{1}\alpha _{2}^{l})=\emptyset .$ 
We fix an $l\in \NN $ with $l\geq \max \{ l_{0},l_{1}\}.$ 
We now show the following claim.\\ 
\noindent Claim 1. 
The semigroup $H_{0}:= \langle \alpha _{2}\alpha _{1}^{l}, 
\alpha _{1}\alpha _{2}^{l}\rangle $ is hyperbolic, and 
for the skew product $\tilde{f}:\GN _{0}\times \CCI 
\rightarrow \GN _{0}\times \CCI $ associated with 
$\G _{0}=\{ \alpha _{2}\alpha _{1}^{l}, \alpha _{1}\alpha _{2}^{l}\} $, 
there exists a constant $K\geq 1$ such that 
for any $\g \in \GN _{0}$, $J_{\g }(\tilde{f})$ is a $K$-quasicircle. 

 To show this claim, 
applying Theorem~\ref{mainth3}-\ref{mainth3-3}
with $\G =\{ \alpha _{1},\alpha _{2}\} , 
S=\G \setminus \G _{\min }, $ and $p=2l+1$, 
we see that 
the polynomial skew product 
$\overline{f}:W_{S, 2l+1}\times \CCI \rightarrow W_{S,2l+1}\times 
\CCI $ over $\sigma :W_{S,2l+1}\rightarrow W_{S,2l+1}$ 
is hyperbolic, and that 
there exists a constant $K\geq 1$ such that 
for each $\g \in W_{S,2l+1}$, $J_{\g }(\overline{f})$ 
is a $K$-quasicircle. 
Moreover, combining the hyperbolicity of $\overline{f}$ above and 
Remark~\ref{hypskewsemigrrem}, 
we see that the semigroup 
$H_{1}$ generated by the family 
$\{ \alpha_{j_{1}}\circ \cdots \circ \alpha _{j_{l+1}}\mid 
1\leq \exists k_{1}\leq l+1 \mbox{ with } j_{k_{1}}=1, 
1\leq \exists k_{2}\leq l+1 \mbox{ with }j_{k_{2}}=2\} $ 
is hyperbolic. 
Hence, the semigroup $H_{0}$, which is a subsemigroup of $H_{1}$, 
is hyperbolic. 
 Therefore, Claim 1 holds.

 We now show the following claim.\\ 
Claim 2. We have either 
$J(\alpha _{2}\alpha _{1}^{l})<J(\alpha _{1}\alpha _{2}^{l})$,  
or $J(\alpha _{1}\alpha _{2}^{l})<J(\alpha _{2}\alpha _{1}^{l}).$

 To show this claim, since 
$J(\alpha _{2}\alpha _{1}^{l})\cap J(\alpha _{1}\alpha _{2}^{l})   
=\emptyset $ and $H_{0}\in {\cal G}$, 
combined with Lemma~\ref{fiborder}, we obtain Claim 2. 

 By Claim 2, we have the following two cases.\\ 
Case 1. $J(\alpha _{2}\alpha _{1}^{l})<J(\alpha _{1}\alpha _{2}^{l}).$\\ 
Case 2. $J(\alpha _{1}\alpha _{2}^{l})<J(\alpha _{2}\alpha _{1}^{l}).$  

 We may assume that we have Case 1 (when we have Case 2, we can show 
 all statements of our theorem, using the same method as below).  
Let $A:= K(\alpha _{1}\alpha _{2}^{l})\setminus $ 
int$(K(\alpha _{2}\alpha _{1}^{l})).$ 
By Claim 1, we have that 
$J(\alpha _{1}\alpha _{2}^{l})$ and 
$J(\alpha _{2}\alpha _{1}^{l})$ are quasicircles. 
Moreover, 
since $H_{0}\in {\cal G}_{dis }$ and 
$H_{0}$ is hyperbolic, 
we must have 
$P^{\ast }(H_{0})\subset $
int$(K(\alpha _{2}\alpha _{1}^{l})).$ 
Therefore, 
it follows that if we take a small open neighborhood $U$ of 
$A$, then there exists a number $n\in \NN $ such that, 
setting $h_{1}:= (\alpha _{2}\alpha _{1}^{l} )^{n}$ 
and $h_{2}:= (\alpha _{1}\alpha _{2}^{l})^{n}$, 
we have that 
\begin{equation}
\label{cantorqcpfeq01}
h_{1}^{-1}(\overline{U})\cup 
h_{2}^{-1}(\overline{U})\subset U \mbox{ and } 
h_{1}^{-1}(\overline{U})\cap h_{2}^{-1}(\overline{U})=\emptyset . 
\end{equation}
Moreover, combining Lemma~\ref{fibfundlem}-\ref{fibfundlem2} 
and that $J(h_{1})\cap V\neq \emptyset $, 
we get that there exists a $k\in \NN $ such that 
$J(h_{2}h_{1}^{k})\cap V\neq \emptyset .$ 
We set $g_{1}:= h_{1}^{k+1}$ and $g_{2}:= h_{2}h_{1}^{k}.$ 
Moreover, we set $H:= \langle g_{1},g_{2}\rangle .$ Since 
$H$ is a subsemigroup of $H_{0}$ and $H_{0}$ is hyperbolic, 
we have that $H$ is hyperbolic. Moreover, 
(\ref{cantorqcpfeq01}) implies that 
$g_{1}^{-1}(\overline{U})\cup 
g_{2}^{-1}(\overline{U})\subset U$ and 
$g_{1}^{-1}(\overline{U})\cap g_{2}^{-1}(\overline{U})=\emptyset .$ 
Hence, we have shown that for the semigroup $H=\langle g_{1},g_{2}\rangle $, 
statements \ref{cantorqc1},\ref{cantorqc2},
and \ref{cantorqc3} in Theorem~\ref{cantorqc} hold. 

 From statement \ref{cantorqc2} and 
  Lemma~\ref{hmslem}-\ref{backmin}, 
 we obtain $J(H)\subset \overline{U}$ 
 and $g_{1}^{-1}(J(H))\cap g_{2}^{-1}(J(H))=\emptyset .$  
Combining this with Lemma~\ref{hmslem}-\ref{bss} and 
Lemma~\ref{fiblem}-\ref{fibfundlem5},  
it follows that the 
skew product $f:\GN _{1}\times \CCI 
\rightarrow \GN _{1}\times \CCI $ 
associated with $\G _{1}= \{ g_{1},g_{2}\} $ 
satisfies that 
$J(H)$ is equal to the disjoint union of 
the sets $\{ \hat{J}_{\g }(f)\} _{\g \in \GN _{1}}.$ 
Moreover, since $H$ is hyperbolic, 
\cite[Theorem 2.14-(2)]{S1} implies that for each 
$\g \in \G _{1}^{\NN }$, 
$\hat{J}_{\g }(f)=J_{\g }(f).$  
In particular, 
the map $\g \mapsto J_{\g }(f)$ 
from $\GN _{1}$ into the space of 
non-empty compact sets in $\CCI $,  is injective. 
Since $J_{\g }(f)$ is connected for each 
$\g \in \GN _{1}$ (Claim 1), 
it follows that for each connected component $J$ of $J(H)$,  
there exists an element $\g \in \GN _{1}$ such that 
$J=J_{\g }(f).$ 
Furthermore, by Claim 1,  
each connected component $J$ of $J(H)$ is a $K$-quasicircle, 
where $K$ is a constant not depending on $J.$ 
 Moreover, by \cite[Theorem 2.14-(4)]{S1}, the map 
 $\g \mapsto J_{\g }(f)$ from $\GN _{1}$ 
 into the space of non-empty compact sets in $\CCI $, 
 is continuous with respect to 
 the Hausdorff topology. 
Therefore, we have shown that statements 
 \ref{cantorqc4a},\ref{cantorqc4b},\ref{cantorqc4c}, and \ref{cantorqc5} hold for $H=\langle g_{1},g_{2}\rangle $ and 
 $f:\GN _{1}\times \CCI \rightarrow \GN _{1}\times \CCI .$  

 We now show that statement \ref{cantorqc6} holds.
Since we are assuming Case 1, 
Proposition~\ref{bminprop} implies that 
$\{ h_{1},h_{2}\} _{\min }=\{ h_{1}\} .$ 
Hence $J(g_{1})<J(g_{2}).$ 
Combining it with Proposition~\ref{bminprop} and statement \ref{cantorqc4b},   
we obtain  
\begin{equation}
\label{cantorqc6pfeq1}
J(g_{1})= J_{\min }(H) \mbox{ and }
J(g_{2})= J_{\max }(H).
\end{equation} 
Moreover, since 
$J(g_{1})=J(\alpha _{2}\alpha _{1}^{l})$, 
$J(\alpha _{2}\alpha _{1}^{l})\cap V\neq \emptyset $, 
$J(g_{2})=J(h_{2}h_{1}^{k})$, and $J(h_{2}h_{1}^{k})\cap 
V\neq \emptyset $, 
it follows that 
\begin{equation}
\label{cantorqc6pfeq2}
J_{\min }(H)\cap V\neq \emptyset \mbox{ and }
J_{\max }(H)\cap V\neq \emptyset .
\end{equation}
Let $\g \in \GN $ be an element such that 
$J_{\g }(f)\cap (J_{\min }(H)\cup J_{\max }(H))=\emptyset .$ 
By statement \ref{cantorqc4b}, 
we obtain 
\begin{equation}
\label{cantorqc6pfeq3}
J_{\min }(H)<J_{\g }(f)<J_{\max }(H).
\end{equation}
Since we are assuming $V$ is connected, 
combining (\ref{cantorqc6pfeq2}) and 
(\ref{cantorqc6pfeq3}), 
we obtain $J_{\g }(f)\cap V\neq \emptyset .$ 
Therefore, we have proved that 
statement \ref{cantorqc6} holds. 

 We now show that statement \ref{cantorqc7} holds. 
 To show that, 
let $\omega  =(\omega  _{1},\omega  _{2},\ldots )\in \GN _{1}$ be an 
element such that 
$\sharp (\{ j\in \NN 
\mid \omega _{j}=g_{1}\} )=
\sharp (\{ j\in \NN \mid \omega _{j}=g_{2}\} )=\infty .$ 
For each $r\in \NN $, 
let $\omega  ^{r}=(\omega  _{r,1},\omega  _{r,2},\ldots )\in \GN _{1}$ 
be the element such that 
$\begin{cases}
\omega  _{r,j}=\omega  _{j} \ \ (1\leq j\leq r), \\  
\omega  _{r,j}=g_{1}     \ \ (j\geq r+1).
\end{cases} $
 Moreover, 
let $\rho ^{r}=(\rho _{r,1},\rho _{r,2},\ldots )\in \GN _{1}$ 
be the element such that 
 $\begin{cases}
\rho_{r,j}=\omega  _{j} \ \ (1\leq j\leq r), \\  
\rho _{r,j}=g_{2}     \ \ (j\geq r+1).
\end{cases} $
Combining (\ref{cantorqc6pfeq1}), 
statement \ref{cantorqc4a}, and statement \ref{cantorqc4b},   
we see that for each $r\in \NN $, 
$J(g_{1})<J_{\sigma ^{r}(\omega  )}(f)<J(g_{2}).$ 
Hence, by Theorem~\ref{mainth1}-\ref{mainth1-3}, 
we get that for each $r\in \NN $, 
$(f_{\omega  ,r})^{-1}(J(g_{1}))<
(f_{\omega  ,r})^{-1}\left(J_{\sigma ^{r}(\omega  )}(f)\right)
<(f_{\omega  ,r})^{-1}(J(g_{2}))$, 
 where 
$f_{\omega  ,r}(y)=\pi _{\CCI }(f^{r}(\omega ,y)).$   
Since we have 
$ (f_{\omega  ,r})^{-1}(J(g_{1}))=J_{\omega  ^{r}}(f)$, 
$(f_{\omega  ,r})^{-1}\left(J_{\sigma ^{r}(\omega  )}(f)\right)
=J_{\omega  }(f)$, 
and $(f_{\omega  ,r})^{-1}(J(g_{2}))$ $=J_{\rho ^{r}}(f)$, 
it follows that 
\begin{equation}
\label{cantorqcpfeq1}
J_{\omega  ^{r}}(f)<J_{\omega  }(f)<
J_{\rho ^{r}}(f),  
\end{equation}
for each $r\in \NN .$ 
Moreover, since 
$\omega  ^{r}\rightarrow \omega  $ and 
$\rho ^{r}\rightarrow \omega  $ in 
$\GN _{1}$ as $r\rightarrow \infty $, 
statement \ref{cantorqc5}
implies that 
$J_{\omega  ^{r}}(f)\rightarrow 
J_{\omega  }(f)$ and 
$J_{\rho ^{r}}(f)\rightarrow 
J_{\omega  }(f)$ as $r\rightarrow \infty $, 
with respect to the Hausdorff topology. 
Combined with (\ref{cantorqcpfeq1}), statement \ref{cantorqc4b} 
and statement \ref{cantorqc4c}, 
we get that 
for any connected component $W$ of $F(H)$, 
we must have $\partial W\cap J_{\omega  }(f)
=\emptyset .$ Since $F(G)\subset F(H)$, 
it follows that 
for any connected component $W'$ of $F(G)$, 
$\partial W'\cap J_{\omega  }(f)=\emptyset .$
Therefore, we have shown that statement \ref{cantorqc7} holds. 
%
%
 
  Thus, we have proved Theorem~\ref{cantorqc}.
\qed 
\subsection{Proofs of results in \ref{fjjq}}
\label{Proofs of fjjq}
In this section, we demonstrate Theorem~\ref{mainthjbnq}.  
We need the following notations and lemmas.
\begin{df}
Let $h$ be a polynomial with $\deg (h)\geq 2.$ Suppose that 
$J(h)$ is connected. Let $\psi $ be a biholomorphic map 
$\CCI \setminus \overline{D(0,1)}\rightarrow  
F_{\infty }(h)$ with $\psi (\infty )=\infty $ such that 
$\psi ^{-1}\circ h\circ \psi (z)=z^{\deg (h)}$, 
for each $z\in \CCI \setminus \overline{D(0,1)}.$ 
(For the existence of the biholomorphic map $\psi $, see 
\cite[Theorem 9.5]{M}.) For each 
$\theta \in \partial D(0,1)$, we set 
$T(\theta ):= \psi (\{ r\theta \mid 1<r\leq \infty \} ).$ 
This is called the external ray (for $K(h)$) 
with angle $\theta .$ 
\end{df}

\begin{lem}
\label{jbnqlem1}
Let $h$ be a polynomial with $\deg (h)\geq 2.$ 
Suppose that $J(h)$ is connected and locally connected 
and $J(h)$ is not a Jordan curve. Moreover, suppose that 
there exists an attracting periodic point of $h$ in $K(h).$ 
Then, for any $\epsilon >0$, there exist a point $p\in J(h)$  
and elements $\theta _{1}, \theta _{2}\in \partial D(0,1)$ 
with $\theta _{1}\neq \theta _{2}$, such that all of the following hold.
\begin{enumerate}
\item 
For each $i=1,2$,  the external ray $T(\theta _{i})$ lands 
at the point $p.$ 
\item Let $V_{1}$ and $V_{2}$ be the two connected components of  
 $\CCI \setminus (T(\theta _{1})\cup T(\theta _{2})\cup \{ p\} ).$  
 Then, for each $i=1,2,$ 
 $V_{i}\cap J(h)\neq \emptyset .$ Moreover, there exists an $i$ such that 
 diam $(V_{i}\cap K(h))\leq \epsilon .$ 
\end{enumerate}
\end{lem} 
\begin{proof}
Let $\psi : \CCI \setminus \overline{D(0,1)}\rightarrow 
F_{\infty }(h)$ be a biholomorphic map with $\psi (\infty )=\infty $ 
such that for each $z\in 
\CCI \setminus \partial D(0,1),  
\psi ^{-1}\circ h\circ \psi (z)=z^{\deg (h)}.$ 
Since $J(h)$ is locally connected, the map 
$\psi : \CCI \setminus \overline{D(0,1)}\rightarrow 
F_{\infty }(h)$ extends continuously over $\partial D(0,1)$, 
mapping $\partial D(0,1)$ onto $J(h).$ 
Moreover, since $J(h)$ is not a Jordan curve, it follows that 
there exist a point $p_{0}\in J(h)$ and two points $t_{0,1},t_{0,2}\in 
\partial D(0,1)$ with $t_{0,1}\neq t_{0,2}$ such that 
two external rays 
$T(t _{0,1})$ and $T(t _{0,2})$ land at the same 
point $p_{0}.$ 
Considering a higher iterate of $h$ if necessary, 
we may assume that there exists an attracting fixed point of $h$ in 
int$(K(h)).$ 
Let $a\in $ int$(K(h))$ be an attracting fixed point of $h$ 
and let $U$ be the connected component of int$(K(h))$ containing $a.$ 
Then, there exists a critical point $c\in U$ of $h.$ 
Let $V_{0}$ be the connected component of 
$\CCI \setminus (T(t_{1})\cup T(t_{2})\cup \{ p_{0}\} )$ 
containing $a.$ 
Moreover, for each $n\in \NN $, 
let $V_{n}$ be the connected component of $(h^{n})^{-1}(V_{0})$ 
containing $a.$ 
Since $c\in U$, we get that for each $n\in \NN $, 
$c\in V_{n}.$ Hence, setting 
$e_{n}:= \deg (h^{n}:V_{n}\rightarrow V_{0})$, 
it follows that 
\begin{equation}
\label{jbnqlem1eq1}
e_{n}\rightarrow \infty \mbox{ as }n\rightarrow \infty .
\end{equation} 
We fix an $n\in \NN $ satisfying $e_{n}>d$, where 
$d:= \deg (h).$ 
Since 
$\deg (h^{n}:V_{n}\cap F_{\infty }(h)\rightarrow 
V_{0}\cap F_{\infty }(h))=\deg (h^{n}:V_{n}\rightarrow V_{0})$, 
we have that the number of connected components of 
$V_{n}\cap F_{\infty }(h)$ is equal to $e_{n}.$ 
Moreover, every connected component of $V_{n}\cap F_{\infty }(h)$ 
is a connected component of 
$(h^{n})^{-1}(V_{0}\cap F_{\infty }(h)).$ Hence, it follows that 
there exist mutually disjoint arcs $\xi _{1},\xi _{2},\ldots ,\xi _{e_{n}}$ 
in $\CC $ satisfying all of the following.
\begin{enumerate}
\item For each $j$,   
$h^{n}(\xi _{j})=(T(t_{1})\cup T(t_{2})\cup \{ p_{0}\} )\cap \CC .$ 
\item 
For each $j$, 
$\xi _{j}\cup \{ \infty \} $ is the closure of union of 
two external rays and $\xi _{j}\cup \{ \infty \} $ is a Jordan curve.
\item $\partial V_{n}=\xi _{1}\cup \cdots \cup \xi _{e_{n}}\cup \{ \infty \} .$ \end{enumerate} 
For each $j=1,\ldots ,e_{n}$, 
let $W_{j}$ be the connected component 
of $\CCI \setminus (\xi _{j}\cup \{ \infty \} )$ that does not contain 
$V_{n}.$ 
Then, each $W_{j}$ is a connected component of 
$\CCI \setminus \overline{V_{n}}.$ Hence, 
for each $(i,j)$ with $i\neq j$, 
$W_{i}\cap W_{j}=\emptyset .$ 
Since the number of critical values of $h$ in $\CC $ is less than or equal to 
$d-1$, we have that 
$\sharp (\{1\leq j\leq e_{n}\mid W_{j}\cap CV(h)=\emptyset \} )\geq 
e_{n}-(d-1).$ 
Therefore, denoting by $u_{1,j}$ 
the number of well-defined inverse branches of 
$h^{-1}$ on $W_{j}$, we obtain
$$\sum _{j=1}^{e_{n}}u_{1,j}\geq d(e_{n}-(d-1))\geq d.$$ 
Inductively, denoting by $u_{k,j}$ the number of well-defined inverse 
branches of $(h^{k})^{-1}$ on $W_{j}$, we obtain
\begin{equation}
\label{jbnqlem1pfeqa1}
\sum _{j=1}^{e_{n}}u_{k,j}\geq d(d-(d-1))\geq d, \mbox{for each }k\in \NN .
\end{equation}
For each $k\in \NN $, we take a well-defined   
inverse branch $\zeta _{k}$ of $(h^{k})^{-1}$ on a 
domain $W_{j}$, and let $B_{k}:=\zeta _{k}(W_{j}).$ 
Then, $h^{k}:B_{k}\rightarrow W_{j}$ is biholomorphic.
Since $\partial B_{k}$ is the closure of finite union of 
external rays and $h^{n+k}$ maps each connected component of 
$(\partial B_{k})\cap \CC $ 
onto $(T(t_{1})\cup T(t_{2})\cup \{ p_{0}\} )\cap \CC $, 
$B_{k}$ is a Jordan domain. 
Hence, $h^{k}:B_{k}\rightarrow W_{j}$ induces a 
homeomorphism $\partial B_{k}\cong \partial W_{j}.$ 
Therefore, $\partial B_{k}$ is the closure of 
union of two external rays, which implies that 
$B_{k}\cap F_{\infty }(h)$ is a connected component of 
$(h^{k})^{-1}(W_{j}\cap F_{\infty }(h)).$ 
Hence, we obtain  
\begin{equation}
\label{jbnqlem1eq2}
l\left( \overline{\psi ^{-1}(B_{k}\cap F_{\infty }(h))}\cap \partial D(0,1)
\right) \rightarrow 0 \mbox{ as } k\rightarrow \infty ,
\end{equation}  
where $l(\cdot )$ denotes the arc length of a subarc of 
$\partial D(0,1).$ 
Since $\psi :\CCI \setminus \overline{D(0,1)}\rightarrow  F_{\infty }(h)$ 
extends continuously over $\partial D(0,1)$, 
(\ref{jbnqlem1eq2}) implies that 
diam $(B_{k}\cap J(h))\rightarrow 0$ as $k\rightarrow \infty .$ 
Hence, there exists a $k\in \NN $ such that 
diam $(B_{k}\cap K(h))\leq \epsilon .$ 
Let $\theta _{1},\theta _{2}\in \partial D(0,1)$ be 
two elements such that 
$\partial B_{k}=\overline{T(\theta _{1})\cup T(\theta _{2})}.$ 
Then, there exists a point $p\in J(h)$ such that 
each $T(\theta _{i})$ lands at the point $p.$ 
By \cite[Lemma 17.5]{M}, any of two connected components of 
$\CCI \setminus (T(\theta _{1})\cup T(\theta _{2})\cup \{ p\} )$ 
intersects $J(h).$ 

 Thus, we have proved Lemma~\ref{jbnqlem1}.
 \end{proof}
\begin{lem}
\label{060414astlem1}
Let $G$ be a polynomial semigroup generated by a 
compact subset $\G $ of {\em Poly}$_{\deg \geq 2}.$
Let $f:\GNCR $ be the skew product associated with the family 
$\G .$ Suppose $G\in {\cal G}_{dis}.$ 
Let $m\in \NN $ and suppose that there exists an element 
$(h_{1},\ldots ,h_{m})\in \G ^{m}$ such that setting 
$h=h_{m}\circ \cdots \circ h_{1}$, 
$J(h)$ is connected and locally connected, and 
$J(h)$ is not a Jordan curve. Moreover, suppose that 
there exists an attracting periodic point of $h$ in 
$K(h).$ Let $\alpha =(\alpha _{1},\alpha _{2},\ldots )\in 
\GN $ be the element such that for each 
$k,l\in \NN \cup \{ 0\} $ with $1\leq l\leq m$, 
$\alpha _{km+l}=h_{l}.$ 
Let $\rho _{0}\in \G \setminus \G _{\min }$ be an element and 
let 
$\beta =(\rho _{0},\alpha _{1},\alpha _{2},\ldots )\in \GN .$ 
Moreover, let 
$\psi _{\beta }:\CCI \setminus \overline{D(0,1)}
\rightarrow A_{\beta }(f)$ be a biholomorphic map 
with $\psi _{\beta }(\infty )=\infty .$ 
Furthermore, for each $\theta \in 
\partial D(0,1)$, 
let 
$T_{\beta }(\theta )=\psi _{\beta }(\{ r\theta \mid 
1<r\leq \infty \} ).$ 
Then, for any $\epsilon >0$, there exist a point 
$p\in J_{\beta }(f)$ and elements 
$\theta _{1},\theta _{2}\in \partial D(0,1)$ with 
$\theta _{1}\neq \theta _{2}$, such that 
all of the following statements 1 and 2 hold.
\begin{enumerate}
\item For each $i=1,2$, $T_{\beta }(\theta _{i})$ lands at $p.$
\item Let $V_{1}$ and $V_{2}$ be the two connected components 
of 
$\CCI \setminus (T_{\beta }(\theta _{1})\cup T_{\beta }(\theta _{2})
\cup \{ p\} ).$ Then, 
for each $i=1,2$, $V_{i}\cap J_{\beta }(f)\neq \emptyset .$ 
Moreover, there exists an $i$ such that 
diam $(V_{i}\cap K_{\beta }(f))\leq \epsilon $ and 
such that $V_{i}\cap J_{\beta }(f)\subset 
\rho _{0}^{-1}(J(G))\subset \CC \setminus P(G).$ 
\end{enumerate}  
\end{lem}
\begin{proof}
We use the notation and argument in the proof of 
Lemma~\ref{jbnqlem1}. 
Taking a higher iterate of $h$, 
we may assume that 
$d:=\deg (h)>\deg (\rho _{0}).$ 
Then, from (\ref{jbnqlem1pfeqa1}), 
it follows that 
for each $k\in \NN $, 
we can take a well-defined inverse branch 
$\zeta _{k}$ of 
$(h^{k})^{-1}$ on a domain $W_{j}$ such that 
setting $B_{k}:=\zeta _{k}(W_{j})$, 
$B_{k}$ does not contain any critical value of 
$\rho _{0}.$ 
By (\ref{jbnqlem1eq2}), 
there exists a $k\in \NN $ such that 
diam $(B_{k}\cap J(h))\leq \epsilon '$, where 
$\epsilon '>0$ is a small number. 
Let $B$ be a connected component of 
$\rho _{0}^{-1}(B_{k}).$ 
Then, there exist a point 
$p\in J_{\beta }(f)$ and elements 
$\theta _{1},\theta _{2}\in \partial D(0,1)$ with 
$\theta _{1}\neq \theta _{2}$ such that 
for each $i=1,2,$ $T_{\beta }(\theta _{i})$ lands at 
$p$, and such that $B$ is a connected component of 
$\CCI \setminus (T_{\beta }(\theta _{1})\cup 
T_{\beta }(\theta _{2})\cup \{ p\} ).$ Taking $\epsilon '$ so small, 
we obtain diam $(B\cap K_{\beta }(f))=$ diam $(B\cap J_{\beta }(f))
\leq \epsilon .$ Moreover, 
since $\rho _{0}\in \G \setminus \G _{\min }$, 
combining Theorem~\ref{mainth2}-\ref{mainth2-3} and 
Theorem~\ref{mainth2}-\ref{mainth2-4-2}, 
we obtain 
$J_{\beta }(f)
=\rho _{0}^{-1}(J(h))
\subset \rho _{0}^{-1}(J(G))
\subset  
\CC \setminus P(G).$ 
Hence, $B\cap J_{\beta }(f)\subset 
\rho _{0}^{-1}(J(G))\subset 
\CC \setminus P(G).$ 
Therefore, we have proved Lemma~\ref{060414astlem1}.  
\end{proof}
\begin{lem}
\label{060416lema}
Let $f:X\times \CCI \rightarrow X\times \CCI $ be a polynomial skew product over 
$g:X\rightarrow X$ such that 
for each $x\in X$, 
$d(x)\geq 2.$ Let $\g \in X$ be a point. 
Suppose that $J_{\g }(f)$ is a Jordan curve. 
Then, for each $n\in \NN $, 
$J_{g^{n}(\g )}(f)$ is a Jordan curve. 
Moreover, for each $n\in \NN $, 
there exists no critical value of 
$f_{\g ,n}$ in $J_{g^{n}(\g )}(f).$ 
\end{lem}
\begin{proof}
Since 
$(f_{\g ,1})^{-1}(K_{g(\g )}(f))=K_{\g }(f)$, 
it follows that int$(K_{g(\g )}(f))$ is a non-empty 
connected set. 
Moreover, 
$J_{g(\g )}(f)=f_{\g ,1}(J_{\g }(f))$ is locally connected. 
Furthermore, by Lemma~\ref{fibfundlem}-\ref{fibfundlem4} and 
Lemma~\ref{fibfundlem}-\ref{fibfundlema}, 
$\partial (\mbox{int}(K_{g(\g )}(f)))=\partial (A_{g(\g )}(f))= J_{g(\g )}(f).$ 
Combining the above arguments and \cite[Lemma 5.1]{PT}, 
we get that $J_{g(\g )}(f)$ is a Jordan curve. 
Inductively, 
we conclude that for each $n\in \NN $, 
$J_{g^{n}(\g )}(f)$ is a Jordan curve. 

 Furthermore, applying the Riemann-Hurwitz formula 
 to the map 
 $f_{\g, n}: \mbox{int}(K_{\g }(f))\rightarrow 
 \mbox{int}(K_{g^{n}(\g )}(f))$, we obtain 
 $1+p=\deg (f_{\g ,n})$, where $p$ 
 denotes the cardinality of the critical 
 points of 
 $f_{\g, n}: \mbox{int}(K_{\g }(f))\rightarrow 
 \mbox{int}(K_{g^{n}(\g )}(f))$ counting multiplicities. 
Hence, $p=\deg (f_{\g ,n})-1.$ 
It implies that there exists no critical value of 
$f_{\g ,n}$ in 
$J_{g^{n}(\g )}(f).$     
\end{proof}
\begin{lem}
\label{060416lemast}
Let $f:X\times \CCI \rightarrow X\times \CCI  $ be a polynomial skew product over 
$g:X\rightarrow X$ such that 
for each $x\in X$, $d(x)\geq 2.$ 
Let $\mu >0$ be a number. 
Then, there exists a number $\delta >0$ such that 
the following statement holds.
\begin{itemize}
\item 
Let $\omega \in X$ be any point and 
$p\in J_{\omega }(f)$ any point with 
$\min \{ |p-b|\mid (\omega, b)\in P(f), b\in \CC \} >\mu .$  
Suppose that $J_{\omega }(f)$ is connected. 
Let 
$\psi :\CCI \setminus \overline{D(0,1)} 
\rightarrow A_{\omega }(f)$ be a biholomorphic map 
with $\psi (\infty )=\infty .$ For each 
$\theta \in \partial D(0,1)$, 
let $T(\theta )=\psi (\{ r\theta \mid 
1<r\leq \infty \} ).$ 
Suppose that there exist two elements 
$\theta _{1},\theta _{2}\in \partial D(0,1)$ with $\theta _{1}\neq \theta _{2}$ such that for each $i=1,2,$, 
$T(\theta _{i})$ lands at $p.$ 
Moreover, suppose that 
a connected component $V$ of 
$\CCI \setminus 
(T(\theta _{1})\cup T(\theta _{2})\cup \{ p\} )$ 
satisfies that 
diam $(V\cap K_{\omega }(f))\leq \delta .$ 
Furthermore, let 
$\g \in X$ be any point and suppose that 
there exists a sequence $\{ n_{k}\} _{k\in \NN }$ 
of positive integers such that 
$g^{n_{k}}(\g )\rightarrow \omega $ as 
$k\rightarrow \infty .$ Then, 
$J_{\g }(f)$ is not a quasicircle. 
\end{itemize}
\end{lem}
\begin{proof}
Let $\mu >0.$ Let $R>0$ with 
$\pi _{\CCI }(\tilde{J}(f))\subset D(0,R).$ 
Combining Lemma~\ref{invnormal} and 
Lemma~\ref{fibfundlem}-\ref{fibfundlemast}, 
we see that there exists a $\delta _{0}>0$ with\\  
$0<\delta _{0}<\frac{1}{20}\min 
\{ \inf _{x\in X}\mbox{diam }J_{x}(f),\mu \} $ such that the 
following statement holds:
\begin{itemize}
\item 
Let $x\in X$ be any point and $n\in \NN $ any element. 
Let 
$p\in D(0,R)$ be any point with 
$\min \{ |p-b|\mid (g^{n}(x), b)\in P(f), b\in \CC \} >\mu .$ 
Let $\phi :D(p,\mu )\rightarrow \CC $ be any well-defined 
inverse branch of 
$(f_{x,n})^{-1}$ on $D(p,\mu ).$ Let 
$A$ be any subset of $D(p,\frac{\mu }{2})$ with diam $A\leq \delta _{0}.$ 
Then, 
\begin{equation}
\label{060416lemastpfeq1}
\mbox{diam }\phi (A)\leq \frac{1}{10}\inf 
_{x\in X}\mbox{diam }J_{x}(f). 
\end{equation}
\end{itemize}
We set $\delta :=\frac{1}{10}\delta _{0}.$ 
Let $\omega \in X$ and $p\in J_{\omega }(f)$ with 
$\min \{ |p-b|\mid (\omega , b)\in P(f), b\in \CC \} >\mu .$ 
Suppose that $J_{\omega }(f)$ is connected and 
let $\psi :\CCI \setminus \overline{D(0,1)}\rightarrow 
A_{\omega }(f)$ be a biholomorphic map with 
$\psi (\infty )=\infty .$ 
Setting $T(\theta ):=\psi (\{ r\theta 
\mid 1<r\leq \infty \} )$ for each 
$\theta \in \partial D(0,1)$, 
suppose that there exist two elements 
$\theta _{1},\theta _{2}\in \partial D(0,1)$ with $\theta _{1}\neq \theta _{2}$  such that 
for each $i=1,2$, 
$T(\theta _{i})$ lands at $p.$ Moreover, 
suppose that a connected component 
$V$ of 
$\CCI \setminus (T(\theta _{1})\cup T(\theta _{2})\cup 
\{ p\} )$ satisfies that 
\begin{equation}
\label{060416lemastpfeq2}
\mbox{diam}(V\cap K_{\omega }(f))\leq \delta .
\end{equation}
Furthermore, let $\g \in X$ and suppose that 
there exists a sequence $\{ n_{k}\} _{k\in \NN }$ of 
positive integers such that 
$g^{n_{k}}(\g )\rightarrow \omega $ as 
$k\rightarrow \infty .$ 
 We now suppose that $J_{\g }(f)$ is a quasicircle, 
and we will deduce a contradiction. 
Since $g^{n_{k}}(\g )\rightarrow \omega $ as 
$k\rightarrow \infty $, 
we obtain 
\begin{equation}
\label{060416lemastpfeq3}
\max \{ d_{e}(b,K_{\omega }(f))\mid b\in J_{g^{n_{k}}(\g )}(f)\} 
\rightarrow 0 \mbox{ as } k\rightarrow \infty .
\end{equation}  
We take a point $a\in V\cap J_{\omega }(f)$ and fix it. 
By Lemma~\ref{fibfundlem}-\ref{fibfundlem2}, 
there exists a number $k_{0}\in \NN $ such that 
for each $k\geq k_{0}$, 
there exists a point $y_{k}$ satisfying that 
\begin{equation}
\label{060416lemastpfeq4}
y_{k}\in J_{g^{n_{k}}(\g )}(f)\cap D(a,\frac{|a-p|}{10k}).
\end{equation}
Let $V'$ be the connected component 
of $\CCI \setminus (T(\theta _{1})\cup T(\theta _{2})\cup 
\{ p\} )$ with $V'\neq V.$ 
Then, by \cite[Lemma 17.5]{M}, 
\begin{equation}
\label{060416lemastpfeq5}
V'\cap J_{\omega }(f)\neq \emptyset .
\end{equation}
Combining (\ref{060416lemastpfeq5}) and 
Lemma~\ref{fibfundlem}-\ref{fibfundlem2}, 
we see that there exists a $k_{1}(\geq k_{0})\in \NN $
 such that for each $k\geq k_{1}$, 
\begin{equation}
\label{060416lemastpfeq6}
V'\cap J_{g^{n_{k}}(\g )}(f)\neq \emptyset .
\end{equation}  
By assumption and Lemma~\ref{060416lema}, 
for each $k\geq k_{1}$, 
$J_{g^{n_{k}}(\g )}(f)$ is a Jordan curve. 
Combining it with (\ref{060416lemastpfeq4}) and 
(\ref{060416lemastpfeq6}), there exists a 
$k_{2}(\geq k_{1})\in \NN $ satisfying that 
for each $k\geq k_{2}$,  
there exists a smallest closed subarc 
$\xi _{k}$ of 
$J_{g^{n_{k}}(\g )}(f)\cong S^{1}$ such that 
$y_{k}\in \xi _{k}$, $\xi _{k}\subset \overline{V}, $ 
$\sharp (\xi _{k}\cap (T(\theta _{1})\cup T(\theta _{2})\cup 
\{ p\} ))=2,$ and such that 
$\xi _{k}\neq J_{g^{n_{k}}(\g )}(f).$ 
For each $k\geq k_{2}$, let 
$y_{k,1}$ and $y_{k,2}$ be the two points such that 
$\{ y_{k,1},y_{k,2}\} =
\xi _{k}\cap 
(T(\theta _{1})\cup T(\theta _{2})\cup \{ p\} ).$ 
Then, (\ref{060416lemastpfeq3}) implies that 
\begin{equation}
\label{060416lemastpfeq7}
y_{k,i}\rightarrow p \mbox{ as }k\rightarrow \infty , \mbox{ for each }i=1,2.
\end{equation}
Combining that 
$\xi _{k}\subset V\cup \{ y_{k,1},y_{k,2}\} $, 
(\ref{060416lemastpfeq3}), and 
(\ref{060416lemastpfeq2}), we get that 
there exists a $k_{3}(\geq k_{2})\in \NN $ such that 
for each $k\geq k_{3}$, 
\begin{equation}
\label{060416lemastpfeq8}
\mbox{diam }\xi _{k}\leq \frac{\delta _{0}}{2}.
\end{equation}
Moreover, combining (\ref{060416lemastpfeq4}) and 
(\ref{060416lemastpfeq7}), 
we see that there exists a constant 
$C>0$ such that 
\begin{equation}
\label{060416lemastpfeq9}
\mbox{diam }\xi _{k}>C.
\end{equation}
Combining (\ref{060416lemastpfeq7}), (\ref{060416lemastpfeq8}), 
and (\ref{060416lemastpfeq9}), 
we may assume that there exists a constant $C>0$ such that 
for each $k\in \NN $, 
\begin{equation}
\label{060416lemastpfeq10}
C<\mbox{ diam }\xi _{k}\leq \frac{\delta _{0}}{2}
\mbox{ and } 
\xi _{k}\subset D(p,\delta _{0}).
\end{equation}
By Lemma~\ref{060416lema}, each connected component $v$ of 
$(f_{\g ,n_{k}})^{-1}(\xi _{k})$ is a subarc of 
$J_{\g }(f)\cong S^{1}$ and 
$f_{\g ,n_{k}}:v\rightarrow \xi _{k}$ is a homeomorphism.
For each $k\in \NN $, 
let $\lambda _{k}$ be a connected component of 
$(f_{\g ,n_{k}})^{-1}(\xi _{k})$, and let 
$z_{k,1},z_{k,2}\in \lambda _{k}$ be the two endpoints 
of $\lambda _{k}$ such that 
$f_{\g ,n_{k}}(z_{k,1})=y_{k,1}$ and 
$f_{\g ,n_{k}}(z_{k,2})=y_{k,2}.$ Then, combining  
(\ref{060416lemastpfeq1}) and (\ref{060416lemastpfeq10}), 
we obtain 
\begin{equation}
\label{060416lemastpfeq11}
\mbox{ diam} \lambda _{k}<\mbox{ diam }(J_{\g }(f)\setminus 
\lambda _{k}), \mbox{ for each large }k\in \NN .
\end{equation}
Moreover, combining (\ref{060416lemastpfeq7}), 
(\ref{060416lemastpfeq10}), and Koebe distortion theorem, 
it follows that 
\begin{equation}
\label{060416lemastpfeq12}
\frac{\mbox{diam }\lambda _{k}}{|z_{k,1}-z_{k,2}|}\rightarrow \infty 
\mbox{ as }k\rightarrow \infty .
\end{equation}
Combining (\ref{060416lemastpfeq11}) and (\ref{060416lemastpfeq12}), 
we conclude that $J_{\g }(f)$ cannot be a quasicircle, 
since we have the following well-known fact:\\ 
Fact (\cite[Chapter 2]{LV}): 
Let $\xi $ be a Jordan curve in $\CC .$ Then, 
$\xi $ is a quasicircle if and only if 
there exists a constant $K>0 $ such that 
for each $z_{1},z_{2}\in \xi $ with $z_{1}\neq z_{2}$, 
we have $\frac{\mbox{diam }\lambda (z_{1},z_{2})}{|z_{1}-z_{2}|}\leq K$, 
where $\lambda (z_{1},z_{2})$ denotes the smallest closed subarc 
of $\xi $ such that 
$z_{1},z_{2}\in \lambda (z_{1},z_{2})$ and such that 
diam $\lambda (z_{1},z_{2})<$ 
diam $(\xi \setminus \lambda (z_{1},z_{2})).$ 

 Hence, we have proved Lemma~\ref{060416lemast}. 
\end{proof}

 We now demonstrate Theorem~\ref{mainthjbnq}-\ref{mainthjbnq1}.\\ 
{\bf Proof of Theorem~\ref{mainthjbnq}-\ref{mainthjbnq1}:} 
Let $\g $ be as in Theorem~\ref{mainthjbnq}-\ref{mainthjbnq1}. 
Then, by Theorem~\ref{mainth3-2},
$J_{\g }(f)$ is a Jordan curve. 
Moreover, setting $h=h_{m}\circ \cdots \circ h_{1}$, 
since $h$ is hyperbolic and $J(h)$ is not a quasicircle, 
$J(h)$ is not a Jordan curve.
Combining it with 
Lemma~\ref{060416lemast} and Lemma~\ref{jbnqlem1}, it follows that 
$J_{\g }(f)$ is not a quasicircle. 
Moreover, $A_{\g }(f)$ is a John domain 
(cf. \cite[Theorem 1.12]{S4}). 
Combining the above arguments with \cite[Theorem 9.3]{NV}, 
we conclude that the bounded component $U_{\g }$ of 
$F_{\g }(f)$ is not a John domain. 

 Thus, we have proved Theorem~\ref{mainthjbnq}-\ref{mainthjbnq1}.
\qed 

\ 

 We now demonstrate Theorem~\ref{mainthjbnq}-\ref{mainthjbnq2}.\\ 
{\bf Proof of Theorem~\ref{mainthjbnq}-\ref{mainthjbnq2}:}
Let $\rho _{0}, \beta ,\g $ be as in 
Theorem~\ref{mainthjbnq}-\ref{mainthjbnq2}.
By Theorem~\ref{mainth3-2}, 
$J_{\g }(f)$ is a Jordan curve. 
By Theorem~\ref{mainth2}-\ref{mainth2-4}, 
we have 
$\emptyset \neq $ int$(\hat{K}(G))\subset $ int$(K(h)).$ 
Moreover, $h$ is semi-hyperbolic. Hence, 
$h$ has an attracting periodic point in $K(h).$ 
Combining Lemma~\ref{060416lemast} and 
Lemma~\ref{060414astlem1}, we get that 
$J_{\g }(f)$ is not a quasicircle. 
Combining it with the argument in the proof of 
Theorem~\ref{mainthjbnq}-\ref{mainthjbnq1}, 
it follows that 
$A_{\g }(f)$ is a John domain, but the bounded component 
$U_{\g }$ of $F_{\g }(f)$ is not a John domain. 

 Thus, we have proved Theorem~\ref{mainthjbnq}-\ref{mainthjbnq2}.
\qed 

\subsection{Proofs of results in \ref{random}}
\label{Proofs of random}
In this subsection, we will demonstrate results in 
Section~\ref{random}. 

 we now prove Corollary~\ref{rancor1}.\\ 
{\bf Proof of Corollary~\ref{rancor1}:} 
By Remark~\ref{jminrem}, there exists 
a compact subset 
$S$ of $\G \setminus \G _{\min }$ such that 
the interior of $S$ with respect to the 
space $\G $ is not empty.  
Let ${\cal U}:=R(\G ,S).$ Then, 
it is easy to see that 
${\cal U}$ is residual in $\GN $, and that 
for each Borel probability measure 
$\tau $ on Poly$_{\deg \geq 2}$ with $\G _{\tau }=\G $, 
we have $\tilde{\tau }({\cal U})=1.$ Moreover, 
by Theorem~\ref{mainth3}-\ref{mainth3-0} and 
Theorem~\ref{mainth3}-\ref{mainth3-1}, 
each $\g \in {\cal U}$ satisfies properties 
\ref{rancor1-1},\ref{rancor1-2},\ref{rancor1-3}, and 
\ref{rancor1-4} in Corollary~\ref{rancor1}. Hence, 
we have proved Corollary~\ref{rancor1}.
\qed  

\ 

To demonstrate Theorem~\ref{mainthran1}, 
we need several lemmas.
\begin{lem}
\label{mainthranlem1}
Let $\G $ be a compact set in {\em Poly}$_{\deg \geq 2}.$ 
Let $f:\GN \times \CCI \rightarrow \GN \times \CCI $ be the  
skew product associated with the family $\G .$ 
Let $G$ be the polynomial semigroup generated by $\G .$ 
Suppose that $G\in {\cal G}$ and that $G$ is semi-hyperbolic. 
Moreover, suppose that there exist two elements 
$\alpha ,\beta \in \GN $ such that 
$J_{\beta }(f)<J_{\alpha }(f).$ 
Let $\g \in \GN $ and suppose that there exists a strictly increasing sequence 
$\{ n_{k}\} _{k\in \NN }$ of positive integers such that 
$\sigma ^{n_{k}}(\g )\rightarrow \alpha $ as $k\rightarrow \infty .$ 
Then, $J_{\g }(f)$ is a Jordan curve.  
\end{lem}
\begin{proof}
Since $G$ is semi-hyperbolic, \cite[Theorem 2.14-(4)]{S1} implies that 
\begin{equation}
\label{mainthranlem1eq1}
J_{\sigma ^{n_{k}}(\g )}(f)\rightarrow J_{\alpha }(f) 
\mbox{ as } k\rightarrow \infty, 
\end{equation}
with respect to the Hausdorff topology in the space of non-empty 
compact subsets of $\CCI .$ 
Combining it with Lemma~\ref{fiborder}, 
we see that there exists a number $k_{0}\in \NN $ such that 
for each $k\geq k_{0}$, 
\begin{equation}
\label{mainthranlem1eq2}
J_{\beta }(f)<J_{\sigma ^{n_{k}}(\g )}(f).
\end{equation}
We will show the following claim.\\ 
Claim: int$(K_{\g }(f))$ is connected.
 
 To show this claim, suppose that there exist two distinct components 
 $U_{1}$ and $U_{2}$ of 
 int$(K_{\g }(f)).$ 
Let $y_{i}\in U_{i}$ be a point, for each $i=1,2.$ 
Let $\epsilon >0$ be a number such that 
$\overline{D(K_{\beta }(f),\epsilon )}$ 
is included in a connected component $U$ of 
int$(K_{\alpha }(f)).$ 
Then, combining \cite[Theorem 2.14-(5)]{S1} and Lemma~\ref{constlimlem}, 
we get that there exists a number $k_{1}\in \NN $ with 
$k_{1}\geq k_{0}$ such that for each $k\geq k_{1}$ and each $i=1,2$, 
\begin{equation}
\label{mainthranlem1eq3}
f_{\g ,n_{k}}(y_{i})\in D(P^{\ast }(G),\epsilon )\subset 
\overline{D(K_{\beta }(f),\epsilon )}\subset U.
\end{equation}
Combining (\ref{mainthranlem1eq3}), (\ref{mainthranlem1eq1})
 and (\ref{mainthranlem1eq2}), 
we get that there exists a number $k_{2}\in \NN $ with 
$k_{2}\geq k_{1}$ such that for each $k\geq k_{2}$, 
\begin{equation}
\label{mainthranlem1eq4}
f_{\g ,n_{k}}(U_{1})=f_{\g ,n_{k}}(U_{2})=V_{k},
\end{equation}
where $V_{k}$ denotes the connected component of 
int$(K_{\sigma ^{n_{k}}(\g )}(f))$ containing 
$J_{\beta }(f).$ 
From (\ref{mainthranlem1eq2}) and (\ref{mainthranlem1eq4}), 
it follows that 
\begin{equation}
\label{mainthranlem1eq5}
(f_{\g ,n_{k}})^{-1}(J_{\beta }(f))
\subset \mbox{int}(K_{\g }(f)) 
\mbox{ and }
(f_{\g ,n_{k}})^{-1}(J_{\beta }(f))\cap 
U_{i}\neq \emptyset \ (i=1,2), 
\end{equation}
which implies that 
\begin{equation}
\label{mainthranlem1eq6}
(f_{\g ,n_{k}})^{-1}(J_{\beta }(f))
\mbox{ is disconnected.}
\end{equation}
For each $k\geq k_{2}$, let 
$\omega ^{k}:= (\g _{1},\ldots ,\g _{n_{k}},\beta _{1},\beta _{2},\ldots )
\in \GN .$ Then for each $k\geq k_{2}$, 
\begin{equation}
\label{mainthranlem1eq7}
(f_{\g ,n_{k}})^{-1}(J_{\beta }(f))=
J_{\omega ^{k}}(f).
\end{equation}
Since $G\in {\cal G}$, combining (\ref{mainthranlem1eq6}), 
(\ref{mainthranlem1eq7}) and Lemma~\ref{fibconnlem}
yields a contradiction. 
Hence, we have proved the claim.

 From the above claim and Proposition~\ref{shonecomp}, 
 it follows that 
$J_{\g }(f)$ is a Jordan curve.  
\end{proof}
\begin{lem}
\label{mainthranlem2}
Let $\G $ be a non-empty compact subset of 
{\em Poly}$_{\deg \geq 2}.$ 
Let $f:\GN \times \CCI \rightarrow \GN \times \CCI $
 be the skew product associated with the family $\G $ of 
 polynomials. 
 Let $G$ be the polynomial semigroup generated by $\G .$ 
Let $\alpha ,\rho \in \GN $ be two elements. 
Suppose that $G\in {\cal G}$, that $G$ is semi-hyperbolic, 
that $\alpha $ is a periodic point of $\sigma :\GN \rightarrow \GN $, 
that $J_{\alpha }(f)$ is a quasicircle, 
and that $J_{\rho }(f)$ is not a Jordan curve. 
Then, for each $\epsilon >0$, there exist $n\in \NN $ and 
two elements $\theta _{1},\theta _{2}\in \partial D(0,1)$ 
with $\theta _{1}\neq \theta _{2}$ 
satisfying all of the following.
\begin{enumerate}
\item \label{mainthranlem2-1}
Let $\omega =(\alpha _{1},\ldots ,\alpha _{n},\rho _{1},\rho _{2},\ldots )
\in \GN $ and let $\psi :\CCI \setminus 
\overline{D(0,1)}\cong A_{\omega }(f)$ be a biholomorphic map 
with $\psi (\infty )=\infty .$ Moreover, 
for each $i=1,2$, let 
$T(\theta _{i}):=\psi (\{ r\theta _{i}\mid 1<r\leq \infty \} ).$ 
Then, there exists a point $p\in J_{\omega }(f)$ 
such that for each $i=1,2$, $T(\theta _{i})$ lands at $p.$ 
\item \label{mainthranlem2-2}
Let $V_{1}$ and $V_{2}$ be the two connected components of 
$\CCI \setminus (T(\theta _{1})\cup T(\theta _{2})\cup \{ p\} ).$ 
Then, for each $i=1,2$, 
$V_{i}\cap J_{\omega }(f)\neq \emptyset .$ 
Moreover, there exists an $i\in \{ 1,2\} $ such that 
{\em diam} $(V_{i}\cap K_{\omega }(f))\leq \epsilon $, 
and such that $V_{i}\cap J_{\omega }(f)\subset 
D(J_{\alpha }(f),\epsilon ).$  
\end{enumerate}   
\end{lem} 
\begin{proof}
For each $\g \in \GN $, let 
$\psi _{\g }:\CCI \setminus \overline{D(0,1)}\cong 
A_{\g }(f)$ be a biholomorphic map with 
$\psi _{\g }(\infty )=\infty .$ Moreover, 
for each $\theta \in \partial D(0,1)$, 
let 
$T_{\g }(\theta ):= \psi _{\g }(\{ r\theta \mid 1<r\leq \infty \} ).$ 
Since $G$ is semi-hyperbolic, 
combining \cite[Theorem 1.12]{S4}, 
Lemma~\ref{fibconnlem}, and \cite[page 26]{NV}, 
we see that for each $\g \in \GN $, 
$J_{\g }(f)$ is locally connected. 
Hence, for each $\g \in \GN $, 
$\psi _{\g }$ extends continuously over 
$\CCI \setminus D(0,1)$ such that 
$\psi _{\g }(\partial D(0,1))=J_{\g }(f).$ 
Moreover, since $G\in {\cal G}$, 
it is easy to see that for each $\g \in \GN $, 
there exists a number $a_{\g }\in \CC $ with $|a_{\g }|=1$ such that 
for each $z\in \CCI \setminus \overline{D(0,1)}$, 
we have $\psi _{\sigma (\g )}^{-1}\circ f_{\g ,1}\circ 
\psi _{\g }(z)=a_{\g }z^{d(\g )}.$

 Let $m\in \NN $ be an integer such that $\sigma ^{m}(\alpha )=
\alpha $ and let $h:= \alpha _{m}\circ \cdots \circ \alpha _{1}.$ 
Moreover, for each $n\in \NN $, 
we set $\omega ^{n}:= (\alpha _{1},\ldots ,\alpha _{mn},
\rho _{1},\rho _{2},\ldots )\in \GN .$ 
Then, $\omega ^{n}\rightarrow \alpha $ in $\GN $ as 
$n\rightarrow \infty .$ 
Combining it with \cite[Theorem 2.14-(4)]{S1}, 
we obtain  
\begin{equation}
\label{mainthranlem2eq1}
J_{\omega ^{n}}(f)\rightarrow 
J_{\alpha }(f) \mbox{ as } n\rightarrow \infty , 
\end{equation} 
with respect to the Hausdorff topology.
Let $\xi $ be a Jordan curve in int$(K(h))$ 
such that $P^{\ast }(\langle h\rangle )$ is included 
in the bounded component $B$ of $\CC \setminus \xi .$ 
By (\ref{mainthranlem2eq1}), 
there exists a $k\in \NN $ such that 
$J_{\omega ^{k}}(f)\cap (\xi \cup B) =\emptyset .$ 
We now show the following claim.\\ 
Claim 1: $\xi \subset $ int$(K_{\omega ^{k}}(f)).$ 

 To show this claim, suppose that $\xi $ is included in 
$A_{\omega ^{k}}(f)=
\CCI \setminus (K_{\omega ^{k}}(f)).$ 
Then, it implies that 
$f_{\omega ^{k},u}\rightarrow \infty $ on 
$P^{\ast }(\langle h\rangle )$ as $u\rightarrow \infty .$ 
 However, this is a 
contradiction, since $G\in {\cal G}.$ 
Hence, we have shown Claim 1.

 By Claim 1, we see that $P^{\ast }(\langle h\rangle )$ is 
 included in a bounded component $B_{0}$ of 
 int$(K_{\omega ^{k}}(f)).$ 
We now show the following claim.\\ 
Claim 2: $J_{\omega ^{k}}(f)$ is not a 
Jordan curve. 

 To show this claim, suppose that 
$J_{\omega ^{k}}(f)$ is a Jordan curve. Then, 
Lemma~\ref{060416lema} implies 
that $J_{\rho }(f)$ is a Jordan curve. 
However, this is a contradiction. Hence, we have shown Claim 2.

By Claim 2, 
there exist two distinct elements $t_{1},t_{2}\in \partial D(0,1)$ 
 and a point $p_{0}\in J_{\omega ^{k} }(f)$ 
 such that for each $i=1,2$, 
 $T_{\rho }(t_{i})$ lands at the point $p_{0}.$ 
Let $W_{0}$ be the connected component of  
$\CCI \setminus (T_{\rho }(t_{1})\cup T_{\rho }(t_{2})\cup 
\{ p_{0}\} )$ such that 
$W_{0}$ does not contain $B_{0}.$ 
Then, we have 
\begin{equation}
\label{mainthranlem2eq3}
\overline{W_{0}}\cap P^{\ast }(\langle h\rangle )=\emptyset .
\end{equation}
For each $j\in \NN $, we take a connected component 
$W_{j}$ of $(h^{j})^{-1}(W_{0}).$ 
Then, $h^{j}:W_{j}\rightarrow W_{0}$ is biholomorphic.
We set $\zeta _{j}:= (h^{j}|_{W_{j}})^{-1} $ on 
$W_{0}.$ By (\ref{mainthranlem2eq3}), 
 there exists a number $R>0$ and a number $a >0$ 
 such that 
for each $j$, $\zeta _{j}$ is analytically continued to a univalent function 
$\tilde{\zeta _{j}}:B(\overline{W_{0}\cap D(0,R)}, a )\rightarrow 
\CCI $ and 
$W_{j}\cap (J_{\omega ^{k+j}}(f))\subset 
\tilde{\zeta _{j}}(W_{0}\cap  
D(0,R)).$ 
Hence, we obtain
\begin{equation}
\label{mainthranlem2eq4}
\mbox{diam }(W_{j}\cap K_{\omega ^{k+j}}(f))=
\mbox{diam }(W_{j}\cap J_{\omega ^{k+j}}(f))\rightarrow 
0 \mbox{ as } j\rightarrow \infty .
\end{equation}
Combining (\ref{mainthranlem2eq1}) and (\ref{mainthranlem2eq4}), 
there exists an $s\in \NN $
 such that 
diam $(W_{s}\cap K_{\omega ^{k+s}}(f))\leq \epsilon $, and such that 
$W_{s}\cap J_{\omega ^{k+s}}(f)\subset D(J_{\alpha }(f),\epsilon ).$ 

 Each connected component 
of $(\partial W_{s})\cap \CC $ is a connected component of 
\\ $(h^{s})^{-1}((T_{\omega ^{k}}(t_{1})\cup T_{\omega ^{k}}(t_{2})
\cup \{ p_{0}\} )\cap \CC )$, and there are some 
$u_{1},\ldots ,u_{v}\in \partial D(0,1)$ such that 
$\partial W_{s}=\overline{\cup _{i=1}^{v}T_{\omega ^{k+s}}(u_{i})}.$ 
Hence, $W_{s}$ is a Jordan domain. Therefore, 
$h^{s}:\overline{W_{s}}\rightarrow \overline{W_{0}}$ is a 
homeomorphism. 
Thus, $h^{s}:(\partial W_{s})\cap \CC \rightarrow 
(\partial W_{0})\cap \CC $ is a homeomorphism. 
Hence, $(\partial W_{s})\cap \CC $ is connected. 
It follows that there exist 
two elements $\theta _{1},\theta _{2}\in \partial D(0,1)$ 
with $\theta _{1}\neq \theta _{2}$  and 
a point $p\in J_{\omega ^{k+s}}(f)$ 
such that 
$\partial W_{s}=T_{\omega ^{k+s}}(\theta _{1})\cup 
T_{\omega ^{k+s}}(\theta _{2})\cup \{ p\} $, and such that 
for each $i=1,2$, $T_{\omega ^{k+s}}(\theta _{i})$ 
lands at the point $p.$ 
 By \cite[Lemma 17.5]{M}, 
each of two connected components of 
$\CCI \setminus 
( T_{\omega ^{k+s}}(\theta _{1})\cup 
T_{\omega ^{k+s}}(\theta _{2})\cup \{ p\} )$ intersects 
$J_{\omega ^{k+s}}(f).$ 

 Hence, we have proved Lemma~\ref{mainthranlem2}.
\end{proof}
\begin{lem}
\label{mainthranlem3}
Let $\G $ be a non-empty compact subset of 
{\em Poly}$_{\deg \geq 2}.$ 
Let $f:\GN \times \CCI \rightarrow \GN \times \CCI $
 be the skew product associated with the family $\G $ of 
 polynomials. 
 Let $G$ be the polynomial semigroup generated by $\G .$ 
Let $\alpha ,\beta ,\rho \in \GN $ be three elements. 
Suppose that $G\in {\cal G}$, that $G$ is semi-hyperbolic, 
that $\alpha $ is a periodic point of $\sigma :\GN \rightarrow \GN $, 
that $J_{\beta }(f)<J_{\alpha }(f)$, 
and that $J_{\rho }(f)$ is not a Jordan curve. 
Then, there exists an $n\in \NN $ such that 
setting $\omega := (\alpha _{1},\ldots ,\alpha _{n},
\rho _{1},\rho _{2},\ldots )\in \GN $ and 
${\cal U}:= \{ \g \in \GN 
\mid \exists \{ m_{j}\} _{j\in \NN }, \exists \{ n_{k}\} _{k\in \NN }, 
\sigma ^{m_{j}}(\g )\rightarrow \alpha ,\ 
\sigma ^{n_{k}}(\g )\rightarrow \omega \} $, 
we have that for each $\g \in {\cal U},$ 
$J_{\g }(f)$ is a Jordan curve but not a quasicircle,  
$A_{\g }(f)$ is a John domain, and the bounded component 
$U_{\g }$ of $F_{\g }(f)$ is not a John domain. 
\end{lem}
\begin{proof}

Let $p\in \NN $ be a number such that 
$\sigma ^{p}(\alpha )=\alpha $ and let 
$u :=\alpha _{p}\circ \cdots \circ \alpha _{1}.$ 
We show the following claim.\\ 
Claim 1: $J(u)$ is a quasicircle.

 To show this claim, by assumption,  
we have 
$J_{\beta }(f)<J(u).$ 
Let $\zeta := (\alpha _{1},\ldots $ $,\alpha _{p},\beta _{1},\beta _{2},\ldots 
)\in \GN .$ Then, 
we have $J_{\zeta }(f)=
u^{-1}(J_{\beta }(f)).$ Moreover, 
since $G\in {\cal G}$, we have that 
$J_{\zeta }(f)$ is connected. 
Hence, it follows that 
$u^{-1}(J_{\beta }(f))$ is connected. 
Let $U$ be a connected component of int$(K(u))$ containing 
$ J_{\beta }(f)$ and $V$ a connected component 
of int$(K(u))$ containing $u^{-1}(J_{\beta }(f)).$ 
By Lemma~\ref{fiborder}, it must hold that $U=V.$ 
Therefore, we obtain $u^{-1}(U)=U.$ Thus,  
int$(K(u))=U.$ Since $G$ is semi-hyperbolic, 
it follows that $J(u)$ is a quasicircle.  Hence, we have proved Claim 1.

 Let 
$\mu :=\frac{1}{3}\min 
 \{ |b-c|\mid b\in J_{\alpha }(f), c\in P^{\ast }(G)\} .$
Since $J_{\beta }(f)<J_{\alpha }(f)$, we have 
$P^{\ast }(G)\subset K_{\beta }(f).$ Hence, 
$\mu >0.$  
Applying Lemma~\ref{060416lemast} to the above 
$(f,\mu )$, let $\delta $ be the number 
in the statement of Lemma~\ref{060416lemast}. 
We set $\epsilon := \min \{ \delta ,\mu \} (>0).$
Applying Lemma~\ref{mainthranlem2} to 
the above $(\G ,\alpha , \rho , \epsilon )$, 
let $(n,\theta _{1},\theta _{2},\omega )$ be the 
element in the statement of Lemma~\ref{mainthranlem2}.
We set 
${\cal U}:= 
\{ \g \in \GN \mid 
\exists \{ m_{j}\} _{j\in \NN }, 
\exists \{ n_{k}\} _{k\in \NN }, 
\sigma ^{m_{j}}(\g )\rightarrow \alpha , 
\sigma ^{n_{k}}(\g )\rightarrow \omega \} .$ 
Then, combining 
the statement Lemma~\ref{060416lemast} and that of 
Lemma~\ref{mainthranlem2}, 
it follows that 
for any $\g \in {\cal U}$, 
$J_{\g }(f)$ is not a quasicircle. 
Moreover, by Lemma~\ref{mainthranlem1}, 
we see that 
for any $\g \in {\cal U}$, 
$J_{\g }(f)$ is a Jordan curve. 
Furthermore, combining the above argument, 
\cite[Theorem 1.12]{S4}, Lemma~\ref{fibconnlem}, 
and \cite[Theorem 9.3]{NV}, 
we see that for any $\g \in {\cal U}$, 
$A_{\g }(f)$ is a John domain, and 
the bounded component $U_{\g }$ of $F_{\g }(f)$ is not a John domain.
Therefore, we have proved Lemma~\ref{mainthranlem3}. 
\end{proof}

 We now demonstrate Theorem~\ref{mainthran1}.\\ 
{\bf Proof of Theorem~\ref{mainthran1}:} 
We suppose the assumption of Theorem~\ref{mainthran1}. 
We will consider several cases. 
First, we show the following claim.\\ 
Claim 1: If $J_{\g }(f)$ is a Jordan curve 
 for each $\g \in \GN $, then statement \ref{mainthran1-1} 
in Theorem~\ref{mainthran1} holds.

 To show this claim, 
 Lemma~\ref{060416lema} implies that  for each $\g \in X$, 
any critical point $v\in \pi ^{-1}(\{ \g \} )$ 
of $f_{\g }:\pi ^{-1}(\{ \g \} )\rightarrow 
\pi ^{-1}(\{ \sigma (\g )\} )$ 
(under the canonical identification $\pi ^{-1}(\{ \g \} )
\cong \pi ^{-1}(\{ \sigma (\g )\} )\cong \CCI $) 
belongs to 
$F^{\g }(f).$
Moreover, by \cite[Theorem 2.14-(2)]{S1}, 
$\tilde{J}(f)=\cup _{\g \in \GN }J^{\g }(f).$ 
Hence, it follows that 
$C(f)\subset \tilde{F}(f).$ Therefore, 
$C(f)$ is a compact subset of $\tilde{F}(f).$  
Since $f$ is semi-hyperbolic, 
\cite[Theorem 2.14-(5)]{S1} implies that 
$P(f)=\overline{\bigcup _{n\in \NN  }f^{n}(C(f))}\subset \tilde{F}(f).$ 
Hence, $f:\GN \times \CCI \rightarrow \GN \times  \CCI $ 
is hyperbolic. 
Combining it with Remark~\ref{hypskewsemigrrem}, 
we conclude that 
 $G$ is hyperbolic. 
  Moreover, Theorem~\ref{hypskewqc} implies that 
  there exists a constant $K\geq 1$ such that for each $\g \in \GN $, 
  $J_{\g }(f)$ is a $K$-quasicircle.  
Hence, we have proved Claim 1.

 Next, we will show the following claim.\\ 
Claim 2: If $J_{\alpha }(f)\cap J_{\beta }(f)
\neq \emptyset $ 
for each $(\alpha ,\beta )\in \GN \times \GN $, 
then $J(G)$ is arcwise connected. 

 To show this claim, since $G$ is semi-hyperbolic, 
combining \cite[Theorem 1.12]{S4},
Lemma~\ref{fibconnlem}, and \cite[page 26]{NV},   
we get that for each $\g \in \GN $, 
$A_{\g }(f)$ is a John domain and 
$J_{\g }(f)$ is locally connected. 
In particular, for each $\g \in \GN $, 
\begin{equation}
\label{mainthran1pfeq2}
J_{\g }(f) \mbox{ is arcwise connected.} 
\end{equation}
Moreover, by \cite[Theorem 2.14-(2)]{S1}, 
we have 
\begin{equation}
\label{mainthran1pfeq3}
\tilde{J}(f)=\cup _{\g \in \GN }J^{\g }(f).
\end{equation} 
 Combining (\ref{mainthran1pfeq2}), (\ref{mainthran1pfeq3}) and 
Lemma~\ref{fiblem}-\ref{pic}, 
we conclude that $J(G)$ is arcwise connected. Hence, we have proved 
Claim 2. 

  Next, we will show the following claim. \\ 
Claim 3: If $J_{\alpha }(f)\cap J_{\beta }(f)
\neq \emptyset $ 
for each $(\alpha ,\beta )\in \GN \times \GN $, and if 
there exists an element $\rho \in \GN $ such that 
$J_{\rho }(f)$ is not a Jordan curve, 
then statement \ref{mainthran1-3} in Theorem~\ref{mainthran1} holds. 

 To show this claim, let 
 ${\cal V}:= \cup _{n\in \NN }(\sigma ^{n})^{-1}(\{ \rho \} ).$ 
 Then, ${\cal V}$ is a dense subset of $\GN .$ 
From Lemma~\ref{060416lema}, it follows that 
 for each $\g \in {\cal V}$, $J_{\g }(f)$ is not a Jordan 
 curve. Combining this result with Claim 2, we conclude that statement 
 \ref{mainthran1-3} in Theorem~\ref{mainthran1} holds. Hence, 
 we have proved Claim 3. 

 We now show the following claim. \\ 
Claim 4: If there exist two elements $\alpha ,\beta \in \GN $ such that 
$J_{\alpha }(f)\cap J_{\beta }(f)= 
\emptyset $, and if there exists an element $\rho \in \GN $ such that 
$J_{\rho }(f)$ is not a Jordan curve, then 
statement \ref{mainthran1-2} in Theorem~\ref{mainthran1} holds.

 To show this claim, using Lemma~\ref{fiborder}, 
 We may assume that $J_{\beta }(f)<
 J_{\alpha }(f).$ 
Combining this, Lemma~\ref{fiborder}, 
\cite[Theorem 2.14-(4)]{S1}, and that 
the set of all periodic points of $\sigma $ in 
$\GN $ is dense in $\GN $, 
we may assume further that 
$\alpha $ is a periodic point of $\sigma .$ 
Applying Lemma~\ref{mainthranlem3} to 
$(\G, \alpha ,\beta , \rho )$ above, 
let $n\in \NN $ be the element 
in the statement of Lemma~\ref{mainthranlem3}, and  
we set  
$\omega =(\alpha _{1},\ldots ,\alpha _{n}, \rho _{1},\rho _{2},\ldots 
)\in \GN $ and 
${\cal U}:= 
\{ \g \in \GN \mid 
\exists (m_{j}), \exists (n_{k}), 
\sigma ^{m_{j}}(\g )\rightarrow \alpha ,
\sigma ^{n_{k}}(\g )\rightarrow \omega \} .$ 
Then, by the statement of Lemma~\ref{mainthranlem3}, 
we have that for each $\g \in {\cal U}$, 
$J_{\g }(f)$ is a Jordan curve but not a quasicircle, 
$A_{\g }(f)$ is a John domain, and the 
bounded component $U_{\g }$ of $F_{\g }(f)$ is not 
a John domain.   
Moreover, ${\cal U} $ is residual in $\GN $, and 
for any Borel probability measure 
$\tau $ on Poly$_{\deg \geq 2}$ with $\G _{\tau }=\G$, 
we have $\tilde{\tau }({\cal U})=1.$ 
Furthermore, let ${\cal 
V}:= \cup _{n\in \NN }(\sigma ^{n})^{-1}(\{ \rho \} ).$ Then, 
${\cal V}$ is a dense subset of $\GN $, 
and the argument in the proof 
of Claim 3 implies that for each $\g \in {\cal V}$, 
$J_{\g }(f)$ is not a Jordan curve. 
Hence, we have proved Claim 4.

 Combining Claims 1,2,3 and 4, Theorem~\ref{mainthran1} follows.    
\qed 

\ 

We now demonstrate Corollary~\ref{rancor2}.\\ 
{\bf Proof of Corollary~\ref{rancor2}:} 
From Theorem~\ref{mainthran1}, 
Corollary~\ref{rancor2} immediately follows.  
\qed 

\ 

 To demonstrate Theorem~\ref{mainthran2}, we need several lemmas.

\noindent {\bf Notation:} 
For a subset $A$ of $\CCI $, we denote by ${\cal C}(A)$ the set of 
all connected components of $A.$ 
\begin{lem}
\label{mainthran2lem1}
Let $f:X\times \CCI \rightarrow X\times \CCI $ be a polynomial 
skew product over $g:X\rightarrow X$ such that for each 
$x\in X$, $d(x)\geq 2.$ Let $\alpha \in X$ be a point. 
Suppose that $2\leq \sharp \left( {\cal C}(  
\mbox{{\em int}}(K_{\alpha }(f)))\right)<\infty .$ 
Then, 
$\sharp \left( {\cal C}(  
\mbox{{\em int}}(K_{g(\alpha )}(f)))\right)
<\sharp \left( {\cal C}(  
\mbox{{\em int}}(K_{\alpha }(f)))\right).$ 
In particular, there exists an $n\in \NN $ such that 
{\em int}$(K_{g^{n}(\alpha )}(f))$ is a non-empty connected set.
\end{lem} 
\begin{proof}
Suppose that 
$2\leq \sharp ({\cal C}(\mbox{int}(K_{g(\alpha )}(f))))=
\sharp ({\cal C}(\mbox{int}(K_{\alpha }(f))))<\infty .$ We 
will deduce a contradiction. 
Let $\{ V_{j}\} _{j=1}^{r}={\cal C}(\mbox{int}(K_{g(\alpha )}(f)))$, 
where $2\leq r<\infty .$  
Then, by the assumption above, 
we have that 
${\cal C}(\mbox{int}(K_{g(\alpha )}(f)))
=\{ f_{\alpha ,1}(V_{j})\} _{j=1}^{r}.$ 
For each $j=1,\ldots ,r$, let $p_{j}$ be the 
number of critical points of 
$f_{\alpha ,1}:V_{j}\rightarrow 
f_{\alpha ,1}(V_{j})$ counting multiplicities. 
Then, by the Riemann-Hurwitz formula, 
we have that for each $j=1,\ldots ,r$, 
$\chi (V_{j})+p_{j}=d\chi (f_{\alpha ,1}(V_{j}))$, 
where $\chi (\cdot) $ denotes the Euler number and 
$d:=\deg (f_{\alpha ,1}).$ 
Since $\chi (V_{j})=\chi (f_{\alpha ,1}(V_{j}))=1$ for each $j$, 
we obtain $r+\sum _{j=1}^{r}p_{j}=rd.$ 
Since $\sum _{j=1}^{r}p_{j}\leq d-1$, 
it follows that $rd-r\leq d-1.$ 
Therefore, we obtain $r\leq 1$, which is a contradiction. 
Thus, we have proved Lemma~\ref{mainthran2lem1}.
\end{proof}
\begin{lem}
\label{mainthran2lem2}
Let $f:X\times \CCI \rightarrow X\times \CCI $ be a 
polynomial skew product over $g:X\rightarrow X$ such that 
for each $x\in X$, $d(x)\geq 2.$ Let 
$\omega \in X$ be a point. 
Suppose that $f$ is hyperbolic, 
that $\pi _{\CCI }(P(f))\cap \CC $ is bounded in $\CC $, 
and that {\em int}$(K_{\omega }(f))$ is not connected. Then, 
there exist infinitely many connected components of 
{\em int}$(K_{\omega }(f)).$ 
\end{lem} 
\begin{proof}
Suppose that 
$2\leq \sharp ({\cal C}(\mbox{int}(K_{\omega }(f))))<\infty .$ 
Then, by Lemma~\ref{mainthran2lem1}, 
there exists an $n\in \NN $ such that 
int$(K_{g^{n}(\omega )}(f))$ is connected. 
We set $U:=$ int$(K_{g^{n}(\omega )}(f)).$ 
Let $\{ V_{j}\} _{j=1}^{r}$ be the set of 
all connected components 
of $(f_{\omega ,n})^{-1}(U).$ Since 
int$(K_{\omega }(f))$ is not connected, 
we have $r\geq 2.$   
For each $j=1,\ldots ,r$, we set 
$d_{j}:=\deg (f_{\omega ,n}:V_{j}\rightarrow 
U).$ Moreover, 
we denote by $p_{j}$ the number of critical points 
of $f_{\omega ,n}:V_{j}\rightarrow U$ 
counting multiplicities. 
Then, by the Riemann-Hurwitz formula, we see that for each 
$j=1,\ldots ,r$, 
$\chi (V_{j})+p_{j}=d_{j}\chi (U).$
Since $\chi (V_{j})=\chi (U)=1$ for each $j=1,\ldots ,r$, 
it follows that 
\begin{equation}
\label{mainthran2lem2eq2}
r+\sum _{j=1}^{r}p_{j}=d,
\end{equation}
where $d:= \deg (f_{\omega ,n}).$ 
Since $f$ is hyperbolic and $\pi _{\CCI }(P(f))\cap \CC $ is bounded in 
$\CC $, we have $\sum _{j=1}^{r}p_{j}=d-1.$ Combining it with 
(\ref{mainthran2lem2eq2}), we obtain $r=1$, 
which is a contradiction. Hence, we have proved Lemma~\ref{mainthran2lem2}. 
\end{proof}
\begin{lem}
\label{mainthran2lem3}
Let 
$f:X \times \CCI \rightarrow X \times \CCI $ be a polynomial skew 
product over $g:X\rightarrow X.$ 
Let 
$\alpha \in X $ be an element.  
Suppose that $\pi _{\CCI }(P(f))\cap \CC $ is bounded in $\CC $, that $f$ is hyperbolic, and 
that {\em int}$(K_{\alpha }(f)))$ is connected. Then, there exists a 
neighborhood ${\cal U}_{0}$ of $\alpha $ in $X$ satisfying the 
following.
\begin{itemize}
\item Let $\g \in X $ and suppose that there exists a sequence 
$\{ m_{j}\} _{j\in \NN }\subset \NN , m_{j}\rightarrow \infty $ such that 
for each $j\in \NN $, $g ^{m_{j}}(\g )\in {\cal U}_{0}.$ Then, 
$J_{\g }(f)$ is a Jordan curve.  
\end{itemize}  
\end{lem}
\begin{proof}
Let $P^{\ast }(f):= P(f)\setminus \pi _{\CCI }^{-1}(\{ \infty \} ).$ 
By assumption, we have 
$\pi _{\CCI }(P^{\ast }(f)\cap \pi ^{-1}(\{ \alpha \} ))\subset \mbox{ int}(K_{\alpha }(f)).$ 
Since int$(K_{\alpha }(f))$ is simply connected, 
there exists a Jordan curve $\xi $ in 
int$(K_{\alpha }(f))$ such that 
$\pi _{\CCI }(P^{\ast }(f)\cap \pi ^{-1}(\{ \alpha \} ))$ is included in the bounded component $B$ of 
$\CC \setminus \xi .$ 
Since $f$ is hyperbolic, \cite[Theorem 2.14-(4)]{S1} implies that 
the map $x \mapsto J_{x }(f)$ is continuous with respect to 
the Hausdorff topology. 
Hence, there exists a neighborhood ${\cal U}_{0}$ of $\alpha $ 
in $X $ such that for each $\beta \in {\cal U}_{0}$, 
$J_{\beta }(f)\cap (\xi \cup B)=\emptyset .$ Moreover, since $P(f) $ is compact, 
shrinking ${\cal U}_{0}$ if necessary, we may assume that for each $\beta \in {\cal U}_{0}$, 
$\pi _{\CCI }(P^{\ast }(f)\cap \pi ^{-1}(\{ \beta \} ))\subset B.$   
Since $\pi _{\CCI }(P(f))\cap \CC $ is bounded in $\CC $, it follows that for each $\beta \in {\cal U}_{0}$, 
$\xi <J_{\beta }(f).$ Hence, for each $\beta \in {\cal U}_{0}$, 
there exists a connected component $V_{\beta }$ of 
int$(K_{\beta }(f))$ such that 
\begin{equation}
\label{mainthran2lem3eq1}
\pi _{\CCI }(P^{\ast }(f)\cap \pi ^{-1}(\{ \beta \} ))\subset V_{\beta }.
\end{equation} 
Let $\g \in X$ be an element and suppose that 
there exists a sequence $\{ m_{j}\} _{j\in \NN }\subset \NN , 
m_{j}\rightarrow \infty $
 such that for each $j\in \NN $, $g ^{m_{j}}(\g )\in 
 {\cal U}_{0}.$ 
We will show that int$(K_{\g }(f))$ is connected. 
Suppose that there exist two distinct connected components  
$W_{1}$ and $W_{2}$ of int$(K_{\g }(f)).$ 
Then, combining 
\cite[Corollary 2.7]{S4} and 
(\ref{mainthran2lem3eq1}), we get that 
there exists a $j\in \NN $ such that 
\begin{equation}
\label{mainthran2lem3eq2}
\pi _{\CCI }(P^{\ast }(f)\cap \pi ^{-1}(\{ \beta \} ))\subset f_{\g ,m_{j}}(W_{1})=f_{\g ,m_{j}}(W_{2}).
\end{equation}
We set $W=f_{\g ,m_{j}}(W_{1})=f_{\g ,m_{j}}(W_{2}).$ 
Let $\{ V_{i}\} _{i=1}^{r}$ be the set of all connected components of 
$(f_{\g ,m_{j}})^{-1}(W).$ 
Since $W_{1}\neq W_{2}$, we have $r\geq 2.$ 
For each $i=1,\ldots ,r$, we denote by $p_{i}$ the 
number of critical points of $f_{\g ,m_{j}}:V_{i}\rightarrow W$ 
counting multiplicities. Moreover, we set 
$d_{i}:= \deg (f_{\g ,m_{j}}:V_{i}\rightarrow W).$ 
Then, by the Riemann-Hurwitz formula, we see that for each 
$i=1,\ldots ,r$, 
$\chi (V_{i})+p_{i}=d_{i}\chi (W).$ 
Since $\chi (V_{i})=\chi (W)=1$, it follows that 
\begin{equation}
\label{mainthran2lem3eq3}
r+\sum _{i=1}^{r}p_{i}=d, \mbox{ where }d:= \deg (f_{\g ,m_{j}}).
\end{equation} 
By (\ref{mainthran2lem3eq2}), we have $\sum _{i=1}^{r}p_{i}=d-1.$ 
Hence, (\ref{mainthran2lem3eq3}) implies 
$r=1$, which is a contradiction. Therefore, 
int$(K_{\g }(f))$ is a non-empty connected set. 
Combining it with Proposition~\ref{shonecomp}, 
we conclude that $J_{\g }(f)$ is a Jordan curve.   

 Thus, we have proved Lemma~\ref{mainthran2lem3}.
\end{proof}
We now demonstrate Theorem~\ref{mainthran2}.\\ 
{\bf Proof of Theorem~\ref{mainthran2}:} 
We suppose the assumption of Theorem~\ref{mainthran2}. 
We consider the following three cases. 

 \noindent Case 1: For each $\g \in \GN $, int$(K_{\g }(f))$ is connected. \\ 
 Case 2: For each $\g \in \GN $, int$(K_{\g }(f))$ is disconnected.\\ 
 Case 3: There exist two elements $\alpha \in \GN $ and $\beta \in \GN $ 
 such that int$(K_{\alpha }(f))$ is connected and such that 
 int$(K_{\beta }(f))$ is 
 disconnected. 

 Suppose that we have Case 1. 
 Then, by Theorem~\ref{hypskewqc}, there exists a constant $K\geq 1$ such that 
 for each $\g \in \GN $, $J_{\g }(f)$ is a $K$-quasicircle. 

 Suppose that we have  Case 2. Then, by 
 Lemma~\ref{mainthran2lem2}, we get that 
 for each $\g \in \GN $, there exist infinitely many 
 connected components of int$(K_{\g }(f)).$ 
 Moreover, by Theorem~\ref{mainthran1}, 
 we see that statement \ref{mainthran1-3} in Theorem~\ref{mainthran1} 
 holds. Hence, statement \ref{mainthran2-3} in Theorem~\ref{mainthran2} holds.

 Suppose that we have Case 3. 
 By Lemma~\ref{mainthran2lem2}, there exist infinitely many 
 connected components of int$(K_{\beta }(f)).$ 
 Let ${\cal W}:= \cup _{n\in \NN }(\sigma ^{n})^{-1}(\{ \beta \} ).$ 
 Then, for each $\g \in {\cal W}$, there exist infinitely many 
 connected components of int$(K_{\g }(f)).$ Moreover, 
 ${\cal W}$ is dense in $\GN .$ 
 
Next, combining Lemma~\ref{mainthran2lem3} and that the set of all periodic points of $\sigma :\GN \rightarrow 
  \GN $ is dense in $\GN $, we may assume that 
  the above $\alpha $ is a periodic point of $\sigma .$ 
  Then, $J_{\alpha }(f)$ is a quasicircle. 
 We set ${\cal V}:= \cup _{n\in \NN }(\sigma ^{n})^{-1}(\{ \alpha \} ).$ 
 Then ${\cal V}$ is dense in $\GN .$ 
 Let $\g \in {\cal V}$ be an element. 
Then there exists an $n\in \NN $ such that 
$\sigma ^{n}(\g )=\alpha .$ 
Since $(f_{\g ,n})^{-1}(K_{\alpha }(f))=
K_{\g }(f)$, it follows that 
$\sharp ({\cal C}($int$(K_{\g }(f))))<\infty .$ 
Combining it with Lemma~\ref{mainthran2lem2} 
and Proposition~\ref{shonecomp}, 
we get that $J_{\g }(f)$ is a Jordan curve. 
Combining it with 
 that $J_{\alpha }(f)$ is a quasicircle, 
 it follows that $J_{\g }(f)$ is a quasicircle. 

Next, let 
$\mu := \frac{1}{3}
\min \{ |b-c| \mid b\in J(G),\ c\in P^{\ast }(G)\} (>0).$ 
Applying Lemma~\ref{060416lemast} to 
$(f, \mu )$ above, 
let $\delta $ be the number in the statement of 
Lemma~\ref{060416lemast}. 
We set $\epsilon := \min \{ \delta ,\mu \} $ and 
$\rho := \beta .$  
Applying Lemma~\ref{mainthranlem2} to $(\G ,\alpha ,\rho ,\epsilon )$ 
above, let 
$(n,\theta _{1},\theta _{2},\omega )$ be the element 
in the statement of Lemma~\ref{mainthranlem2}. 
Let 
${\cal U}:= 
\{ \g \in \GN \mid \exists \{ m_{j}\} _{j\in \NN }, 
\exists \{ n_{k}\} _{k\in \NN },  
\sigma ^{m_{j}}(\g )\rightarrow \alpha , 
\sigma ^{n_{k}}(\g )\rightarrow \omega \} .$ 
Then, combining the statement of Lemma~\ref{060416lemast} 
and that of Lemma~\ref{mainthranlem2}, 
it follows that for any $\g \in {\cal U}$, 
$J_{\g }(f)$ is not a quasicircle. 
Moreover, by Lemma~\ref{mainthran2lem3}, 
we get that for any $\g \in {\cal U}$, 
$J_{\g }(f)$ is a Jordan curve. 
Combining the above argument, \cite[Theorem 1.12]{S4}, 
Lemma~\ref{fibconnlem}, and \cite[Theorem 9.3]{NV}, 
we see that for any $\g \in {\cal U}$, 
$A_{\g }(f)$ is a John domain, and the bounded component 
$U_{\g }$ of $F_{\g }(f)$ is not a John domain. 
Furthermore, 
it is easy to see that 
${\cal U}$ is residual in $\GN $, and that 
for any 
Borel probability measure $\tau $ on 
Poly$_{\deg \geq 2}$ with $\G _{\tau }=\G $, 
$\tilde{\tau }({\cal U})=1.$ 
Thus, we have proved Theorem~\ref{mainthran2}. 
%
%
\qed 
\begin{rem}
Using the above method (especially, using Lemma~\ref{jbnqlem1}, Lemma~\ref{060416lemast} and Lemma~\ref{mainthran2lem3}), 
we can also construct an example of a polynomial skew product $f:\CC ^{2}\rightarrow \CC ^{2}, 
f(x,y)=(p(x), q_{x}(y))$, where $p:\CC \rightarrow \CC $ is a polynomial with $\deg (p)\geq 2$, $q_{x}: \CC \rightarrow 
\CC $ is a monic polynomial with $\deg (q_{x})\geq 2$ for each $x\in \CC $, and $(x,y)\rightarrow q_{x}(y)$ is a polynomial of 
$(x,y)$,  
such that all of the following hold: 
\begin{itemize}
\item $f$ satisfies the Axiom A; and 
\item for almost every $x\in J(p)$ with respect to the maximal entropy measure of $p:\CC \rightarrow \CC $, 
the fiberwise Julia set $J_{x}(f)$ is a Jordan curve but not a quasicircle,  the fiberwise basin 
$A_{x}(f)$ of $\infty $ is a John domain, and the bounded component of $F_{x}(f)$ is not a John domain.  
\end{itemize}
For the related topics of Axiom A polynomial skew products on $\CC ^{2}$, see \cite{DH}. 
\end{rem}

We now demonstrate Proposition~\ref{ranprop1}.\\ 
{\bf Proof of Proposition~\ref{ranprop1}:}
Since $P^{\ast }(G)\subset $ int$(\hat{K}(G))\subset F(G)$, 
$G$ is hyperbolic.  
Let $\g \in \GN $ be any element. We will show the following claim.\\ 
Claim: int$(K_{\g }(f))$ is a non-empty connected set. 

 To show this claim, since $G$ is hyperbolic, int$(K_{\g }(f))$ is non-empty. 
Suppose that there exist two distinct connected components 
$W_{1}$ and $W_{2}$ of int$(K_{\g }(f)).$ 
Since $P^{\ast }(G)$ is included in a 
connected component $U$ of int$(\hat{K}(G))$ $
\subset F(G)$, 
\cite[Corollary 2.7]{S4} implies that 
there exists an $n\in \NN $ such that 
$P^{\ast }(G)\subset f_{\g ,n}(W_{1})=f_{\g ,n}(W_{2}).$ 
Let $W:= f_{\g ,n}(W_{1})=f_{\g ,n}(W_{2}).$ Then, 
any critical value of $f_{\g ,n}$ in $\CC $ is included in  $W.$ 
Using the method in the proof of Lemma~\ref{mainthran2lem3}, 
we see that $(f_{\g ,n})^{-1}(W)$ is connected. However, 
this is a contradiction, since $W_{1}\neq W_{2}.$ Hence, we have 
proved the above claim. 

 From Claim above and Theorem~\ref{hypskewqc}, 
 it follows that there exists a constant $K\geq 1 $ such that 
 for each $\g \in \GN $, 
 $J_{\g }(f)$ is a $K$-quasicircle. 

 Hence, we have proved Proposition~\ref{ranprop1}. 
\qed  
\subsection{Proofs of results in \ref{Const}}
We now demonstrate Proposition~\ref{Constprop}.\\ 
{\bf Proof of Proposition~\ref{Constprop}:}
Conjugating $G$ by $z\mapsto z+b$, we may assume that $b=0.$ 
For each $h\in \G$, we set $a_{h}:=a(h)$ 
and $d_{h}:=\deg (h).$ Let $r>0$ be a number 
such that $\overline{D(0,r)}\subset 
$ int$(\hat{K}(G)).$ 

Let $h\in \G $ and let $\alpha >0$ be a number. 
Since $d\geq 2$ and $(d,d_{h})\neq (2,2)$, 
it is easy to see that 
$(\frac{r}{\alpha })^{\frac{1}{d}}>
2\left(\frac{2}{|a_{h}|}(\frac{1}{\alpha })
^{\frac{1}{d-1}}\right)^{\frac{1}{d_{h}}}
$ if and only if 
\begin{equation}
\label{Contproppfeq1}
\log \alpha <
\frac{d(d-1)d_{h}}{d+d_{h}-d_{h}d}
( \log 2-\frac{1}{d_{h}}\log \frac{|a_{h}|}{2}-\frac{1}{d}\log r) .
\end{equation} 
We set 
\begin{equation}
\label{Contproppfeq2}
c_{0}:=\min _{h\in \G }\exp \left(\frac{d(d-1)d_{h}}{d+d_{h}-d_{h}d}
( \log 2-\frac{1}{d_{h}}\log \frac{|a_{h}|}{2}-\frac{1}{d}\log r) \right)
\in (0,\infty ).
\end{equation}
Let $0<c<c_{0}$ be a small number and let $a\in \CC $ 
be a number with $0<|a|<c.$ 
Let $g_{a}(z)=az^{d}.$ 
Then, we obtain $K(g_{a})=\{ z\in \CC \mid 
|z|\leq (\frac{1}{|a|})^{\frac{1}{d-1}}\} $ and 
$g_{a}^{-1}(\{ z\in \CC \mid  |z|=r\} )=
\{ z\in \CC \mid |z|=(\frac{r}{|a|})^{\frac{1}{d}}\} .$
Let 
$D_{a}:=\overline{D(0,2(\frac{1}{|a|})^{\frac{1}{d-1}})}.$ 
Since $h(z)=a_{h}z^{d_{h}}(1+o(1))\ (z\rightarrow \infty )$ uniformly on 
$\G $, 
it follows that if $c$ is small enough, then 
for any $a\in \CC $ with $0<|a|<c$ and for any $h\in \G $, 
$h^{-1}(D_{a})\subset 
\left\{ z\in \CC \mid 
|z|\leq 2\left( \frac{2}{|a_{h}|}(\frac{1}{|a|})^{d-1}\right) 
^{\frac{1}{d_{h}}}\right\} .$  
This implies that for each $h\in \G $, 
\begin{equation}
\label{Contproppfeq3}
h^{-1}(D_{a})\subset g_{a}^{-1}(\{ z\in \CC \mid |z|<r\} ).
\end{equation} 
Moreover, if $c$ is small enough, then for any $a\in \CC $ with 
$0<|a|<c$ and any $h\in \G $,  
\begin{equation}
\label{Contproppfeq4}
\hat{K}(G)\subset g_{a}^{-1}(\{ z\in \CC \mid |z|<r\} ),\ 
\overline{h(\CCI \setminus D_{a})}\subset 
\CCI \setminus D_{a}.
\end{equation}
Let $a\in \CC $ with $0<|a|<c.$  
By (\ref{Contproppfeq3}) and (\ref{Contproppfeq4}), 
there exists a compact neighborhood $V$ of $g_{a}$ in Poly$_{\deg \geq 2}$,  
such that 
\begin{equation}
\label{Contproppfeq5} 
\hat{K}(G)\cup \bigcup _{h\in \G }h^{-1}
(D_{a})
\subset 
\mbox{int}\left( \cap _{g\in V}
g^{-1}(\{ z\in \CC \mid |z|<r\} )\right), \mbox{ and } 
\end{equation}
\begin{equation}
\label{Contproppfeq5-a}
\bigcup _{h\in \G \cup V}
\overline{h(\CCI \setminus D_{a})}
\subset \CCI \setminus D_{a},
\end{equation}
which implies that 
\begin{equation}
\label{Contproppfeq5-1}
\mbox{int}(\hat{K}(G))\cup 
(\CCI \setminus D_{a})
\subset F(H_{\G ,V}),
\end{equation}
where $H_{\G ,V}$ denotes the polynomial semigroup generated by 
the family $\G \cup V.$   

 By (\ref{Contproppfeq5}), we obtain that for any non-empty subset 
 $V'$ of $V$, 
\begin{equation}
\label{Contproppfeq7}
\hat{K}(G)=\hat{K}(H_{\G, V'} ),
\end{equation} 
where $H_{\G ,V'}$ denotes the polynomial semigroup generated by 
the family $\G \cup V'.$
If the compact neighborhood $V$ of $g_{a}$ is so small, then  
\begin{equation}
\label{Contproppfeq8}
\bigcup _{g\in V} CV^{\ast }(g)  \subset  \mbox{int}(\hat{K}(G)).
\end{equation} 
Since $P^{\ast }(G)\subset \hat{K}(G)$, 
combining it with (\ref{Contproppfeq7}) and (\ref{Contproppfeq8}),
we get  
that for any non-empty subset $V'$ of $V$, 
$P^{\ast }(H_{\G ,V'} )\subset 
\hat{K}(H_{\G ,V'} ).$ 
 Therefore, for any non-empty subset $V'$ of $V$,  
$H_{\G ,V'} \in {\cal G}.$ 

We now show that for any non-empty subset $V'$ of $V$, $J(H_{\G,V'})$ is disconnected and 
$(\Gamma \cup V')_{\min }\subset \Gamma.$  
Let $$U:=\left(\mbox{int}(\bigcap _{g\in V}g^{-1}(\{ z\in \CC \mid |z|<r\} ))\right)
\setminus 
\bigcup _{h\in \G }h^{-1}(D_{a}).$$ 
Then, for any $h\in \G $, 
\begin{equation}
\label{Contproppfeq6}
h(U)\subset \CCI \setminus D_{a}.
\end{equation} 
Moreover, for any $g\in V$, $g(U)\subset 
$ int$(\hat{K}(G)).$ 
Combining it with (\ref{Contproppfeq5-1}), 
(\ref{Contproppfeq6}), and Lemma~\ref{hmslem}-\ref{bss}, 
it follows that $U\subset 
F(H_{\G ,V} ).$ 
If the neighborhood $V$ of $g_{a}$ is  so small, then 
there exists an annulus $A$ in $U$ such that for any $g\in V$,  
$A$ separates $J(g)$ and $\cup _{h\in \G }h^{-1}(J(g)).$  
Hence, it follows that for any non-empty subset $V'$ of $V$, 
the polynomial semigroup $H_{\G, V'}$ generated by 
the family $\G \cup V'$ satisfies that $J(H_{\G ,V'})$ is 
disconnected and $(\G \cup V')_{\min }\subset \G .$  

 We now suppose that in addition to the assumption, 
 $G$ is semi-hyperbolic. Let $V'$ be any non-empty subset of 
 $V .$ Since $(\G \cup \overline{V'})_{\min }\subset \G $, 
 Theorem~\ref{shshprop} implies that the above 
 $H_{\G ,V'}$ is semi-hyperbolic. 

 We now suppose that in addition to the assumption, 
$G$ is hyperbolic.  Let $V'$ be any non-empty subset of 
$V.$ By (\ref{Contproppfeq7}) and 
(\ref{Contproppfeq8}), 
we have 
\begin{equation}
\label{Contproppfeq9}
\bigcup _{g\in \Gamma \cup \overline{V'}}CV^{\ast }(g)\subset 
\mbox{int}(\hat{K}(H_{\G ,\overline{V'}})).
\end{equation}  
Since $(\G \cup \overline{V'})_{\min }\subset \G $, 
combining it with (\ref{Contproppfeq9}) and 
Theorem~\ref{hhprop}, 
we obtain that $H_{\G ,V'}$ is hyperbolic. 

 Thus, we have proved Proposition~\ref{Constprop}. 
\qed  

\ 

 We now demonstrate Theorem~\ref{shshfinprop}.\\ 
 {\bf Proof of Theorem~\ref{shshfinprop}:} 
First, we show \ref{shshfinprop1}. 
Let $r>0$ be a number such that 
$D(b_{j},2r)\subset \mbox{int}(K(h_{1}))$ for each 
$j=1,\ldots ,m.$ 
If we take $c>0$ so small, then 
for each $(a_{2},\ldots ,a_{m})\in \CC ^{m-1}$ 
such that $0<|a_{j}|<c$ for each $j=2,\ldots ,m$, 
setting $h_{j}(z)=a_{j}(z-b_{j})^{d_{j}}+b_{j}$ 
($j=2,\ldots ,m$), we have 
\begin{equation}
\label{shshfinpropeq1}
h_{j}(K(h_{1}))\subset D(b_{j},r)\subset 
\mbox{int}(K(h_{1}))\ (j=2,\ldots ,m). 
\end{equation} 
Hence, $K(h_{1})=\hat{K}(G)$, 
where $G=\langle h_{1},\ldots ,h_{m}\rangle .$ 
Moreover, by (\ref{shshfinpropeq1}), 
we have $P^{\ast }(G)\subset K(h_{1}).$ 
Hence, $G\in {\cal G}.$ 

 If $\langle h_{1}\rangle $ is semi-hyperbolic, 
then using the same method as that of Case 1 in the proof of 
Theorem~\ref{shshprop}, we obtain that $G$ is semi-hyperbolic. 

 We now suppose that $\langle h_{1}\rangle $ is hyperbolic. 
By (\ref{shshfinpropeq1}), we have 
$\cup _{j=2}^{m}CV^{\ast }(h_{j})\subset 
\mbox{int}(\hat{K}(G)).$ Combining it with 
the same method as that in the proof of Theorem~\ref{hhprop}, 
we obtain that $G$ is hyperbolic. 
Hence, we have proved statement \ref{shshfinprop1}. 

 We now show statement \ref{shshfinprop2}. 
Suppose we have case (i). 
We may assume $d_{m}\geq 3.$ 
Then, by statement \ref{shshfinprop1}, 
there exists an element $a>0$ such that 
setting $h_{j}(z)=a(z-b_{j})^{d_{j}}+b_{j}$ ($j=2,\ldots ,m-1$), 
$G_{0}=\langle h_{1},\ldots ,h_{m-1}\rangle $ satisfies 
that $G_{0}\in {\cal G}$ and $\hat{K}(G_{0})=$ $K(h_{1})$ 
and if $\langle h_{1}\rangle $ is semi-hyperbolic (resp. hyperbolic), 
then $G_{0}$ is semi-hyperbolic (resp. hyperbolic). 
Combining it with Proposition~\ref{Constprop}, 
it follows that there exists an $a_{m}>0$ such that 
setting $h_{m}(z)=a_{m}(z-b_{m})^{d_{m}}+b_{m}$, 
$G=\langle h_{1},\ldots ,h_{m}\rangle $ satisfies that 
$G\in {\cal G}_{dis}$ and $\hat{K}(G)=\hat{K}(G_{0})=K(h_{1})$ and if 
$G_{0}$ is semi-hyperbolic (resp. hyperbolic), 
then $G$ is semi-hyperbolic (resp. hyperbolic).  

 Suppose now we have case (ii). 
 Then by Proposition~\ref{Constprop}, 
 there exists an $a_{2}>0$ such that 
setting $h_{j}(z)=a_{2}(z-b_{j})^{2}+b_{j}$ $(j=2,\ldots ,m)$, 
$G=\langle h_{1},\ldots ,h_{m}\rangle =\langle h_{1}, h_{2}\rangle $ 
satisfies that $G\in {\cal G}_{dis}$ and $\hat{K}(G)=K(h_{1})$ and if $\langle h_{1}\rangle $ 
is semi-hyperbolic (resp. hyperbolic), then 
$G$ is semi-hyperbolic (resp. hyperbolic). 

 Thus, we have proved Theorem~\ref{shshfinprop}.  
\qed  

We now demonstrate Theorem~\ref{sphypopen}.\\ 
{\bf Proof of Theorem~\ref{sphypopen}:} 
Statements \ref{sphypopen2} and \ref{sphypopen3} follow 
from Theorem~\ref{shshfinprop}. 

 We now show statement \ref{sphypopen1}. 
By \cite[Theorem 2.4.1]{S5}, 
${\cal H}_{m}$ and ${\cal H}_{m}\cap {\cal D}_{m}$ are open. 

 We now show that ${\cal H}_{m}\cap {\cal B}_{m}$ is open. 
In order to do that, 
let 
$(h_{1},\ldots ,h_{m})\in {\cal H}_{m}\cap {\cal B}_{m}.$ 
Let $\epsilon >0$ such that 
$D(P^{\ast }(\langle h_{1},\ldots ,h_{m}\rangle ),\ 3\epsilon )
\subset F(\langle h_{1},\ldots ,h_{m}\rangle ).$ 
By \cite[Theorem 1.35]{S1}, 
there exists an $n\in \NN $ such that 
for each $(i_{1},\ldots ,i_{n})\in 
\{ 1,\ldots ,m\} ^{n}$, 
$$h_{i_{n}}\cdots h_{i_{1}}
(D(P^{\ast }(\langle h_{1},\ldots ,h_{m}\rangle ),\ 2\epsilon ))
\subset D(P^{\ast }(\langle h_{1},\ldots ,h_{m}\rangle ),\ \epsilon /2).$$ 
Hence, there exists a neighborhood $U$ of $(h_{1},\ldots ,h_{m})$ 
in (Poly$_{\deg \geq 2})^{m}$  such that for each 
$(g_{1},\ldots ,g_{m})\in U$ and each 
$(i_{1},\ldots ,i_{n})\in \{ 1,\ldots ,m\} ^{n}$, 
$$g_{i_{n}}\cdots g_{i_{1}}
(D(P^{\ast }(\langle h_{1},\ldots ,h_{m}\rangle ),\ 2\epsilon ))
\subset D(P^{\ast }(\langle h_{1},\ldots ,h_{m}\rangle ),\ \epsilon ).$$ 
If $U$ is small, then for each 
$(g_{1},\ldots ,g_{m})\in U$, 
$\cup _{j=1}^{m}CV^{\ast }(g_{j})\subset 
D(P^{\ast }(\langle h_{1},\ldots ,h_{m}\rangle ),\ \epsilon ).$ 
Hence, if $U$ is small enough, 
then for each $(g_{1},\ldots ,g_{m})\in U$, 
$P^{\ast }(\langle g_{1},\ldots ,g_{m}\rangle )
\subset D(P^{\ast }(\langle h_{1},\ldots ,h_{m}\rangle ), \epsilon ).$ 
Hence, for each $(g_{1},\ldots ,g_{m})\in U$, 
$\langle g_{1},\ldots ,g_{m}\rangle \in {\cal G}.$ 
Therefore, ${\cal H}_{m}\cap {\cal B}_{m}$ is open. 

 Thus, we have proved Theorem~\ref{sphypopen}.   
\qed

\end{document}